\documentclass[11pt,reqno]{amsart}
\usepackage[margin=1in]{geometry}
\usepackage{amsmath,amssymb,mathrsfs,amsthm,amsfonts,stmaryrd,enumitem}
\usepackage{bbm}

\usepackage[dvipsnames]{xcolor}

\usepackage{orcidlink}

\usepackage{caption}
\usepackage{subcaption}

\usepackage{graphicx}
\usepackage{pgfplots}
\pgfplotsset{compat=1.18}

\usepackage{tikz}
\usetikzlibrary{decorations.markings}
\usetikzlibrary{arrows}
\usetikzlibrary{calc}

\usepackage{float}

\usepackage{aligned-overset}

\usepackage{aliascnt}

\usepackage{hyperref}
\hypersetup{
     colorlinks=true,
     linkcolor=blue,
     filecolor=black,
     citecolor = blue,      
     urlcolor=blue,
     }

\usepackage{mathrsfs}
\usepackage{bm}
\usepackage{xspace}

\usepackage{scalerel}
\newcommand\smallbullet{{\raisebox{1.1pt}{\hspace{0.6pt}\scaleobj{0.55}{\bullet}\hspace{0.75pt}}}}

\newcommand{\nin}{\noindent}

\renewcommand{\ge}{\geqslant}
\renewcommand{\le}{\leqslant}

\renewcommand{\leq}{\leqslant}

\newcommand{\Rect}{\mathcal{R}}

\newcommand{\Cptwise}{\mathfrak{C}_{\ref{p:dov-pointwise}}}
\newcommand{\CHolder}{\mathfrak{C}_{\ref{p:dov-holder-general}}}
\newcommand{\CTFloc}{\mathfrak{C}_{\ref{p:DOV-geodesic-holder}}}
\newcommand{\CTF}{\mathfrak{C}_{\ref{p:GZ-uniform-TF}}}

\newcommand{\Cresample}{\mathfrak{C}_{\ref{l:uniform-bound-as}}}

\newcommand{\MeshGeoDeloc}{\mathsf{MeshGeoDeloc}}
\newcommand{\Inside}{\mathsf{MeshGeoInside}}
\newcommand{\BothIn}{\mathsf{BothInside}}
\newcommand{\BothDeloc}{\mathsf{BothDeloc}}
\newcommand{\NoBigTF}{{\mathsf{BigTFCostly}^{\delta}}}

\newcommand{\far}{{\mathrm{far}}}
\newcommand{\near}{{\mathrm{near}}}
\newcommand{\bdry}{{\mathrm{bdry}}}

\newcommand{\bL}{\bm{\mathcal{L}}}

\newcommand{\LBdryComp}[1]{{\cL}^{\bdry}_{#1}}
\newcommand{\LBdryCompRe}[1]{\widetilde{\cL}^{\bdry}_{#1}}
\newcommand{\LFarComp}[1]{{\cL}^{\far}_{#1}}
\newcommand{\LFarCompRe}[1]{\widetilde{\cL}^{\far}_{#1}}

\newcommand{\blue}{{\color[HTML]{648FFF}\textbf{blue}}\xspace}
\newcommand{\pink}{{\color[HTML]{DC267F}\textbf{pink}}\xspace}
\newcommand{\orange}{{\color[HTML]{FE6100}\textbf{orange}}\xspace}

\newcommand{\HReg}{\mathsf{PolymerTF}}
\newcommand{\CondReg}{\mathsf{PolymerCondTF}}
\newcommand{\HHol}{\textup{\textsf{EnergyH\"{o}lder}}}
\newcommand{\HOnePt}{\mathsf{EnergyConc}}
\newcommand{\Good}{\mathsf{Good}}
\newcommand{\Max}{\mathsf{PolymerDeloc}}

\newcommand{\XTF}{\mathrm{TF}}

\newcommand{\reg}{{\,\XTF_k}}

\newcommand{\Pann}{\overline{\mathbb{P}}}
\newcommand{\CDRP}{\mathrm{CDRP}}

\newcommand{\Rup}{\mathbb{R}^{4}_{\uparrow}}

\newcommand{\Fprod}{\mathcal{F}^{\mathrm{prod}}}

\newcommand{\res}{\mathrm{re}}
\newcommand{\Xires}{\Xi_\res}
\newcommand{\Fres}{\mathcal{F}_\res}
\newcommand{\Pres}{\P_\res}
\newcommand{\Eres}{\E_\res}

\newcommand{\Res}{\mathrm{res}} 

\newcommand{\Prob}{\mathrm{Prob}}

\newcommand{\TT}{T_{\delta}}

\newcommand{\transaff}{\mathtt{T}}

\newcommand{\shearaff}{\mathtt{S}}
\newcommand{\dilaff}{\mathtt{D}}

\newcommand{\autaff}{\mathtt{Q}}

\newcommand{\trans}{\mathrm{T}}

\newcommand{\shear}{\mathrm{S}}
\newcommand{\dil}{\mathrm{D}}

\newcommand{\aut}{\mathrm{Q}}

\newcommand{\Trans}{\mathcal{T}}

\newcommand{\Shear}{\mathcal{S}}
\newcommand{\Dil}{\mathcal{D}}

\newcommand{\Aut}{\mathcal{Q}}

\newcommand{\Pbr}{\mathrm{P}_{\mathrm{BB}}}
\newcommand{\Ebr}{\mathrm{E}_{\mathrm{BB}}}
\DeclareMathOperator{\VarXi}{\mathbb{V}ar}

\newcommand{\ee}{\mathrm{e}}

\renewcommand{\P}{\mathbb{P}}
\newcommand{\E}{\mathbb{E}}

\renewcommand{\th}{\textsuperscript{th}}

\newcommand{\norm}[1]{\left\lVert#1\right\rVert}

\DeclarePairedDelimiter\floor{\lfloor}{\rfloor}
\newcommand{\lb}{\llbracket}
\newcommand{\rb}{\rrbracket}

\newcommand{\wh}{\widehat}
\newcommand{\wt}{\widetilde}

\newcommand{\R}{\mathbb{R}}  
\newcommand{\Z}{\mathbb{Z}}
\newcommand{\Zpos}{{\mathbb{Z}_{\ge0}}}
\newcommand{\N}{\mathbb{N}}

\newcommand{\e}{\varepsilon}
\renewcommand{\d}{\delta}

\newcommand{\ls}{\lesssim}
\newcommand{\gs}{\gtrsim}

\DeclareMathOperator{\Corr}{Corr}
\DeclareMathOperator{\Cov}{Cov}
\DeclareMathOperator{\Var}{Var}

\newcommand{\law}{\overset{d}{=}}

\newcommand{\dto}{\xrightarrow{\,d\,}}

\newcommand{\pto}{\xrightarrow{\,p\,}}

\newcommand{\1}{\mathbf{1}}

\newcommand{\mf}{\mathfrak}

\newcommand{\mrm}{\mathrm}

\newcommand{\cB}{\mathcal{B}}

\newcommand{\cE}{\mathcal{E}}
\newcommand{\cF}{\mathcal{F}}
\newcommand{\cG}{\mathcal{G}}
\newcommand{\cH}{\mathcal{H}}
\newcommand{\cI}{\mathcal{I}}

\newcommand{\cK}{\mathcal{K}}
\newcommand{\cL}{\mathcal{L}}
\newcommand{\cM}{\mathcal{M}}
\newcommand{\cN}{\mathcal{N}}

\newcommand{\cR}{\mathcal{R}}

\newcommand{\cZ}{\mathcal{Z}}

\newcommand{\sA}{\mathsf{A}}
\newcommand{\sB}{\mathsf{B}}

\newcommand{\sE}{\mathsf{E}}
\newcommand{\sF}{\mathsf{F}}
\newcommand{\sG}{\mathsf{G}}

\newcommand{\sN}{\mathsf{N}}

\newcommand{\sS}{\mathsf{S}}

\newcommand{\fC}{\mathfrak{C}}

\newcommand{\fL}{\mathfrak{L}}

\newcommand{\bfX}{\mathbf{X}}

\newcommand{\bfs}{\mathbf{s}}

\newcommand{\bfu}{\mathbf{u}}

\newcommand{\bfx}{\mathbf{x}}
\newcommand{\bfy}{\mathbf{y}}
\newcommand{\bfz}{\mathbf{z}}

\DeclareRobustCommand{\SkipTocEntry}[5]{} 

\usepackage[capitalise]{cleveref}

\theoremstyle{plain}
\newtheorem*{informal-thm}{Informal theorem}

\newtheorem{thm}{Theorem}[section]

\newaliascnt{lemm}{thm}
\newtheorem{lemm}[lemm]{Lemma}
\aliascntresetthe{lemm}

\crefname{lemm}{Lemma}{Lemmas}
\Crefname{lemm}{Lemma}{Lemmas}

\newaliascnt{prop}{thm}
\newtheorem{prop}[prop]{Proposition}
\aliascntresetthe{prop}

\crefname{prop}{Proposition}{Propositions}
\Crefname{prop}{Proposition}{Propositions}

\newaliascnt{cor}{thm}
\newtheorem{cor}[cor]{Corollary}
\aliascntresetthe{cor}

\crefname{cor}{Corollary}{Corollaries}
\Crefname{cor}{Corollary}{Corollaries}

\newaliascnt{conj}{thm}
\newtheorem{conj}[conj]{Conjecture}
\aliascntresetthe{conj}

\crefname{conj}{Conjecture}{Conjectures}
\Crefname{conj}{Conjecture}{Conjectures}

\newaliascnt{claim}{thm}

\aliascntresetthe{claim}

\crefname{claim}{Claim}{Claims}
\Crefname{claim}{Claim}{Claims}

\newaliascnt{obs}{thm}

\aliascntresetthe{obs}

\crefname{obs}{Observation}{Observations}
\Crefname{obs}{Observation}{Observations}

\theoremstyle{definition}

\newaliascnt{defn}{thm}
\newtheorem{defn}[defn]{Definition}
\aliascntresetthe{defn}

\crefname{defn}{Definition}{Definitions}
\Crefname{defn}{Definition}{Definitions}

\newaliascnt{qn}{thm}

\aliascntresetthe{qn}

\crefname{qn}{Question}{Questions}
\Crefname{qn}{Question}{Questions}

\theoremstyle{remark}

\newaliascnt{rem}{thm}
\newtheorem{rem}[rem]{Remark}
\aliascntresetthe{rem}

\crefname{rem}{Remark}{Remarks}
\Crefname{rem}{Remark}{Remarks}

\numberwithin{equation}{section}

\theoremstyle{plain}

\theoremstyle{definition}

\theoremstyle{remark}

\numberwithin{equation}{section}
\numberwithin{thm}{section}

\title{Temperature chaos in directed polymers}

\author{Shirshendu Ganguly \and Victor Ginsburg \and Zoe Himwich}

\address{Shirshendu Ganguly\\
Department of Statistics,
University of California, Berkeley, CA, USA} 
\email{sganguly@berkeley.edu}

\address{Victor Ginsburg \orcidlink{0000-0001-9399-6748}\\
Department of Mathematics, University of California, Berkeley, CA, USA}
\email{victor@math.berkeley.edu}

\address{Zoe Himwich\\
Department of Statistics, University of California, Berkeley, CA, USA}
\email{himwich@berkeley.edu}

\begin{document}

\begin{abstract}
    Disordered systems such as spin glasses and polymers characteristically exhibit random energy landscapes with many macroscopically separated energetic valleys corresponding to near-ground states.
    This high complexity renders these systems extremely sensitive to perturbations of external parameters, such as the temperature or the random disorder itself.
    For instance, the support of associated Gibbs measures may change macroscopically under such perturbations, a phenomenon known as \emph{chaos} in the literature. 
    In the experimental literature, chaotic phenomena are typically studied via temperature perturbations owing to the ease of tuning the temperature in laboratory experiments.
    Several non-rigorous studies of disorder and temperature chaos for low-dimensional polymers have predicted exponents governing the transition from stability to chaos \cite{FV88,FH91,SY02,dSB04}, 
    whereas mathematical studies have been primarily restricted to disorder chaos, relying on the additional randomness imparted by perturbing the disorder \cite{Cha14,GH20}.

    In this article, we initiate the rigorous study of temperature-chaotic properties of the continuum directed random polymer (CDRP), a canonical model in the KPZ universality class.
    The CDRP is driven by white noise and is parametrized by inverse temperature $\beta$, and is known \cite{W23,DZ24} to converge in the zero-temperature limit $\beta\to\infty$ to the directed landscape  constructed in \cite{dov}, the putative universal scaling limit of models in the KPZ universality class.
    The main result of this article considers the  CDRP free energies coupled through the same white noise at a pair of inverse temperatures $(\beta_1,\beta_2)$, and shows that they decouple in the limit $\beta_2\gg \beta_1\gg 1$, converging to a pair of independent directed landscapes.
    This is the first such ``energetic de-correlation across temperatures'' result.
    Our key estimate measures the ``pivotality'' or ``influence'' of spatially thin strips in models of last passage percolation.
    As an interesting byproduct, the proof strategy also allows to show that the directed landscape is a two-dimensional black noise (in the sense of Tsirelson and Vershik \cite{Tsir}),
    previously conjectured by  Vir\'ag.
    This provides the third known example of a two-dimensional black noise after critical planar percolation \cite{SS11} and the Brownian web \cite{EF16}.
    Furthermore, the fragility of Gibbs measures on high-dimensional spaces such as the path space is exhibited by proving the mutual singularity of the CDRP measures across different temperatures, regardless of the free energy correlation behavior. 
    This relies on the Gaussian multiplicative chaos perspective on CDRP put forth in \cite{QRV25}.
    Finally, a precise conjecture for the sharp exponent governing temperature chaos is formulated, based on a heuristic coarse-graining argument linking the notions of temperature chaos and disorder chaos.
    Such a connection was speculated earlier in the physics literature.
\end{abstract}

\maketitle

\clearpage

\begin{figure}[p]
\centering
  \begin{subfigure}[c]{0.7\textwidth}
    \centering
    \includegraphics[width=\linewidth]{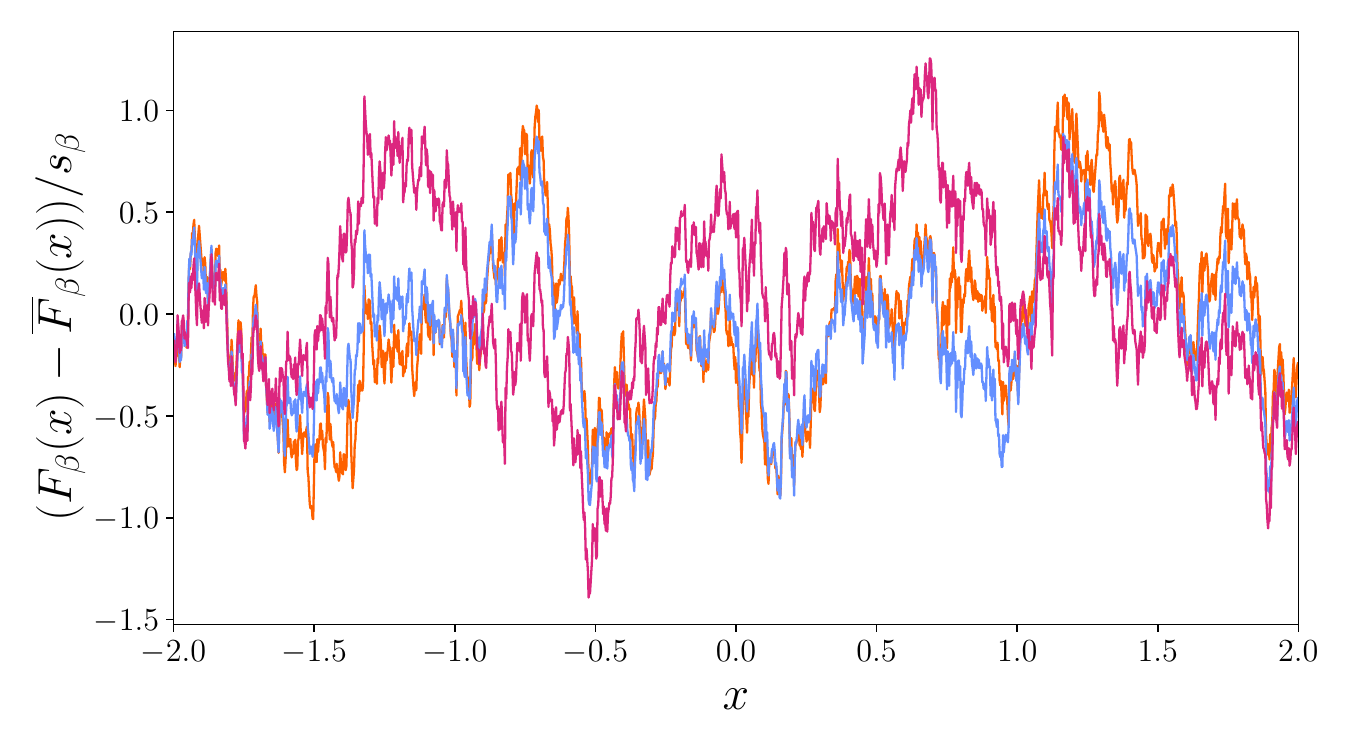}
  \end{subfigure}
  \hfill
  \begin{subfigure}[c]{0.28\textwidth}
    \centering
    \includegraphics[width=\linewidth]{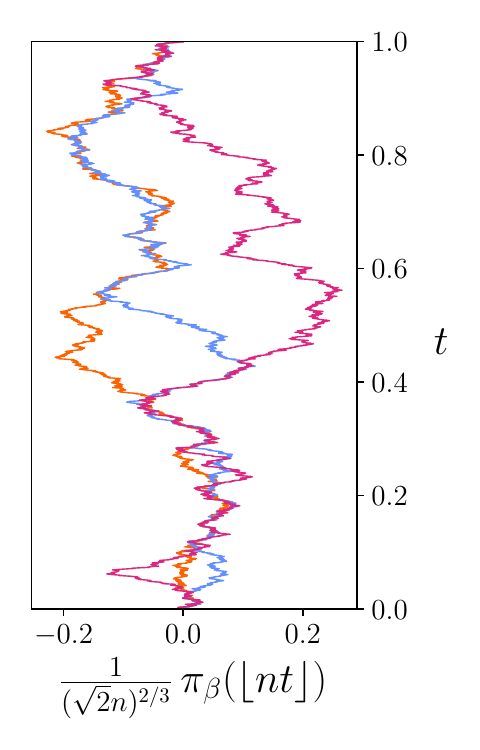}
  \end{subfigure}
\caption{ 
    A simulation of lattice directed polymers on $\Z^{2}$ (rotated by $45^{\circ}$) and their free energies at three different inverse temperatures, coupled through a single instance of an i.i.d Gaussian environment.
    The polymers all start at $(0,0)$ and travel upwards to height $n$ with $n = 2^{14}$ taking north-east $(1,1)$ or north-west $(-1,1)$ steps.
    In the plots, the colors orange, blue, pink  correspond respectively to inverse temperatures $\beta_1,\beta_2,\beta_3$,
    where $\beta_1=0.5$, and where $\beta_2,\beta_3$ are related to $\beta_1$ by 
    $\beta_2-\beta_1 = 0.79 \Delta_n \approx 0.157$ and $\beta_3-\beta_1= 2.38 \Delta_n \approx 0.472$.
    Here $\Delta_n\coloneqq n^{-1/6}$ is the predicted critical temperature perturbation scale marking the onset of temperature-chaotic behavior for small perturbations
    of $\beta_1$ (\cref{c:main} elaborates on this).
    \\
    \textbf{Left (free energies):}
    For each inverse temperature $\beta=\beta_1,\beta_2,\beta_3$, we simulate the corresponding point-to-point free energy from $(0,0)$ to $(\floor{x n^{2/3}}, n)$, denoted $F_{\beta}(x)$, for $x\in[-2,2]$.
    Plotted are $x \mapsto \frac{1}{s_\beta}(F_\beta(x) - \overline{F}_\beta(x))$ for $\beta=\beta_1,\beta_2,\beta_3$, where $\overline{F}_\beta(x)$ is the empirical mean of $F_\beta(x)$, and $s_\beta$ is the empirical standard deviation of $F_\beta(0)$ (computed using $N=50$ samples).
    This centering and scaling ensures that each free energy profile approximates the stationary Airy$_2$ process normalized to have unit one-point variance.
    As is evident from the simulation, and in agreement with the predicted criticality of $\Delta_n$, the three profiles lie in the critical window, 
    with the $\beta_1,\beta_2$ pair exhibiting a stronger correlation than the $\beta_1,\beta_3$ pair.
    \\
    \textbf{Right (polymers):}
    For each $\beta=\beta_1,\beta_2,\beta_3$, a sample from the corresponding directed polymer measure between $(0,0)$ and $(n,0)$ given the same random environment as in the left panel is displayed.
    Denoting the polymer by $\pi_\beta : [0,n]\cap\Z \to \Z$,
    plotted are the paths $[0,1]\ni t\mapsto \frac{1}{(\sqrt{2}n)^{2/3}}\pi_\beta(\floor{nt})$ for $\beta=\beta_1,\beta_2,\beta_3$.
    Again, the scaling ensures that each polymer approximates the directed landscape geodesic from $(0,0)$ to $(0,1)$.
    Consistent with the free energy plot,  the polymers at inverse temperatures $\beta_1,\beta_2$ have a significantly larger overlap than the polymers at $\beta_1,\beta_3$.
}\label{fig:simulation-free-energy}
\end{figure}

\clearpage

\setcounter{tocdepth}{2}{
  \hypersetup{linkcolor=black}
  \tableofcontents
}

\section{Introduction}\label{s:intro}

\subsection{Preface}\label{ss:preface}
Many models of central importance in probability and statistical mechanics, particularly in the study of many-body interactions, are formulated in terms of Gibbs measures.
A Gibbs measure is a probability measure which assigns a weight to each state through a Hamiltonian.
For a finite, high-dimensional set $S$ equipped with a measure $\mu$ (for instance, $\{-1,1\}^n$ with the uniform measure) and a Hamiltonian $H : S \to \R$,
the Gibbs weight assigned to $\sigma \in S$  is proportional to $\exp(-\beta H(\sigma)) \mu(\{\sigma\})$, where $\beta=1/T$ denotes inverse temperature, and where $H(\sigma)$ is interpreted as the \emph{energy} of the state $\sigma$.
Viewed together with the underlying geometry of $S$, the Hamiltonian defines an energy landscape on $S$.
We refer to a state as a \emph{ground state} if it minimizes $H$. 
In many models, including polymers and spin glasses, the Hamiltonian is driven by external \emph{disorder} given by independent randomness.
Generically, such random Hamiltonians have many qualitatively distinct near-ground states, corresponding to an energy landscape with many macroscopically separated deep valleys.
This complexity of the static picture has important dynamical implications, including the extreme fragility of such systems under perturbations of external parameters such as the disorder or temperature.
For instance, such perturbations can change which valley dominates and thereby produce macroscopic changes in observables such as the ground state, a phenomenon termed \emph{chaos} in the literature.
This is closely related to the notion of \emph{noise sensitivity} first developed in the context of Boolean functions in theoretical computer science \cite{BKS99}.

\addtocontents{toc}{\SkipTocEntry}
\subsection*{Dynamical models}
In the physics literature, perturbations of the parameters have often involved varying the system temperature or the underlying disorder. We briefly expand on these commonly adopted dynamical models of perturbation.

\addtocontents{toc}{\SkipTocEntry}
\subsubsection*{Disorder perturbation}
The random disorder in a statistical mechanical model may be discrete or continuous depending on the definition of the model.
A discrete disorder may be a countable collection
of independent random variables (for example i.i.d. $\mathrm{Bernoulli}(1/2)$),
whereas a continuum disorder might take the form of a Gaussian white noise.
There is a canonical notion of perturbation in each case.
In the discrete case, each random variable is independently resampled at the ring times of independent Poisson clocks of rate one.
Typically, under suitable rescaling the discrete disorder converges to Gaussian white noise $\xi$, and the corresponding resampling dynamics converges
to stationary Ornstein--Uhlenbeck dynamics
$\{\xi_t\}_{t\ge 0}$, where each $\xi_t$ is marginally a Gaussian white noise.
It is illuminating to note that the Ornstein--Uhlenbeck dynamics has the following two-point description:
\begin{align*}
    \bigl(\xi_0,\, \xi_t\bigr) \overset{d}{=} \bigl(\xi,\; e^{-t}\xi + \sqrt{1-e^{-2t}}\xi'\bigr)
\end{align*}
where $\xi'$ is an independent copy of $\xi$.
For details on Ornstein--Uhlenbeck dynamics on Gaussian white noise, see \cite[Section 1.4]{nualart}.

\addtocontents{toc}{\SkipTocEntry}
\subsubsection*{Temperature perturbation}
For positive-temperature models, with inverse temperature $\beta\in(0,\infty)$, varying $\beta$ results in a reweighting of the same realization of the disordered energy landscape.
This produces a family of Gibbs measures indexed by $\beta$ and coupled through a shared random disorder.\\

We refer to chaotic behavior of observables under these dynamical models, respectively, as \emph{disorder chaos} and \emph{temperature chaos}.
In experiments, it is common to study a temperature perturbation, since temperature is typically far easier to control experimentally than the disorder.
However, as might not be surprising, most mathematical developments have focused on disorder perturbations, profiting from the infusion of extra independent randomness.\\

While there have been exciting recent developments around the notion of disorder chaos (to be reviewed shortly), the central aim of this article is to make the first rigorous progress in probing temperature chaotic behavior of models in the Kardar--Parisi--Zhang (KPZ) universality class.
For parallel developments in mean-field spin glasses, see \cite{Cp17,Che14,AC16,BSZ20}.
Formally stating the main result requires significant preparation, so for now we present an informal version to act as a signpost for the overarching goal.

\begin{informal-thm}
    For the continuum directed random polymer, a canonical example in the KPZ universality class,
    its free energies at different inverse temperatures $\beta_1$ and $\beta_2$ become asymptotically independent when $\beta_1, \beta_2 \to \infty$ such that $\beta_1=o(\beta_2)$.
\end{informal-thm}

\nin

A conjecture for a sharp version of the above result is also formulated, and justified heuristically by drawing a connection between disorder and temperature chaos.
Such a connection was previously speculated in \cite{SY02,dSB04}.
Sharp predictions can also be found in the non-rigorous literature \cite{FH91,SY02,dSB04} based on numerical experiments and scaling arguments.
Formal versions of the above theorem and conjecture appear later in \cref{ss:results}. 
\\

\nin
The last three decades have seen substantial progress in the rigorous study of chaos and the companion notion of \emph{noise sensitivity}, which studies the degree of disorder perturbation needed for an observable of interest to substantially lose correlation with its initial value. 
The aforementioned work \cite{BKS99} developed a general theory of noise sensitivity for Boolean functions based on Fourier analysis on the discrete hypercube. 
Informally, a function $f$ is \emph{noise sensitive} iff its projection on the Fourier modes of low frequency is negligible.
A particular example of focus in \cite{BKS99} was critical planar percolation. In the subsequent study \cite{GPS}, a remarkably refined picture of the Fourier spectrum in this case was obtained relying on the geometry of pivotal points (see also \cite{GPS13,GPS18}).

\addtocontents{toc}{\SkipTocEntry}
\subsection*{Black noise}

At around the same time as the above developments in discrete noise sensitivity, Tsirelson and Vershik \cite{TV} introduced a related  notion of a \emph{black noise} as part of an axiomatic framework for the notion of a \emph{noise} (see \cref{d:noise}).
Their framework centers around a vast generalization of the Wiener chaos decomposition, a fundamental tool in the study of Gaussian white noise.
In their terminology, a \emph{classical} noise is generated by its first chaos.
This class includes Gaussian and Poissonian white noises and their combinations.
In contrast, a black noise has trivial first chaos, which informally means that for any observable, its Fourier expansion is supported on modes of ``infinite frequency.''
A celebrated result of Schramm and Smirnov \cite{SS11} proved that scaling limits of critical planar percolation are black noises, despite the pre-limiting model being driven by a discrete version of white noise (namely i.i.d. Bernoulli random variables).
In fact, they proved that these scaling limits are \textit{two-dimensional} black noises (the definition is somewhat involved and appears later in  \cref{def:black-noise}).
Subsequently, Ellis and Feldheim \cite{EF16} proved that the Brownian web is also a two-dimensional black noise, making these the only two examples of two-dimensional black noises in the literature to date.

Closer to the theme of the present paper, a conjecture of B\'alint Vir\'ag (recorded in \cite[Conjecture 1.17]{HP24}) asserts that the directed landscape, the universal scaling limit of models of directed geometry in the KPZ universality class constructed in \cite{dov}, to be defined shortly, 
is a two-dimensional black noise.
As an application of the methods developed to prove our main result on temperature chaotic properties of polymers,
we establish this conjecture, 
thereby providing a third example of a two-dimensional black noise. \\

This is an apt moment to introduce the KPZ universality class and recent developments in this topic along with the predictions about its dynamical properties to set up the necessary foundation for the statements of our results.

\addtocontents{toc}{\SkipTocEntry}
\subsection*{Polymer models and KPZ universality}

Kardar, Parisi, and Zhang \cite{KPZ86} introduced the Kardar--Parisi--Zhang (KPZ) equation, the nonlinear stochastic PDE given formally by
\begin{align}\label{e:kpz-1}
    \partial_{t}\mathcal{H} = \frac{1}{2}\partial_{x}^{2}\mathcal{H}  + \frac{1}{2}\left(\partial_{x}\mathcal{H}\right)^{2} + \xi
\end{align}
(here $\xi$ is a space-time white noise, see \cref{s:meas-short} for background),
as a universal model for random interface growth.
The KPZ equation is the central object in the \emph{KPZ universality class}, a vast family of models described by random evolving height functions that are predicted to exhibit the same large-scale fluctuation behavior.
Non-rigorous scaling arguments going back to \cite{KPZ86} predict exponents of $1/3$ for height fluctuations and $2/3$ for the spatial correlation scale of all such models.
These observables also follow universal non-Gaussian limiting laws, such as the Tracy--Widom distributions that govern fluctuations of edge eigenvalues of large random matrices.
We discuss the KPZ equation itself in \cref{ss:KPZ}. 
For background on KPZ universality see the surveys \cite{C16,Q11}, with some of its connections to random geometry reviewed in \cite{G22}.

The central principle of KPZ universality is that the large-scale statistics depend only on coarse features such as dimension, local growth asymmetry,
and noise structure, rather than on microscopic details.
Systems believed to lie in the KPZ class include stochastic interface growth models, asymmetric exclusion processes, directed polymers in random environments, last passage percolation, and many more. 
This article focuses on the continuum directed random polymer (CDRP), a canonical continuous directed polymer model associated with the KPZ equation.
We defer the formal discussion to \cref{ss:CDRP}, but as a reference point for the current discussion, we introduce the discrete directed polymer model as a random Gibbs measure on up-right, or directed, lattice paths $\pi$ in $\mathbb{Z}^2$ between $(0,0)$ and $(n,n)$, such that the probability of a path $\pi$ is proportional to 
\begin{equation}\label{strong23}
    \exp\left(\beta \sum_{v\in \pi}X_v\right),
\end{equation}
where $\{X_v\}_{v\in \Z^2}$ is an independent, identically distributed field of disorder and $\beta=1/T$ is inverse temperature. 
Thus, for any directed path $\pi$, its \emph{energy} is the sum of disorder variables along it.

Last passage percolation (LPP), a central model for local random growth, is the zero-temperature limit ($\beta\to\infty$) of the directed polymer \eqref{strong23}.
Namely, in lattice LPP, one considers the \emph{maximum} energy of any up-right path between $(0,0)$ and $(n,n)$.
Paths with maximal energy are called \emph{geodesics}, and correspond to ground states of the directed polymer model.
Natural observables include the geodesics themselves and their energies.

Although the fluctuation behavior of LPP is expected to be universal under mild conditions on the disorder, for instance for i.i.d. weights with sufficiently light tails, 
rigorous verification of this behavior has largely been restricted to a handful of \emph{integrable} models with additional algebraic structure, a particular example being when the disorder variables are given by i.i.d. Exponential variables.
In LPP and polymer models, such integrability often arises via the Robinson--Schensted--Knuth (RSK) correspondence.
This goes back to the early one-point fluctuation results \cite{BDJ99,Joh00}.
A continuous analogue of the RSK correspondence was developed by O'Connell and Yor \cite{OCY02} and O'Connell \cite{OC03}.
In a breakthrough development, Dauvergne, Ortmann, and Vir\'ag \cite{dov} used the continuous RSK correspondence to construct the \emph{directed landscape}, a full space-time scaling limit of last passage percolation models (see \cref{d:dl} below).
Their construction relied on a remarkable generalization of the RSK correspondence, versions of which had previously featured in \cite{BBOC05,NY04}.
Around the same time as \cite{dov}, Matetski, Quastel, and Remenik \cite{MQR21} constructed the KPZ fixed point, the universal Markov process governing the evolution of KPZ-scaled height profiles.
The reader is referred to \cite{Zyg22} for a more comprehensive review of the role of the RSK correspondence in KPZ universality.

The construction of the directed landscape in \cite{dov} relied crucially on connections between Brownian LPP (an integrable LPP model) and the \emph{Airy line ensemble} \cite{CH}, an infinite family of non-intersecting curves whose top curve is the scaling limit of LPP energy profiles.
A closely-related set of tools compares Airy/KPZ energy profiles to Brownian motion, beginning with the Brownian Gibbs property for the Airy line ensemble \cite{CH} and refined in increasingly strong forms in \cite{Ham22,CHH,duncan2}.
These Brownian comparison estimates will also feature centrally in our arguments.

\addtocontents{toc}{\SkipTocEntry}
\subsection*{Dynamical aspects of KPZ universality}

Since the 1980s, there have been various non-rigorous, both numerical and heuristic, approaches to studying dynamical perturbations of the disorder and temperature of polymer models.
In \cite{Zha87} Zhang  considered the geometry of near maximizers and their clustering properties and how the ground state jumps from one cluster to the other upon perturbation.
M\'ezard \cite{Mz90} studied positive-temperature models and properties of samples with interactions penalizing overlap between them.
Feigel'man and Vinokur \cite{FV88} and Fisher and Huse \cite{FH91} developed scaling arguments predicting the critical exponents governing decorrelation/relaxation behavior of dynamical $1+1$-dimensional polymer models in the KPZ universality class.
More recently, this critical behavior was studied in greater depth in \cite{SY02,dSB04} via numerical experiments, the replica Bethe ansatz \cite{K87}, and a $1+\e$-expansion.

A particularly interesting idea put forth in \cite{SY02,dSB04}, and described in \cite{dSB04} as somewhat non-intuitive, is that changing temperature should have an effect comparable to perturbing the disorder.
This claim is further supported by numerical findings in \cite{SY02}, where related arguments based on the replica method \cite{K87} are also discussed.
Though we will not pursue this rigorously, we will sketch an argument in support of this later in this article (see \cref{s:heuristic}).
\\

Subsequently, there have been important developments which have succeeded in establishing a part of the predicted picture on a mathematically sound footing.

\addtocontents{toc}{\SkipTocEntry}
\subsection*{Chaos, superconcentration, and multiple valleys}
In \cite{Cha14}, Chatterjee put forth a unified framework  
in which he established an equivalence between the dynamical notion of disorder chaos, and the static notions of superconcentration of the free energy (exhibiting stronger concentration than what the Gaussian Poincar\'e inequality predicts) and a multiple-valleyed nature of the associated energy landscape (existence of many near-optimal states that are well-separated in the configuration space).
Related and refined results for spin glass models were later obtained in \cite{DEZ15,CHL18}.

Within this general framework, \cite{Cha14} considered a random matrix version of the disorder chaos problem.
Given an $n\times n$ Hermitian matrix, its quadratic form defines an energy landscape on the unit sphere, with the ground state being the eigenvector associated to the largest eigenvalue.
A natural model for a dynamical $n\times n$ random matrix is the Gaussian unitary ensemble, with its disorder perturbed by running an Ornstein--Uhlenbeck flow for time $t>0$ independently in each Gaussian entry.
Chatterjee \cite{Cha14} found that chaos---the near-orthogonality of the time-$0$ and time-$t$ ground state eigenvectors---must occur by any time $t\gg n^{-1/3}$. 
Subsequently,  Bordenave, Lugosi, and Zhivotovskiy \cite{BLZ20} studied a discrete dynamics in which $k$ random entries in a random $n\times n $ Wigner matrix are resampled, establishing the sharpness---and hence universality---of the $-1/3$ dynamical exponent for random matrices. 
A dynamical approach has also been central in proofs of random matrix universality for spectral statistics and quantum ergodicity \cite{ESY11,ESYY12}.

\addtocontents{toc}{\SkipTocEntry}
\subsection*{Noise-sensitivity and chaos in polymers}

As indicated already, in critical planar percolation, the success of the dynamical studies \cite{GPS,GPS18} relied heavily on a refined understanding of the static picture involving the geometry of pivotal points which exhibits a rich fractal structure \cite{SW01,GPS13}.
Analogous fractal properties of LPP and the directed landscape have only recently been explored in a series of works.
This includes Hausdorff dimension results for exceptional sets of disjoint or non-unique geodesics \cite{BGH22,BGH21,GZ22,Bha24,Dau25}, as well as the structure and coalescence of semi-infinite geodesics and competition interfaces \cite{RV25,BSS24,bus25} (see also the survey \cite{GH24}).
Such developments have made it plausible to initiate a program to study rigorously how models in the KPZ class react to natural dynamics.

Towards this, Ganguly and Hammond \cite{GH20,GH23} undertook a systematic study of the transition from stability to chaos in Brownian last passage percolation (LPP), a central semi-discrete model in the KPZ universality class. In this model, the environment consists of independent white noises (for a thorough definition, see \cite[Section 1.2]{GH20}). 
An up-right path accumulates energy by integrating this noise along its horizontal segments, and a geodesic is a path of maximal energy between its endpoints.
The dynamical model is obtained by evolving the white noise environment under an Ornstein--Uhlenbeck flow.
\cite{GH20} built on a collection of static estimates developed in \cite{GH23} to show that geodesics at times $0$ and $t$ have non-trivial overlap when $t\ll n^{-1/3}$, while their overlap is negligible when $t\gg n^{-1/3}$.
Thus $n^{-1/3}$ is the critical time scale for the onset of geodesic chaos in Brownian LPP,
matching the counterpart results for random matrices discussed above.
While Chatterjee's theory \cite{Cha14} implies that the energies at times $0$ and $t$ exhibit strong correlation when $t\ll n^{-1/3}$, proving energetic decorrelation in the super-critical regime $t\gg n^{-1/3}$ remains an interesting open problem.
Subsequent works of Bhatia \cite{bha25,bhatia2025} studied the related problems of switching of geodesics under such an Ornstein--Uhlenbeck flow and the potential existence of exceptional times admitting bi-geodesics passing through the origin.

On the continuum side, \cite{HP24} showed that the directed landscape is a one-dimensional black noise with respect to its time-indexed noise structure. We will discuss one-dimensional and two-dimensional black noise in \cref{s:black-noise}. 
As already indicated, one of the results of this paper will strengthen this and establish that the directed landscape is a \emph{two-dimensional} black noise.
\\
 
The main focus of this article, however, is \emph{temperature chaos} for polymers in the KPZ universality class.
In dimension $1+1$, directed polymers (as in \eqref{strong23}) are expected to lie in the strong disorder regime for every fixed inverse temperature $\beta>0$, 
exhibiting similar qualitative behavior as the zero-temperature model (see the survey \cite[Section 3]{Z24} and the references therein).
Like in the zero-temperature case, 
rigorous understanding of fixed-temperature polymer models remains limited outside of special integrable examples, such as the log-gamma polymer \cite{Sep12,COSZ14,BCR13}. 
In contrast, in the \emph{intermediate disorder regime}, where $\beta$ is scaled like $n^{-1/4}$, a universal continuum theory has emerged through the construction of the continuum directed random polymer (CDRP) \cite{AKQ14,AKQ14JSP}.
This is the model that we study in the present article. 
The CDRP is also closely connected to the KPZ line ensemble, and this connection will play a crucial role in our arguments.\\

In the next two subsections, we present the definitions required to state our main results which are then presented in \cref{ss:results}.

\subsection{The Continuum Directed Random Polymer (CDRP)}\label{ss:CDRP}

Although the CDRP was first constructed in \cite{AKQ14JSP}, for our purposes we rely on the more recent comprehensive treatment in \cite{AJRS22}. 
Let
\begin{align*}
    \Rup \coloneqq  \left\{(s,x;t,y)\in \mathbb{R}^{4} : s<t\right\},
\end{align*}
and let $C(\Rup)$ be the space of continuous real-valued functions on $\Rup$ equipped with the topology of uniform convergence on compact sets. 
Let $\xi$ be a (space-time) white noise on $\R^2$, that is, the centered Gaussian process on $L^2(\R^2)$ with covariances $\E[\xi(f)\xi(g)]=\langle f,g\rangle_{L^2(\R^2)}$.
We denote the law of $\xi$ by $\P$, which is a probability measure on the product space $(\Xi, \Fprod) \coloneqq  (\R, \cB(\R))^{\otimes L^2(\R^2)}$, where $\cB(\R)$ is the Borel $\sigma$-algebra of $\R$.
A detailed discussion of white noise appears in Sections \ref{s:meas-short} and \ref{s:meas-long}.

Fix $(s,x;t,y)\in\Rup$ and $\beta\ge0$, and a realization of the white noise $\xi$.
At a heuristic level, the CDRP (from $(s,x)$ to $(t,y)$, at inverse temperature $\beta$)
is the  probability measure on $C([s,t])$ given by the following formal tilt of Brownian bridge:
\begin{align}\label{wienertilt}
   \frac{{d}\CDRP^\xi_\beta(s,x;t,y)}{{d}{\mathrm{BB}_{(s,x),(t,y)}}}(X)
    \coloneqq 
    ``\;
        \frac{p(t-s, y-x)}{\cZ_\beta^\xi(s,x;t,y)}
        \exp\left(
            \beta \int_s^t \xi(r, X(r))\, dr
        \right),"
\end{align}
where the normalizing constant $\cZ_\beta^\xi(s,x;t,y)$ is the partition function, and $p(r,z)\coloneqq \frac{1}{\sqrt{2\pi r}}e^{-z^2/2r}$ is the heat kernel.
This expression can be made rigorous by interpreting the partition function and finite-dimensional distributions of the CDRP in terms of a formal Feynman--Kac representation for the \emph{stochastic heat equation (SHE) with multiplicative noise}:
\begin{align}\label{e:she-intro}
    \partial_{t}\mathcal{Z}^{\xi}_{\beta}(s,x;t,y) = \frac{1}{2}\partial_{y}^{2}\mathcal{Z}^{\xi}_{\beta}(s,x;t,y) + \beta\mathcal{Z}^{\xi}_{\beta}(s,x;t,y)\xi(t,y),
    \qquad\quad t>s,\; y\in\R,
    \tag{SHE}
\end{align}
with initial data $\lim_{t\downarrow s}\mathcal{Z}^{\xi}_{\beta}(s,x;t,y) =\delta(x-y)$.
The aforementioned \cite{AJRS22} constructs a coupling across $(s,x),\beta$ of the solutions to \eqref{e:she-intro} and the corresponding CDRP measures, and we turn now to formulating this precisely.
The following definitions are somewhat involved, so the reader may prefer to simply glance at the finite-dimensional distribution formula \eqref{e:CDRP-fdd1} and then skip ahead to \cref{ss:KPZ}, keeping in mind the above heuristic picture.

The following definition paraphrases \cite[Theorem 2.2]{AJRS22} (a more detailed version is recorded later in \cref{t:AJRAS-Z-existence}).
\begin{defn}[CDRP partition function process]\label{def:intro-SHE}
    By \cite[Theorem 2.2]{AJRS22}, there exists a measurable map 
    \begin{align*}
        \begin{array}{ccc}
            \cZ :(\Xi,\Fprod,\P) &\longrightarrow &(C(\Rup\times \R_{\ge 0}),\; \cB(C(\Rup\times \R_{\ge 0})))\\
            \;\;\xi &\longmapsto & \cZ^{\xi}_{\smallbullet}(\smallbullet,\smallbullet;\smallbullet,\smallbullet)
        \end{array}
    \end{align*}
    satisfying the following properties.
    \begin{itemize}
        \item \textup{(Coupling of SHEs).} 
        For any fixed $(s,x)\in\R^2$ and $\beta\ge 0$, the process $(t,y)\mapsto \cZ_{\beta}^{\xi}(s,x;t,y)$ coincides a.s. with the unique continuous adapted mild solution of \eqref{e:she-intro}.
        \item\textup{(Strict positivity).} Almost surely, $\cZ^\xi_{\beta}(s,x;t,y)>0$ for all $(s,x;t,y)\in\Rup$ and $\beta\ge 0$.
        \item\textup{(Chapman--Kolmogorov, \cite[Lemma 3.12]{AJRS22}).} 
        Almost surely, for all $(s,x;t,y)\in\Rup$, all $r\in (s,t)$, and all $\beta\ge 0$, we have
        \begin{align}\label{e:ck}
            \mathcal{Z}^{\xi}_{\beta}(s,x;t,y) & = \int_{\mathbb{R}}
            \mathcal{Z}^{\xi}_{\beta}(s,x;r,z)\mathcal{Z}^{\xi}_{\beta}(r,z;t,y)dz.
        \end{align}
    \end{itemize}
\end{defn}

Given \cref{def:intro-SHE}, one can build the CDRP by treating the partition function process as its heat kernel. 
Namely, given a realization of the white noise $\xi$,
and given $(s,x;t,y)\in\Rup$ and $\beta \ge 0$, the multi-point densities of $X\sim \CDRP^{\xi}_{\beta}(s,x;t,y)$ are defined as follows: for $k\in\N$ and $s<s_1<\cdots<s_k<t$,
\begin{align}\label{e:CDRP-fdd1}
    \P^\xi_{\beta,(s,x),(t,y)}(X(s_1)\in dx_1,\dots, X(s_k)\in dx_k)
    &\coloneqq  \frac{\prod_{i=0}^k \cZ_\beta^\xi(s_i,x_i;s_{i+1},x_{i+1})}{\cZ_\beta^\xi(s,x;t,y)}
    dx_{1}\dots dx_{k},
\end{align}
where $(s_0,x_0)\coloneqq (s,x)$ and $(s_{k+1},x_{k+1})\coloneqq (t,y)$.
The RHS is well-defined since $\cZ_{\beta}^{\xi}>0$ by \cref{def:intro-SHE}.
Further, by the Chapman--Kolmogorov equation \eqref{e:ck}, the formula \eqref{e:CDRP-fdd1} defines a consistent family of finite-dimensional distributions.
This along with the Kolmogorov extension theorem is the starting point for the following result.

\begin{thm}[{CDRP measure, \cite[Theorem 2.14]{AJRS22}}]\label{t:AJRAS-polymer-existence}
    The following holds for $\P$-a.e. realization of the white noise $\xi$.
    For every $(s,x;t,y)\in\Rup$ and every $\beta\ge0$, there exists a unique probability measure $\P^\xi_{\beta,(s,x),(t,y)}$ on $C([s,t])$ with finite-dimensional marginals given by \eqref{e:CDRP-fdd1}.

    Given a realization of the white noise $\xi$, we refer to $\P^\xi_{\beta,(s,x),(t,y)}$ as the \emph{(quenched) polymer measure}.
    We will often write $X\sim \P^\xi_{\beta,(s,x),(t,y)}$ or $X\sim \CDRP^\xi_\beta(s,x;t,y)$ to mean that $X$ is sampled from $\P^\xi_{\beta,(s,x),(t,y)}$ (conditional on $\xi$).
\end{thm}

The free energy profile of the CDRP is given by the KPZ equation. 
This is recorded next, along with some recent developments about its zero-temperature limit, all of which will serve as key inputs for us.

\subsection{The KPZ Equation and the Directed Landscape}\label{ss:KPZ}

We recall the KPZ equation from \eqref{e:kpz-1}, now rewritten to match the notation of \eqref{e:she-intro}.
For $(s,x)\in\R^2$ and $\beta \ge 0$, the KPZ equation is the stochastic PDE
\begin{align}\label{e:KPZ-eq}
    \partial_{t} \mathcal{H}^\xi_\beta(s,x;t,y) & = \frac{1}{2}\partial_{y}^{2}\mathcal{H}^\xi_\beta(s,x;t,y) + \frac{1}{2}\left(\partial_{y}\mathcal{H}^\xi_\beta(s,x;t,y)\right)^{2} + \beta\xi(t,y),
    \qquad
    t>s,\, y\in\R,
    \tag{KPZ}
\end{align}
where $\xi$ is space-time white noise on $\R^2$.
The above equation is ill-posed, but there is a standard notion of solution given by formally applying the Cole--Hopf transformation to \eqref{e:she-intro}.
That is, the Cole--Hopf solution to \eqref{e:KPZ-eq} is \emph{defined} as
\begin{align}\label{e:KPZ-cole-hopf}
    \cH^{\xi}_\beta(s,x;t,y) \coloneqq  \log \cZ^\xi_\beta(s,x;t,y)
    \qquad\qquad\text{for } (s,x;t,y)\in\Rup,\; \beta>0,
\end{align}
with $\cZ^\xi_\beta$ as in \cref{def:intro-SHE}
(the logarithm is well-defined since $\cZ^\xi > 0$ a.s. by definition).
The above is commonly known as the ``narrow wedge'' solution to \eqref{e:KPZ-eq}, due to the $\delta(x-y)$ initial data for \eqref{e:she-intro}.

We next define a scaled version of \eqref{e:KPZ-cole-hopf} that, as we will explain momentarily, converges in the zero-temperature limit $\beta\to\infty$ to the directed landscape.

\begin{defn}[CDRP free energy profile]
    \label{d:narrow-wedge}
    The \emph{CDRP free energy profile} is the random continuous function on $\Rup\times \R_{>0}$ given by
    \begin{align}\label{sheet12}
        \mathfrak{H}^\xi_\beta(s,x;t,y) \coloneqq  
        \frac{2^{1/3}}{\beta^{4/3}}
        \left(
            \cH_\beta^\xi(s, x2^{1/3}\beta^{2/3}; t, y2^{1/3}\beta^{2/3}) 
            - 2\log\beta 
            + \frac{(t-s)\beta^{4}}{24}
        \right),
    \end{align}
    where $\cH_\beta^\xi(s,x;t,y) \coloneqq  \log \cZ_\beta^\xi(s,x;t,y)$ as in \eqref{e:KPZ-cole-hopf}.

    In other works such as \cite{W23}, $\mf{H}^\xi_\beta$ is referred to as the \emph{scaled narrow wedge solution to the KPZ equation}. 
    Since our focus is on the CDRP, we use the above terminology instead.
\end{defn}
 
In a major development, Dauvergne, Ortmann, and Vir\'ag \cite{dov} constructed the \emph{directed landscape} as the full space-time scaling limit of Brownian last passage percolation. 
The directed landscape is expected to be the universal $1:2:3$-scaling limit of directed random geometry in the KPZ universality class. 
Subsequent work of Dauvergne and Vir\'ag \cite{dv2} established convergence to the same object for several exactly solvable LPP models, including geometric, exponential, and Poissonian LPP. 
We now recall the definition of the directed landscape from \cite[Definition 10.1]{dov}.
\begin{defn}[Directed landscape]\label{d:dl}
The \emph{directed landscape} $\cL$ is the unique $C(\Rup)$-valued random variable satisfying the following properties. 
\begin{enumerate}
    \item 
    For any finite collection of disjoint time intervals $\{(s_i,t_i)\}_{i=1,\dots,k}$, the random functions $\{\cL(s_i,\smallbullet; t_i, \smallbullet)\}_{i=1,\dots, k}$ are independent.
    \item 
    Almost surely, for all $(s,x;t,y)\in\Rup$ and all $r\in (s,t)$, we have
    \begin{align*}
        \cL(s,x;t,y) = \sup_{z\in \R}\left(\cL(s,x;r,z) + \cL(r,z;t,y)\right).
    \end{align*}
    \item 
    For all $s<t$, the random function $(x,y)\mapsto \mathcal{L}(s,x;t,y)$ has the same law as $(x,y)\mapsto (t-s)^{1/3}\mathcal{L}(0,(t-s)^{-2/3}x;1,(t-s)^{-2/3}y)$.
    \item $\mathcal{L}(0,\smallbullet;1,\smallbullet)$ is distributed as the Airy sheet (see \cite[Definition 1.2]{dov}).
\end{enumerate}
\end{defn}

More recently, and more pertinent for our application, the positive-temperature model of the CDRP free energy profile \eqref{sheet12} has also been shown to converge to the directed landscape, with an earlier one-point convergence shown in \cite{ACQ11}. 

\begin{thm}[{\cite[Theorem 1.6]{W23}}]\label{t:KPZ-to-landscape}
    As $\beta\to\infty$, the CDRP free energy profile $\mf{H}^\xi_{\beta}$ converges in distribution to the directed landscape $\cL$.
    Here convergence is with respect to the uniform-on-compact topology on $C(\Rup)$.
\end{thm}

A crucial ingredient in the proof of \cref{t:KPZ-to-landscape}
is the convergence of the KPZ line ensemble to the Airy line ensemble. 
While we will not need the latter objects for our purposes, for the uninitiated reader we simply mention that they are both infinite sequences of random curves with a natural Gibbsian structure.
For us, the only relevant observables will be the top curves, which are the functions $\mf{H}^{\xi}_{\beta}(0,0;1,\smallbullet)$ and $\cL(0,0;1,\smallbullet)$.
The convergence of the KPZ line ensemble to the Airy line ensemble can be delivered by combining the recent works \cite{V20,DM,W23LE,AH23}, as noted in \cite[Proposition 5.1]{W23}. 
A self-contained argument can also be found in \cite[Corollary 25.1]{AH23}.

\subsection{Main result and a conjecture}\label{ss:results}
Our main result on temperature chaos concerns the coupling structure of the processes $\mf{H}^{\xi}_{\beta}(\smallbullet, \smallbullet; \smallbullet,\smallbullet)$ across different inverse temperatures $\beta$, and sharing the same white noise $\xi$.
As mentioned above \cref{ss:CDRP}, 
$\mf{H}^\xi_\beta$ is the universal scaling limit of lattice polymer free energies in the intermediate disorder regime, wherein inverse temperature $\beta$ is scaled by $n^{-1/4}$ so that the Gibbs weight of a directed lattice path $\pi$ from $(0,0)$ to $(n,n)$ is proportional to
\begin{align}\label{e:intermediate-disorder}
    \exp\left(
        \beta n^{-1/4} \sum_{v\in\pi} X_v
    \right),
\end{align}
where $\beta>0$ is constant.
To motivate temperature-chaotic properties in the intermediate disorder setting, we first discuss counterpart notions for lattice polymers in the strong disorder regime \eqref{strong23}.

For the length-$n$ lattice polymer in the strong disorder regime (fixed inverse temperature $\beta>0$), temperature chaos predicts the existence of a \emph{critical temperature perturbation scale} $\Delta\beta_n=n^{-\alpha}$ for some $\alpha>0$ such that two polymers with different inverse temperatures are decorrelated \emph{if and only if} their inverse temperatures differ by an amount $\gg \Delta \beta_n$.
In fact, we now formulate a conjecture for the critical exponent $\alpha$ marking the onset of temperature chaos with a justification presented later in the article. The  prediction in fact can be traced back to the important work \cite{FH91} in the physics literature.   
Let us recall more precisely the lattice polymer model defined in \eqref{strong23}. 
As before, $\mathbf{X}=\{X_v\}_{v\in \Z^2}$ is a field of i.i.d. disorder, say with a finite exponential moment.
Given a realization of $\mathbf{X}$, for any inverse temperature $\beta>0$ let $\P_{\beta}$ be the corresponding polymer measure as in \eqref{strong23}, let $Z_\beta$ be its partition function, and let $F_{\beta}=\log Z_\beta$ be its free energy. 
For two inverse temperatures $\beta_1,\beta_2>0$, let $\pi_1\sim \P_{\beta_1}$ and $\pi_2\sim \P_{\beta_2}$ be polymers sampled conditionally independently given the environment $\mathbf{X}$.

\begin{conj}[Temperature chaos exponent]\label{c:main}
    Let $\beta>0$. 
    Setting $\beta_1=\beta$ and letting $\beta_2=\beta_2(n) > \beta_1$ depend on $n\ge 1$,
    as $n\to \infty$ the following hold.

    \begin{itemize}
        \item \textbf{\emph{Energetic decorrelation:}}
        \begin{align*}
            \Corr(F_{\beta_1},F_{\beta_2})
            =
            \begin{cases}
                1-o(1), & \beta_2-\beta_1\ll n^{-1/6},\\[2mm]
                o(1),   & \beta_2-\beta_1\gg n^{-1/6}.
            \end{cases}
        \end{align*}
        See Figure \ref{fig:simulation-free-energy}.
        \item \textbf{\emph{Polymer chaos:}}
        For the normalized polymer overlap,
        \begin{align}\label{chaos2}
            \E\left[
                \frac{|\pi_1\cap \pi_2|}{n}
            \right]
            =
            \begin{cases}
                \Theta(1),
                & \beta_2-\beta_1\ll n^{-1/6},\\[2mm]
                o(1),
                & \beta_2-\beta_1\gg n^{-1/6},
            \end{cases}
        \end{align}
        where the expectation is over the disorder and the conditionally independent polymer measures.
        See Figure \ref{fig:simulation-free-energy}.
    \end{itemize}

    The corresponding conjecture for the CDRP is obtained by an intermediate disorder change of variables, replacing $\beta$ by $\beta n^{-1/4}$, so that in the $n\to\infty$ limit one obtains the CDRP at inverse temperature $\beta$.
    After this change of variables, to recover behavior analogous to the fixed-temperature lattice polymer, one must send the CDRP inverse temperature $\beta$ to infinity.
    Under this change of variables, the conjectural lattice temperature chaos scale $n^{-1/6}$ becomes the CDRP scale $\beta^{1/3}$.
    \footnote{
        By a space-time-temperature scaling, the CDRP of height $1$ and inverse temperature $\beta_2$ is comparable to a lattice model of height $n=\beta_1^{4}$ and inverse temperature $\frac{\beta_2}{\beta_1}$.
        Thus, under this scaling, the additive change in inverse temperature $1\to 1+ n^{-1/6}$ corresponds to the multiplicative change $\beta_1 \to \beta_1 (1+ \beta_1^{-2/3}) = \beta_1 + \beta_1^{1/3}$.
        This will be expanded on further in \cref{s:heuristic}.
        }
    Thus, as $\beta_1,\beta_2\to\infty$ with $\beta_2 > \beta_1$, the CDRP free energies satisfy
    \begin{align}\label{energy1} 
        (\mathfrak{H}^\xi_{\beta_1},\;\mathfrak{H}^\xi_{\beta_2})
        &\dto
        \begin{cases} 
            (\cL_1, \cL_1), &  \beta_2-\beta_1\ll \beta_1^{1/3},
                \\[2mm]
            (\cL_1, \cL_2),  &  \beta_2-\beta_1\gg \beta_1^{1/3},
        \end{cases}  
    \end{align}
    where $\cL_1$ and $\cL_2$ are independent directed landscapes, 
    and convergence is in the uniform-on-compact topology on $C(\Rup)$ in each coordinate.
    The marginal convergence is known by \cref{t:KPZ-to-landscape}.
\end{conj}

Note that in \eqref{chaos2} the overlap is conjectured to be $\Theta(1)$, not $1-o(1)$.
This is indeed the best one can hope for, since even when $\beta_1=\beta_2$, two independent polymers have a constant chance of not overlapping at a given height, conditional on their entire remaining trajectories.
This is a feature of the positive-temperature nature of the model.
We will explain the source of the exponents $-1/6$ and $1/3$ and why we expect \cref{c:main} to be true in \cref{s:heuristic} using a heuristic argument connecting temperature chaos and disorder chaos. 
As already indicated, the exponent in fact can be traced back to \cite{FH91} where a different heuristic argument based on relating the fluctuation of the path entropy and the free energy fluctuations was presented as a justification (see also \cite{SY02,dSB04}).
\\

We are now in a position to state our main result, which is the first rigorous manifestation of \emph{asymptotic decoupling} of the free energy profiles across different temperatures coupled through the same noise.

\begin{thm}[Polymer free energy decoupling]\label{t:main} 
    Let $(\mf{H}^{\xi}_{\beta_1}, \mf{H}^{\xi}_{\beta_2})$ be as in \cref{d:narrow-wedge}, coupled through a shared space-time white noise $\xi$.
    Then as $\beta_1,\beta_2\to\infty$ such that $\beta_1 = o(\beta_2)$,   
    \begin{align*}
        (\mf{H}^{\xi}_{\beta_1},\, \mf{H}^{\xi}_{\beta_2})
        \dto 
        \left(\mathcal{L}_{1},\mathcal{L}_{2}\right),
    \end{align*}
    where $\mathcal{L}_{1},\mathcal{L}_{2}$ are two independent copies of the directed landscape.
    Here, convergence is as $C(\Rup)\times C(\Rup)$-valued random variables. 
\end{thm}

The assumption $\beta_1 = o(\beta_2)$ is far from sharp, as is evident from \cref{c:main} which remains elusive.
The main difficulty in proving decoupling results of this kind is that most available estimates for the CDRP concern the moments of the partition function $\cZ_\beta^\xi$, which do not capture the typical behavior of the free energy.
To illustrate this, observe that by \cref{t:KPZ-to-landscape}, $\log \cZ_\beta^\xi(0,0;1,0)= -c_0\beta^4+c_1\beta^{4/3} Y_\beta + O(\log\beta),$ where $c_0,c_1>0$ are constants and $Y_\beta$ converges to a Tracy--Widom GUE random variable $Y$ as $\beta\to\infty$.
In particular, the upper tail of the limiting free energy fluctuation satisfies $\P(Y > m)=\exp(-\Theta(m^{3/2}))$ as $m\to\infty$.
The partition function is thus heavy-tailed, and hence its moments are dominated by rare events and reveal little about the free energy's typical behavior.
We revisit this issue in \cref{s:heuristic}.

As for the stability regime in \eqref{energy1}, it is not too difficult to prove that if $\beta_2$ is chosen polynomially close to $\beta_1$, that is  $\beta_2-\beta_1 \le \beta_1^{-\alpha}$ for some large enough $\alpha>0$, then the corresponding polymer free energies are indeed strongly correlated, with $(\mf{H}^{\xi}_{\beta_1},\, \mf{H}^{\xi}_{\beta_2}) \dto (\cL_1,\cL_1)$ as $\beta_2 > \beta_1\to \infty$.
This is a straightforward consequence of the H\"older continuity properties of  $\cZ_\smallbullet^\xi(\smallbullet, \smallbullet; \smallbullet, \smallbullet)$ established in \cite[Proposition 3.8]{AJRS22}, and a  rigorous analogue of the heuristic coarse-graining construction underlying \cref{c:main}.
We elaborate on this in \cref{ss:polynomial-stability}.
Establishing energetic stability up to the optimal threshold of $\alpha=1/3$ remains an important open problem.

In contrast to the above discussion, our next result shows that regardless of whether the polymer free energies are strongly correlated or decoupled, the polymer measures are mutually singular as soon as $\beta_1\ne \beta_2$.
This indicates that the notion of singularity of random measures in a high-dimensional space such as the path space is rather fragile and does not reveal any information about the correlation structure. 

\begin{thm}[Mutual singularity of polymer measures]\label{t:pathsing} 
    For any $(s,x;t,y)\in\Rup$ and any $\beta_2>\beta_1\ge 0$, almost surely the polymer measures $\mathbb{P}_{\beta_{1},(s,x),(t,y)}^{\xi}$ and $\mathbb{P}_{\beta_{2},(s,x),(t,y)}^{\xi}$ are mutually singular.
\end{thm}
The case $\beta_1=0$ recovers \cite[Theorem 4.5]{AKQ14JSP}, which proved singularity of the CDRP measure with respect to the Wiener measure (see also \cite{COTX26}).
Our proof takes a different route than theirs and treats all pairs of distinct temperatures simultaneously by showing that each CDRP measure is supported on a different \emph{level set} of the Hamiltonian $X\mapsto \int_s^t \xi(r,X(r))\,dr$ in \eqref{wienertilt}.
The latter is only an informal expression, and to make this strategy rigorous we rely on the Gaussian multiplicative chaos framework for CDRP recently developed in \cite{QRV25}.
This result is strongly reminiscent of the well-known fact that almost surely, the $\gamma$-Liouville quantum gravity measure is supported on $\gamma$-thick points of the Gaussian free field and hence the measures are singular across different values of $\gamma$ (see for instance \cite[Proposition 3.4]{DS11}).
Similar results for a related directed polymer model appeared previously in \cite{BM20,BLM24} (see also \cite{BM19}).

Finally, it is worth commenting that the approach of viewing polymer measures as GMC measures on path space allows one to go beyond white noise to study polymers driven by correlated space-time disorder.
In the forthcoming work \cite{GG}, this will be explored in detail in the log-correlated setting with connections to Liouville quantum gravity.

\subsection{The directed landscape is a two-dimensional black noise} 
As a byproduct of our methods for proving \cref{t:main}, we establish the following zero-temperature result.
Recall the notion of black noise introduced by Tsirelson and Vershik \cite{TV} alluded to above (the formal definition appears in \cref{def:black-noise}, see also \cref{r:informal-def-black-noise}).

\begin{thm}\label{t:black-noise}
    The directed landscape is a two-dimensional black noise.    
\end{thm}

As already mentioned, this was previously conjectured by B\'alint Vir\'ag and recorded as \cite[Conjecture 1.17]{HP24}. The above result makes the directed landscape the third known instance of a two-dimensional black noise (following critical planar percolation \cite{SS11} and the Brownian web \cite{EF16}). 
In \cite{HP24} it was shown that the directed landscape is a one-dimensional black noise (in the time direction) and as a consequence in \cite[Theorem 1.11]{HP24} that the directed landscape must be independent of any white noise on the same filtered probability space. 
On the other hand, \cref{t:black-noise} shows that the directed landscape is nevertheless driven by a two-dimensional noise---albeit a non-classical one.\\

Finally, though not the topic of this paper, it is worth pointing out that 
a result analogous to \cite[Theorem 1.10]{HP24} was proven in \cite{GT25} 
for the \emph{critical two-dimensional stochastic heat flow} (SHF), 
the universal scaling limit of critical $2+1$-dimensional directed polymer partition functions constructed in \cite{CSZ23}. Namely, it was shown that the 2D SHF is a one-dimensional black noise in the temporal direction. In a similar vein, \cite{CD25} recently established a notion of \emph{enhanced noise sensitivity} for broad classes of functions of independent random variables, including such critical $2+1$-dimensional polymer partition functions.
Though we do not pursue this in full detail, we sketch briefly in \cref{ss:SHF-black-noise} how the technique for the proof of \cref{t:black-noise} may be used to show that the SHF is a 3D black noise. Curiously, there are no other known examples of black noises in dimension three or higher, and finding one has appeared as \cite[Open Problem 3]{EF16}.
\\

Having stated the results, we now proceed with the remainder of the article.
Since it is rather long, we have attempted to guide the reader by providing frequent roadmaps both at the beginning of each section and at key points within them. We have also included explanatory remarks throughout the paper to motivate how the results developed in each section fit into the overall goal.
We first record a brief outline of the rest of the paper:
\subsection{Organization of the paper}
\begin{itemize}
    \item Before any formal development, we start with a brief discussion of some of the key proof ideas in the upcoming \cref{ss:idea}.
    \item Notational conventions are laid out in \cref{s:notation}.
    \item The measure-theoretic framework required for the proof of \cref{t:main} is outlined in  Section \ref{s:meas-short} and developed in Section \ref{s:meas-long}.
    \item Technical inputs concerning the directed landscape, the CDRP free energy profile, and polymer geometry appear in \cref{s:LE}.
    \item The free energy decoupling result \cref{t:main} is proved in Sections \ref{s:mainproof} and \ref{s:es}.
    \item The polymer measure singularity result \cref{t:pathsing} is proved in \cref{s:mutual-singularity}.
    \item \cref{t:black-noise} stating that the directed landscape is a 2D black noise is proved in \cref{s:black-noise}.
    Further, a brief discussion towards proving that the critical 2D stochastic heat flow is a  3D black noise is presented in \cref{ss:SHF-black-noise}.
    \item Heuristic discussions of temperature chaos including the argument justifying \cref{c:main} are presented in \cref{s:heuristic}, along with a rigorous argument proving energetic stability up to a polynomial temperature gap.
    \item Basic properties of conditional variance are recorded in Appendix \ref{s:conditional-variance}.
\end{itemize}

\subsection{Proof ideas}\label{ss:idea}

Let us first mention our coordinate conventions.
Note that in Definitions \ref{d:narrow-wedge} and \ref{d:dl}, in a space-time pair $(s,x)\in\R^2$, the first coordinate denotes time and the second coordinate denotes space.
This is done to be consistent with \cite{AJRS22}, and we maintain this convention throughout the paper.
However, most of our arguments involve drawing parallels to LPP, so to match the LPP literature, in our figures we use the opposite coordinate convention: the temporal axis is vertical and the spatial axis is horizontal.

We start with the proof of \cref{t:main}.
However, the same principle guiding the argument applied in the zero-temperature setting will yield \cref{t:black-noise}.

The first key observation is that when $\beta_{1}\ll \beta_2$,
the two polymer measures $\CDRP^\xi_{\beta_1}(0,0;1,0)$
and $\CDRP^\xi_{\beta_2}(0,0;1,0)$ fluctuate on different scales, making the former significantly more \emph{localized} than the latter.
It is well-known that in models of LPP, and in particular in the directed landscape, geodesics traveling from $(0,0)$ to $(t,0)$ fluctuate transversally at scale $t^{2/3}$. 
The same is true for $\CDRP^\xi_1(0,0;t,0)$ for large $t$ by \cite[Corollary 3.5]{DZ24}, which by the CDRP height-temperature scaling relation (\cref{p:CDRP-scaling}\ref{property-CDRP-scaling}) implies that $\CDRP^\xi_\beta(0,0;1,0)$ fluctuates at scale $\beta^{2/3}$ for large $\beta$ (this is recorded in \cref{p:polymer-TF-DZ}).
Thus, when $\beta_1\ll \beta_2$, with high probability $\xi$ is such that 
with high probability under $\CDRP^\xi_{\beta_1}(0,0;1,0)$, the polymer is localized to a strip $S$ of width 
$\beta_1^{2/3}=\e \beta_{2}^{2/3}$ for some small $\e = o(1)$. 
As a consequence, the polymer free energy $\mf{H}^{\xi}_{\beta_1}(0,0;1,0)$ is ``almost'' a measurable function of the restriction of $\xi$ to $S$.
To prove the required decoupling we will employ a resampling strategy.
We will resample the white noise $\xi$ on $S$ to obtain a new white noise $\eta$ that agrees with $\xi$ outside $S$ and is independent of $\xi$ on $S$.
The above ``almost-measurability'' reasoning then implies that $\mf{H}^{\xi}_{\beta_1}$ is essentially independent of $\eta$. 
The next key observation is that spatially thin strips are not \emph{pivotal} or \emph{influential}.
That is, this resampling alters the free energy $\mf{H}^{\xi}_{\beta_2}$ negligibly:
\begin{align}\label{resampling-effect}
    \mf{H}^{\xi}_{\beta_2}\approx \mf{H}^{\eta}_{\beta_2},
\end{align}
where $\approx$ will not be defined but a quantitative statement will be presented shortly. As alluded to in the introduction, the topic of this paper has natural connections to the theory of noise sensitivity of Boolean functions \cite{BKS99}.
To help draw this parallel, we will often borrow terms like ``influence'' and ``pivotality'' from the noise sensitivity literature.
In the Boolean setting, the \emph{influence} of a Bernoulli bit on a function $f$ is the probability of changing the value of $f$ upon flipping the bit's value, i.e. the probability that the bit is \emph{pivotal} for $f$.
See for instance \cite{BKS99}.

\begin{figure}[htb]
\centering
    \includegraphics[height = 0.37\textheight]{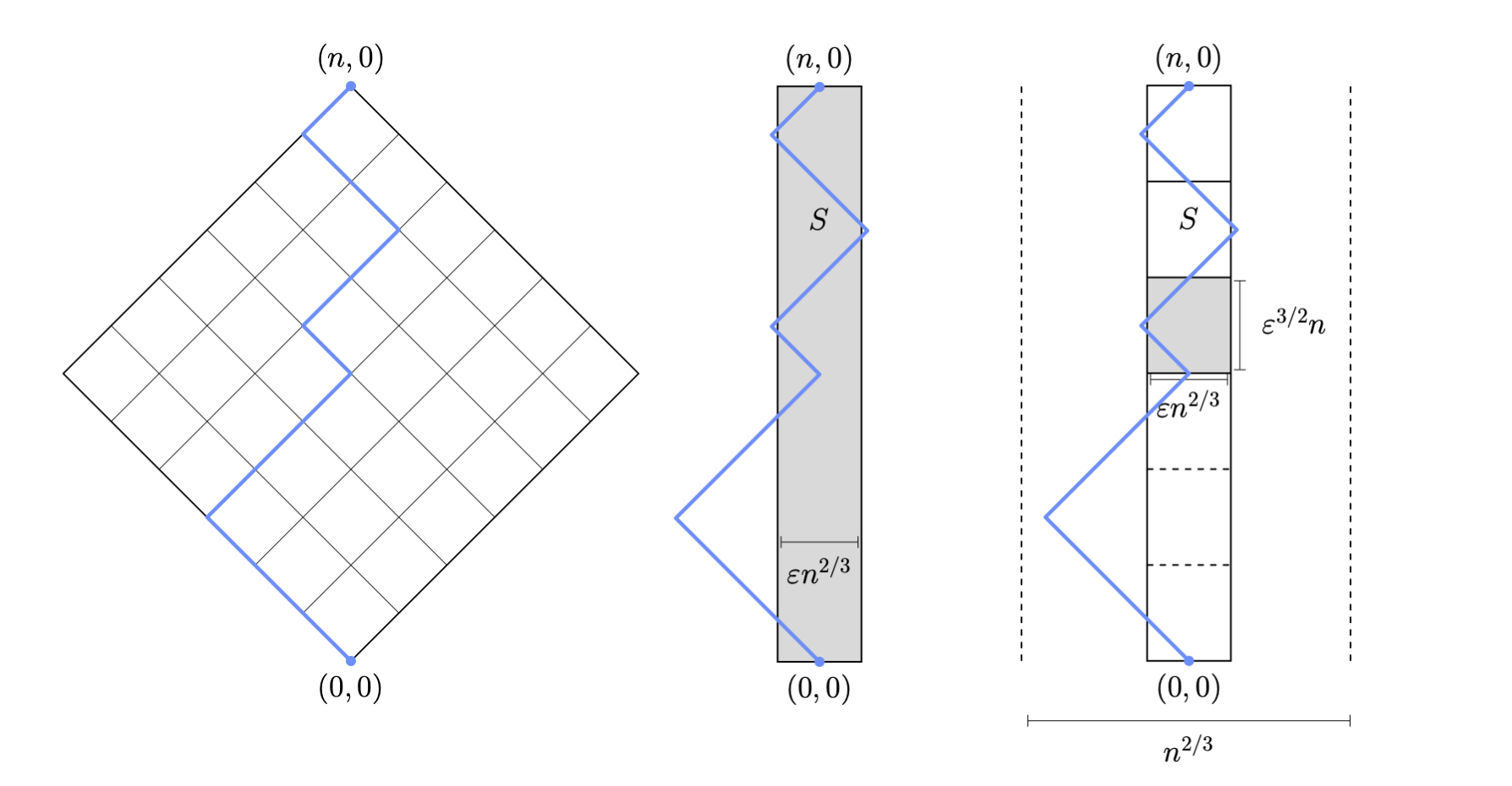}
\caption{
    \textbf{Left:} Exponential LPP on $\Z^2$.
    The vertices of $\Z^2$ carry i.i.d. $\mathrm{Exponential}(1)$ weights.
    The \emph{last passage time} is the maximum sum of weights along any directed (north-east/north-west) lattice path from $(0,0)$ to $(n,0)$.
    The blue path depicts the \emph{geodesic} (maximizing path).
    \\
    \textbf{Middle:} The geodesic from the left panel is redrawn on top of a thin strip $S$ (shaded) of height $n$ and width $\e n^{2/3}$.
    The key bound \eqref{efron-stein} asserts that the $L^2$ influence of the noise in $S$ on the last passage time is $O(\e^{1/2} n^{2/3})$.\\
    \textbf{Right:} We tile $S$ with boxes of size $\e^{3/2}n\times \e n^{2/3}$.
    By the Efron--Stein inequality, the influence of $S$ is at most the sum of the influences of the boxes.
    Each box has influence $O(\e^2 n^{2/3})$, which summed over $\e^{-3/2}$-many boxes yields the $O(\e^{1/2}n^{2/3})$ bound claimed in  \eqref{efron-stein}.
    The per-box influence bound follows from KPZ fluctuation considerations. 
    For a given box $B$ (shaded), resampling it only changes the last passage time on the rare event that the geodesic lands in $B$ before or after resampling.
    This occurs with probability $\Theta(\e)$, since the geodesic explores a strip of width $\Theta(n^{2/3})$.
    On this rare event, resampling $B$ changes the last passage time by $O(\e^{1/2}n^{1/3})$, and hence the $L^2$-influence is at most $\e\cdot (\e^{1/2}n^{1/3})^2 = \e^2n^{2/3}$. 
    }\label{fig:iop-LPP}
\end{figure}

To explain why \eqref{resampling-effect} should be true, it is slightly clearer to provide a heuristic sketch in the zero-temperature setting.
The actual implementation of this argument for the CDRP uses the same idea, relying on the weak convergence statement in \cref{t:KPZ-to-landscape} to import the relevant zero-temperature estimates. 
As a further simplification, we explain the mechanism in a discrete zero-temperature model.
Consider Exponential LPP on $\Z^2$: the vertices carry i.i.d. $\mathrm{Exponential}(1)$ weights, and the last passage time between two vertices is the maximum total weight over all directed paths connecting them, as depicted in Figure \ref{fig:iop-LPP} (left).
Let $S$ be the strip of width $\e n^{2/3}$ depicted in Figure \ref{fig:iop-LPP} (middle).
Let $\omega=\{\omega(v)\}_{v\in\Z^2}$ and $\omega^S=\{\omega^S(v)\}_{v\in\Z^2}$ be a pair of coupled exponential environments that agree outside $S$, and are independently resampled inside $S$: that is, $\omega(v)=\omega^S(v)$ for $v\not\in S$, while inside $S$ the two collections are independent copies of the exponential disorder.
Let $L_n$ and $L_n^S$ denote the corresponding last passage times from $(0,0)$ to $(0,n)$.
The analogous statement to \eqref{resampling-effect} in this setting is that resampling the thin strip $S$ changes the last passage time by less than its natural fluctuation scale of $\Theta(n^{1/3})$.
The key quantitative assertion is the following influence bound:
\begin{align}\label{efron-stein}
    \mathbb{E}\bigl[(L_n-L_n^S)^2\bigr]\lesssim \e^{1/2}\, n^{2/3}.
\end{align}

The proof of \eqref{efron-stein} employs the Efron--Stein inequality, resampling $S$ not all in one go, but rather subdividing $S$ into boxes of size $\e^{3/2}n \times \e n^{2/3}$ 
as in  Figure \ref{fig:iop-LPP} (right) and resampling them one by one.
This yields the bound
\begin{align*}
    \mathbb{E}\bigl[(L_n-L_n^S)^2\bigr] \leq \sum_{i=1}^{\e^{-3/2}}\E\bigl[(L_{n}-L_{n}^{B_i})^{2}\bigr]
\end{align*}
where $L_{n}$ and $L_{n}^{B_{i}}$ are the last passage times in environments $\omega$ and $\omega^{B_i}$, with $\omega^{B_i}$ obtained from $\omega$ by resampling the noise in the box $B_i$.
A key estimate then states that each summand in the right hand side is approximately $\e^{2} n^{2/3}$, which finishes the proof.
This estimate comes from the observation that resampling $B_i$ typically does not change $L_n$ at all, i.e. $L_n = L_n^{B_i}$ with high probability.
Indeed, a change occurs only if the geodesic intersects $B_i$ in either the environment $\omega$ or the resampled environment $\omega^{B_i}$.
The probability of this is approximately $\e$, since the geodesic is roughly uniformly distributed on an interval of width $n^{2/3}$ by KPZ considerations.
On this rare event, the change in last passage time is at most $\e^{1/2}n^{1/3}$, again by KPZ considerations.
Thus the second moment of $L_n-L_n^{B_i}$ is of order $\e\cdot (\e^{1/2}n^{1/3})^2 = \e^2 n^{2/3}$.

The bound on the probability of the geodesic passing through a small box is obtained by translating this into a statement about the location of the argmax of the top line of the Airy line ensemble, and relying on the strong comparisons of the latter to Brownian motion developed in \cite{CH,Ham22,CHH,duncan2} (see \cref{ss:locbrown}).
The same approach applies for the CDRP, but using the KPZ line ensemble instead of Airy.

The above argument in the positive-temperature CDRP setting is carried out in Sections \ref{s:mainproof} and \ref{s:es}.
In the actual argument, we will prove a slightly weaker version of \eqref{efron-stein}.
The precise resampling framework for white noise needed to implement the above strategy in the continuum is discussed in \cref{s:resampling-framework-short}.
It is worth mentioning that the latter demands significant measure-theoretic preparations, which we outline in Section \ref{s:meas-short} and detail in Section \ref{s:meas-long}.
While primarily technical in nature, we hope that this finds future use in the study of the CDRP.

\subsubsection{Directed landscape is a 2D black noise}\label{s:iop-black-noise}
We end this section with a quick discussion on how the above proof strategy can be used to prove \cref{t:black-noise}. 
As previously mentioned, \cite{HP24} showed that the directed landscape is a one-dimensional black noise. 
At this point an abstract result allows one to reduce the task of proving \cref{t:black-noise} to proving that the directed landscape is a \emph{two-dimensional noise} in the sense of Tsirelson and Vershik \cite{Tsir}. 
The notion of a two-dimensional noise involves a collection of $\sigma$-algebras indexed by open rectangles in $\R^2$ satisfying a few axioms.

Let $\cL$ be the directed landscape, as defined in \cref{d:dl}.
Let $\Rect$ be the collection of open rectangles in $\R^2$.
For any $R\in\Rect$, we define the \emph{restricted length}
\begin{align*}
    \cL_R(s,x;t,y) \coloneqq  
    \sup\left\{
        \int_s^t d\cL\circ \pi : 
        \pi\in C([s,t]),
        \quad
        \pi(s)=x, \pi(t)=y,
        \quad (r,\pi(r))\in R\;\forall r\in[s,t]
    \right\},
\end{align*}
where $\int_s^t d\cL\circ \pi$ is the \emph{length} of the path $\pi$ with respect to $\cL$ (defined formally later in \cref{def:length-dl}).
Informally, the restricted length $\cL_R$ corresponds to the ``internal length metric'' on $R$ induced by $\cL$.

We define $\cF_R$ to be the $\sigma$-algebra generated by $\cL_R$,
\begin{align*}
    \cF_R \coloneqq  \sigma\left(
        \cL_R(s,x;t,y) : (s,x),(t,y)\in R, s<t
    \right).
\end{align*}
Beyond certain natural axioms that hold trivially,
for the collection of $\sigma$-algebras  $\{\cF_R\}_{R \in \cR}$ to be a two-dimensional noise, the following two central properties must be satisfied.
\begin{itemize}
    \item \textbf{Disjoint independence property.} If $R_1,R_2\in\cR$ satisfy $R_1\cap R_2=\varnothing$, then $\cF_{R_1}$ and $\cF_{R_2}$ are independent.
    \item \textbf{Join property.} 
    If $R_1,R_2,R_3\in\Rect$ satisfy $R_1\cap R_2=\varnothing$ and $\overline{R_1\cup R_2}=\overline{R_3}$, then $\cF_{R_1}\vee \cF_{R_2} = \cF_{R_3}$ up to null sets, where $\vee$ denotes the join of $\sigma$-algebras.
\end{itemize}
The disjoint independence property is easier to prove. If the rectangles are separated in time, it follows from the independent temporal increments of the directed landscape. If they overlap in time but are separated in space, it follows from the spatial mixing property of the directed landscape as recorded in \cite[Lemma 3.2]{HP24}, which itself is a straightforward consequence of geodesic transversal fluctuation bounds.

\begin{figure}[htb]
\centering
    \begin{subfigure}[t]{0.35\textwidth}
        \includegraphics[width=\linewidth]{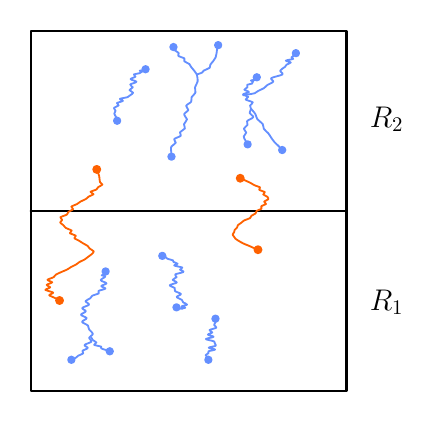}
    \end{subfigure}
    \hfill
    \begin{subfigure}[t]{0.6\textwidth}
        \includegraphics[width=\linewidth]{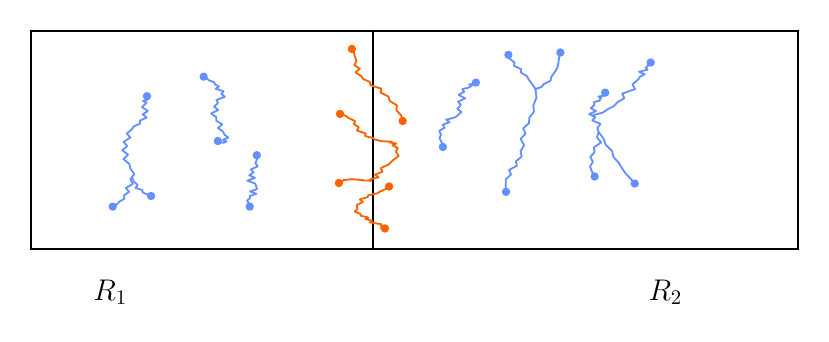}
    \end{subfigure}
\caption{
    \textbf{Left:} Temporally adjacent boxes.
    \textbf{Right:} Spatially adjacent boxes.\\
    In either case, let $R_3$ be the rectangle comprised of $R_1,R_2$, and the interface between them.
    The task is to reconstruct the restricted lengths in $R_3$ from the restricted lengths in $R_1,R_2$.
    These restricted lengths are respectively depicted using orange and blue ``restricted geodesics.''
    }\label{fig:iop-noise-1}
\end{figure}

The main challenge of proving \cref{t:black-noise} is verifying the join property, and we will break this argument into two cases.
The first case, where the two rectangles are temporally adjacent (Figure \ref{fig:iop-noise-1} (left)), is straightforward.
In this case the join property follows from the metric composition property of the directed landscape (see \cref{d:dl}) together with the continuity of restricted lengths. 
The second case, where the two rectangles are spatially adjacent (Figure \ref{fig:iop-noise-1} (right)), is really the heart of the matter, and most of the proof of \cref{t:main} is devoted to it.

\begin{figure}[tbh]
\centering
    \includegraphics[width=\linewidth]{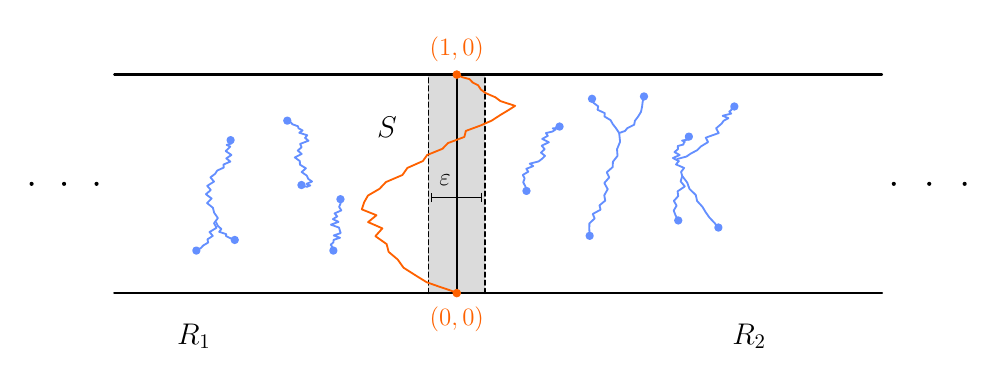}
\caption{
    The limiting case of spatially adjacent half-infinite strips.
    Given the data of the restricted lengths in $R_1$ and $R_2$ (depicted via blue ``restricted geodesics''), we seek to reconstruct the length of the orange geodesic, which can cross the interface between $R_1$ and $R_2$ (vertical black line through $(0,0)$).
    To do this, we consider a thin strip $S$ (shaded) of width $\e$ around the interface, and show that its influence on the orange geodesic's length is $O(\e^{1/2})$, which is negligible on the fluctuation scale.
    This uses an argument analogous to that in \eqref{efron-stein} (see also Figure \ref{fig:iop-LPP}).
    }\label{fig:iop-noise-2}
\end{figure}

To convey the strategy, consider the simplified limiting picture in which
\begin{align}\label{semi1}
    R_1=(0,1)\times(-\infty,0),
    & & 
    R_2=(0,1)\times(0,\infty),
\end{align}
so that the two rectangles are infinite half-strips lying on opposite sides of the line $x=0$, and their closures cover $\overline{R_3}=[0,1]\times\R$, as depicted in Figure \ref{fig:iop-noise-2}.
The goal is to show that the information of the restricted lengths in the two half-strips determines the restricted lengths in the full strip.
The resampling estimate enters by giving a quantitative way to remove a thin strip $S$ around the interface $x=0$: the length (last passage time) across the full strip can be reconstructed, up to an error negligible on the fluctuation scale, from the information on the two sides of the thin strip $S$.

Equivalently, the same thin-strip resampling principle can be expressed as a conditional variance bound. 
Let $S=(0,1)\times[-\e/2,\e/2]$ be the strip of width $\e$ around the interface between $R_1$ and $R_2$.
Let $\cG_\e$ be the $\sigma$-algebra generated by restricted lengths of paths constrained to avoid $S$, i.e.
\begin{align*}
    \cG_\e = \cF_{(0,1)\times(-\infty,-\e/2)} \vee \cF_{(0,1)\times (\e/2, \infty)}.
\end{align*}
The estimate analogous to \eqref{efron-stein} then reads
\begin{align}\label{equiv}
    \E[\Var(\cL(0,0;1,0) \mid \mathcal G_\varepsilon)] \ls \e^{1/2}.
\end{align}
In words, the information of the restricted lengths on either side of $S$ determines $\cL(0,0;1,0)$ up to an error that is negligible on the fluctuation scale $\Theta(1)$.
We prove \cref{t:black-noise} by implementing a version of this argument.
This appears in \cref{s:black-noise}.

There is a further complication to this argument that merits a brief discussion here.
As previously mentioned, we obtained a key bound on the probability of a geodesic passing through a given interval by recasting the event in terms of the Airy line ensemble.
In the simplified setting we have just discussed, where the union of the half-infinite rectangles is the entirety of the directed landscape on $(0,1)\times\R$, this connection to the Airy line ensemble persists.
However, for finite rectangles such a connection is lost.
To address this, we resort to a further discretization of the rectangle into short rectangles.
Each short rectangle can be individually coupled to a half-infinite problem \eqref{semi1} with high probability, allowing us to carry out the above proof strategy.
\\

The proof of Theorem \ref{t:pathsing} relies on a different set of ideas based on the analysis of Gaussian multiplicative chaos measures and viewing the CDRP measure as one, a perspective first put forth in \cite{QRV25}.
We refrain from discussing further details here.

\subsection{Notation}\label{s:notation}

Throughout the paper, $C,C',C'',c,c',c''$ etc. denote positive deterministic constants that do not depend on any parameters unless otherwise indicated, and whose values may change from one line to the next.
    We also adopt the usual Landau asymptotic notation:
    \begin{itemize}
        \item $A=O(B)$ (or $A\ls B$) means that $|A|\le CB$ for some $C>0$ that does not depend on any relevant asymptotic parameters.
        \item $A=\Theta(B)$ (or $A\asymp B$) means that $A=O(B)$ and $B=O(A)$.
        \item $A=o(B)$ (or $A\ll B$) means that $A/B$ converges to zero as the limit is taken in some asymptotic parameter.
    \end{itemize} 
    We write $A=O_{x,y,z}(B)$ (or $A\ls_{x,y,z}B$) to mean that $A=O(B)$ with an implied constant $C$ that depends on parameters $x,y,z$.

Below are other general notational conventions.
\begin{itemize}
    \item For real numbers $x,y$, we write $x\wedge y\coloneqq \min\{x,y\}$ and $x\vee y\coloneqq \max\{x,y\}$.
    \item $\N = \Z_{\ge 1} \coloneqq  \{1,2,\dots\}$ and $\Zpos \coloneqq  \{0,1,\dots\}$.
    \item For real numbers $x<y$, we write $\lb x,y\rb \coloneqq  [x,y]\cap\Z$.
    \item $\Rup \coloneqq  \{(s,x;t,y)\in\R^4 : s<t\}$.
    \item For a topological space $S$, we denote by $\cB(S)$ its Borel $\sigma$-algebra.
    \item For a topological space $S$, we denote by $C(S) = C(S,\R)$ the space of real-valued continuous functions on $S$ equipped with the topology of uniform convergence on compact sets.
    For $I\subset \R$, we denote by $C(S, I)$ the set of functions in $C(S,\R)$ taking values in $I$, equipped with the subspace topology.
    \item For $B\in \cB(\R^2)$, we define $L^2_0(B)\coloneqq \{f\in L^2(\R^2) : f\text{ vanishes a.e. on } B^c\}$.
    \item We denote by $\norm{\cdot}_2$ the $\ell^2$-norm on $\R^k$ for any $k\ge 1$ (the value of $k$ will always be clear from context).
    \item We reserve \textsf{sans serif} font to denote events about the white noise or the directed landscape.
    \item The complement of an event $A$ is denoted $\neg A$.
    \item Random constants are written in fraktur font, for example $\mf{C}$.
\end{itemize}

\addtocontents{toc}{\SkipTocEntry}
\subsection*{Acknowledgments}
SG thanks Sourav Chatterjee for introducing him to the topic of temperature and disorder chaos and for valuable discussions over the years. He also thanks Francesco Caravenna for discussions on \cite{CD25} as well as Chris Janjigian for discussions on \cite{AJRS22}. He was supported by a Miller Research Professorship at the Miller Institute for Basic Research in Science and NSF Career grant-1945172.
VG thanks Jeremy Clark and Vilas Winstein for helpful discussions on Gaussian multiplicative chaos, and was partially supported by the NSF
Graduate Research Fellowship Program under Grant No. DGE-2146752.
ZH was supported by NSF MSPRF 2503374.
A large language model was used to convert hand-drawn figures into TikZ code, and for a final round of proofreading.

\section{White noise, restricted CDRP partition function, and Efron--Stein}\label{s:meas-short}

As indicated in Section \ref{ss:idea}, our arguments based on resampling the white noise in strips demand a measure-theoretic framework. 
While the theory of white noise has a rich history (see for instance \cite{nualart,J97}), some of the statements featuring in our arguments, particularly the ones pertaining to the properties of the CDRP, are novel with non-trivial proofs.
To ease readability as well as highlight the results we believe could be of independent interest and likely to find future use, we have split the entire discussion on white noise and its properties into two sections.
The present section is short and records the statements of two central inputs.
A comprehensive discussion of white noise, CDRP, and resampling along with all the proofs of the statements in this section as well as more straightforward ones appear in the forthcoming  \cref{s:meas-long}.

\subsection{White noise and associated filtration}\label{ss:whitenoise-short}
For a topological space $S$, we denote by $\cB(S)$ its Borel $\sigma$-algebra.
Define
\begin{align}\label{e:def-Xi}
    (\Xi, \Fprod)  \coloneqq  \left(\R^{L^2(\R^2)},\; \cB(\R)^{\otimes L^2(\R^2)}\right).
\end{align}
We denote by $\xi$ a generic element of $\Xi$, i.e. $\xi$ is a function $L^2(\R^2)\to\R$.
By the Kolmogorov extension theorem \cite[Theorem 8.23]{Kal21},
there exists a unique probability measure $\P$ on the product space $(\Xi, \Fprod)$ 
such that for all $n\in\Z_{\ge 1}$ and all $f_1,\dots,f_n\in L^2(\R^2)$, under $\P$ the finite-dimensional marginal $(\xi(f_1),\dots,\xi(f_n))$ follows a multivariate Gaussian distribution with mean zero and covariances $\E[\xi(f_i)\xi(f_j)] = \langle f_i,f_j\rangle_{L^2(\R^2)}$ for $i,j\in\lb 1,n\rb$.

\begin{defn}[White noise]\label{def:white-noise}
    A \emph{white noise} on $\R^2$ is a $(\Xi,\Fprod)$-valued random variable with law $\P$.

    We realize white noise on the canonical complete probability space.
    Let $\cF$ be the completion of $\Fprod$ with respect to $\P$.
    We slightly abuse notation by writing $\xi$ for the identity map $\Xi\to\Xi$, so that $\xi : (\Xi,\cF,\P) \to (\Xi,\Fprod)$ is a white noise.
\end{defn}

The notion of \emph{restricting} white noise to subsets of $\R^2$ will be important for us, and this is formulated next.
For any Borel set $B\in\cB(\R^2)$, define
\begin{align}
    \label{e:L20}
    L^2_0(B) \coloneqq  \left\{f\in L^2(\R^2) : f\text{ vanishes a.e. on $B^c$}\right\},
\end{align}
where $B^c \coloneqq  \R^2\setminus B$.
Note that $L^2_0(B)$ is a closed subspace of $L^2(\R^2)$.
We define the product measurable space
\begin{align}\label{e:def-XiB}
    (\Xi_B, \Fprod_B) \coloneqq  \left(\R^{L^2_0(B)},\; \cB(\R)^{\otimes L^2_0(B)}\right).
\end{align}
Consider the \emph{restriction map}
\begin{equation}\label{e:def-Res}
    \begin{array}{ccc}
        \Res_B: (\Xi,\Fprod) &\longrightarrow &(\Xi_B,\Fprod_B)\\
        \xi &\longmapsto &\xi\big|_{L^2_0(B)}
    \end{array}.
\end{equation}
Note that $\Res_B$ is measurable.
\footnote{
    Since $\Fprod$ is a product $\sigma$-algebra, measurability of $\Res_B$ follows from the fact that for all $f\in L^2_0(B)$, the map $\xi \mapsto (\Res_B(\xi))(f) = \xi(f)$ is $\Fprod$-measurable.
    }
We remark that the $\sigma$-algebra generated by $\Res_B$ is naturally identified with $\Fprod_B = \cB(\R)^{\otimes L^2_0(B)}$.
Let $\cN\subset\cF$ be the $\sigma$-algebra generated by all $\P$-null subsets in $\cF$.
We denote
\begin{equation}\label{e:def-FB}
    \begin{split}
        \cF_B &\coloneqq  
        \sigma(\Res_B) \vee \cN\\
        &= \sigma\left(
            \xi(f) : f\in L^2_0(B)
        \right) \vee \cN,
    \end{split}
\end{equation}
where $\vee$ denotes the join of $\sigma$-algebras.
Note that $\cF_B \subset \cF$ is a $\sigma$-algebra on $\Xi$, and it depends on $B$ only up to Lebesgue-null sets.

\subsection{Restricting the polymer}\label{s:restricting-polymer-short}

An obvious difference between zero- and positive-temperature polymer models is as follows. 
In last passage percolation, a given geodesic is likely contained in a sufficiently wide deterministic box due to transversal fluctuation bounds, and hence its energy is ``almost'' measurable with respect to the noise inside that box.
At positive temperature, however, the polymer partition function always receives a contribution, though diminishing in value, from paths that wander arbitrarily far away.
In the proof of \cref{t:main} (occupying Sections \ref{s:mainproof} and \ref{s:es}), we will need a similar ``almost measurability'' statement for the partition function with respect to the local noise.
Towards this, in the present subsection we consider the  partition function \emph{restricted to those paths contained in a box}, and assert that it is measurable with respect to the noise in the box.

Fix $(s,x;t,y)\in\Rup$ and $s'<t'$ with $[s',t']\subset [s,t]$, as well as $\beta>0$.
For $a<b$, denote
\begin{align}
    B(a,b) \coloneqq  [s',t']\times [a,b].
\end{align}
Fix any Borel set $S\in \cB(C([s,t]))$.
Given a realization of the white noise $\xi\in \Xi$, set
\begin{align}\label{e:def-ZB}
    \cZ_{B(a,b)}^S &= \cZ_{B(a,b)}^S(s,x;t,y) \nonumber\\
    &\coloneqq  
    \cZ_\beta^\xi(s,x;t,y)\cdot \P_{\beta,(s,x),(t,y)}^\xi\left(
        S\cap\{X(r)\in [a,b]\text{ for all $r\in[s',t']$}\}
    \right),
\end{align}
where $\cZ_\beta^\xi(s,x;t,y)$ and $\P^\xi_{\beta,(s,x),(t,y)}$ are defined in \cref{def:intro-SHE} and \cref{t:AJRAS-polymer-existence} respectively.
Like in \eqref{wienertilt}, we can formally express $\cZ_{B(a,b)}^S(s,x;t,y)$ as an integral with respect to Brownian bridge:
\begin{multline*}
    \cZ_{B(a,b)}^S(s,x;t,y)
    \\
    =
    ``\;
    p(t-s,y-x)
    \mathrm{E}^{\mathrm{BB}}_{(s,x),(t,y)}
    \left[
        \exp\left(
            \beta \int_s^t \xi(r, X(r))\, dr
        \right)
        \1_{X\in S}
        \1_{X(r)\in [a,b]\text{ for all $r\in[s',t']$}}
    \right],\;"
\end{multline*}
where $p$ is the heat kernel.
Notice that,  heuristically, the expression on the RHS ought to depend only on the restriction of the white noise $\xi$ to $(([s,s']\cup[t',t])\times \R) \cup B(a,b)$, which suggests that the LHS---which is well-defined via \eqref{e:def-ZB}---is a measurable function of the same restriction of $\xi$.
This is confirmed by the following proposition.
\begin{prop}\label{p:Z-measurable-real}
    Let all notation be as above.
    Then
    $\cZ_{B(a,b)}^S$ is 
    measurable with respect to
    \begin{align*}
        \cF_{(([s,s']\cup[t',t])\times \R) \cup B(a,b)}
    \end{align*}
    where the latter was defined in \eqref{e:def-FB}.
\end{prop}
While an intuitive statement, the formal proof of \cref{p:Z-measurable-real} will take some work and is deferred to Section \ref{s:restricting-polymer}.

\subsection{Efron--Stein inequality, and resampling white noise}\label{s:resampling-framework-short}

The main result of this subsection (\cref{p:efron-stein-ineq}) is a version of the Efron--Stein inequality for white noise, which will play a fundamental role in the proof of \cref{t:main} (see the proof sketch in \cref{ss:idea}).
The statement is exactly what one would expect, but requires first setting up a precise framework for \emph{resampling} the white noise inside any Borel subset of $\R^2$.
As a signpost, we first recall the classical Efron--Stein inequality  (e.g. \cite[Theorem 3.1]{Bou13}), which itself will find use in the proof of \cref{t:black-noise}.

\begin{prop}[Efron--Stein inequality]\label{p:efron-stein-classical}
    Let $X_0,X_1,\dots,X_n, Y$ be random vectors (possibly with different and/or infinite lengths) on the same probability space, such that $X_0,\dots,X_n$ are conditionally independent given $Y$.
    Write $\bfX \coloneqq  (X_0,\dots,X_n)$.
    Let $\wt{\bfX} = (\wt{X}_0,\dots,\wt{X}_n)$ be a conditionally independent copy of $\bfX$ given $Y$.
    Write
    \begin{align*}
        \wt{\bfX}^{(k)} \coloneqq  (X_0,\dots,X_{k-1},\wt{X}_k, X_{k+1},\dots,X_n).
    \end{align*}
    Note that we have the equality of unconditional joint laws
    \begin{align*}
        (\bfX,Y)\law (\wt{\bfX},Y) \law (\wt{\bfX}^{(0)}, Y)\law \cdots \law (\wt{\bfX}^{(n)},Y).
    \end{align*}

    Let $Z=Z(\bfX,Y)$ be a square-integrable function of $(\bfX,Y)$.
    Then
    \begin{align*}
        \E[\Var(Z \mid Y)]
        = \frac12
        \E\left[
            \left(
                Z(\bfX,Y)
                - Z(\wt{\bfX},Y)
            \right)^2
        \right]
        &\le
        \frac12
        \sum_{k=0}^n 
        \E\left[
            \left(
                Z(\bfX,Y)
                - Z(\wt{\bfX}^{(k)},Y)
            \right)^2
        \right].
    \end{align*}
    In fact, the above holds almost surely with $\E[\,\cdot\,|\,Y]$ in place of $\E[\,\cdot\,]$.
\end{prop}

In \cref{p:efron-stein-classical}, the pair $(\wt{\bfX},Y)$ should be interpreted as being obtained from $(\bfX,Y)$ by \emph{resampling} $\bfX$.
Similarly, $(\wt{\bfX}^{(k)},Y)$ is obtained from $(\bfX,Y)$ by resampling $X_k$.\\

We now set up the aforementioned framework for resampling white noise inside any Borel subset of $\R^2$.
Define the probability space
\begin{align}\label{e:def-Xi-resample}
    (\Xires, \Fres, \Pres)
    \coloneqq  
    (\Xi\times \Xi,
    \overline{\cF\otimes \cF},
    \P\otimes \P),
\end{align}
where $(\Xi,\cF,\P)$ is the canonical space of white noise on $\R^2$ (see \cref{def:white-noise}), and
where $\overline{\cF\otimes \cF}$ is the completion of $\cF\otimes \cF$ with respect to $\P\otimes \P$.
The subscript ``$\res$'' is short for ``resample.''
Observe that the coordinate projections $(\Xires,\Fres,\Pres) \to (\Xi,\Fprod)$ given by $(\xi_1,\xi_2)\mapsto \xi_1$ and $(\xi_1,\xi_2)\mapsto \xi_2$ comprise a pair of independent white noises on $\R^2$.
Further, since $\Fprod\otimes\Fprod \subset \cF\otimes\cF \subset \Fres$, and $\P$ is the pushforward of $\Pres$ by each of the above two coordinate projections, the latter extend to measurable maps $(\Xires,\Fres,\Pres) \to (\Xi,\cF,\P)$.

Fix any Borel set $B\in\cB(\R^2)$ and let $L^2_0(B)$ be as in \eqref{e:L20}.
Let $P_{B} : L^2(\R^2) \to L^2_0(B)$ be the orthogonal projection onto $L^2_0(B)$, i.e. $P_B f \coloneqq  f\1_B$.
We define two measurable maps $\xi,\eta_B : (\Xires, \Fres,\Pres) \to (\Xi,\cF,\P)$ by
\begin{equation}\label{e:def-xi-eta}
    \begin{split}
        \xi(\xi_1, \xi_2) &\coloneqq  \xi_1,\\
        \eta_B(\xi_1,\xi_2) &\coloneqq  \xi_1 \circ P_{B^c} + \xi_2 \circ P_{B}.
    \end{split}
\end{equation} 
To help parse this, note that given $(\xi_1,\xi_2)\in \Xires = \Xi\times\Xi$, the expressions on the RHS of \eqref{e:def-xi-eta} are elements of $\Xi = \R^{L^2(\R^2)}$ which act on $f\in L^2(\R^2)$ by
\begin{align*}
    \xi(\xi_1,\xi_2)(f) &= \xi_1(f),\\
    \eta_B(\xi_1,\xi_2)(f) &= \xi_1(f\1_{B^c}) + \xi_2(f\1_{B}).
\end{align*}
When $B$ is clear from context, we will sometimes abbreviate $\eta = \eta_B$.

We now verify that $\xi,\eta_B$ are marginally white noises.
\begin{lemm}\label{l:coupled-white-noises}
    With the above notation,
    $\xi$ and $\eta_B$ are white noises on $\R^2$.
\end{lemm}
\begin{proof}
This is clear for $\xi$ since it is just a coordinate projection.
To show that $\eta=\eta_B$ is a white noise, by the $\pi$-$\lambda$ theorem it suffices to compute its finite-dimensional marginals.
For any $n\in\N$ and $f_1,\dots,f_n\in L^2(\R^2)$, we have
\begin{align*}
    (\eta(f_1), \dots,\eta(f_n))
    = (\xi_1(f_1\1_{B^c}),\dots,\xi_1(f_n\1_{B^c}))
    + (\xi_2(f_1\1_{B}),\dots,\xi_2(f_n\1_{B})).
\end{align*}
Since $\xi_1,\xi_2$ are independent white noises, the RHS is the sum of two independent centered multivariate Gaussians, and thus the LHS is a centered multivariate Gaussian.
By independence, the entries have covariances
\begin{align*}
    \Eres\left[
        \eta(f_i)\eta(f_j)
    \right]
    &= \Eres\left[
        \xi_1(f_i\1_{B^c}) \xi_1(f_j\1_{B^c})
    \right]
    + \Eres\left[
        \xi_2(f_i\1_{B})\xi_2(f_j\1_{B})
    \right]\\
    &= \langle f_i\1_{B^c}, f_j\1_{B^c}\rangle_{L^2(\R^2)}
    + \langle f_i\1_{B}, f_j\1_{B}\rangle_{L^2(\R^2)}\\
    &= \langle f_i,f_j\rangle_{L^2(\R^2)}.
\end{align*}
So $\eta$ has the finite-dimensional distributions of a white noise (recall \cref{def:white-noise}).
\end{proof}

In summary, by \cref{l:coupled-white-noises}, we have a probability space $(\Xires,\Fres,\Pres)$ supporting a white noise $\xi$, such that given any Borel set $B\in\cB(\R^2)$, we can construct another white noise $\eta_B$ on the same probability space that can be viewed as obtained from $\xi$ by \emph{resampling} the restriction of $\xi$ to $B$.
Indeed, it may be helpful to note that for all $f,g\in L^2(\R^2)$,
\begin{align*}
    \Eres\left[\xi(f)\eta_B(g)\right] 
    = \E\left[\xi(f\1_{B^c})\xi(g\1_{B^c})\right].
\end{align*}

With the above preparations, we can now formulate the analogue of the Efron--Stein inequality \cref{p:efron-stein-classical} in the white noise context.
The proof appears in Section \ref{s:resampling-framework}.

\begin{prop}[Efron--Stein inequality for white noise]\label{p:efron-stein-ineq}
    Let $\{B_k\}_{k\in \Z}\subset\cB(\R^2)$ be a bi-infinite sequence of pairwise essentially disjoint Borel sets (meaning $|B_k\cap B_{k'}|=0$ for all $k\ne k'$).
    Denote $B \coloneqq  \bigcup_{k\in\Z} B_k$.

    For $A\in\cB(\R^2)$ let $\xi,\eta_A : (\Xires,\Fres,\Pres)\to(\Xi,\cF,\P)$ be coupled white noises as in \eqref{e:def-xi-eta} (so $\eta_A$ is obtained from $\xi$ by resampling $\xi|_A$).
    Then for any $Z\in L^2(\Xi,\cF,\P)$, we have
    \begin{align*}
        \Eres\left[
            \bigl(Z(\xi)  - Z(\eta_{B})\bigr)^2
        \right]
        \le \sum_{k\in\Z}\Eres\left[
            \bigl(Z(\xi) - Z(\eta_{B_k})\bigr)^2
        \right].
    \end{align*}
\end{prop}

Finally, we reiterate that additional basic properties of the resampling framework \eqref{e:def-xi-eta} beyond those presented above will be needed for \cref{t:main}.
These are developed in the forthcoming section (see Section \ref{s:resampling-framework}).

\section{Measure-theoretic details}\label{s:meas-long}

This section discusses white noise, the CDRP, and the resampling framework \eqref{e:def-Xi-resample}--\eqref{e:def-xi-eta}, picking up from the overview in \cref{s:meas-short}.
While we have attempted to provide a comprehensive treatment and hope that this finds use in future studies of the CDRP, this section is nonetheless long and technical, and the reader may therefore be well served to skip it on first reading and consult it subsequently as a reference when needed.
At a high level, this section accomplishes three goals:
\begin{enumerate}[label={\rm(\arabic*)},ref={\rm \arabic*}]
    \item\label{meas-goal-1} Prove \cref{p:Z-measurable-real}.
    \item\label{meas-goal-2} Collect inputs from \cite{AJRS22}, including almost-sure scaling properties of the CDRP.
    \item\label{meas-goal-3} Establish basic properties of the resampling framework \eqref{e:def-Xi-resample}--\eqref{e:def-xi-eta}, and prove the Efron--Stein inequality for white noise (\cref{p:efron-stein-ineq}).
\end{enumerate}
Goal \ref{meas-goal-1} will take by far the most work.
Our proof of \cref{p:Z-measurable-real} relies on analyzing the \emph{(Wiener) chaos expansion} of the stochastic heat equation \eqref{e:she-intro}.
With the above three goals in mind, we now present a more detailed outline of the section.
\begin{itemize}
    \item In Section \ref{ss:white-noise-filtration} we prove some basic properties of the white noise ``filtration'' defined in \eqref{e:def-FB}.
    These are only used as inputs for other results in \cref{s:meas-long}.
    \item In Section \ref{ss:stocint}, we recall the basic theory of Wiener chaos.
    This is only used in the proof of \cref{p:Z-measurable-real}.
    \item In Sections \ref{s:automorphisms} and \ref{ss:stocheatpolymer} we record inputs from \cite{AJRS22}, achieving Goal \ref{meas-goal-2}.
    The results of these subsections are used in the proof of \cref{p:Z-measurable-real} and throughout Sections \ref{s:LE}--\ref{s:es}.
    \item In \cref{s:restricting-polymer}, we prove \cref{p:Z-measurable-real}, achieving Goal \ref{meas-goal-1}.
    \item In \cref{s:resampling-framework}, we expand on the resampling framework \eqref{e:def-Xi-resample}--\eqref{e:def-xi-eta}, achieving Goal \ref{meas-goal-3}.
\end{itemize}

\subsection{Properties of the white noise filtration}\label{ss:white-noise-filtration}

The following lemma records some basic properties of the white noise ``filtration'' $\{\cF_B\}_{B\in\cB(\R^2)}$ defined in \eqref{e:def-FB}.

\begin{lemm}[Properties of white noise filtration]\label{l:independence-disjoint}
    Fix $k\ge 2$ and Borel sets $B_1,\dots,B_k\in\cB(\R^2)$.
    Then the following hold.
    \begin{enumerate}[label={\rm(\roman*)}]
        \item\label{generation} We have 
        $\cF_{B_1}\vee\cdots\vee\cF_{B_k}=\cF_{B_1\cup\cdots\cup B_k}$.
        \item\label{independence} If $B_i\cap B_j$ has measure zero for all $i\ne j$, then $\cF_{B_1},\dots,\cF_{B_k}$  are independent.
    \end{enumerate}
\end{lemm}
\begin{proof}
    \ref{generation}:
    For $B\in\cB(\R^2)$, denote $S_B \coloneqq  \{\xi(f) : f\in L_0^2(B)\}$.
    Clearly,
    \begin{align*}
        S_{B_1}\cup\cdots\cup S_{B_k} \subset S_{B_1\cup\cdots\cup B_k}.
    \end{align*}
    Moreover, given $f \in L^2_0(B_1\cup\cdots\cup B_k)$, we can write
    \begin{align*}
        \xi(f) = \xi(f\1_{B_1}) + \sum_{i=2}^k
        \xi(f\1_{B_i\setminus(B_1\cup\cdots\cup B_{i-1})})
        \qquad\text{almost surely.}
    \end{align*}
    Since the $i\textsuperscript{th}$ summand above belongs to $S_{B_i}$, it follows that 
    $S_{B_1}\cup\cdots\cup S_{B_k}$ and $S_{B_1\cup\cdots\cup B_k}$ span the same linear subspace of $L^2(\Xi,\cF,\P)$,
    and hence $\sigma(S_{B_1}\cup\cdots\cup S_{B_k})= \sigma(S_{B_1\cup\cdots\cup B_k})$ up to null sets.
{Adding in all $\P$-null sets proves \ref{generation} (see \eqref{e:def-FB}).}

    \ref{independence}:
    Fix $f_{i,1},\dots,f_{i,n_i}\in L^2_0(B_i)$ for  $i\in\lb 1,k\rb$.
    For all $i\ne j$, since $B_i\cap B_j$ has measure zero, it follows that for all $p\in \lb 1,n_i\rb$ and all $q\in \lb 1,n_j\rb$,
    \begin{align*}
        \E\left[\xi(f_{i,p})\xi(f_{j,q})\right] = \langle f_{i,p},f_{j,q}\rangle_{L^2(\R^2)}
        = \int_{B_i\cap B_j} f_{i,p}(a)f_{j,q}(a)\,da = 0.
    \end{align*} 
    Since $\{\xi(f_{i,p}) : i\in\lb 1,k\rb, p \in \lb 1,n_i\rb\}$ are jointly Gaussian with mean zero, it follows that the collections $\{\xi(f_{i,p}):p\in \lb 1,n_i\rb\}$ are independent across $i\in\lb 1,k\rb$.
    The claim now follows by the $\pi$-$\lambda$ theorem.
\end{proof}

\subsubsection{Conditional expectation}\label{sss:CE-general-noise}

Conditional expectations will act as a useful device for us in several arguments involving white noise, and in this subsection we record properties of the conditional expectation operator $\E[\,\smallbullet\,|\,\cF_B]$ for $B\in\cB(\R^2)$ and $\cF_B$ defined in \eqref{e:def-FB}.
These results may seem unmotivated for the moment but will find applications in Sections \ref{s:restricting-polymer} and \ref{s:resampling-framework}.

We first record a standard measure-theoretic fact that will be used to prove the upcoming Lemmas \ref{l:FAFBisFAB} and \ref{l:conditionally-independent}.

\begin{lemm}\label{l:independent-sub-sigma-alg}
    Let $\cF_1,\cF_2 \subset \cF$ be independent sub-$\sigma$-algebras on $(\Xi,\cF,\P)$, and for $i=1,2$ let $X_i$ be an $\cF_i$-measurable random variable with finite mean.
    Then for any sub-$\sigma$-algebras $\cF_1'\subset \cF_1$ and $\cF_2'\subset \cF_2$, we have
    \begin{align*}
        \E\left[
            X_1X_2\,\middle|\,\cF_1'\vee \cF_2'
        \right]
        &= \E[X_1\,|\,\cF_1']\cdot\E[X_2\,|\,\cF_2'].
    \end{align*}
\end{lemm}
\begin{proof}
    For any events $F_1\in\cF_1'$ and $F_2\in\cF_2'$, we have by independence
    \begin{align*}
        \E\left[
            \1_{F_1}\1_{F_2} \E[X_1\,|\,\cF_1']\cdot
            \E[X_2\,|\,\cF_2']
        \right]
        &= \E\left[
            \1_{F_1}\E[X_1\,|\,\cF_1']
            \right]
        \cdot
        \E\left[
            \1_{F_2}
            \E[X_2\,|\,\cF_2']
        \right]\\
        &=\E\left[
            \1_{F_1}X_1
        \right]
        \cdot
        \E\left[
            \1_{F_2}X_2
        \right]\\
        &= \E\left[
            \1_{F_1}\1_{F_2}X_1X_2
        \right].
    \end{align*}
    The lemma now follows from the $\pi$-$\lambda$ theorem.
\end{proof}

\begin{lemm}[Conditional expectation preserves measurability]\label{l:FAFBisFAB}
    For any Borel sets $A,B\in\cB(\R^2)$ and any $X\in L^2(\Xi,\cF,\P)$, if $X$ is $\cF_B$-measurable, then
    \begin{align*}
        \E[X|\cF_A] = \E[X|\cF_{A\cap B}].
    \end{align*}
\end{lemm}
\begin{proof}
    By \cref{l:independence-disjoint}, we have $\cF_A = \cF_{A\cap B} \vee \cF_{A\cap B^c}$,
    and $\cF_{A\cap B^c}$ is independent of $\cF_B$.
    Thus since $X$ is $\cF_B$-measurable, by \cref{l:independent-sub-sigma-alg} (take $X_1=X$ and $X_2=1$), we get
    \begin{align*}
        \E[X|\cF_A]
        &= \E[X|\cF_{A\cap B} \vee \cF_{A\cap B^c}]
        =\E[X|\cF_{A\cap B}]
    \end{align*}
    as claimed.
\end{proof}

By \cref{l:independence-disjoint}\ref{independence}, the restrictions of the white noise to essentially disjoint subsets of $\R^2$ are independent.
The next lemma, a consequence of the same argument used to prove \cref{l:FAFBisFAB}, states that this independence property persists even after revealing part of the white noise.

\begin{lemm}\label{l:conditionally-independent}
    Fix an integer $k\ge 2$ and Borel sets $B_1,\dots,B_k\in\cB(\R^2)$ such that $B_i\cap B_j$ has measure zero for all $i\ne j$.
    Let $X_1,\dots,X_k\in L^2(\Xi,\cF,\P)$ be random variables such that $X_i$ is $\cF_{B_i}$-measurable for every $i\in\lb 1,k\rb$.
    Then for any Borel set $A\in\cB(\R^2)$, we have
    \begin{align}\label{e:210}
        \E\left[
            \prod_{i=1}^k X_i
            \;\middle|\;\cF_A
        \right]
        &= \prod_{i=1}^k \E\left[
            X_i\,\middle|\,\cF_A
        \right].
    \end{align}
\end{lemm}
\begin{proof}
    By \cref{l:independence-disjoint}\ref{generation}, the product $\prod_{i=1}^kX_i$ is $\cF_{B_1\cup\cdots\cup B_k}$-measurable.
    Therefore, by \cref{l:FAFBisFAB},
    the identity \eqref{e:210} is equivalent to
    \begin{align*}
        \E\left[
            \prod_{i=1}^k X_i\;\middle|\;
            \cF_{(A\cap B_1)\cup\cdots\cup(A\cap B_k)}
        \right]
        =
        \prod_{i=1}^k \E\left[
            X_i\,\middle|\,\cF_{A\cap B_i}
        \right].
    \end{align*}
    To establish the latter, we note that $\cF_{A\cap B_1},\dots,\cF_{A\cap B_k}$ are independent by \cref{l:independence-disjoint}\ref{independence}.
    Therefore, by \cref{l:independence-disjoint}\ref{generation} and repeated applications of \cref{l:independent-sub-sigma-alg}, we get
    \begin{align*}
        \E\left[
            \prod_{i=1}^k X_i\;\middle|\;
            \cF_{(A\cap B_1)\cup\cdots\cup(A\cap B_k)}
        \right]
        &=
        \E\left[
            \prod_{i=1}^k X_i\;\middle|\;
            \bigvee_{i=1}^k \cF_{A\cap B_i}
        \right]
        \\
        &= \E\left[X_1\,\middle|\,\cF_{A\cap B_1}\right]
        \E\left[
            \prod_{i=2}^k X_i\;\middle|\;
            \bigvee_{i=2}^k \cF_{A\cap B_i}
        \right]
        \\
        &\;\,\vdots\\
        &= \prod_{i=1}^k \E\left[
            X_i\,\middle|\,\cF_{A\cap B_i}
        \right]
    \end{align*}
    as desired.
\end{proof}

\subsection{Wiener chaos}
\label{ss:stocint}

In this subsection we recall basic definitions and properties of Wiener chaos for white noise, which will only be used in the proof of \cref{p:Z-measurable-real}.
The content of this subsection is largely quoted from \cite{nualart}, and the reader is referred there for a more systematic treatment (see also \cite{J97}).
It is straightforward to show that $f\mapsto \xi(f)$ defines a linear isometric embedding $L^2(\R^2)\hookrightarrow L^2(\Xi,\cF,\P)$  (see \cite[p.~4]{nualart}),
and we adopt this perspective for the present discussion.

Let $\cB_0(\R^2) \coloneqq  \{A\in\cB(\R^2):|A|<\infty\}$ be the collection of Borel sets with finite Lebesgue measure, the latter denoted by $|\cdot|$.
For $n\in\Z_{\ge1}$, we denote by $\cE_n$ the following real vector space of \emph{elementary functions}:
\begin{align}\label{e:elementary-functions}
    \cE_n \coloneqq  \mathrm{Span}\left(\left\{\1_{A_1\times\cdots\times A_n}:A_1,\dots,A_n\in\cB_0(\R^2) \text{ pairwise disjoint}\right\}\right).
\end{align}
We define the \emph{stochastic integral} map $I_n:\cE_n\to L^2(\Xi,\cF,\P)$ by setting
\begin{align}\label{e:def-I}
    I_n(\1_{A_1\times\cdots\times A_n}) \coloneqq  \xi(\1_{A_1})\cdots\xi(\1_{A_n})
    \qquad\text{for all pairwise disjoint $A_1,\dots,A_n\in\cB_0(\R^2)$},
\end{align}
and extending by linearity to all of $\cE_n$.
It can be shown that $I_n$ is a continuous linear map, and that $\cE_n$ is dense in $L^2(\R^{2n})$ (see \cite[Section 1.1.2]{nualart}), so $I_n$ extends to a continuous linear map $I_n:L^2(\R^{2n})\to L^2(\Xi,\cF,\P)$.
It will be convenient to extend this notation to the case $n=0$: we write $L^2(\R^0)\coloneqq \R$ for the space of constant functions, and write $I_0$ for the identity map on $L^2(\R^0)$.

The following basic properties of $I_n$ will be used throughout this section.
\begin{lemm}[Properties of $I_n$]\label{l:properties-of-I}
    The following properties hold.
    \begin{enumerate}[label={\rm(\roman*)}]
        \item\label{symmetry} \textup{(Symmetry)}
        For all $n\in\Z_{\ge1}$ and all $f\in L^2(\R^{2n})$, we have
        $I_n(f) = I_n(\wt{f})$, 
        where $\wt{f}$ denotes the \emph{symmetrization} of $f$:
        \begin{align}\label{e:def-symmetrization}
            \wt{f}(a_1;\cdots;a_n) \coloneqq  \frac{1}{n!}
            \sum_{\pi\in \mathfrak{S}_n} f(a_{\pi(1)};\cdots;a_{\pi(n)}),
            \qquad\qquad (a_1;\cdots;a_n)\in \R^{2n},
        \end{align}
        where $\mathfrak{S}_n$ denotes the group of permutations of $\{1,\dots,n\}$.

        For $n=0$ and $f\in L^2(\R^0)$, we denote $\wt{f}\coloneqq f$, so that $I_0(\wt{f})=I_0(f)$.
        \item\label{orthogonality} \textup{(Orthogonality)}
        For all $n,m\in\Zpos$, all $f\in L^2(\R^{2n})$, and all $g\in L^2(\R^{2m})$,
        \begin{align*}
            \E\left[I_n(f)I_m(g)\right]&= \begin{cases}
                0,&\quad n\ne m\\
                n!\langle \wt{f},\wt{g}\rangle_{L^2(\R^{2n})},&\quad n=m.
            \end{cases}
        \end{align*}
        Moreover,
        \begin{align}\label{e:symm-inner-product}
            n!\langle \wt{f},\wt{g}\rangle_{L^2(\R^{2n})}
            &= \sum_{\pi\in\mathfrak{S}_n} \int_{\R^{2n}}f(a_1;\cdots;a_n)\,g(a_{\pi(1)};\cdots;a_{\pi(n)})\,da_1\dots da_n.
        \end{align}
    \end{enumerate}
\end{lemm}
\begin{proof}
    Proofs can be found in \cite[Section 1.1.2, pp. 8--10]{nualart}, with the exception of the identity \eqref{e:symm-inner-product} which follows from a straightforward calculation. 
\end{proof}

The following fundamental result says that stochastic integrals span $L^2(\Xi,\cF,\P)$
(see \cite[Theorem 1.1.2]{nualart} for a proof). 
\begin{prop}[Chaos expansion]\label{p:chaos-expansion}
    Every square-integrable random variable $X\in L^2(\Xi,\cF,\P)$ admits an $L^2(\Xi,\cF,\P)$-orthogonal expansion (called a \emph{(Wiener) chaos expansion}) of the form
    \begin{align*}
        X = \sum_{n=0}^\infty I_n(f_n),
        \qquad\qquad f_n\in L^2(\R^{2n}),
    \end{align*}
    where the equality is as elements of $L^2(\Xi,\cF,\P)$.
    
    For $X,Y\in L^2(\Xi,\cF,\P)$ with chaos expansions $X=\sum_{n=0}^\infty I_n(f_n)$ and $Y=\sum_{n=0}^\infty I_n(g_n)$, 
    we have
    \begin{align}\label{e:parseval}
        \E[XY] = \sum_{n=0}^\infty n! \langle\wt{f}_n, \wt{g}_n\rangle_{L^2(\R^{2n})},
    \end{align}
    by Parseval's identity and \cref{l:properties-of-I}\ref{orthogonality}.
\end{prop}

\subsubsection{Conditional expectation and Wiener chaos}\label{ss:condexp}

We now continue the discussion of conditional expectation from \cref{sss:CE-general-noise} in the Wiener chaos setting.
We start by recording a formula for $\E[X\,|\,\cF_B]$ in terms of the chaos expansion of $X\in L^2(\Xi,\cF,\P)$, where $B\in\cB(\R^2)$ and $\cF_B$ was defined in \eqref{e:def-FB}.
A proof can be found in \cite[Lemma 1.2.5]{nualart}.
\footnote{
    Strictly speaking, \cite[Lemma 1.2.5]{nualart} considers the conditional expectation given $\sigma(\xi(\1_A) : A\subset B, A\in \cB_0(\R^2))$, but it is straightforward to check that the latter equals $\cF_B$ modulo null sets.
}
\begin{lemm}[Conditional expectations of chaos expansions]\label{l:I-given-B}
    Fix any Borel set $B\in\cB(\R^2)$ and any square-integrable random variable $X\in L^2(\Xi,\cF,\P)$, and denote its chaos expansion by 
    $X = \sum_{n=0}^\infty I_n(f_n)$ (see \cref{p:chaos-expansion}).
    Then we have
    \begin{align*}
        \E[X|\cF_B]
        = \sum_{n=0}^\infty I_n(f_n \1_{B}^{\otimes n})
        \qquad\text{in $L^2(\Xi,\cF,\P)$,}
    \end{align*}
    where for $n=0$, the $0\textsuperscript{th}$ tensor power of any function is identically $1$.
    Note that for $n\ge 1$, we have $\1_B^{\otimes n} = \1_{B^n}$, where $B^n\coloneqq B\times\cdots\times B$ is the $n$-fold Cartesian product.

    It follows by \eqref{e:parseval} that for any $X,Y\in L^2(\Xi,\cF,\P)$ with chaos expansions $X=\sum_{n=0}^\infty I_n(f_n)$ and $Y=\sum_{n=0}^\infty I_n(g_n)$, we have
    \begin{align}\label{e:parseval-CE}
        \E\Bigl[
            \E[X|\cF_B]\cdot\E[Y|\cF_B]
        \Bigr]
        &=
        \sum_{n=0}^\infty \E\left[
            I_n(f_n\1_B^{\otimes n})
            I_n(g_n\1_B^{\otimes n})
        \right].
    \end{align}
\end{lemm}

\cref{l:I-given-B} implies continuity of the filtration $\{\cF_B\}_{B\in\cB(\R^2)}$, as we show in the next lemma.
A similar statement for $\{\cF_{[s,t]\times\R}\}_{s<t}$ appears in \cite[Proposition 3d3 and Corollary 3d5]{Tsir}.
\begin{lemm}[Continuity of the white noise filtration]\label{l:nested}
    Let $\{B_r\}_{r>0}\subset \cB(\R^2)$ be an increasing collection of Borel sets: $B_{r}\subset B_{r'}$ for all $0<r<r'$.
    Then
    \begin{align}\label{e:decreasing}
        \cF_{\bigcap_{r>0}B_r} = \bigcap_{r>0}\cF_{B_r}
    \end{align}
    and
    \begin{align}\label{e:increasing}
        \cF_{\bigcup_{r>0}B_r} = 
        \bigvee_{r>0}\cF_{B_r}.
    \end{align}
\end{lemm}
\begin{proof}
    Note that $\bigcap_{r>0}B_r=\bigcap_{m\in\N}B_{1/m}$ and $\bigcup_{r>0}B_r = \bigcup_{m\in\N}B_m$ are Borel sets.
    We first prove \eqref{e:decreasing}, and afterwards explain how to modify the argument to establish \eqref{e:increasing}.

    Fix $X\in L^2(\Xi,\cF,\P)$.
    By \cref{p:chaos-expansion} we have the chaos expansion
    $X = \sum_{n=0}^\infty I_n(f_n)$
    for suitable $f_n\in L^2(\R^{2n})$, and by \cref{l:properties-of-I}\ref{symmetry} 
    we can assume $f_n=\wt{f_n}$ for all $n$, where $\wt{f_n}$ is the symmetrization defined in \eqref{e:def-symmetrization}.

    Denote $B_0\coloneqq \bigcap_{r>0}B_r$.
    We now show that $\E[X|\cF_{B_r}]\to \E[X|\cF_{B_0}]$ in $L^2(\Xi,\cF,\P)$ as $r\to 0$.
    First, using \cref{l:I-given-B}, and applying the second moment formula \eqref{e:parseval},
    we get
    \begin{align}\label{e:parseval1}
        \lim_{r\to 0}
        \norm{\E[X|\cF_{B_r}] - \E[X|\cF_{B_0}]}_{L^2(\Xi,\cF,\P)}^2
        &= 
        \lim_{r\to 0}
        \norm{
            \sum_{n=1}^\infty 
            I_n\left(f_n(\1_{B_r}^{\otimes n} - \1_{B_0}^{\otimes n})\right)
        }_{L^2(\Xi,\cF,\P)}^2 \nonumber\\
        \overset{\eqref{e:parseval}}&{=} 
        \lim_{r\to 0}
        \sum_{n=1}^\infty 
        n!
        \norm{f_n \1_{B^n_r\setminus B^n_0}}^2_{L^2(\R^{2n})}.
    \end{align}
    By the dominated convergence theorem, since $n! \lVert f_n \1_{B^n_r\setminus B^n_0} \rVert_{L^2(\R^{2n})}^2 \le n! \norm{f_n}_{L^2(\R^{2n})}^2 < \infty$ and $\sum_{n=1}^\infty n! \norm{f_n}_{L^2(\R^{2n})}^2 = \Var(X) <\infty$,
    we can move the $\lim\limits_{r\to 0}$ inside the infinite sum and the $L^2(\R^{2n})$-norms:
    \begin{align*}
        \eqref{e:parseval1}
        &=
        \sum_{n=1}^\infty 
        \lim_{r\to 0}
        n! \norm{f_n \1_{B^n_r\setminus B^n_0}}_{L^2(\R^{2n})}^2
        =
        \sum_{n=1}^\infty 
        n! \norm{f_n \lim_{r\to 0}\1_{B^n_r\setminus B^n_0}}_{L^2(\R^{2n})}^2
        =0.
    \end{align*}
    
    We next compare the $\sigma$-algebras $\cF_{B_0}$ and $\bigcap_{r>0}\cF_{B_r}$.
    Since $\cF_{B_r} \subset \cF_{B_{r'}}$ for $r<r'$, standard backwards martingale theory implies that
    \begin{align}\label{e:backwards-mg}
        \E[X|\cF_{B_r}] \xrightarrow{r\to 0} \E\left[X\,\middle|\,\bigcap_{r>0}\cF_{B_r}\right]
        \qquad\text{in $L^2(\Xi,\cF,\P)$}.
    \end{align}
    The above two displays imply
    $\E[X|\cF_{B_0}]=\E[X|\bigcap_{r>0}\cF_{B_r}]$.
    Applying this with $X=\1_{\sA}$ for any event $\sA\in \cF_{B_0}$ or $\sA\in \bigcap_{r>0}\cF_{B_r}$, we deduce that $\cF_{B_0} = \bigcap_{r>0}\cF_{B_r}$ modulo null sets, hence by completeness $\cF_{B_0} = \bigcap_{r>0}\cF_{B_r}$.
    This proves \eqref{e:decreasing}.

    The above argument applies with straightforward modifications to yield \eqref{e:increasing}.
    The main difference is that instead of using backwards martingale theory as in \eqref{e:backwards-mg}, one instead notes that by usual martingale theory, since $X\in L^2(\Xi,\cF,\P)$,
    \begin{align*}
        \E[X|\cF_{B_r}] \xrightarrow{r\to\infty} \E\left[X\,\middle|\,\bigvee_{r>0}\cF_{B_r}\right]
        \qquad\text{in $L^2(\Xi,\cF,\P)$}.
    \end{align*}
    We omit the details.
\end{proof}

We close this subsection with the following lemma stating that conditioning contracts inner products of stochastic integrals of non-negative functions:

\begin{lemm}\label{l:crossterms} 
    Fix a Borel set $B\in\cB(\R^2)$ and $n\in\Zpos$.
    For all $f,g\in L^{2}(\mathbb{R}^{2n})$ with $f,g\ge 0$ almost everywhere,
    \begin{align*}
        \mathbb{E}\left[I_{n}(f\1_{B}^{\otimes n})I_{n}(g\1_{B}^{\otimes n})\right] 
        &\le \mathbb{E}\left[I_{n}(f)I_{n}(g)\right].
    \end{align*}
\end{lemm}
\begin{proof}
    For $n=0$ we have $\1_B^{\otimes 0} = 1$, so the inequality is actually an equality.
    Assume $n\ge 1$.
    Since $\1_{B}^{\otimes n}$ is a symmetric function, the symmetrization of $f\1_{B}^{\otimes n}$ equals $\wt{f}\cdot \1_{B}^{\otimes n}$ (see \eqref{e:def-symmetrization}).
    Then by \cref{l:properties-of-I}\ref{orthogonality} and non-negativity of $\wt{f},\wt{g}$, we get
    \begin{align*}
        \mathbb{E}\left[I_{n}(f\1_{B}^{\otimes n})I_{n}(g\1_{B}^{\otimes n})\right]
        &=
        n!\langle \wt{f}\1_{B}^{\otimes n}, \wt{g}\1_{B}^{\otimes n}\rangle_{L^2(\R^{2n})} 
        \le 
        n!\langle \wt{f}, \wt{g}\rangle_{L^2(\R^{2n})} 
        =
        \mathbb{E}\left[I_{n}(f)I_{n}(g)\right]
    \end{align*}
    as desired.
\end{proof}

\subsection{Automorphisms of white noise}\label{s:automorphisms}

Due to its covariance structure, the law of white noise is invariant under unitary automorphisms of $L^2(\R^2)$.
This along with the CDRP's underlying Brownian structure (recall the formal expression \eqref{wienertilt}) implies that the law of the CDRP is covariant with respect to several natural transformations on path space, such as translation and diffusive scaling.
Before getting to the CDRP, in this short subsection we define the relevant automorphisms of the white noise space $(\Xi,\cF,\P)$.

Following \cite{AJRS22},
\footnote{The framework of \cite{AJRS22} includes additional families of automorphisms beyond those presented here, but for brevity we only recorded what we need for our analysis.}
we define the following families of invertible affine maps $\R^2\to\R^2$:
\begin{equation}\label{e:affine-maps}
    \begin{alignedat}{3}
        \transaff_{u,z}(t,x) &\coloneqq  (t,x) + (u,z) 
        &&\quad\text{for}\quad u,z\in\R;\\
        \shearaff_{r,\nu} (t,x) &\coloneqq  (t, x + \nu(t-r))
        &&\quad\text{for}\quad r,\nu\in\R;\\
        \dilaff_{\alpha,\lambda} (t,x) &\coloneqq  (\alpha t, \lambda x)
        &&\quad\text{for}\quad \alpha,\lambda>0.
    \end{alignedat}
\end{equation}
The above affine maps induce the following unitary automorphisms of $L^2(\R^2)$:
\begin{equation}\label{e:automorphisms}
    \begin{split}
        \trans_{u,z}\,f &\coloneqq  f \circ \transaff_{u,z}\\
        \shear_{r,\nu}\, f &\coloneqq  f\circ \shearaff_{r,\nu}\\
        \dil_{\alpha,\lambda}\, f &\coloneqq  \sqrt{\alpha \lambda}  f\circ \dilaff_{\alpha,\lambda}
    \end{split}.
\end{equation}
These in turn induce the following maps $(\Xi,\cF,\P)\to (\Xi,\cF,\P)$ 
(recall from \eqref{e:def-Xi} that $\Xi$ is the set of functions $L^2(\R^2)\to\R$):
\begin{equation}\label{e:automorphism-induced}
    \begin{split}
        \Trans_{u,z}\,\xi &\coloneqq  \xi \circ \trans_{u,z}^{-1}\\
        \Shear_{r,\nu}\, \xi &\coloneqq  \xi \circ \shear_{r,\nu}^{-1}\\
        \Dil_{\alpha,\lambda}\, \xi &\coloneqq  \xi \circ \dil_{\alpha,\lambda}^{-1}
    \end{split}.
\end{equation}
The maps in \eqref{e:automorphisms} are continuous bijections, so the maps in \eqref{e:automorphism-induced} are measurable bijections.
Their inverses are given by $\Trans_{u,z}^{-1} = \Trans_{-u,-z},\;\; \Shear_{r,\nu}^{-1} = \Shear_{r, -\nu}$, and $\Dil_{\alpha,\lambda}^{-1} = \Dil_{\alpha^{-1}, \lambda^{-1}}$.
Note that  $\Aut \xi$ is a proxy for ``$\aut\xi$,'' since if $\xi$ were actually an element of $L^2(\R^2)^* \cong L^2(\R^2)$, then for $f\in L^2(\R^2)$ we would have by unitarity $\langle \aut\xi, f\rangle_{L^2(\R^2)} = \langle \xi, \aut^{-1}f\rangle_{L^2(\R^2)}$.

In what follows we use $\autaff$ as a placeholder for any of the maps in \eqref{e:affine-maps}, and we write $\aut$ and $\Aut$ for the corresponding maps in \eqref{e:automorphisms} and \eqref{e:automorphism-induced} respectively.

\begin{lemm}\label{l:Q-measure-preserving}
    The map $\Aut$ is a measure-preserving automorphism of $(\Xi,\cF,\P)$, i.e. 
    a measure-preserving bijection with measure-preserving inverse.
\end{lemm}
\begin{proof}
We argued above that $\Aut$ is a measurable bijection.
We claim that $\Aut$ preserves finite-dimensional distributions.
For all $n\in\N$ and all $f_1,\dots,f_n\in L^2(\R^2)$, the vector 
$(\Aut\xi(f_1),\dots,\Aut\xi(f_n))$
is a centered multivariate Gaussian with covariances given by
\begin{align*}
    \E[\Aut\xi(f_i)\Aut\xi(f_j)]
    =
    \E[\xi(\aut^{-1}f_i)\xi(\aut^{-1}f_j)]
    =\langle \aut^{-1}f_i, \aut^{-1}f_j\rangle_{L^2(\R^2)}
    =\langle f_i,f_j\rangle_{L^2(\R^2)}
    =\E\left[
        \xi(f_i)\xi(f_j)
    \right],
\end{align*}
where we used the unitarity of $\aut^{-1}$.
Thus by Gaussianity, $(\Aut\xi(f_1),\dots,\Aut\xi(f_n))$ is equal in distribution to $(\xi(f_1),\dots,\xi(f_n))$.
Since finite-dimensional cylinder sets generate $\Fprod$, it follows
that $\Aut$ is a measure-preserving bijection $(\Xi, \Fprod,\P)\to(\Xi, \Fprod,\P)$, which then extends to a measure-preserving bijection of the completion $(\Xi,\cF,\P)\to (\Xi,\cF,\P)$.
Repeating this argument with $\Aut^{-1}$ in place of $\Aut$ completes the proof.
\end{proof}

\subsection{Continuum directed random polymer (CDRP)}
\label{ss:stocheatpolymer}

In this section we record the basics of the CDRP model as realized in \cite{AJRS22}.
An abbreviated version of this section appeared already in \cref{ss:CDRP}.

Fix $(s,x)\in\R^2$ and $\beta\ge 0$.
The $1+1$-dimensional stochastic heat equation (SHE) with multiplicative noise is the stochastic partial differential equation
\begin{align}\label{e:she}
    \partial_{t}\mathcal{Z}^{\mrm{SHE}}_{\beta}(s,x;t,y)
    = \frac{1}{2}\partial_{y}^2 \mathcal{Z}^\mrm{SHE}_{\beta}(s,x;t,y)
    + \beta\mathcal{Z}^\mrm{SHE}_{\beta}(s,x;t,y)\xi(t,y),
    \qquad t>s,\; y\in\R
    \tag{SHE}
\end{align}
with initial data $\lim_{t\downarrow s}\mathcal{Z}^{\mrm{SHE}}_{\beta}(s,x;t,y) =\delta(x-y)$.
The solution theory for \eqref{e:she} is classical, and we will not detail it here.
It is known that there exists a unique continuous and $\{\cF_{[s,t]\times\R}\}_{t\ge s}$-adapted
mild solution to \eqref{e:she}; we refer the reader to \cite[Appendix A]{AJRS22} (see especially Lemma A.3 therein) for the precise formulation of the mild problem and a review of the existence and uniqueness theory.
\footnote{
    Note that \cite[Lemma A.3]{AJRS22} is formulated in terms of the augmented filtration $\cF_{s,t}^W \coloneqq  \bigcap_{a<s,\;b>t}\cF_{[a,b]\times\R}$.
    However, by \cref{l:nested} we have $\cF_{s,t}^W = \cF_{[s,t]\times\R}$.
}

Let $\xi:(\Xi,\cF,\P)\to (\Xi,\Fprod)$ be the white noise on the canonical complete probability space defined in \cref{def:white-noise}.
By \cref{l:Q-measure-preserving}, $\xi$ and $(\Xi,\cF,\P)$ together satisfy the hypotheses assumed in \cite{AJRS22} (see the beginning of Section 2 therein) for their results to apply and hence we will freely do so.

The chaos expansion of the partition function will play a fundamental role in the coming arguments.
To prepare for this we introduce the following notation.
We denote the heat kernel by $p(s,x)\coloneqq \frac{1}{\sqrt{2\pi s}}e^{-x^2/2s}$
for $x\in\R$ and $s>0$.
For  $(s,x;t,y)\in\Rup$ and $n\in\N$, we define, for $(s_1,z_1;\cdots;s_n,z_n)\in \R^{2n}$,
\begin{align}\label{e:def-pn}
    p^{(n)}_{(s,x),(t,y)}(s_{1},z_{1};\cdots;s_{n},z_{n})
    \coloneqq   
    \mathbf{1}_{s<s_{1}<\cdots <s_{n}<t}
    \prod_{i=0}^{n}p(s_{i+1}-s_{i},z_{i+1}-z_{i}),
\end{align}
where $(s_{0},z_{0})\coloneqq (s,x)$ and $(s_{n+1},z_{n+1})\coloneqq (t,y)$.
We also denote $p^{(0)}_{(s,x),(t,y)} \coloneqq  p(t-s, y-x)$.
Thus $(z_1,\dots,z_n)\mapsto p^{(n)}_{(s,x),(t,y)}(s_1,z_1;\cdots;s_n,z_n)$ is (up to a normalization) the density of $(X(s_1),\dots,X(s_n))$, where $X$ is a Brownian bridge from $(s,x)$ to $(t,y)$.
For future reference, we record the following analogue of the It\^o isometry, which is an immediate consequence of \cref{l:properties-of-I}\ref{orthogonality} and the presence of the indicator $\1_{s<s_1<\cdots<s_n<t}$ in the definition of $p^{(n)}_{(s,x),(t,y)}$.
\begin{lemm}[It\^o isometry for Brownian multi-point densities]\label{l:ito-isometry}
    For every $s<t$ and $x,x',y,y'\in\R$,
    every $n\in\Z_{\ge1}$, and
    every $A\in \cB(\R^{2n})$, we have
    \begin{align*}
        \E\left[
            I_n\left(\1_A\, p^{(n)}_{(s,x),(t,y)}\right)
            I_n\left(\1_A\, p^{(n)}_{(s,x'),(t,y')}\right)
        \right]
        &=
        \left\langle \1_{A}\, p^{(n)}_{(s,x),(t,y)},\;
        \1_{A}\, p^{(n)}_{(s,x'),(t,y')}
        \right\rangle_{L^2(\R^{2n})}.
    \end{align*}
\end{lemm}

We proceed now to discussing the CDRP model, starting with the partition function.
The following result is a more comprehensive version of \cref{def:intro-SHE}.
This is quoted from \cite[Theorem 2.2]{AJRS22}, but some of their statements are weakened or omitted to suit our applications.

\begin{thm}[CDRP partition function, {\cite[Theorem 2.2]{AJRS22}}]\label{t:AJRAS-Z-existence}
    There exists a measurable map
    $\cZ : (\Xi,\cF,\P) \to C(\Rup\times\R_{\ge 0}, \R)$, denoted $\xi\mapsto \cZ^\xi_{\smallbullet}(\smallbullet,\smallbullet;\smallbullet,\smallbullet)$,
    such that
    \begin{enumerate}[label={\rm(\roman*)}]
        \item\textup{(Chaos expansion).}\label{property-chaos}        
        For each $(s,x;t,y)\in\Rup$ 
        and $\beta\ge 0$, 
        the random variable $\xi\mapsto \cZ_\beta^\xi(s,x;t,y)$ belongs to $L^2(\Xi,\cF,\P)$,
        and its chaos expansion (\cref{p:chaos-expansion}) is given by 
        \begin{align}\label{e:chaos-Z}
            \cZ_\beta^\xi(s,x;t,y) = \sum_{n=0}^\infty \beta^n I_n\left(
                p^{(n)}_{(s,x),(t,y)}
            \right),
        \end{align}
        where $p^{(n)}_{(s,x),(t,y)}$ is defined in \eqref{e:def-pn}.
        The above equality is as elements of $L^2(\Xi,\cF,\P)$.

        \item\textup{(Strict positivity).}\label{property-strict-positivity}
        For $\P$-a.e. $\xi\in\Xi$, we have $\cZ_\beta^\xi(s,x;t,y) > 0$
        for all $(s,x;t,y)\in\Rup$ and all $\beta\ge 0$.

        We henceforth replace $\cZ$ with a strictly positive version, i.e. we assume $\cZ_\beta^\xi(s,x;t,y)>0$ for \emph{all} $\xi\in\Xi$ and $(s,x;t,y)\in\Rup, \beta\ge0$.

        \item\textup{(Adapted).}\label{property-adapted} 
        For all $s<t$ and all $\beta\ge 0$, 
        the random function 
        $\xi\mapsto \cZ_\beta^\xi(s,\smallbullet;t,\smallbullet)$ is $\cF_{[s,t]\times\R}$-measurable.
        \footnote{
            As mentioned earlier, $\cF_{[s,t]\times\R}$ is indeed the same $\sigma$-algebra as that in \cite[Theorem 2.2(v)]{AJRS22}, by \cref{l:nested}.
        }  

        \item\textup{(Coupling of SHEs).}\label{property-solution-SHE}
        For each $(s,x)\in\R^2$ and $\beta\ge 0$, we have
        \begin{align*}
            \P\Bigl(
                \cZ_\beta^\xi(s,x;t,y) = \cZ^{\mathrm{SHE}}_{\beta}(s,x;t,y)
                \quad\textup{for all $t>s$ and all $y\in \R$}
            \Bigr)=1,
        \end{align*}
        where $\cZ^{\mathrm{SHE}}_\beta(s,x;\smallbullet,\smallbullet)$ is the unique (up to indistinguishability) continuous and $\{\cF_{[s,t]\times\R}\}_{t\ge s}$-adapted mild solution to \eqref{e:she} (see \cite[Appendix A]{AJRS22} for details).
    \end{enumerate}
\end{thm}

The following result on symmetries of the partition function was alluded to earlier in \cref{s:automorphisms}.

\begin{prop}[{Partition function symmetries, \cite[Proposition 2.3]{AJRS22}}]\label{t:AJRS-Z-scaling}
    Recall
    the families of measure-preserving automorphisms of $(\Xi,\cF,\P)$ defined in \eqref{e:automorphism-induced}.
    The following statements are true.
    \begin{enumerate}[label={\rm(\roman*)}]
        \item\textup{(Translation).}\label{property-translation}
        For each $(u,z)\in\R^2$, there exists an event $\mathsf{T}(u,z) \subset \Xi$ with 
        $\P(\mathsf{T}(u,z))=1$ such that on $\mathsf{T}(u,z)$, for all $(s,x;t,y)\in\Rup$ and all $\beta\ge0$,
        \begin{align*}
            \cZ_\beta(s+u,x+z;t+u,y+z) \circ \Trans_{-u, -z}
            = \cZ_\beta(s,x;t,y).
        \end{align*}
        \item\textup{(Shearing).}\label{property-shear}
        For each $(r,\nu)\in\R^2$, there exists an event $\mathsf{S}(r,\nu) \subset \Xi$ with $\P(\mathsf{S}(r,\nu))=1$ such that on $\mathsf{S}(r,\nu)$, for all $(s,x;t,y)\in\Rup$ and all $\beta\ge0$,
        \begin{align*}
            e^{\nu(y-x) + \frac{\nu^2}{2}(t-s)}\cZ_\beta(s,x+\nu(s-r); t, y+\nu(t-r))
            \circ \Shear_{r,-\nu}
            = \cZ_\beta(s,x;t,y).
        \end{align*}

        \item\textup{(Diffusive scaling).}\label{property-scaling}
            For each $\lambda>0$, there exists an event $\mathsf{D}(\lambda) \subset \Xi$ with $\P(\mathsf{D}(\lambda))=1$ such that on $\mathsf{D}(\lambda)$,
            for all $(s,x;t,y)\in\Rup$ and all $\beta\ge0$,
            \begin{align*}
                \lambda 
                \cZ_{\beta/\sqrt{\lambda}}(\lambda^2s,\lambda x; \lambda^2 t, \lambda y)
                \circ \Dil_{\lambda^{-2}, \lambda^{-1}} = 
                \cZ_\beta(s,x;t,y).
            \end{align*}
    \end{enumerate}
    In the above three equalities, $\cZ_\beta(s,x;t,y)$ is regarded as a real-valued function of the white noise, $\xi \mapsto \cZ^\xi_\beta(s,x;t,y)$.
    This is why we use the composition notation $\circ$.
\end{prop}

We turn now to the CDRP measure.
First let us recap the discussion from \cref{ss:CDRP}.
By \cref{t:AJRAS-polymer-existence}, for $\P$-a.e. realization of the white noise $\xi$, the following holds.
For every $(s,x;t,y)\in\Rup$ and every $\beta\ge 0$, there exists a unique probability measure $\P^\xi_{\beta,(s,x),(t,y)}$ on $C([s,t])$ with finite-dimensional distributions
\begin{align}\label{e:CDRP-fdd}
    \P^\xi_{\beta,(s,x),(t,y)}(X(s_1)\in dx_1,\dots, X(s_k)\in dx_k)
    &\coloneqq  \frac{\prod_{i=0}^k \cZ_\beta^\xi(s_i,x_i;s_{i+1},x_{i+1})}{\cZ_\beta^\xi(s,x;t,y)} dx_1 \dots dx_k,
\end{align}
where $(s_0,x_0)\coloneqq(s,x)$ and $(s_{k+1},x_{k+1})\coloneqq(t,y)$.
We use $X\sim \P^\xi_{\beta,(s,x),(t,y)}$ or $X\sim \CDRP^\xi_\beta(s,x;t,y)$ to mean that $X$ is sampled from $\P^\xi_{\beta,(s,x),(t,y)}$ (conditional on the white noise $\xi$).

Throughout the paper we will consider events on the white noise space defined in terms of the polymer measure.
To ensure that these are well-defined, we record for posterity the (intuitively obvious) fact that the polymer measure is itself a measurable function of the white noise $\xi$.

\begin{lemm}[Measurability of the polymer measure]\label{l:polymer-measure-is-measurable}
    For each $(s,x;t,y)\in\Rup$ and $\beta>0$, the map $\xi\mapsto \P_{\beta, (s,x),(t,y)}^\xi$ is a measurable map $(\Xi,\cF_{[s,t]\times \R},\P) \to \Prob(C([s,t]))$,
    where $\Prob(C([s,t]))$ is the space of Borel probability measures on $C([s,t])$ equipped with the weak topology and corresponding Borel $\sigma$-algebra.

    Equivalently, for every Borel set $S\in\cB(C([s,t]))$, the function $\xi\mapsto \P_{\beta,(s,x),(t,y)}^\xi(S)$ is $\cF_{[s,t]\times\R}$-measurable.
\end{lemm}
\begin{proof}
    That the first and second claims are equivalent follows from the fact that since $C([s,t])$ is a separable metric space, 
    the Borel $\sigma$-algebra on $\Prob(C([s,t]))$ coincides with the $\sigma$-algebra generated by the evaluation maps $\mu\mapsto \mu(S)$ for $S\in \cB(C([s,t]))$ (e.g. \cite[Lemma 2.3]{Varad63}).
    So it suffices to prove that $\xi\mapsto \P^\xi_{\beta,(s,x),(t,y)}(S)$ is measurable for all $S\in\cB(C([s,t]))$.
    We first handle the case when $S$ is a cylinder set, then extend to arbitrary $S$ using the $\pi$-$\lambda$ theorem.

    Let $S = \bigcap_{i=1}^n \{X(s_i) \in S_i\}$  for some $n\in\N$ and $s < s_1<\cdots<s_n < t$ and $S_1,\dots,S_n\in\cB(\R)$.
    Consider the map $\xi\mapsto K^\xi \in C(\R^n,(0,\infty))$ defined by
    \begin{align}
        K^\xi(x_1,\dots,x_n) \coloneqq  
        \frac{\prod_{i=0}^n
            \cZ^\xi_\beta(s_i,x_i; s_{i+1}, x_{i+1})}{\cZ^\xi_\beta(s,x;t,y)},
        \qquad (x_1,\dots,x_n)\in\R^n,
    \end{align}
    where $(s_0,x_0)\coloneqq (s,x)$ and $(s_{n+1},x_{n+1})\coloneqq (t,y)$, and where the strict positivity of the above map is by \cref{t:AJRAS-Z-existence}\ref{property-strict-positivity}.
    By \cref{t:AJRAS-Z-existence}\ref{property-adapted}, $\xi\mapsto K^\xi$ is $\cF_{[s,t]\times\R}$-measurable.
    Moreover, the functional $\Psi:C(\R^{n},(0,\infty)) \to [0,\infty]$ given by 
    \begin{align*}
        \Psi(G) \coloneqq  \int_{S_1\times\cdots\times S_n} G(x_1,\dots,x_n)\,dx_1\cdots dx_n
    \end{align*}
    is measurable, since the input $G$ is a positive function.
    It follows that
    \begin{align*}
        \xi\mapsto \Psi(K^\xi) 
        = \int_{S_1\times\cdots\times S_n}
        \frac{\prod_{i=0}^n \cZ_\beta^\xi(s_i,x_i;s_{i+1},x_{i+1})}{\cZ^\xi_\beta(s,x;t,y)}
        \,dx_1\cdots dx_n
        = \P^\xi_{\beta,(s,x),(t,y)}(S)
    \end{align*} 
    is $\cF_{[s,t]\times\R}$-measurable.

    Having shown that $\xi\mapsto \P^\xi_{\beta,(s,x),(t,y)}(S)$ is measurable for any cylinder set $S$, we now deduce the same statement for arbitrary $S\in \cB(C([s,t]))$.
    We abbreviate $\P^\xi\coloneqq \P^\xi_{\beta,(s,x),(t,y)}$ and define
    \begin{align}\label{e:pi-lambda}
        \mathscr{L} \coloneqq  \left\{
            S\in\cB(C([s,t])) : \xi\mapsto \P^\xi(S)
            \text{ is $\cF_{[s,t]\times\R}$-measurable}
        \right\}.
    \end{align}
    Since $\cB(C([s,t]))$ is generated by cylinder sets,
    by the $\pi$-$\lambda$ theorem it suffices to show that $\mathscr{L}$ is a $\lambda$-system.
    We have $C([s,t])\in\mathscr{L}$ since $C([s,t])$ is a cylinder set.
    For $S,S'\in \mathscr{L}$ with $S\subset S'$, we have $\P^\xi(S'\setminus S) = \P^\xi(S') - \P^\xi(S)$,
    which is the difference of two $\cF_{[s,t]\times \R}$-measurable functions and hence is $\cF_{[s,t]\times \R}$-measurable.
    So $S'\setminus S \in \mathscr{L}$. 
    Next, for $S_1,S_2,\dots\in\cB(C([s,t]))$ with $S_1\subset S_2\subset\cdots$, we have $\P^\xi(\bigcup_{n=1}^\infty S_n) = \lim_{n\to\infty}\P^\xi(S_n)$, and a limit of $\cF_{[s,t]\times\R}$-measurable functions is $\cF_{[s,t]\times\R}$-measurable by completeness.
    So $\bigcup_{n=1}^\infty S_n\in \mathscr{L}$.
    This shows that $\mathscr{L}$ is a $\lambda$-system, completing the proof of \cref{l:polymer-measure-is-measurable}.
\end{proof}

The following lemma is the CDRP analogue of \cref{t:AJRS-Z-scaling} alluded to in \cref{s:automorphisms}.
\begin{prop}[Quenched CDRP symmetries]\label{p:CDRP-scaling}
    There exists an event $\Xi_0\subset \Xi$ with $\P(\Xi_0)=1$ such that
    \begin{enumerate}[label={\rm(\roman*)}]
        \item\textup{(Translation).}\label{property-CDRP-translation} For each $(u,z) \in \R^2$, on the event $\Xi_0\cap\mathsf{T}(u,z)$, it holds for all $(s,x;t,y)\in\Rup$, all $\beta\ge 0$, and all Borel sets $S\in\cB(C([s,t]))$ that
        \begin{align*}
            \P_{\beta,(s+u, x+z),(t+u, y+z)}\left(
                X(\cdot + u) - z
                \in S\right)
            \circ \Trans_{-u, -z}
            =
            \P_{\beta,(s,x),(t,y)}(X\in S).
        \end{align*}
        \item\textup{(Shearing).}\label{property-CDRP-shear} For each $(r,\nu)\in\R^2$, on the event $\Xi_0\cap\mathsf{S}(r,\nu)$, it holds for all $(s,x;t,y)\in\Rup$, all $\beta\ge 0$, and all Borel sets $S\in\cB(C([s,t]))$ that
        \begin{align*}
            \P_{\beta,(s, x+\nu(s-r)),(t, y+\nu(t-r))}\left(
                X(\cdot) - \nu(\cdot - r)
                \in S
            \right)
            \circ \Shear_{r, -\nu}
            =
            \P_{\beta,(s,x),(t,y)}(X\in S).
        \end{align*}
        \item\textup{(Diffusive scaling).}\label{property-CDRP-scaling}
        For each $\lambda>0$, on the event $\Xi_0\cap \mathsf{D}(\lambda)$, it holds for all $(s,x;t,y)\in\Rup$, all $\beta\ge 0$, and all Borel sets $S\in\cB(C([s,t]))$ that
        \begin{align*}
            \P_{\beta\lambda^{-1/2}, (\lambda^2 s, \lambda x),(\lambda^2 t,\lambda y)}\left(
                \lambda^{-1}X(\lambda^2\cdot)
                \in S
            \right)
            \circ \Dil_{\lambda^{-2},\lambda^{-1}}
            =
            \P_{\beta,(s,x),(t,y)}\left(X\in S\right).
        \end{align*}
    \end{enumerate}
    The events $\mathsf{T}(u,z), \mathsf{S}(r,\nu), \mathsf{D}(\lambda)$ are defined in \cref{t:AJRS-Z-scaling}.
    Like in \cref{t:AJRS-Z-scaling}, the composition notation $\circ$ reflects that the above polymer probabilities are regarded as real-valued functions of the white noise.
 \end{prop}
 \begin{proof} 
    Let $\Xi_0$ be the event in \cref{t:AJRAS-polymer-existence} on which all the polymer measures exist. 
    Since $\cB(C([s,t]))$ is generated by finite-dimensional cylinder sets, we just need to check that the measures in \ref{property-CDRP-translation}--\ref{property-CDRP-scaling} have the same finite-dimensional distributions.
    We do this using \cref{t:AJRS-Z-scaling} along with \eqref{e:CDRP-fdd}.

    Fix $(s,x;t,y)\in\Rup$ and $\beta\ge 0$.
    Fix $k\in\N$ and $s<s_1<\cdots<s_k<t$, and write $(s_0,x_0)\coloneqq (s,x)$ and $(s_{k+1},x_{k+1})\coloneqq (t,y)$.
    
    \ref{property-CDRP-translation}:
    Write $\wt{X}(\cdot) \coloneqq  X(\cdot + u) - z$.
    Then by \eqref{e:CDRP-fdd} and \cref{t:AJRS-Z-scaling}\ref{property-translation},
    \begin{multline*}
        \P_{\beta,(s+u, x+z),(t+u, y+z)}\left(\wt{X}(s_1)\in dx_1,\dots, \wt{X}(s_k) \in dx_k\right)
        \circ \Trans_{-u,-z}\\
        \begin{aligned}
            \overset{\eqref{e:CDRP-fdd}}&{=}
            \frac{\prod_{i=0}^k\cZ_\beta(s_i+u, x_i+z; s_{i+1}+u, x_{i+1}+z)}{\cZ_\beta(s+u, x+z; t+u, y+z)}
            \circ \Trans_{-u, -z}\\
            &= \frac{\prod_{i=0}^k \cZ_\beta(s_i,x_i; s_{i+1},x_{i+1})}{\cZ_\beta(s,x;t,y)}
            \overset{\eqref{e:CDRP-fdd}}{=}
            \P_{\beta,(s,x),(t,y)}(X(s_1)\in dx_1,\dots, X(s_k)\in dx_k).
        \end{aligned}
    \end{multline*}

    \ref{property-CDRP-shear}:
    Write $\wt{X}(\cdot) \coloneqq  X(\cdot) - \nu(\cdot-r)$.
    Then by \eqref{e:CDRP-fdd} and \cref{t:AJRS-Z-scaling}\ref{property-shear},
    \begin{multline*}
        \P_{\beta,(s, x+\nu(s-r)),(t, y+\nu(t-r))}\left(\wt{X}(s_1)\in dx_1,\dots, \wt{X}(s_k) \in dx_k\right)
        \circ \Shear_{r,-\nu}
        \\
        \begin{aligned}
        \overset{\eqref{e:CDRP-fdd}}&{=}
        \frac{\prod_{i=0}^k\cZ_\beta(s_i, x_i+\nu(s_i-r); s_{i+1}, x_{i+1}+\nu(s_{i+1}-r))}
        {\cZ_\beta(s, x+\nu(s-r); t, y+\nu(t-r))}
        \circ \Shear_{r,-\nu}\\
        &=\frac{\prod_{i=0}^k \cZ_\beta(s_i, x_i; s_{i+1}, x_{i+1})}
        {\cZ_\beta(s, x; t, y)}
        \cdot \frac{\prod_{i=0}^k e^{-\nu(x_{i+1}-x_i) - \frac{\nu^2}{2}(s_{i+1}-s_i)}}{e^{-\nu(y-x) - \frac{\nu^2}{2}(t-s)}}\\
        &=\frac{\prod_{i=0}^k \cZ_\beta(s_i, x_i; s_{i+1}, x_{i+1})}
        {\cZ_\beta(s, x; t, y)}
        \overset{\eqref{e:CDRP-fdd}}{=}
        \P_{\beta,(s,x),(t,y)}(X(s_1)\in dx_1,\dots, X(s_k)\in dx_k).
        \end{aligned}
    \end{multline*}

    \ref{property-CDRP-scaling}:
    Write $\wt{X}(\cdot) \coloneqq  \lambda^{-1}X(\lambda^2 \cdot)$.
    Then by \eqref{e:CDRP-fdd} and \cref{t:AJRS-Z-scaling}\ref{property-scaling},
    \begin{multline*}
        \P_{\beta\lambda^{-1/2}, (\lambda^2 s, \lambda x),(\lambda^2 t,\lambda y)}\left(
                \wt{X}(s_1)\in dx_1,\dots, \wt{X}(s_k) \in dx_k
            \right)
            \circ \Dil_{\lambda^{-2},\lambda^{-1}}
        \\
        \begin{aligned}
        \overset{\eqref{e:CDRP-fdd}}&{=}
        \lambda^{k}
        \frac{\prod_{i=0}^k \cZ_{\beta\lambda^{-1/2}}(\lambda^2 s_i, \lambda x_i; \lambda^2 s_{i+1}, \lambda x_{i+1})}
        {\cZ_{\beta\lambda^{-1/2}}(\lambda^2 s, \lambda x;\lambda^2 t,\lambda y)}
        \circ \Dil_{\lambda^{-2},\lambda^{-1}} \\
        &=
        \lambda^k
        \frac{\prod_{i=0}^k \lambda^{-1} \cZ_\beta(s_i,x_i;s_{i+1},x_{i+1})}
        {\lambda^{-1} \cZ_\beta(s,x;t,y)}\\
        &=\frac{\prod_{i=0}^k\cZ_\beta(s_i,x_i;s_{i+1},x_{i+1})}
        {\cZ_\beta(s,x;t,y)}
        \overset{\eqref{e:CDRP-fdd}}{=}
        \P_{\beta,(s,x),(t,y)}(X(s_1)\in dx_1,\dots, X(s_k)\in dx_k).
        \end{aligned}
    \end{multline*}
    We are done.
\end{proof}

The next result is (a special case of) the Markov property for the polymer.
We will only need the case $\beta=1$, but a similar statement holds for any $\beta\ge 0$.
\begin{prop}[{Markov property, \cite[Theorem 2.14(iii)]{AJRS22}}]\label{p:markov-property-CDRP}
    Fix $(s,x;t,y)\in\Rup$ and fix $s',t'$ satisfying $s<s'<t'<t$.
    Define the following $\sigma$-algebras on path space $C([s,t])$:
    \begin{align*}
        \cG_{s,s'} \coloneqq  \sigma\left(
            X|_{[s,s']}:X\in C([s,t])
        \right),
        \qquad
        \cG_{t',t} \coloneqq  \sigma\left(
            X|_{[t',t]}:X\in C([s,t])
        \right).
    \end{align*}
    Then almost surely,
    for all bounded  measurable functions $F:C([s',t'])\to\R$,  we have
    \begin{align*}
        \E^\xi_{1,(s,x),(t,y)}\left[
            F(X|_{[s',t']})\,\middle|\,\cG_{s,s'},\, \cG_{t',t}
        \right]
        &=
        \E^\xi_{1,(s',X(s')), (t', X(t'))}\left[
            F(X)
        \right]
        \qquad\quad
        \textup{$\P^\xi_{1,(s,x),(t,y)}$-a.s.,}
    \end{align*} 
    where $\E^\xi_{1,(s,x),(t,y)}$ denotes expectation with respect to $\P^\xi_{1,(s,x),(t,y)}$.
\end{prop}

The final result in this subsection states that pairs of polymers with ordered endpoints are stochastically ordered.
This phenomenon, known as \emph{polymer ordering} in the literature, is a consequence of the planarity of the model. This has primarily featured in the literature in the context of zero-temperature geodesics in last passage percolation, e.g., \cite[Lemma 11.2]{BSS14}.
As in the previous result, we will only need the case $\beta=1$.
\begin{prop}[{Polymer ordering, \cite[Proposition 2.18]{AJRS22}}]\label{p:polymer-ordering}
    Almost surely, 
    for all $s<t$ and all $x_{1}\leq x_{2}$ and $y_{1}\leq y_{2}$,
    we have the stochastic dominance  
    \footnote{
        For probability measures $\mu,\nu$ on $C([s,t])$, we say that $\nu$ stochastically dominates $\mu$ (denoted $\mu\preceq_{\mathrm{sd}}\nu$) if $\int F\,d\mu \le \int F\,d\nu$ for all continuous bounded increasing functionals $F: C([s,t])\to\R$, where $F : C([s,t]) \to \R$ is called increasing if for all $X,Y\in C([s,t])$ satisfying $X(r)\le Y(r)$ for all $r\in[s,t]$, we have $F(X)\le F(Y)$.
    }
    \begin{align*}
        \P^\xi_{1,(s,x_1),(t,y_1)}
        \preceq_{\mathrm{sd}}
        \P^\xi_{1,(s,x_2),(t,y_2)}.
    \end{align*}
\end{prop}

Propositions \ref{p:markov-property-CDRP} and \ref{p:polymer-ordering} are only used in the proof of \cref{l:cond-reg} in \cref{s:es}.

\subsection{Restricting the polymer}\label{s:restricting-polymer}

In this subsection we prove \cref{p:Z-measurable-real}.
Let us recall the relevant notation.
Throughout this subsection we fix $(s,x;t,y)\in\Rup$ and $s'<t'$ with $[s',t']\subset [s,t]$, as well as $\beta>0$.
For $a<b$, denote
\begin{align}
    B(a,b) \coloneqq  [s',t']\times [a,b].
\end{align}
Fix any Borel set $S\in \cB(C([s,t]))$.
Recall from \eqref{e:def-ZB} that given a realization of the white noise $\xi\in \Xi$, we write
\begin{align*}
    \cZ_{B(a,b)}^S &= \cZ_{B(a,b)}^S(s,x;t,y)\\
    &\coloneqq  
    \cZ_\beta^\xi(s,x;t,y)\cdot \P_{\beta,(s,x),(t,y)}^\xi\left(
        S\cap\{X(r)\in [a,b]\text{ for all $r\in[s',t']$}\}
    \right).
\end{align*}
\cref{p:Z-measurable-real} asserts that $\cZ^S_{B(a,b)}$ is measurable with respect to 
\begin{align*}
    \cF_{(([s,s']\cup[t',t])\times \R) \cup B(a,b)}
    =
    \cF_{([s,s']\cup[t',t])\times \R}\vee \cF_{B(a,b)},
\end{align*}
where the equality is by \cref{l:independence-disjoint}\ref{generation}.

Before getting into the proof of \cref{p:Z-measurable-real}, we first deduce a simple corollary that will be useful for our later applications.
\begin{cor}\label{l:meas3}
    Fix $(s,x;t,y)\in\Rup$ and $s'<t'$ with $[s',t']\subset[s,t]$.
    For any Borel set $S\in \cB(C([s,t]))$ and any $z<w$, define
    \begin{align*}
        \cZ^S_{B(z,w)^c} &= \cZ^S_{B(z,w)^c}(s,x;t,y)\\
        &\coloneqq  \cZ_\beta^\xi(s,x;t,y)\cdot \P_{\beta,(s,x),(t,y)}^\xi\left(
            S\cap \{X(r)\in [z,w]^c\text{ for all }r\in[s',t']\}
        \right).
    \end{align*}
    Then $\cZ^S_{B(z,w)^c}$ is $\cF_{B(z,w)^c}$-measurable.
\end{cor}
\begin{proof}[Proof of \cref{l:meas3}]
    We suppress the dependence on $(s,x;t,y)$ from the notation.
    By path continuity of the polymer and upwards continuity of measures, we have
    \begin{multline*}
        \cZ_\beta^\xi \P^\xi_\beta(S\cap\{
            X(r)\in [z,w]^c\text{ for all }r\in[s',t']
            \})
        \\
        \begin{aligned} 
            &=\cZ_\beta^\xi \P^\xi_\beta(S\cap\{
                X(r)\in (-\infty,z)\text{ for all }r\in[s',t']
            \})\\
            &\quad+ \cZ_\beta^\xi\P^\xi_\beta(S\cap\{
                X(r)\in (w,\infty)\text{ for all }r\in[s',t']
            \})\\
            &= 
            \lim_{n\to \infty}
            \cZ_\beta^\xi\P^\xi_\beta(S\cap\{
                X(r)\in [-n, z- \tfrac{1}{n}]\text{ for all }r\in[s',t']
            \})\\
            &\quad+ 
            \lim_{n\to \infty}
            \cZ_\beta^\xi
            \P^\xi_\beta(S\cap\{
                X(r)\in [w+\tfrac{1}{n}, n]\text{ for all }r\in[s',t']
            \}).
        \end{aligned}
    \end{multline*}
    By \cref{p:Z-measurable-real},
    for fixed $n$, the quantity inside the first limit is measurable with respect to
    \begin{align*}
            \cF_{([s,s']\cup[t',t])\times\R} \vee \cF_{[s',t']\times[-n,z-\frac{1}{n}]}
            \quad\subset\quad
            \cF_{([s,s']\cup[t',t])\times\R} \vee \cF_{[s',t']\times(-\infty, z)}.
    \end{align*}
    It follows that the limit is also measurable with respect to $\cF_{([s,s']\cup[t',t])\times\R} \vee \cF_{[s',t']\times(-\infty, z)}$, since the above $\sigma$-algebras contain all $\P$-null sets in $\cF$ by definition (recall \eqref{e:def-FB}).
    The same argument implies that the second limit is measurable with respect to
    $\cF_{([s,s']\cup[t',t])\times\R} \vee \cF_{[s',t']\times(w, \infty)}$.
    Therefore, by \cref{l:independence-disjoint}\ref{generation}, the sum is measurable with respect to
    \begin{align*}
        \left(\cF_{([s,s']\cup[t',t])\times\R} \vee \cF_{[s',t']\times(-\infty, z)}\right)\vee 
        \left(
            \cF_{([s,s']\cup[t',t])\times\R} \vee \cF_{[s',t']\times(w,\infty)}
        \right)
        &=
        \cF_{([s,s']\cup[t',t])\times\R} \vee \cF_{[s',t']\times[z,w]^c}\\
        &\subset 
        \cF_{B(z,w)^c}\,.
    \end{align*}
    This completes the proof.
\end{proof}

\begin{figure}[htb]
    \centering
    \includegraphics[width=0.7\textwidth]{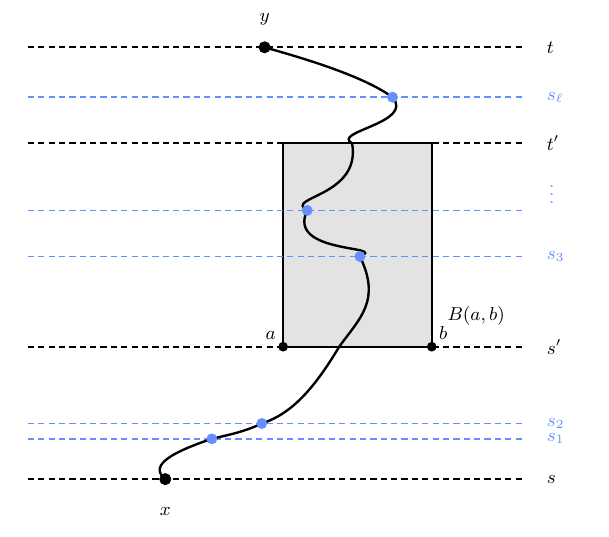}
\caption{
    The setting of \cref{p:ZB-measurable}.
    A continuous path from $(s,x)$ to $(t,y)$ is drawn in black.
    The path belongs to some fixed cylinder set $S = \bigcap_{j=1}^\ell\{X(s_j) \in S_j\} \subset C([s,t])$ (this constrains the horizontal locations of the blue points, but these constraints are not depicted),
    and it stays inside the shaded box $B(a,b)=[s',t']\times[a,b]$ throughout the time interval $[s',t']\subset [s,t]$.
    \cref{p:ZB-measurable} asserts that $\cZ_B^S$---the total contribution made by such paths to the polymer partition function $\cZ_\beta(s,x;t,y)$---is a measurable function of the white noise inside $B(a,b)$ and the infinite horizontal strips $[s,s']\times\R$ and $[t',t]\times\R$.
    }\label{fig:meas-1}
\end{figure}

We turn to the proof of \cref{p:Z-measurable-real}.
From now on we fix $a<b$ and abbreviate $B\coloneqq B(a,b)$.
We also abbreviate $\cZ_B^S\coloneqq \cZ_B^S(s,x;t,y)$.
We will first prove \cref{p:Z-measurable-real} for cylinder sets $S$, and then deduce the full result via a straightforward $\pi$-$\lambda$ argument.
For clarity, we record the first step in the following proposition whose statement is also depicted in Figure \ref{fig:meas-1}.

\begin{prop}\label{p:ZB-measurable}
    Fix any $\ell\in\Z_{\ge 1}$ and any cylinder set $S\in\cB(C([s,t]))$ of the form
    \begin{align}\label{e:cylinder-goal}
        S = \bigcap_{j=1}^\ell
        \{X(s_j) \in S_j\},
        \qquad 
        s_1,\dots,s_\ell \in (s,t)\setminus\{s',t'\},
        \quad s_1<\cdots<s_\ell,
        \quad S_1,\dots,S_\ell \in \cB(\R).
    \end{align}
    Then 
    $\cZ_{B(a,b)}^S$ is measurable with respect to
    $\cF_{([s,s']\cup[t',t])\times\R} \vee \cF_{B(a,b)}$.
\end{prop}

As mentioned, \cref{p:Z-measurable-real} follows from \cref{p:ZB-measurable} by a simple $\pi$-$\lambda$ argument:

\begin{proof}[Proof of \cref{p:Z-measurable-real}]
    We abbreviate  $\cF_0\coloneqq  \cF_{([s,s']\cup[t',t])\times\R} \vee \cF_{B(a,b)}$ and $B\coloneqq B(a,b)$.
    It is straightforward to check that the collection of cylinder sets of the form \eqref{e:cylinder-goal} forms a $\pi$-system that generates the Borel $\sigma$-algebra $\cB(C([s,t]))$ (e.g. \cite[Example 1.3]{bil}).
    By \cref{p:ZB-measurable}, the cylinder sets are contained in the collection
    \begin{align*}
        \mathscr{L} 
        \coloneqq  \left\{
            S \in \cB(C([s,t])) : \cZ_{B}^S \text{ is $\cF_0$-measurable}
        \right\}.
    \end{align*}
    By the $\pi$-$\lambda$ theorem, it now suffices to show that $\mathscr{L}$ forms a $\lambda$-system.
    A similar statement was established earlier in the proof of \cref{l:polymer-measure-is-measurable} (see below \eqref{e:pi-lambda}), and the same argument goes through essentially verbatim here.
    We omit the details.
\end{proof}

The rest of this subsection is devoted to proving \cref{p:ZB-measurable}, which will take several pages.
We begin with an informal overview of the proof strategy.

\subsubsection{Proof strategy}\label{sss:meas-proof-overview}
Write $B\coloneqq B(a,b)$ and $\cF_0 \coloneqq  \cF_{([s,s']\cup[t',t])\times\R} \vee \cF_{B}.$
Since the CDRP is defined in terms of finite-dimensional distributions, the condition that the polymer lies inside $B$ throughout the time interval $[s',t']$ is somewhat unwieldy.
We will therefore work with a finite-dimensional approximation of $\cZ^S_B$ given by discretizing the time coordinate.
Namely, for $\d>0$, we will consider a proxy $\cZ^S_{B,\d}$ given (roughly) by the polymer partition function restricted to Brownian paths $X\in S$ that lie inside $[a,b]$ at times $s',s'+\d, \dots,t'-\d, t'$ (these times correspond roughly to the horizontal dashed orange lines in Figure \ref{fig:meas-2} below).
Since $S$ is a cylinder set, $\cZ^S_{B,\d}$ is just an integral of a CDRP finite-dimensional density,
and it converges to $\cZ^S_B$ as $\d\to 0$.
So it suffices to prove that $\cZ^S_{B,\d}$ is ``asymptotically $\cF_0$-measurable'' as $\d\to0$.
To accomplish this, we make a further approximation.
Fixing $\e>0$, let $B_\e$ be the box given by enlarging $B$ by $\e$ in the spatial direction, i.e.
\begin{align}\label{e:def-B-eps}
    B_\e \coloneqq  [s',t']\times[a-\e,b+\e].
\end{align}
The box $B_\e$ is depicted in Figure \ref{fig:meas-2} below.
Define the $\sigma$-algebra
\begin{align}\label{e:def-F-eps}
    \cF_\e \coloneqq  \cF_{(([s,s']\cup[t',t])\times\R)\cup B_\e}.
\end{align}
Since $\cF_\e\downarrow \cF_0$ as $\e\to 0$ (\cref{l:nested}), it suffices to prove that $\cZ^S_{B}$ is $\cF_\e$-measurable for any fixed $\e>0$.
Moreover, since $\cZ^S_{B,\d}\to \cZ^S_B$ as $\d\to 0$, it suffices to prove ``asymptotic $\cF_\e$-measurability'' of $\cZ^S_{B,\d}$ as $\d\to0$.
To make this precise, we decompose $\cZ^S_{B,\d}$ as
\begin{align}\label{e:290}
    \cZ^S_{B,\d}
    = \E[\cZ^S_{B,\d}\mid \cF_{\e}]
    + \left(\cZ^S_{B,\d} - \E[\cZ^S_{B,\d}\mid \cF_{\e}]\right).
\end{align}
The first term is $\cF_\e$-measurable by definition, and we will show that the second term converges to $0$ in probability as $\d\to 0$. 
Let us heuristically explain this convergence.
Intuitively, $\cZ^S_{B,\d} - \E[\cZ^S_{B,\d}\mid \cF_{\e}]$ can be thought of as measuring the contribution to $\cZ^S_{B,\d}$ made by Brownian paths that exit $[a-\e,b+\e]$ during the time interval $[s',t']$.
By definition of $\cZ^S_{B,\d}$, the Brownian paths contributing to $\cZ^S_{B,\d}$ are required to be inside $[a,b]$ at the times $s'+\d, s'+2\d,\dots,t'-\d,t'$,
which by diffusivity implies that all but a vanishing fraction of them stay inside $[a - \d^{0.49}, b + \d^{0.49}]$ throughout the time interval $[s',t']$. 
Thus when $\d\to 0$ and in particular $\d^{0.49}\ll\e$, the contribution to $\cZ^S_{B,\d}$ from paths that exit $[a-\e,b+\e]$ is negligible, simply because these paths are negligible under Brownian measure.

In the actual proof of \cref{p:ZB-measurable}, we will show that the error term in \eqref{e:290} converges to $0$ in $L^2$ using the Wiener chaos machinery developed in the previous subsections.
The key estimate (\cref{l:moment-of-chaos-Bc}) is indeed a chaos-level version of the above Brownian diffusivity heuristic.
With this broad strategy in mind, we now present a more detailed outline of the proof.
\begin{itemize}
    \item In \cref{sss:defining-ZBdelta} we define the proxy $\cZ^S_{B,\d}$ and establish some of its basic properties, including that $\cZ^S_{B,\d}\to\cZ^S_B$ as $\d\to 0$.
    We also set up notation for the finite-dimensional CDRP density whose integral equals $\cZ^S_{B,\d}$.

    \item In \cref{sss:proof-modulo-chaos} we implement a chaos-level version of the above Brownian heuristic and complete the proof of \cref{p:ZB-measurable} (modulo some calculations that we defer to \cref{sss:L2-computations}).
    To do this, we study how the chaos expansion of $\cZ^S_{B,\d}$ transforms under orthogonal projections in $L^2(\Xi,\cF,\P)$ (note that the second term on the RHS of \eqref{e:290} is the orthogonal projection onto $L^2(\Xi,\cF_\e,\P)^\perp$).
    \item In \cref{sss:L2-computations} we supply the calculations required to make the argument of \cref{sss:proof-modulo-chaos} fully rigorous. 
\end{itemize}
While the above proof strategy is natural, making it rigorous  demands lengthy calculations with the CDRP finite-dimensional density \eqref{e:CDRP-fdd} and the chaos expansion \eqref{e:chaos-Z}. Unfortunately, these calculations are unavoidable, and they make the forthcoming proofs technically involved and somewhat difficult to parse.

\addtocontents{toc}{\SkipTocEntry}
\subsection*{Proof of Proposition \ref{p:ZB-measurable}}

From now on we fix a cylinder set $S=\bigcap_{j=1}^\ell \{X(s_j)\in S_j\}$ as in \eqref{e:cylinder-goal}.
Note that if $s_j\in (s',t')$, then we have the following equality of subsets of $C([s,t])$:
\begin{multline}\label{e:WLOG-subset-ab}
    \{X(s_j)\in S_j\}\cap
        \{X(r)\in[a,b]\text{ for all } r\in [s',t']\}\\
    = \{X(s_j)\in S_j\cap [a,b]\}\cap
        \{X(r)\in[a,b]\text{ for all } r\in [s',t']\}.
\end{multline}
So we can assume that $S_j\subset [a,b]$ whenever $s_j\in (s',t')$.

For notational convenience we assume that $\{s_1,\dots,s_\ell\}\cap [s,s']\ne\varnothing$ (we can always intersect $S$ with $\{X(s)\in \R\}=C([s,t])$), and similarly that $\{s_1,\dots,s_\ell\}\cap [t',t]\ne \varnothing$.

\subsubsection{Finite-dimensional proxy}\label{sss:defining-ZBdelta}

\begin{figure}[htbp]
    \centering
    \includegraphics[width=0.958\linewidth]{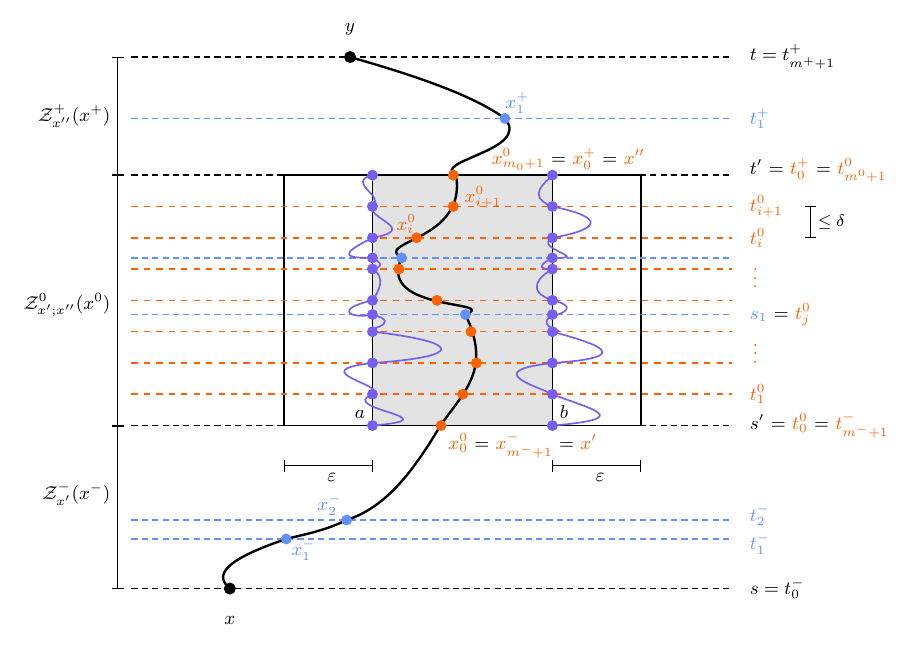}
\caption{
    The proof of \cref{p:ZB-measurable}.
    Several times are depicted as horizontal dashed lines:
    black lines for $s\le s'<t'\le t$ (drawn here as $s<s'<t'<t$), blue lines for $s_1<\cdots<s_\ell$ from \cref{p:ZB-measurable}, and orange lines for the  $\d$-mesh $R_\d\coloneqq (s',t')\cap\d\Z$.
    The blue/orange times are collected into three groups, depending on whether they land in $(s,s'), (s',t'),$ or $(t',t)$.
    The times in $(s,s')$ are denoted $t_1^- < \cdots < t_{m^-}^-$, the times in $(s',t')$ are denoted $t_1^0<\cdots<t_{m^0}^0$, and the times in $(t',t)$ are denoted $t_1^+<\cdots<t_{m^+}^+$ (as drawn, $m^- = 2$, $m^0 = 2$, and $m^+=1$).
    We also write $t^-_0\coloneqq s$ and $t^+_{m^+ +1} \coloneqq  t$, and $t_0^0 = t^-_{m^- +1} \coloneqq  s'$ and $t^0_{m^0+1} = t^+_0 \coloneqq  t'$.
    A similar notational convention is used for the spatial coordinates (horizontal).
    \\[0.8\baselineskip]
    Drawn in black is a Brownian path from $(s,x)$ to $(t,y)$.
    It belongs to the cylinder set $S$, and it lands in the shaded box $B=B(a,b)=[s',t']\times[a,b]$ at times $t^0_0,\dots,t^0_{m^0+1}$.
    We study the contribution $\cZ_{B,\d}^S$ to the polymer partition function from all such paths.
    In \eqref{e:ZBdelta-final} we express $\cZ_{B,\d}^S$ as an integral over the dashed horizontal lines of a product of three (un-normalized) CDRP finite-dimensional densities $\cZ^-_{x'}(\bfx^-), \cZ^0_{x';x''}(\bfx^0),$ and $\cZ^+_{x''}(\bfx^+)$, corresponding respectively to the polymer's behavior before, during, and after the time interval $[s',t']$.
    By directedness, these three densities depend on disjoint parts of the white noise and hence are independent.
    \\[0.8\baselineskip]
    Since each time interval $[t_i^0, t_{i+1}^0]$ has length $\le\d$, for a Brownian bridge on the same time interval with endpoints in $B$, 
    it exits $B_\e$ with probability $O(\exp(-c \e^2/\d))$.
    Therefore, by a union bound, the concatenation of $O(\d^{-3/2})$-many such paths stays inside $B_\e$ with probability $1-o(1)$ as $\d\to 0$, as exhibited by the purple paths traveling along the boundary of $B$.
    }\label{fig:meas-2}
\end{figure}

As indicated in the proof sketch in \cref{sss:meas-proof-overview}, we begin by approximating $\cZ_B^S$ in terms of the polymer's behavior on a fine mesh of times.
To this end, 
given $\d>0$, we define the following $\d$-mesh (depicted as orange dashed lines in Figure \ref{fig:meas-2}):
\begin{align*}
    R_{\d} = R_{\d}(s',t') 
    &\coloneqq  
    (s',t')\cap\d\Z
\end{align*}
and define 
\begin{align}\label{e:def-ZBdelta}
    \cZ_{B,\d}^S
    &\coloneqq  \cZ_\beta^\xi(s,x;t,y) 
    \cdot
    \P_{\beta,(s,x),(t,y)}^\xi\left(
        S\cap\{X(r)\in[a,b]\text{ for all } r\in R_\d
        \cup\{s',t'\}\}
    \right).
\end{align}

\begin{lemm}\label{l:partitionalmostsure}
    For $h\in\Z_{\ge 1}$ let $\d_h\coloneqq 2^{-h}$.
    Then $\cZ_{B,\d_h}^S \xrightarrow{h\to\infty} \cZ_B^S$ almost surely.
\end{lemm}
\begin{proof}
    Fix a realization of the white noise $\xi$.
    Since 
    $R_{\d_1}\subset R_{\d_2}\subset\cdots$ and $\bigcup_{h=1}^\infty R_{\d_h}$ is dense in $[s',t']$,
    it follows by path continuity of the polymer that
    \begin{align*}
        \P^\xi_\beta
        \left(S\cap \{X(r)\in [a,b]\text{ for all }r\in 
        R_{\d_h}\cup\{s',t'\}
        \}
        \right)
        \xrightarrow{h\to\infty}
        \P^\xi_\beta
        \left(S\cap\{X(r)\in [a,b]\text{ for all }r\in [s',t']\}\right),
    \end{align*}
    where we suppressed $(s,x),(t,y)$ from the notation.
    This implies the lemma.
\end{proof}

We will also need the following simple $L^2$ bound:
\begin{lemm}\label{l:l2bound} 
    We have $\sup_{\d>0} \E\left[(\cZ_{B,\d}^S)^2\right]< \infty$.
\end{lemm}
\begin{proof}
    By \cref{t:AJRAS-Z-existence}\ref{property-strict-positivity} we have $\cZ_{B,\d}^S\ge 0$.
    We also have the simple upper bound
    \begin{align*}
        \cZ_{B,\d}^S
        &= \cZ_\beta^\xi(s,x;t,y) \cdot \P^\xi_{\beta,(s,x),(t,y)}\left(S\cap\{X(r) \in [a,b]\text{ for all $r\in R_\d\cup\{s',t'\}$}\}\right)\\
        &\le \cZ_\beta^\xi(s,x;t,y).
    \end{align*} 
    This together with \cref{t:AJRAS-Z-existence}\ref{property-chaos} yields
    \begin{align*}
        \sup_{\d>0}\E\left[
            (\cZ_{B,\d}^S)^2
        \right]
        \le \E\left[
            \cZ_\beta^\xi(s,x;t,y)^2
        \right]
        <\infty
    \end{align*}
    as claimed.
\end{proof}

We now fix $\d>0$ and derive a  formula for $\cZ_{B,\d}^S$ in terms of the polymer finite-dimensional distributions.
For notational simplicity we will only consider the case where $s'>s$ and $t'<t$, i.e. $[s',t']\subset (s,t)$.
The general case follows from a straightforward modification of our arguments.
\footnote{For example, in the case where $s'=s$ and/or $t'=t$, the upcoming proof applies verbatim if one formally defines $\cZ_\beta^\xi(r,x;r,y)\coloneqq \delta_0(x-y)$ for $r,x,y\in\R$.}
We write
\begin{equation}\label{e:delta-mesh-combined}
    \begin{split}
        \TT &\coloneqq  R_\d \cup \{s_1,\dots,s_\ell\}\\
        &= 
        ((s',t')\cap\d\Z)
        \cup\{s_1,\dots,s_\ell\}.
    \end{split}
\end{equation}
Note that by construction, $\TT\cap \{s',t'\}=\varnothing$.
We partition $\TT$ into three sets:
\begin{align}\label{e:T-partitioned}
    \TT^- \coloneqq  \TT\cap(s,s'),
    \qquad \TT^0 \coloneqq  \TT \cap (s',t'),
    \qquad \TT^+ \coloneqq  \TT\cap(t',t).
\end{align}
That is, we group the times in $\TT$ according to their positions relative to $[s',t']$, as shown in Figure \ref{fig:meas-2}.
For each of the above three sets, we label its elements in increasing order: 
\begin{align*}
    \TT^- = \{t_1^-<\cdots<t_{m^-}^-\},
    \qquad
    \TT^0 = \{t_1^0<\cdots<t_{m^0}^0\},
    \qquad
    \TT^+ = \{t_1^+ <\cdots < t_{m^+}^+\}.
\end{align*}
This notation is demonstrated in Figure \ref{fig:meas-2}.

We now relabel the sets $S_j$ and $[a,b]$ to help write the event $S\cap\left\{X(r)\in [a,b],\;\forall r\in R_\d\cup\{s',t'\}\right\}$
as an intersection of five cylinder sets: three corresponding to the times in $\TT^-, \TT^0, \TT^+$, and one each for $s'$ and $t'$.
Fix $\ast\in\{-,0,+\}$.
For $t_i^\ast \in \TT^\ast$, if $t_i^\ast = s_j$ for some $j\in\lb 1,\ell\rb$, then we define $S^\ast_i \coloneqq  S_j$.
Otherwise, if $t_i^\ast \not\in\{s_1,\dots,s_\ell\}$, then we define $S^\ast_i \coloneqq  [a,b]$.
We form the products
\begin{align}\label{e:def-S-ast}
    S^\ast \coloneqq  \prod_{i=1}^{m^\ast} S^\ast_i
    \qquad\text{for }\ast\in\{-,0,+\}.
\end{align}
With all this notation, we can write
\begin{multline}\label{e:def-ZBdelta-mixed}
    S\cap\{X(r)\in [a,b]\text{ for all }r\in R_\d\cup\{s',t'\}\}\\
    \begin{aligned}
        = 
        \left\{
            (X(t_1^-),\dots, X(t_{m^-}^-))
            \in S^-
        \right\}
        &\cap
        \{X(s')\in [a,b]\}
        \cap
        \left\{
            (X(t_1^0),\dots, X(t_{m^0}^0))
            \in S^0
        \right\}\\
        &\cap
        \{X(t')\in [a,b]\}
        \cap
        \left\{
            (X(t_1^+),\dots, X(t_{m^+}^+))
            \in S^+
        \right\},
    \end{aligned}
\end{multline}
which accomplishes the goal described at the start of this paragraph.
We next use \eqref{e:def-ZBdelta-mixed} to write $\cZ^S_{B,\d}$ in terms of the CDRP finite-dimensional density \eqref{e:CDRP-fdd}.
Fix a realization of the white noise $\xi$, and let $X\sim \P_{\beta,(s,x),(t,y)}^\xi$ be sampled from the polymer measure.

Consider the following finite-dimensional marginal of $X$, depicted in Figure \ref{fig:meas-2} via marked points:
\begin{multline*}
    \mathbf{X}\coloneqq 
    \left(
        X(t_1^-), \dots, X(t_{m^-}^-),\;
        X(s'),\;
        X(t_1^0),\dots, X(t_{m^0}^0),\;
        X(t'),\;
        X(t_1^+),\dots,
        X(t_{m^+}^+)
    \right)\\
    \in \R^{m^-}\times \R\times \R^{m^0} \times \R\times \R^{m^+}.
\end{multline*}
By \eqref{e:CDRP-fdd}, the density 
of $\mathbf{X}$
with respect to Lebesgue measure on $\R^{m^-}\times \R\times \R^{m^0} \times \R\times \R^{m^+}$ is given by (suppressing the dependence on $\beta,\xi$ from the notation)
\begin{multline}\label{e:200}
    \frac{1}{\cZ(s,x;t,y)}
    \prod_{i=0}^{m^-}\cZ(t_i^-, x_i^-; t_{i+1}^-, x_{i+1}^-)
    \cdot
    \prod_{i=0}^{m^0}\cZ(t_i^0, x_i^0;t_{i+1}^0, x_{i+1}^0)
    \cdot
    \prod_{i=0}^{m^+} \cZ(t_i^+, x_i^+; t_{i+1}^+, x_{i+1}^+)\\
    =:\frac{1}{\cZ(s,x;t,y)}
    \cZ^-_{x'}(\bfx^-)
    \cdot
    \cZ^0_{x';x''}(\bfx^0)
    \cdot
    \cZ^+_{x''}(\bfx^+)
\end{multline}
where $(\bfx^-, x', \bfx^0, x'', \bfx^+)\in \R^{m^-}\times\R\times \R^{m^0}\times \R \times \R^{m^+}$,
and where we write
\begin{align*}
    &(t_0^-, x_0^-) \coloneqq  (s,x),\\
    &(t_{m^- +1}^-, x_{m^- +1}^-) = (t_0^0, x_0^0) \coloneqq  (s',x'),\\
    &(t_{m^0 +1}^0, x_{m^0 +1}^0) = (t_{0}^+, x_0^+) \coloneqq  (t',x''),\\
    &(t_{m^+ +1}^+, x_{m^+ + 1}^+) \coloneqq  (t,y).
\end{align*}
In particular, by \eqref{e:def-ZBdelta-mixed}, we have 
\begin{align}\label{e:ZBdelta-final}
    \cZ_{B,\d}^S
    =
    \int\limits_{S^- \times [a,b] \times S^0 \times [a,b] \times S^+}
    \cZ^-_{x'}(\bfx^-)
    \cdot
    \cZ^0_{x';x''}(\bfx^0)
    \cdot
    \cZ^+_{x''}(\bfx^+)
    \,d\bfx^-\,dx'\,d\bfx^0\,dx''\,d\bfx^+.
\end{align}
The above representation for $\cZ^S_{B,\d}$ is the starting point for the upcoming analysis.
For future reference, we make the following technical remark:
\begin{rem}[Conditional expectations factorize]\label{r:independence-inside-outside-B}
    Fix $x',x'',y',y''\in[a,b]$, and  fix $\bfx^\ast,\bfy^\ast\in\R^{m^\ast}$ for each $\ast\in\{-,+\}$.
    By \cref{t:AJRAS-Z-existence}\ref{property-adapted} and \eqref{e:T-partitioned}, 
    the random variables $\cZ^{-}_{x'}(\bfx^-)$ and $\cZ^-_{y'}(\bfy^-)$ are $\cF_{[s,s']\times\R}$-measurable, and the random variables
    $\cZ^+_{x''}(\bfx^+)$ and $\cZ^+_{y''}(\bfy^+)$ are $\cF_{[t',t]\times\R}$-measurable (see the labels on the left side of Figure \ref{fig:meas-2}).
    In particular, the product
    \begin{align*}
        \cZ^{-}_{x'}(\bfx^-)
        \cZ^+_{x''}(\bfx^+)
        \cZ^-_{y'}(\bfy^-)
        \cZ^+_{y''}(\bfy^+)
    \end{align*}
    is measurable with respect to $\cF_{[s,s']\times\R}\vee \cF_{[t',t]\times\R} = \cF_{([s,s']\cup[t',t])\times\R}$ (the equality of $\sigma$-algebras is by \cref{l:independence-disjoint}\ref{generation}).
    Similar reasoning shows that for all $\bfx^0,\bfy^0\in \R^{m^0}$, the product
    \[
        \cZ^0_{x';x''}(\bfx^0)\cZ^0_{y';y''}(\bfy^0)
    \]
    is  $\cF_{[s',t']\times \R}$-measurable.
    Since $[s',t']\times \R$ and $([s,s']\cup[t',t])\times\R$ are essentially disjoint, 
    it follows by \cref{l:conditionally-independent} that for any Borel set $A\in\cB(\R^2)$
    (writing $\E_{\cF_A}[\,\cdot\,]\coloneqq \E[\,\cdot\,|\,\cF_A]$)
    \begin{multline}\label{e:factorize-prod-Z}
        \E_{\cF_A}\left[
            \cZ^{-}_{x'}(\bfx^-)
            \cZ^+_{x''}(\bfx^+)
            \cZ^-_{y'}(\bfy^-)
            \cZ^+_{y''}(\bfy^+)
            \;\cdot\; \cZ^0_{x';x''}(\bfx^0)\cZ^0_{y';y''}(\bfy^0)
        \right]\\
        =
        \E_{\cF_A}\left[
            \cZ^{-}_{x'}(\bfx^-)
            \cZ^+_{x''}(\bfx^+)
            \cZ^-_{y'}(\bfy^-)
            \cZ^+_{y''}(\bfy^+)
        \right]
        \cdot
        \E_{\cF_A}\left[
            \cZ^0_{x';x''}(\bfx^0)\cZ^0_{y';y''}(\bfy^0)
        \right].
    \end{multline}
    In fact, the same reasoning implies a stronger statement: the two factors on the RHS are independent (by combining  Lemmas \ref{l:FAFBisFAB} and \ref{l:independence-disjoint}\ref{independence}).
\end{rem}

\subsubsection{Carrying out the proof of \cref{p:ZB-measurable}}\label{sss:proof-modulo-chaos}

This will be done modulo some inputs involving chaos expansions that will be formulated along the way and proved afterwards.

We will need the following elementary fact about measurability of limits of random variables:
\begin{lemm}\label{l:meas1}
    Let $\{X_n\}_{n\ge 1}$ and $\{Y_n\}_{n\ge 1}$ be sequences of random variables on $(\Xi,\cF,\P)$.
    Let $\cG\subset \cF$ be a sub-$\sigma$-algebra that contains $\cN$ (defined in \cref{ss:whitenoise-short}).
    Assume that $X_n$ is $\cG$-measurable for all $n\ge 1$.
    If $X_n + Y_n\xrightarrow{n\to\infty} X$ almost surely and $Y_n\xrightarrow{n\to\infty} 0$ in probability, then $X$ is $\cG$-measurable.
\end{lemm}
\begin{proof}
    Since $Y_n\to0$ in probability, we can extract a subsequence $Y_{n_k}$ converging to $0$ a.s., and hence $X_{n_k}\to X$ a.s.
    Since $\cG$ contains all $\P$-null subsets in $\cF$, it follows that $X$ is $\cG$-measurable.
\end{proof}

We now prove \cref{p:ZB-measurable} following the strategy outlined in \cref{sss:meas-proof-overview}.
Fix any $\e>0$.
Recall from \eqref{e:def-B-eps} that $B_\e\coloneqq [s',t']\times[a-\e,b+\e]$ is the $\e$-enlargement of $B=[s',t']\times[a,b]$ in the spatial direction.
As in \eqref{e:def-F-eps}, we let $\cF_\e \coloneqq  \cF_{(([s,s']\cup[t',t])\times\R)\cup B_\e}$.
To ease notation, we also write
\begin{align*}
    \E_{\cF_\e}[\,\cdot\,] \coloneqq  \E[\,\cdot\,|\,\cF_{\e}].
\end{align*}
For $\d>0$ let $\cZ_{B,\d}^S$ be as in \eqref{e:def-ZBdelta}.
We will prove the following $L^2$-convergence:
\begin{align}\label{e:L2-goal}
    \lim_{\d\to 0} \E\left[
        \Bigl(
            \cZ_{B,\d}^S - \E_{\cF_\e}[\cZ_{B,\d}^S]
        \Bigr)^2
    \right]
    =0.
\end{align}
Let us quickly deduce \cref{p:ZB-measurable} from the above.
Decompose $\cZ^S_{B,\d}$ as
\begin{align*}
    \cZ^S_{B,\d} = 
    \E_{\cF_\e}[\cZ_{B,\d}^S]
    + 
    \bigl(\cZ_{B,\d}^S - \E_{\cF_\e}[\cZ_{B,\d}^S]\bigr).
\end{align*}
By Lemma \ref{l:partitionalmostsure}, as $\d\to 0$ along the sequence $\d_h = 2^{-h}$, we have $\cZ^S_{B,\d_h} \to \cZ^S_B$ almost surely, and by \eqref{e:L2-goal} we have $\cZ_{B,\d_h}^S - \E_{\cF_\e}[\cZ_{B,\d_h}^S] \to 0$ in probability.
Since $\E_{\cF_\e}[\cZ_{B,\d}^S]$ is by definition $\cF_\e$-measurable, it follows by \cref{l:meas1} that $\cZ^S_{B}$ is $\cF_\e$-measurable.
Taking $\e\to0$ and using that $\cF_\e \downarrow \cF_{(([s,s']\cup[t',t])\times\R)\cup B}$  by \cref{l:nested} finishes the proof of \cref{p:ZB-measurable}.

To prove \eqref{e:L2-goal}, we first note that
\begin{align}\label{e:2000}
    \E\left[
        \bigl(
            \cZ_{B,\d}^S - \E_{\cF_\e}[\cZ_{B,\d}^S]
        \bigr)^2
    \right]
    =
    \E\left[
        (\cZ_{B,\d}^S)^2
    \right]
    -
    \E\left[
        (\E_{\cF_\e}[\cZ_{B,\d}^S])^2
    \right].
\end{align}
To analyze the difference on the RHS, we rely on the following formula for the two terms (proved later in \cref{sss:L2-computations}):
\begin{lemm}\label{l:ZBdelta-moment-formula}
    For any Borel set $A\in\cB(\R^2)$, 
    we have (writing $\E_{\cF_A}[\,\cdot\,]\coloneqq \E[\,\cdot\,|\,\cF_A]$)
    \begin{multline}
        \E\left[
            \E_{\cF_A}[\cZ_{B,\d}^S]^2
        \right]\\
        =
        \int_{[a,b]^4}
        \left(
            \int_{(S^-\times S^+)^2}
            \E\left[
                \E_{\cF_A}[\cZ^-_{x'}(\bfx^-)
                \cZ^+_{x''}(\bfx^+)]
                \E_{\cF_A}[\cZ^-_{y'}(\bfy^-)
                \cZ^+_{y''}(\bfy^+)]
            \right]
            d\bfx^- d\bfx^+ d\bfy^- d\bfy^+
        \right)\times\\
        \times
        \left(
            \int_{(S^0)^2}
            \left\{
                \sum_{n_0,n_1,\dots,n_{m^0} \ge 0}
                \beta^{2(n_0+\cdots+n_{m^0})}
                \prod_{i=0}^{m^0}
                \E\left[
                    I_{n_i}\left(
                        \1_{A}^{\otimes n_i}\,
                        p^{(n_i)}_{(t_i^0, x_i^0), (t_{i+1}^0, x_{i+1}^0)}
                    \right)
                    I_{n_i}\left(
                        \1_{A}^{\otimes n_i}\,
                        p^{(n_i)}_{(t_i^0, y_i^0), (t_{i+1}^0, y_{i+1}^0)}
                    \right)
                \right]
                \right\}
            d\bfx^0 d\bfy^0
        \right)\\
        dx'\,dx''\,dy'\,dy''.
    \end{multline}
\end{lemm}
We will expand the difference in \eqref{e:2000} using \cref{l:ZBdelta-moment-formula} with $A=\R^2$ for the first term, and $A=(([s,s']\cup[t',t])\times \R)\cup B_\e$ for the second term.
Towards this, observe that for $A=\R^2$ or $A=(([s,s']\cup[t',t])\times \R)\cup B_\e$,
by \cref{t:AJRAS-Z-existence}\ref{property-adapted} and \eqref{e:T-partitioned} we have that
\begin{align*}
    \cZ^-_{x'}(\bfx^-)\cZ^+_{x''}(\bfx^+)
    \quad\text{is $\cF_{A}$-measurable}
    \qquad \textup{for all } x',x''\in\R,\; \bfx^-\in\R^{m^-},\; \bfx^+\in\R^{m^+}.
\end{align*}
In particular, for fixed $x',x'',\bfx^-,\bfx^+,y',y'',\bfy^-,\bfy^+$, we have
\begin{align*}
    \E_{\cF_A}\left[
        \cZ^-_{x'}(\bfx^-)
            \cZ^+_{x''}(\bfx^+)
    \right]
    \E_{\cF_A}\left[
        \cZ^-_{y'}(\bfy^-)
        \cZ^+_{y''}(\bfy^+)
    \right]
    = 
    \cZ^-_{x'}(\bfx^-)
            \cZ^+_{x''}(\bfx^+)
            \cZ^-_{y'}(\bfy^-)
            \cZ^+_{y''}(\bfy^+)
    \qquad\text{a.s.}
\end{align*}
Therefore, using \cref{l:ZBdelta-moment-formula}, we can rewrite the RHS of \eqref{e:2000}  as:
\begin{multline}\label{e:201}
    \eqref{e:2000}
    =
    \int_{[a,b]^4}
    \left(
        \int_{(S^-\times S^+)^2}
        \E\left[
            \cZ^-_{x'}(\bfx^-)
            \cZ^+_{x''}(\bfx^+)
            \cZ^-_{y'}(\bfy^-)
            \cZ^+_{y''}(\bfy^+)
        \right]
        d\bfx^- d\bfx^+ d\bfy^- d\bfy^+
    \right)\times\\
    \times\Biggl(
        \int_{(S^0)^2}
        \sum_{n_0,\dots,n_{m^0} \ge 0}
        \beta^{2(n_0+\cdots+n_{m^0})}
        \Biggl\{
            \prod_{i=0}^{m^0}
            \E\left[
                I_{n_i}\left(
                    p^{(n_i)}_{(t_i^0, x_i^0), (t_{i+1}^0, x_{i+1}^0)}
                \right)
                I_{n_i}\left(
                    p^{(n_i)}_{(t_i^0, y_i^0), (t_{i+1}^0, y_{i+1}^0)}
                \right)
            \right]\\
            -
            \prod_{i=0}^{m^0}
            \E\left[
                I_{n_i}\left(
                    \1_{B_\e}^{\otimes n_i}\,
                    p^{(n_i)}_{(t_i^0, x_i^0), (t_{i+1}^0, x_{i+1}^0)}
                \right)
                I_{n_i}\left(
                    \1_{B_\e}^{\otimes n_i}\,
                    p^{(n_i)}_{(t_i^0, y_i^0), (t_{i+1}^0, y_{i+1}^0)}
                \right)
            \right]
            \Biggr\}\,
        d\bfx^0 d\bfy^0
    \Biggr)\\
    dx'\,dx''\,dy'\,dy''.
\end{multline}
We now analyze the difference of products inside the $d\bfx^0 d\bfy^0$ integral.
To ease notation, we write
\begin{align*}
    p^{(n_i)}_{z_i^0} \coloneqq  p^{(n_i)}_{(t_i^0, z_i^0), (t_{i+1}^0, z_{i+1}^0)}
    \qquad\text{for } z\in\{x,y\}.
\end{align*}
Using the elementary identity
\begin{align*}
    \prod_{i=0}^{m^0} u_i - \prod_{i=0}^{m^0} v_i
    = \sum_{\substack{\mathcal{I}\subset\lb 0,m^0\rb\\\mathcal{I}\ne\varnothing}}\,
    \prod_{i\in \mathcal{I}} (u_i-v_i)
    \cdot \prod_{i\in\lb 0,m^0\rb\setminus \mathcal{I}}v_i
    \qquad\text{for } u_i,v_i\in\R,
\end{align*}
a straightforward but tedious calculation lets us rewrite the difference of products inside the curly brackets in \eqref{e:201} as 
\begin{multline}\label{e:202}
    \prod_{i=0}^{m^0}
            \E\left[
                I_{n_i}\left(
                    p^{(n_i)}_{x_i^0}
                \right)
                I_{n_i}\left(
                    p^{(n_i)}_{y_i^0}
                \right)
            \right]
            -
            \prod_{i=0}^{m^0}
            \E\left[
                I_{n_i}\left(
                    \1_{B_\e}^{\otimes n_i}\,
                    p^{(n_i)}_{x_i^0}
                \right)
                I_{n_i}\left(
                    \1_{B_\e}^{\otimes n_i}\,
                    p^{(n_i)}_{y_i^0}
                \right)
            \right]\\
    =
    \sum_{\substack{\mathcal{I}\subset\lb 0,m^0\rb\\\mathcal{I}\ne\varnothing}}
    \prod_{i\in\mathcal{I}}
    \E\left[
                I_{n_i}\left(
                    (1-\1_{B_\e}^{\otimes n_i})\,
                    p^{(n_i)}_{x_i^0}
                \right)
                I_{n_i}\left(
                    (1-\1_{B_\e}^{\otimes n_i})\,
                    p^{(n_i)}_{y_i^0}
                \right)
            \right]\times \\
        \times \prod_{i\in\lb 0,m^0\rb\setminus\mathcal{I}}
        \E\left[
                I_{n_i}\left(
                    \1_{B_\e}^{\otimes n_i}\,
                    p^{(n_i)}_{x_i^0}
                \right)
                I_{n_i}\left(
                    \1_{B_\e}^{\otimes n_i}\,
                    p^{(n_i)}_{y_i^0}
                \right)
            \right],
\end{multline}
where we used  the fact that for all $X,Y\in L^2(\Xi,\cF,\P)$ we have
$\E[XY] - \E[\E_{\cF_\e}[X] \E_{\cF_\e}[Y]] = \E[(X-\E_{\cF_\e}[X])(Y-\E_{\cF_\e}[Y])]$, together with the conditional expectation formula \cref{l:I-given-B}.

Next, to bound the RHS of \eqref{e:202}, we take advantage of the underlying Brownianity (here manifested through heat kernels), thereby making rigorous the heuristic presented in \cref{sss:meas-proof-overview}.
Roughly, since each summand above contains a factor with a projection $1-\1^{\otimes n_i}_{B_\e}$, 
the entire sum can be interpreted as the contribution to the second moment from Brownian paths that exit $B_\e$ yet stay inside $B$ at times $t_0^0,\dots,t_{m^0+1}^0$.
By Brownian fluctuations, this set of paths has Wiener measure approximately $\exp(-c\e^2/\d)$, and hence the overall sum in \eqref{e:202} should be of the same order (see also the caption of Figure \ref{fig:meas-2}).
The next lemma, proved in \cref{sss:L2-computations}, states a version of this for a single chaos covariance.
\begin{lemm}\label{l:moment-of-chaos-Bc}
    As above, let $B_\e \coloneqq  [s',t']\times [a-\e,b+\e]$.
    Fix $t_1<t_2$ such that $[t_1,t_2]\subset [s',t']$.
    Then for all $x_1,x_2, y_1,y_2\in[a,b]$ and all $n\in \Zpos$,
    \begin{multline*}
        \E\left[
            I_n\left(
                (1-\1_{B_\e}^{\otimes n})
                \, p^{(n)}_{(t_1,x_1),(t_2,x_2)}
            \right)
            I_n\left(
                (1-\1_{B_\e}^{\otimes n})
                \,p^{(n)}_{(t_1,y_1),(t_2,y_2)}
            \right)
        \right]\\
        \le
        C\exp\left(
            -\frac{c\e^2}{t_2-t_1}
        \right)
        \E\left[
            I_n\left(
                p^{(n)}_{(t_1,x_1),(t_2,x_2)}
            \right)
            I_n\left(
                p^{(n)}_{(t_1,y_1),(t_2,y_2)}
            \right)
        \right],
    \end{multline*}
    where $C,c>0$ are universal constants.
\end{lemm}

    We now use \cref{l:moment-of-chaos-Bc} to bound each factor of the $\prod_{i\in\mathcal{I}}$ product in \eqref{e:202}.
    To do so, we must check that $x_i^0,x_{i+1}^0,y_i^0,y_{i+1}^0 \in [a,b]$ for every $i\in\lb 0,m^0\rb$.
    This is indeed the case: in \eqref{e:201} we are integrating over $\bfx^0,\bfy^0\in S^0\subset [a,b]^{m^0}$ and $x',x'',y',y''\in[a,b]$, where the inclusion $S^0\subset [a,b]^{m^0}$ is by the discussion below \eqref{e:WLOG-subset-ab}.
    Therefore, by \cref{l:moment-of-chaos-Bc}, for any $i\in\lb 0,m^0\rb$, we have
    \begin{multline*}
        \E\left[I_{n_i}\left(
                        (1-\1_{B_\e}^{\otimes n_i})\,
                        p^{(n_i)}_{x_i^0}
                    \right)
                    I_{n_i}\left(
                        (1-\1_{B_\e}^{\otimes n_i})\,
                        p^{(n_i)}_{y_i^0}
                    \right)
                \right]\\
        \begin{aligned}
        &\le
        \exp\left\{-\left(\frac{c\e^2}{t_{i+1}^0-t_i^0}-C\right)\right\}
        \E\left[I_{n_i}\left(
                        p^{(n_i)}_{x_i^0}
                    \right)
                    I_{n_i}\left(
                       p^{(n_i)}_{y_i^0}
                    \right)
                \right]\\
        &\le
        \exp\left\{-\left(\frac{c\e^2}{\d}-C\right)\right\}
        \E\left[I_{n_i}\left(
                        p^{(n_i)}_{x_i^0}
                    \right)
                    I_{n_i}\left(
                       p^{(n_i)}_{y_i^0}
                    \right)
                \right],
        \end{aligned}
    \end{multline*}
    since $t_{i+1}^0-t_i^0 \le \d$ by \eqref{e:delta-mesh-combined}--\eqref{e:T-partitioned}.
    Applying the above in \eqref{e:202} yields the following upper bound:
    \begin{multline}\label{e:2020}
        \prod_{i=0}^{m^0}
                \E\left[
                    I_{n_i}\left(
                        p^{(n_i)}_{x_i^0}
                    \right)
                    I_{n_i}\left(
                        p^{(n_i)}_{y_i^0}
                    \right)
                \right]
                -
                \prod_{i=0}^{m^0}
                \E\left[
                    I_{n_i}\left(
                        \1_{B_\e}^{\otimes n_i}\,
                        p^{(n_i)}_{x_i^0}
                    \right)
                    I_{n_i}\left(
                        \1_{B_\e}^{\otimes n_i}\,
                        p^{(n_i)}_{y_i^0}
                    \right)
                \right]\\
        \le
        \sum_{\substack{\mathcal{I}\subset\lb 0,m^0\rb\\\mathcal{I}\ne\varnothing}}
        \exp\left\{-|\mathcal{I}|\left(\frac{c\e^2}{\d} - C\right)\right\}
         \prod_{i\in\mathcal{I}}
        \E\left[
            I_{n_i}\left(
                p^{(n_i)}_{x_i^0}
            \right)
            I_{n_i}\left(
                p^{(n_i)}_{y_i^0}
            \right)
        \right]\\
        \times \prod_{i\in\lb 0, m^0\rb\setminus \mathcal{I}}
        \E\left[
                    I_{n_i}\left(
                        \1_{B_\e^{n_i}}\,
                        p^{(n_i)}_{x_i^0}
                    \right)
                    I_{n_i}\left(
                        \1_{B_\e^{n_i}}\,
                        p^{(n_i)}_{y_i^0}
                    \right)
                \right].
    \end{multline}
    Continuing, we bound the RHS using \cref{l:crossterms} to drop the indicators $\1_{B_\e^{n_i}}$ from the $\prod_{i\in \lb 0,m^0\rb \setminus \cI}$ factors,
    followed by the inequality
    $\binom{m^0+1}{k} \le (\frac{(m^0+1) e}{k})^k$:
    \begin{align}\label{e:203}
            \text{RHS of }\eqref{e:2020}
            &\le
            \left(\sum_{\substack{\mathcal{I}\subset\lb 0,m^0\rb\\\mathcal{I}\ne\varnothing}}
            \exp\left\{-|\mathcal{I}|\left(\frac{c\e^2}{\d} - C\right)\right\}
            \right)
            \prod_{i=0}^{m^0}
            \E\left[
                I_{n_i}\left(
                    p^{(n_i)}_{x_i^0}
                \right)
                I_{n_i}\left(
                    p^{(n_i)}_{y_i^0}
                \right)
            \right]\nonumber\\
            &\le
            \left(
                \sum_{k=1}^{m^0}
                \exp\left\{
                    -k\left(\frac{c\e^2}{\d} -C 
                    - \log \left(\frac{(m^0+1)e}{k}\right)
                    \right)
                \right\}
            \right)
            \prod_{i=0}^{m^0}
            \E\left[
                I_{n_i}\left(
                    p^{(n_i)}_{x_i^0}
                \right)
                I_{n_i}\left(
                    p^{(n_i)}_{y_i^0}
                \right)
            \right]\nonumber\\
            &\le
            C'
            e^{-c'\e^2/\d}
            \prod_{i=0}^{m^0}
            \E\left[
                I_{n_i}\left(
                    p^{(n_i)}_{x_i^0}
                \right)
                I_{n_i}\left(
                    p^{(n_i)}_{y_i^0}
                \right)
            \right],
    \end{align}
    where the last line holds for all sufficiently small $\d>0$ because $m^0+1 \ls \d^{-1}$.

    Finally, applying \eqref{e:2020}--\eqref{e:203} in \eqref{e:201} and plugging the resulting bound into \eqref{e:2000}, we obtain
    \begin{multline*}
        \E\left[
            \bigl(
                \cZ_{B,\d}^S - \E_{\cF_\e}[\cZ_{B,\d}^S]
            \bigr)^2
        \right]\\
        \le C' e^{-c' \e^2/\d}
        \int_{[a,b]^4}
        \left(
            \int_{(S^-\times S^+)^2}
            \E\left[
                \cZ^-_{x'}(\bfx^-)
                \cZ^+_{x''}(\bfx^+)\cZ^-_{y'}(\bfy^-)
                \cZ^+_{y''}(\bfy^+)
            \right]
            d\bfx^- d\bfx^+ d\bfy^- d\bfy^+
        \right)\times\\
        \times\Biggl(
            \int_{(S^0)^2}
            \sum_{n_0,\dots,n_{m^0} \ge 0}
            \beta^{2(n_0+\cdots+n_{m^0})}
                \prod_{i=0}^{m^0}
                \E\left[
                    I_{n_i}\left(
                        p^{(n_i)}_{(t_i^0, x_i^0), (t_{i+1}^0, x_{i+1}^0)}
                    \right)
                    I_{n_i}\left(
                        p^{(n_i)}_{(t_i^0, y_i^0), (t_{i+1}^0, y_{i+1}^0)}
                    \right)
                \right]
            d\bfx^0 d\bfy^0
        \Biggr)\\
        dx'\,dx''\,dy'\,dy''\\
        =
        C' e^{-c'\e^2/\d}\,
        \E\left[
            (\cZ_{B,\d}^S)^2
        \right],
    \end{multline*}
    where the last equality is by \cref{l:ZBdelta-moment-formula}.
    Taking $\d\to 0$ and using \cref{l:l2bound}, we have
    \begin{align*}
        \limsup_{\d\to 0} \E\left[
            \bigl(
                \cZ_{B,\d}^S - \E_{\cF_\e}[\cZ_{B,\d}^S]
            \bigr)^2
        \right]
        &\le \sup_{\d>0}\E\left[
            (\cZ_{B,\d}^S)^2
        \right]
        \cdot \limsup_{\d\to 0} 
        C' e^{-c'\e^2/\d}
        =0.
    \end{align*}
    This establishes \eqref{e:L2-goal}, thus completing the proof of \cref{p:ZB-measurable}.
    \qed
    \\

It only remains to prove Lemmas \ref{l:ZBdelta-moment-formula} and \ref{l:moment-of-chaos-Bc}.
This is done in the next subsection.

\subsubsection{Chaos computations}\label{sss:L2-computations}

We start with the second moment formula \cref{l:ZBdelta-moment-formula}. 
While the proof may look a bit daunting, it is just a direct calculation relying on the results of \cref{ss:stocint}.
\begin{proof}[Proof of \cref{l:ZBdelta-moment-formula}]
    Starting from the identity \eqref{e:ZBdelta-final}, we apply Fubini's theorem (note that $\cZ_\beta^\xi$ is everywhere positive by \cref{t:AJRAS-Z-existence}\ref{property-strict-positivity})
    and \cref{r:independence-inside-outside-B} (viz. \eqref{e:factorize-prod-Z}) to obtain
    \begin{multline*}
        \E_{\cF_A}[\cZ_{B,\d}^S]
        = 
        \int\limits_{S^- \times [a,b] \times S^0 \times [a,b] \times S^+}
        \E_{\cF_A}\left[
            \cZ^-_{x'}(\bfx^-)
            \cZ^0_{x';x''}(\bfx^0)
            \cZ^+_{x''}(\bfx^+)
        \right]
        d\bfx^-\,dx'\,d\bfx^0\,dx''\,d\bfx^+\\
        =\int\limits_{S^- \times [a,b] \times S^0 \times [a,b] \times S^+}
        \E_{\cF_A}\left[
            \cZ^-_{x'}(\bfx^-)
            \cZ^+_{x''}(\bfx^+)
            \right]
        \cdot
        \E_{\cF_A}\left[
            \cZ^0_{x';x''}(\bfx^0)
        \right]
        d\bfx^-\,dx'\,d\bfx^0\,dx''\,d\bfx^+.
    \end{multline*}
    Using this formula, we compute the second moment by applying Fubini and using the independence of the two factors in the above integrand (by \cref{r:independence-inside-outside-B}):
        \begin{multline}
        \E\left[
            \left(\E_{\cF_A}[\cZ_{B,\d}^S]\right)^2
        \right]
        =
        \int\limits_{(S^- \times [a,b] \times S^0 \times [a,b] \times S^+)^2}
        \Biggl(
        \E\left[
            \E_{\cF_A}[\cZ^-_{x'}(\bfx^-)
            \cZ^+_{x''}(\bfx^+)]
            \E_{\cF_A}[\cZ^-_{y'}(\bfy^-)
            \cZ^+_{y''}(\bfy^+)]
        \right]\times\\
        \times
        \E\left[
                \E_{\cF_A}[\cZ^0_{x';x''}(\bfx^0)]
                \E_{\cF_A}[\cZ^0_{y';y''}(\bfy^0)]
        \right]\Biggr)
        d\bfx^-\,dx'\,d\bfx^0\,dx''\,d\bfx^+
        \;\;
        d\bfy^-\,dy'\,d\bfy^0\,dy''\,d\bfy^+\\
        =
        \int_{[a,b]^4}
        \left(
            \int_{(S^-\times S^+)^2}
        \E\left[
            \E_{\cF_A}[\cZ^-_{x'}(\bfx^-)
            \cZ^+_{x''}(\bfx^+)]
            \E_{\cF_A}[\cZ^-_{y'}(\bfy^-)
            \cZ^+_{y''}(\bfy^+)]
        \right]
            d\bfx^- d\bfx^+ d\bfy^- d\bfy^+
        \right)\times\\
        \times\left(
            \int_{(S^0)^2}
            \E\left[
                \E_{\cF_A}[\cZ^0_{x';x''}(\bfx^0)]
                \E_{\cF_A}[\cZ^0_{y';y''}(\bfy^0)]
            \right]
            d\bfx^0 d\bfy^0
        \right)
        dx'\,dx''\,dy'\,dy''.
    \end{multline}
    It remains to show that
    \begin{multline*}
        \E\left[
                \E_{\cF_A}[\cZ^0_{x';x''}(\bfx^0)]
                \E_{\cF_A}[\cZ^0_{y';y''}(\bfy^0)]
        \right]\\
        =
        \sum_{n_0,n_1,\dots,n_{m^0} \ge 0}
        \beta^{2(n_0+\cdots+n_{m^0})}
        \prod_{i=0}^{m^0}
        \E\left[
            I_{n_i}\left(
                \1_{A}^{\otimes n_i}\,
                p^{(n_i)}_{(t_i^0, x_i^0), (t_{i+1}^0, x_{i+1}^0)}
            \right)
            I_{n_i}\left(
                \1_{A}^{\otimes n_i}\,
                p^{(n_i)}_{(t_i^0, y_i^0), (t_{i+1}^0, y_{i+1}^0)}
            \right)
        \right].
    \end{multline*}
    Towards this,
    recalling the definition of $\cZ^0_{x';x''}(\bfx^0)$ from \eqref{e:200}, we can rewrite the LHS 
    using \cref{l:conditionally-independent} and \cref{t:AJRAS-Z-existence}\ref{property-adapted} as
    \begin{align*}
        \E\left[
                \E_{\cF_A}[\cZ^0_{x';x''}(\bfx^0)]\,
                \E_{\cF_A}[\cZ^0_{y';y''}(\bfy^0)]
        \right]
        =
        \prod_{i=0}^{m^0}
        \E\left[
            \E_{\cF_A}[\cZ_\beta(t_i^0,x_i^0; t_{i+1}^0, x_{i+1}^0)]\,
            \E_{\cF_A}[\cZ_\beta(t_i^0,y_i^0; t_{i+1}^0, y_{i+1}^0)]
        \right].
    \end{align*}
    We rewrite the RHS using the inner product formula \eqref{e:parseval-CE} and the chaos expansion \eqref{e:chaos-Z} of $\cZ_\beta$:
    \begin{multline*}
        \prod_{i=0}^{m^0}
        \E\left[
            \E_{\cF_A}[\cZ_\beta(t_i^0,x_i^0; t_{i+1}^0, x_{i+1}^0)]\,
            \E_{\cF_A}[\cZ_\beta(t_i^0,y_i^0; t_{i+1}^0, y_{i+1}^0)]
        \right]\\
        \begin{aligned}
        &=\prod_{i=0}^{m^0}\left(
            \sum_{n\ge 0} \beta^{2n}
            \E\left[
                I_n\left(
                    \1_{A}^{\otimes n}\,
                    p^{(n)}_{(t_i^0,x_i^0), (t_{i+1}^0, x_{i+1}^0)}
                \right)
                I_n\left(
                    \1_{A}^{\otimes n}\,
                    p^{(n)}_{(t_i^0,y_i^0), (t_{i+1}^0, y_{i+1}^0)}
                \right)
            \right]
        \right)\\
        &=\sum_{n_0,\dots,n_{m^0} \ge 0}
        \beta^{2(n_0+\cdots+n_{m^0})}
        \prod_{i=0}^{m^0}
        \E\left[
            I_{n_i}\left(
                \1_{A}^{\otimes n_i}\,
                p^{(n_i)}_{(t_i^0, x_i^0), (t_{i+1}^0, x_{i+1}^0)}
            \right)
            I_{n_i}\left(
                \1_{A}^{\otimes n_i}\,
                p^{(n_i)}_{(t_i^0, y_i^0), (t_{i+1}^0, y_{i+1}^0)}
            \right)
        \right],
        \end{aligned}
    \end{multline*}
    where we interchanged product and infinite sum using Fubini (every summand is non-negative by the It\^o isometry (\cref{l:ito-isometry}) and the pointwise non-negativity of the function $\1_{A}^{\otimes n} p^{(n)}_{(s,x),(t,y)}$).
    This finishes the proof.
\end{proof}

We next prove \cref{l:moment-of-chaos-Bc}, which bounds the covariance of two chaos terms projected onto $L^2(\Xi,\cF_{B_\e},\P)^{\perp}$.
The proof amounts to straightforward tedious calculations with the heat kernel.

\begin{proof}[Proof of \cref{l:moment-of-chaos-Bc}]
    When $n=0$ the LHS vanishes since $1-\1^{\otimes 0}_{B_\e}=0$, and the RHS equals $C\exp\left(-\frac{c\e^2}{t_2-t_1}\right)p(t_2-t_1, x_2-x_1)p(t_2-t_1,y_2-y_1)$ which is positive.
    So the $n=0$ case holds trivially.
    
    We assume $n\ge 1$ from now on.
    In this case we have $1-\1_{B_\e}^{\otimes n} = \1_{(B_\e^n)^c}$.
    Define the simplex of ordered times
    \begin{align*}
        \Delta^n(t_1,t_2) \coloneqq  \{\bfs=(s_1,\dots,s_n): t_1 =: s_0 < s_1 < \cdots < s_n < s_{n+1} \coloneqq t_2\}.
    \end{align*}
    Fix any Borel set $A\in \cB(\R^n)$ (we will eventually set $A=\R^n$ or $A=\R^n\setminus[a-\e,b+\e]^n$), and write
    \begin{align}\label{e:209}
        A' \coloneqq  \{(s_1,z_1;\cdots;s_n,z_n)\in\R^{2n} : (z_1,\dots,z_n)\in A\}
    \end{align}
    (so $A'$ is just $\R^n\times A$ with the coordinates interleaved).
    Note that we introduce $A'$ because the functions $p^{(n)}_{(s,x),(t,y)}$ are defined on $\R^{2n}$.
    By \cref{l:ito-isometry} (It\^o isometry),
    we have
    \begin{multline}\label{e:chaos-moment1}
        \E\left[
            I_n\left(
                \1_{A'}\,p^{(n)}_{(t_1,x_1),(t_2,x_2)}
            \right)
            I_n\left(
                \1_{A'}\,p^{(n)}_{(t_1,y_1),(t_2,y_2)}
            \right)
        \right]
        = 
        \left\langle
            \1_{A'}\,p^{(n)}_{(t_1,x_1),(t_2,x_2)},\;
            \1_{A'}\,p^{(n)}_{(t_1,y_1),(t_2,y_2)}
        \right\rangle_{L^2(\R^{2n})}
        \\
        = \int_{\Delta^n(t_1,t_2)}
        \int_{A} 
            \Biggl[
            p(s_1-s_0, z_1-x_1)p(s_1-s_0,z_1-y_1)
            \left(\prod_{i=1}^{n-1}
            p(s_{i+1}-s_{i}, z_{i+1}-z_i)
            \right)^2\\
        \times p(s_{n+1}-s_n, x_2-z_n) p(s_{n+1}-s_n, y_2-z_n)
        \Biggr]
            d\bfz\, d\bfs.
    \end{multline}
    To simplify the above integral, we will use the identities
    \begin{align*}
        p(t,x)^2 &= \frac{1}{\sqrt{2\pi t}}p(t,\sqrt{2}x),
        && t>0,\;x\in\R\\
        p(t,z-x)p(t,z-y) &= 
        \frac{1}{\sqrt{2\pi t}}
        \exp\left(
            -\tfrac{(x-y)^2}{4t}
        \right)
        p\left(
            t,\; \sqrt{2}z - \tfrac{x+y}{\sqrt{2}}
        \right),
        &&t>0,\; x,y,z\in\R.
    \end{align*}
    By repeatedly applying the above identities, we can rewrite the integrand in \eqref{e:chaos-moment1} as 
    \begin{multline*}
            p(s_1-s_0, z_1-x_1)p(s_1-s_0,z_1-y_1)
            \left(\prod_{i=1}^{n-1}
            p(s_{i+1}-s_{i}, z_{i+1}-z_i)
            \right)^2\\
            \times p(s_{n+1}-s_n, x_2-z_n) p(s_{n+1}-s_n, y_2-z_n)\\
    \begin{aligned}
        &=
        \frac{1}{\sqrt{2\pi (s_1-s_0)}}
        \exp\left(-\tfrac{(x_1-y_1)^2}{4(s_1-s_0)}\right)
        p\left(s_1-s_0,\; \sqrt{2}z_1 - \tfrac{x_1+y_1}{\sqrt{2}}\right)\\
        &\qquad\times 
        \left(
            \prod_{i=1}^{n-1}\frac{1}{\sqrt{2\pi(s_{i+1}-s_i)}}p(s_{i+1}-s_i, \sqrt{2}(z_{i+1}-z_i))
        \right)\\
        &\qquad\times
        \frac{1}{\sqrt{2\pi (s_{n+1}-s_n)}}
        \exp\left(-\tfrac{(x_2-y_2)^2}{4(s_{n+1}-s_n)}\right)
        p\left(s_{n+1}-s_n,\; \tfrac{x_2+y_2}{\sqrt{2}} - \sqrt{2} z_n\right)\\
        &= \left(\prod_{i=0}^{n}\frac{1}{\sqrt{2\pi(s_{i+1}-s_i)}}\right)
        \exp\left(
            -\tfrac{(x_1-y_1)^2}{4(s_1-s_0)}
            -\tfrac{(x_2-y_2)^2}{4(s_{n+1}-s_n)}
        \right)\\
        &\times
             \underbrace{p\left(s_1-s_0,\; \sqrt{2}z_1 - \tfrac{x_1+y_1}{\sqrt{2}}\right)
            \cdot
            \prod_{i=1}^{n-1}p(s_{i+1}-s_i, \sqrt{2}(z_{i+1}-z_i))
            \cdot 
            p\left(s_{n+1}-s_n,\; \tfrac{x_2+y_2}{\sqrt{2}} - \sqrt{2}z_n\right)}_{=:\,\rho(\bfs,\, \sqrt{2}\bfz)}.
        \end{aligned}
    \end{multline*}
    Observe that up to a normalization, $\rho(\bfs,\cdot)$ is the density of the random vector $(X(s_1),\dots,X(s_n))$, where $X$ is a Brownian bridge from $X(s_0)=\frac{x_1+y_1}{\sqrt{2}}$ to $X(s_{n+1})=\frac{x_2+y_2}{\sqrt{2}}$.
    In particular, for any $\bfs\in\Delta^n(t_1,t_2)$, we have
    \begin{multline*}
        \int_A \rho(\bfs,\sqrt{2}\bfz)\,d\bfz
        = 
        2^{-n/2}
        p\left(s_{n+1}-s_0, \tfrac{x_2+y_2}{\sqrt{2}}-\tfrac{x_1+y_1}{\sqrt{2}}\right)\\
        \times\mathrm{P}^{\mathrm{BB}}
        _{(s_0,\frac{x_1+y_1}{\sqrt{2}}), (s_{n+1},\frac{x_2+y_2}{\sqrt{2}})}
        \left(
            (X(s_1), X(s_2),\dots,X(s_n)) \in \sqrt{2}A
        \right),
    \end{multline*}
    where the factor $2^{-n/2}$ is the inverse Jacobian of $\bfz\mapsto \sqrt{2}\bfz$.
    Finally, combining the above two displays with \eqref{e:chaos-moment1}, we obtain the moment formula
    \begin{multline}\label{e:chaos-moment2}
        \E\left[
            I_n\left(
                \1_{A'}\,p^{(n)}_{(t_1,x_1),(t_2,x_2)}
            \right)
            I_n\left(
                \1_{A'}\,p^{(n)}_{(t_1,y_1),(t_2,y_2)}
            \right)
        \right]
        = 
        2^{-n/2}
        p\left(s_{n+1}-s_0, \tfrac{x_2+y_2}{\sqrt{2}}-\tfrac{x_1+y_1}{\sqrt{2}}\right)
        \times \\
        \times \int_{\Delta^n(t_1,t_2)}
        \Biggl[\exp\left(
            -\tfrac{(x_1-y_1)^2}{4(s_1-s_0)}
            -\tfrac{(x_2-y_2)^2}{4(s_{n+1}-s_n)}
        \right)
        \left(\prod_{i=0}^n \frac{1}{\sqrt{2\pi(s_{i+1}-s_i)}}\right) \times\\
        \times \mathrm{P}^{\mathrm{BB}}
        _{(s_0,\frac{x_1+y_1}{\sqrt{2}}), (s_{n+1},\frac{x_2+y_2}{\sqrt{2}})}
        \left(
            (X(s_1), X(s_2),\dots,X(s_n)) \in \sqrt{2}A
        \right)
        \Biggr]
        d\bfs.
    \end{multline}
    We now specialize to two cases: $A=\R^n$ and $A=\R^n\setminus[a-\e,b+\e]^n$.

    For $A=\R^n$ (i.e. $A'=\R^{2n}$), the Brownian bridge probability equals $1$, and thus \eqref{e:chaos-moment2} yields
    \begin{multline}\label{e:chaos-moment3}
        \E\left[
            I_n\left(
                p^{(n)}_{(t_1,x_1),(t_2,x_2)}
            \right)
            I_n\left(
                p^{(n)}_{(t_1,y_1),(t_2,y_2)}
            \right)
        \right]
        =
        2^{-n/2}
        p\left(s_{n+1}-s_0, \tfrac{x_2+y_2}{\sqrt{2}}-\tfrac{x_1+y_1}{\sqrt{2}}\right)
        \times\\
        \times
        \int_{\Delta^n(t_1,t_2)}
        \exp\left(
            -\tfrac{(x_1-y_1)^2}{4(s_1-s_0)}
            -\tfrac{(x_2-y_2)^2}{4(s_{n+1}-s_n)}
        \right)
        \left(\prod_{i=0}^n \frac{1}{\sqrt{2\pi(s_{i+1}-s_i)}}\right)
        d\bfs.
    \end{multline}

    For $A=\R^n\setminus[a-\e,b+\e]^n$,
    we have the following simple bound for the Brownian bridge probability in \eqref{e:chaos-moment2}: 
    \begin{multline*}
        \mathrm{P}^{\mathrm{BB}}
        _{(s_0,\frac{x_1+y_1}{\sqrt{2}}), (s_{n+1},\frac{x_2+y_2}{\sqrt{2}})}
        \left( (X(s_1), \dots,X(s_n)) \in \sqrt{2}A\right)\\
        \begin{aligned}
            &=
            \mathrm{P}^{\mathrm{BB}}
            _{(s_0,\frac{x_1+y_1}{\sqrt{2}}), (s_{n+1},\frac{x_2+y_2}{\sqrt{2}})}
            \left(
                \exists i \in\lb 1,n\rb : X(s_i) \not\in \sqrt{2}[a-\e,b+\e]
            \right)
            \\
            &\le
            \mathrm{P}^{\mathrm{BB}}
            _{(s_0,\frac{x_1+y_1}{\sqrt{2}}), (s_{n+1},\frac{x_2+y_2}{\sqrt{2}})}
            \left(
                \exists r\in[s_0,s_{n+1}] : X(r)
                \not\in \sqrt{2}[a-\e,b+\e]
            \right)\\
            &\le C\exp\left(
                -\frac{c\e^2}{s_{n+1}-s_0}
            \right),
        \end{aligned}
    \end{multline*}
    where the last line is by standard  sub-Gaussian tail bounds for the extrema of Brownian bridge (observe that $\frac{x_1+y_1}{\sqrt{2}},\frac{x_2+y_2}{\sqrt{2}}\in \sqrt{2}[a,b]$ by assumption, hence $\mathrm{E}^{\mathrm{BB}}[X(r)]\in \sqrt{2}[a,b]$ for all $r\in[s_0,s_{n+1}]$).
    Applying this bound in \eqref{e:chaos-moment2} and comparing the result with \eqref{e:chaos-moment3}, 
    we obtain (recall $s_0\coloneqq t_1$ and $s_{n+1}\coloneqq t_2$)
    \begin{multline*}
        \E\left[
            I_n\left(
                \1_{A'}\, p^{(n)}_{(t_1,x_1),(t_2,x_2)}
            \right)
            I_n\left(
                \1_{A'}\,
                p^{(n)}_{(t_1,y_1),(t_2,y_2)}
            \right)
        \right]\\
        \le 
        C\exp\left(
            -\frac{c\e^2}{t_2-t_1}
        \right)
         \E\left[
            I_n\left(p^{(n)}_{(t_1,x_1),(t_2,x_2)}
            \right)
            I_n\left(
                p^{(n)}_{(t_1,y_1),(t_2,y_2)}
            \right)
        \right],
    \end{multline*}
    where $A'$ is defined as in \eqref{e:209} with respect to $A=\R^n\setminus[a-\e,b+\e]^n$.
    This is almost the inequality stated in \cref{l:moment-of-chaos-Bc}; the only difference is that above, on the LHS we have $\1_{A'}$ instead of $1-\1_{B_\e}^{\otimes n}$.
    To obtain the desired inequality, 
    note that for $w\in\{x,y\}$, the function $p^{(n)}_{(t_1,w_1),(t_2,w_2)}$ is supported in the set $\{(s_1,z_1;\cdots;s_n,z_n) : s_1,\dots,s_n \in [t_1,t_2]\}$.
    Therefore, we actually have the pointwise equality
    \begin{align*}
        \1_{A'}\,
        p^{(n)}_{(t_1,w_1),(t_2,w_2)}
        =
        (1-\1_{B_\e}^{\otimes n})
         p^{(n)}_{(t_1,w_1),(t_2,w_2)},
    \end{align*}
    which when plugged into the above display finishes the proof.
\end{proof}

\subsection{Resampling white noise}\label{s:resampling-framework}

In this subsection we elaborate on the resampling framework presented in \eqref{e:def-Xi-resample}--\eqref{e:def-xi-eta}, and prove the Efron--Stein inequality \cref{p:efron-stein-ineq}.

Recall that the resampling framework consists of a probability space $(\Xires,\Fres,\Pres)$ (defined in \eqref{e:def-Xi-resample}) supporting a white noise $\xi$, such that for any Borel set $B\in\cB(\R^2)$, we can construct another white noise $\eta_B$ on $(\Xires,\Fres,\Pres)$ obtained from $\xi$ by \emph{resampling} the restriction $\xi|_B$.
The following lemma tracks how resampling $\xi|_B$ affects $\xi|_{B}$-measurable and $\xi|_{B^c}$-measurable random variables.

\begin{lemm}\label{l:Bc-meas-resampling-invariant}
    Fix $B\in\cB(\R^2)$ and let $\xi,\eta_B : (\Xires, \Fres, \Pres) \to (\Xi,\cF,\P)$ be white noises as in \eqref{e:def-xi-eta} (so $\eta_B$ is obtained from $\xi$ by resampling $\xi|_B$).
    Let $S$ be a Polish space and let $Z : (\Xi,\cF,\P) \to (S,\cB(S))$ be a measurable map into $S$.
    The following properties hold.
    \begin{enumerate}[label={\rm(\roman*)}]
        \item\label{resample-independent} If $Z$ is $\cF_{B}$-measurable, then $Z\circ \xi$ is independent of $\eta_B$.
        \item\label{resample-invariant} If $Z$ is $\cF_{B^c}$-measurable, then $\Pres(Z\circ \xi = Z \circ \eta_B) = 1$.
    \end{enumerate}
\end{lemm}
\begin{proof}
Fix $A\in\{B,B^c\}$ and assume $Z$ is $\cF_A$-measurable.
Recall from \eqref{e:def-Res} that $\Res_A:\Xi \to \Xi_A$ denotes  the restriction map $\xi \mapsto \xi|_{L^2_0(A)}$.
Since $\cF_{A}$ is defined as $\sigma(\Res_{A})\vee \cN$, where $\cN$ is the $\sigma$-algebra generated by all $\P$-null subsets in $\cF$ (see \eqref{e:def-FB}),
it follows by \cref{l:removing-completion} (stated and proved right after the proof of this lemma) that there exists a $\sigma(\Res_{A})$-measurable version $Y:(\Xi,\sigma(\Res_{A})) \to (S,\cB(S))$ such that $\P(Z=Y)=1$.
Since $S$ is Polish, by the Doob--Dynkin factorization lemma (e.g. \cite[Lemma 1.14]{Kal21})
there exists a measurable $\wt{Y}:(\Xi_{A},\Fprod_{A})\to(S,\cB(S))$ such that
$Y = \wt{Y}\circ \Res_{A}$.
In particular, we have
\begin{align*}
    Z = \wt{Y}\circ \Res_{A}\qquad\text{$\P$-a.s.}
\end{align*}
Since $\P$ is the pushforward of $\Pres$ by each of the maps $\xi,\eta_B : (\Xires,\Fres,\Pres) \to (\Xi,\cF,\P)$, it follows that (up to a $\Pres$-null set) $Z\circ \xi$ is a measurable function of $\Res_A\circ \xi$, and similarly $Z\circ \eta_B$ is a measurable function of $\Res_A\circ \eta_B$.
It therefore suffices to establish the following two properties:
\begin{enumerate}[label={\rm(\roman*$^{\prime}$)}]
    \item\label{resample-independent-prime} $\Res_B\circ \xi$ is independent of $\eta_B$.
    \item\label{resample-invariant-prime} $\Res_{B^c}\circ \xi = \Res_{B^c}\circ \eta_B$ almost surely.
\end{enumerate}

For \ref{resample-independent-prime}, note 
that $\sigma(\eta_B)$ is generated by the coordinate maps $(\xi_1,\xi_2) \mapsto \eta_B(\xi_1,\xi_2)(f)$ for $f\in L^2(\R^2)$.
Similarly, $\sigma(\Res_B\circ \xi)$ is generated by the maps $(\xi_1,\xi_2)\mapsto \xi(\xi_1,\xi_2)(g)$ for $g\in L^2_0(B)$.
Fixing $f\in L^2(\R^2)$ and $g\in L^2_0(B)$, we have
\begin{align*}
    \Eres\left[
        \eta_B(\xi_1,\xi_2)(f)\cdot 
        (\Res_B\circ \xi)(\xi_1,\xi_2)(g)
    \right]
    &= \Eres\left[
        (\xi_1(P_{B^c}f) + \xi_2(P_B f))
        \cdot
        \xi_1(g)
    \right]\\
    &= \langle \1_{B^c} f,g\rangle_{L^2(\R^2)}\\
    &=0,
\end{align*}
where we used that $\xi_1,\xi_2$ are independent under $\Pres$, and that $g$ vanishes on $B^c$.
It follows by Gaussianity that for all $n\in\N$, all $f_1,\dots,f_n\in L^2(\R^2)$, and all $g_1,\dots,g_n \in L^2_0(B)$, the finite-dimensional marginals $(\eta_B(\cdot)(f_1),\dots,\eta_B(\cdot)(f_n))$ and  $((\Res_B\circ \xi)(\cdot)(g_1),\dots,(\Res_B\circ \xi)(\cdot)(g_n))$ are independent.
Claim \ref{resample-independent-prime} follows by the $\pi$-$\lambda$ theorem.

For \ref{resample-invariant-prime}, note that under $\Pres$ the coordinate projection $\xi_2$ satisfies $\xi_2(0)=0$ a.s., since $\xi_2(0)$ is Gaussian with mean and variance zero.
Thus since $L^2_0(B^c)$ and $L^2_0(B)$ are orthogonal subspaces, it holds $\Pres$-a.s. that 
for all $f\in L^2_0(B^c)$,
\begin{align*}
    \eta_B(\xi_1,\xi_2)(f) 
    &= \xi_1(P_{B^c} f) + \xi_2(P_B f)\\
    &= \xi_1(f).
\end{align*}
In other words, $(\Res_{B^c}\circ \eta_B)(\xi_1,\xi_2) = \Res_{B^c}(\xi_1)$ almost surely.
Since by definition $\xi(\xi_1,\xi_2)=\xi_1$,  this implies \ref{resample-invariant-prime}.
\end{proof}

We conclude this discussion with the measure-theoretic lemma used in the above proof.

\begin{lemm}\label{l:removing-completion}
    Let $(\Xi,\cF,\P)$ be as in \cref{def:white-noise}, and let $\cN$ be the $\sigma$-algebra generated by all $\P$-null subsets in $\cF$.
    Fix any sub-$\sigma$-algebra $\cG\subset \cF$ and write $\overline{\cG}\coloneqq \cG\vee\cN$.
    Let $S$ be a Polish space and let $Z:(\Xi,\overline{\cG}) \to (S,\cB(S))$ be a $\overline{\cG}$-measurable $S$-valued random variable.
    Then there exists a $\cG$-measurable $Y:(\Xi,\cG)\to (S,\cB(S))$ such that $\P(Z=Y)=1$, where the event $\{Z=Y\} \in \overline{\cG}$.
\end{lemm}
\begin{proof}
    Applying a Borel isomorphism reduces to the case $S\in\cB(\R)$. 
    We first assume that $Z$ is integrable.
    Set $Y\coloneqq \E[Z|\cG]$.
    Since $Y,Z$ are $\overline{\cG}$-measurable, it suffices to show that $\int_{\sA} (Z-Y)\,d\P=0$ for all $\sA\in\overline{\cG}$.
    In fact, by the $\pi$-$\lambda$ theorem, it suffices to show this for $\sA$ of the form $\sA=\sB\cup\sN$  for $\sB\in\cG$ and $\sN\in\cF$ with $\P(\sN)=0$.
    We have
    \begin{align*}
        \int_{\sB\cup\sN} (Z-Y)\,d\P
        = \int_{\sB} (Z-Y)\,d\P
        + \int_{\sN\setminus \sB} (Z-Y)\,d\P
        = 0,
    \end{align*}
    where we used the definition of conditional expectation, and that $\P(\sN\setminus \sB) \le \P(\sN)=0$.

    Suppose $Z$ is not integrable.
    By the above argument, the truncation $Z_n \coloneqq  \max\{-n, \min\{Z,n\}\}$ satisfies $\E[Z_n|\cG]=Z_n$ a.s. for all $n$.
    Since $Z_n \to Z$ pointwise, we get $Z=\limsup_{n\to\infty}\E[Z_n|\cG]$ a.s., and the RHS is $\cG$-measurable.
\end{proof}

\addtocontents{toc}{\SkipTocEntry}
\subsection*{Efron--Stein inequality for white noise}
We are now in a position to prove the Efron--Stein inequality for white noise (\cref{p:efron-stein-ineq}), which again follows the usual proof but carried out in the current framework.

We start by recalling an elementary property of conditional variance.
Let $(X,Y)$ be a random vector, and let $Z=Z(X,Y)$ be a square-integrable function.
Then
\begin{align*}
    \E[\Var(Z|Y)]
    = \E\left[
        \bigl(Z - \E[Z | Y]\bigr)^2
    \right]
    =
    \frac12 \E\left[
        \bigl(
            Z(X,Y) - Z(X',Y)
        \bigr)^2
    \right],
\end{align*}
where $X'$ is a random variable such that $(X,Y)\law (X',Y)$, and such that $X,X'$ are conditionally independent given $Y$.
The following lemma contains an analogous property for functions of white noise.

\begin{lemm}\label{l:conditional-variance-formula}
    Fix $B\in\cB(\R^2)$ and let $(\xi,\eta_B)$ be coupled white noises as in \eqref{e:def-xi-eta} (so $\eta_B$ is obtained from $\xi$ by resampling $\xi|_B$).
    Then for any square-integrable function 
    $Z:(\Xi,\cF,\P)\to\R$,
    \begin{align*}
        \E\left[
            \Var\bigl(Z|\cF_{B^c}\bigr)
        \right]
        =
        \E\left[
            \bigl(Z - \E[Z|\cF_{B^c}]\bigr)^2
        \right]
        &= \frac12
        \Eres\left[
            \bigl(Z(\xi) - Z(\eta_B)\bigr)^2
        \right].
    \end{align*}
\end{lemm}
\begin{proof}
    The first equality is by definition, so we focus on proving the second equality.
    Since $\cF$ is the completion of the product $\sigma$-algebra $\cB(\R)^{\otimes L^2(\R^2)}$, {by a routine approximation argument} it suffices to consider the case where $Z$ depends on finitely many coordinates, i.e. $Z = Z(\xi(f_1),\dots,\xi(f_n))$ for some $n\in\N$ and $f_1,\dots,f_n\in L^2(\R^2)$, where here $\xi$ denotes a generic element of $\Xi$.
    In this case, by writing $\xi(f_i) = \xi(f_i\1_{B}) + \xi(f_i\1_{B^c})$ (this holds $\P$-a.s. for any fixed $f_1,\dots,f_n\in L^2(\R^2)$, see e.g. \cite[p.~4]{nualart}),
    we may further express $Z$ as a function 
    \begin{align*}
        Z = 
        Z\left(
            \bigl(
                \xi(f_1\1_B),\dots,\xi(f_n\1_B)
            \bigr),\;\;
            \bigl(\xi(f_1\1_{B^c}),\dots,\xi(f_n\1_{B^c})\bigr)
        \right).
    \end{align*}
    For $A\in\{B,B^c\}$, write $\pi_A(\xi) \coloneqq  (\xi(f_1\1_A),\dots,\xi(f_n\1_A))$.
    Note that $\pi_{B^c}$ is $\cF_{B^c}$-measurable, and $\pi_B$ is independent of $\cF_{B^c}$ under $\P$.
    This implies that 
    \begin{align*}
        \E\left[
            \bigl(Z - \E[Z|\cF_{B^c}]\bigr)^2
        \right]
        &= \E\left[
            \Var\left(Z\,\big|\,\cF_{B^c}\right)
        \right]\\
        &= \E\left[
            \Var\left(
                Z(\pi_B, \pi_{B^c})
                \,\big|\,\cF_{B^c}
            \right)
        \right]\\
        &=
        \E\left[
            \Var\left(
                Z(\pi_B, \pi_{B^c})
                \,\big|\,\pi_{B^c}
            \right)
        \right].
    \end{align*}
    Now for $\xi,\eta_B : (\Xires,\Fres,\Pres) \to (\Xi,\cF,\P)$ as in \eqref{e:def-xi-eta}, \cref{l:Bc-meas-resampling-invariant} implies that 
    $\pi_{B^c}(\xi) = \pi_{B^c}(\eta_B)$ almost surely, and that $\pi_B(\xi)$ and $\pi_B(\eta_B)$ are conditionally independent given $\pi_{B^c}(\xi)$.
    Therefore, continuing from the previous display, we have
    \begin{align*}
        \E\left[
            \Var\left(
                Z(\pi_B, \pi_{B^c})
                \,\big|\,\pi_{B^c}
            \right)
        \right]
        &= 
        \Eres\left[
            \Var_\res\left(
                Z(\pi_B(\xi), \pi_{B^c}(\xi))
                \,\big|\,\pi_{B^c}(\xi)
            \right)
        \right]\\
        &= 
        \frac12\Eres\left[
            \left(
                Z(\pi_B(\xi),\pi_{B^c}(\xi))
                - Z(\pi_B(\eta_B), \pi_{B^c}(\xi))
            \right)^2
        \right]\\
        &= 
        \frac12\Eres\left[
            \left(
                Z(\pi_B(\xi),\pi_{B^c}(\xi))
                - Z(\pi_B(\eta_B), \pi_{B^c}(\eta_B))
            \right)^2
        \right]\\
        &=
        \frac12
        \Eres\left[
            \left(
                Z(\xi)
                - Z(\eta_B)
            \right)^2
        \right],
    \end{align*}
    which completes the proof.
\end{proof}

We now prove \cref{p:efron-stein-ineq}.

\begin{proof}[Proof of \cref{p:efron-stein-ineq}]
    By \cref{l:conditional-variance-formula},
    the proposition is equivalent to 
    \begin{align}\label{e:600}
        \E\left[
            \bigl(Z - \E[Z|\cF_{B^c}]\bigr)^2
        \right]
        \le \sum_{k\in\Z} \E\left[
            \bigl(Z - \E[Z|\cF_{B_k^c}]\bigr)^2
        \right].
    \end{align}
    Write
    \begin{align*}
        A_{k} \coloneqq  \bigcap_{k'\in[k+1,\infty)\cap\Z} B_{k'}^c
        \qquad\text{and}\qquad
        \E_{k}[\,\cdot\,] \coloneqq  \E[\,\cdot\,|\,\cF_{A_{k}}].
    \end{align*} 
    Since $B_k$ are pairwise essentially disjoint, there exists a Lebesgue-null set $N\in\cB(\R^2)$ such that
    $A_{k} \uparrow \R^2\setminus N$ as $k\uparrow\infty$, and $A_{k} \downarrow B^c$ as $k\downarrow -\infty$.
    Note that $\cF = \cF_{\R^2\setminus N}$ since $L^2_0(\R^2\setminus N) = L^2(\R^2)$.
    Therefore, by \cref{l:nested} and the  martingale convergence theorem,
    we have
    \begin{align*}
        \E_{k}[Z] \xrightarrow{k\to\infty} Z \qquad\text{in $L^2(\Xi,\cF,\P)$},
    \end{align*}
    and similarly, by the backwards martingale convergence theorem,
    \begin{align*}
        \E_{k}[Z] \xrightarrow{k\to -\infty} 
        \E[Z|\cF_{B^c}]\qquad\text{in $L^2(\Xi,\cF,\P)$.}
    \end{align*}
    It follows by $L^2$-orthogonality of martingale increments that
    \begin{align}\label{e:601}
        \E\left[
            \bigl(Z - \E[Z|\cF_{B^c}]\bigr)^2
        \right]
        &= \sum_{k\in\Z}\E\left[
            \bigl(
                \E_{k}[Z] - \E_{k-1}[Z]
            \bigr)^2
        \right].
    \end{align}
    Observe that by \cref{l:FAFBisFAB} applied to $X=\E[Z|\cF_{B_k^c}]$,
    since $B_k^c \cap A_k = A_{k-1}$, we have
    \begin{align*}
        \E_{k-1}[Z] = \E_{k}\bigl[\E[Z|\cF_{B_k^c}]\bigr].
    \end{align*}
    Using this in \eqref{e:601} and applying Cauchy--Schwarz yields
        \begin{align*}
        \E\left[
            \bigl(Z - \E[Z|\cF_{B^c}]\bigr)^2
        \right]
        &\le \sum_{k\in\Z} \E\left[
            \left(
                \E_{k}\Bigl[
                    Z - \E[Z|\cF_{B_k^c}]
                \Bigr]
            \right)^2
        \right]
        \le \sum_{k\in\Z} \E\left[
            \Bigl(Z - \E[Z|\cF_{B_k^c}]
            \Bigr)^2
        \right],
    \end{align*}
    which is exactly \eqref{e:600}.
\end{proof}

\addtocontents{toc}{\SkipTocEntry}
\subsection*{Automorphism covariance of resampling}
We now lift the automorphisms introduced in \cref{s:automorphisms} to the space $(\Xires,\Fres,\Pres)$ defined in \eqref{e:def-Xi-resample}. This will help us freely apply the scaling invariance enjoyed by CDRP in our arguments using resampling. 
Let $\autaff : \R^2\to\R^2$ be any of the affine maps in \eqref{e:affine-maps}, let $\aut$ be the corresponding unitary automorphism of $L^2(\R^2)$ as in \eqref{e:automorphisms}, and let $\Aut$ be the corresponding measure-preserving automorphism of $(\Xi,\cF,\P)$ as in \eqref{e:automorphism-induced}.
Then $\Aut$ induces a measure-preserving automorphism of $(\Xires,\Fres,\Pres)$ 
by acting coordinate-wise: $(\Aut,\Aut):(\xi_1,\xi_2) \mapsto (\Aut\xi_1, \Aut\xi_2)$.
The following lemma tracks how $(\Aut,\Aut)$ acts on the coupled white noises defined in \eqref{e:def-xi-eta}. 
At the distributional level, resampling $\xi|_B$ and then applying $\Aut$ has the same effect as resampling $\xi|_{\autaff^{-1}(B)}$.
\begin{lemm}\label{l:automorphism-action-resampled}
    Fix $B\in\cB(\R^2)$ and let $\xi,\eta_B : \Xires \to \Xi$ be the white noises defined in \eqref{e:def-xi-eta}.
    Then for $\autaff, \aut,\Aut$ as above, we have
    \begin{align}\label{e:scaled-coupling-1}
        \Aut \circ \xi = 
        \xi \circ (\Aut,\Aut)
        \qquad\text{and}\qquad
        \Aut \circ \eta_B = \eta_{\autaff^{-1}(B)}\circ (\Aut,\Aut),
    \end{align}
    where $\autaff^{-1}(B)$ denotes the image of $B$ under the map $\autaff^{-1}:\R^2\to\R^2$.

    In particular, since $(\Aut,\Aut)$ preserves $\Pres$, we have the equality of joint distributions
    \begin{align}\label{e:scaled-coupling-2}
        (\Aut\circ \xi,\, \Aut\circ \eta_B) \law
        (\xi,\, \eta_{\autaff^{-1}(B)}).
    \end{align}
\end{lemm}
\begin{proof}
It suffices to prove \eqref{e:scaled-coupling-1}.
For the identity involving $\xi$, we have for all $(\xi_1,\xi_2)\in\Xires$ that
\begin{align*}
    (\Aut\circ \xi)(\xi_1,\xi_2) = \Aut\xi_1 = \xi(\Aut\xi_1,\Aut\xi_2).
\end{align*}
For the identity involving $\eta_B$,
we first note that $P_B \aut^{-1} = \aut^{-1} P_{\autaff^{-1}(B)}$ 
as operators on $L^2(\R^2)$, where recall $P_B$ is the orthogonal projection onto $L^2_0(B)$.
Indeed, for any $f\in L^2(\R^2)$, we have
\begin{align*}
    \aut^{-1} P_{\autaff^{-1}(B)} f
    &= \aut^{-1}(f\1_{\autaff^{-1}(B)})\\
    &= 
    \sqrt{\det\autaff^{-1}}
    \cdot
    (f\circ \autaff^{-1}) \cdot 
    (\1_{\autaff^{-1}(B)}\circ \autaff^{-1})\\
    &= (\aut^{-1}f) \cdot \1_{B}\\
    &= P_B\aut^{-1}f.
\end{align*}
Using this, we get
\begin{align*}
    (\Aut\circ \eta_B)(\xi_1,\xi_2)
    &= \Aut\left(
        \xi_1\circ P_{B^c} + \xi_2\circ P_{B}
    \right)\\
    &= \xi_1\circ (P_{B^c} \aut^{-1})
    + \xi_2\circ (P_{B} \aut^{-1})\\
    &= \xi_1\circ  (\aut^{-1} P_{\autaff^{-1}(B^c)})
    + \xi_2 \circ (\aut^{-1} P_{\autaff^{-1}(B)})\\
    &= (\Aut \xi_1)\circ P_{(\autaff^{-1}(B))^c}
    + (\Aut\xi_2)\circ P_{\autaff^{-1}(B)}\\
    &= \eta_{\autaff^{-1}(B)} (\Aut \xi_1, \Aut \xi_2)
\end{align*}
for all $(\xi_1,\xi_2)\in\Xires$ as desired.
\end{proof}

\section{Directed landscape, CDRP free energy profile, and polymer geometry}\label{s:LE}

In this section we import from the literature several technical results on the directed landscape, the CDRP free energy profile, and polymer geometry at zero and positive temperature.
These results will be used throughout the proofs of Theorems \ref{t:main} and \ref{t:black-noise}.
This section is intended as a reference for inputs; for a smoother reading experience, the reader may wish to skip ahead to \cref{s:mainproof} where the main argument is implemented.

To help navigate this section, we begin with an outline of its content.
\begin{itemize}
    \item Basic properties of the directed landscape appear in \cref{p:landscape-symmetries}.
    \item We record uniform pointwise bounds for the directed landscape in \cref{p:dov-pointwise}, and for its spatial modulus of continuity in \cref{p:dov-holder-general}.
    \item \cref{l:fourth-moment} states a fourth moment bound for the CDRP free energy profile. 
    \item We establish an anticoncentration result for the directed landscape geodesic in \cref{ss:locbrown} (see \cref{p:geodesic-deloc}) using the local Brownianity of the Airy line ensemble.
    \item Inputs on polymer geometry are recorded in \cref{s:polymer-geometry}, with a more detailed outline appearing there.
\end{itemize}

We first record several basic properties of the directed landscape $\cL$, some of which appeared already in \cref{d:dl} but are repeated here for ease of referencing.

\begin{prop}[{\cite[Definition 10.1 and Lemma 10.2]{dov}, \cite[Proposition 1.23]{dv2}}]\label{p:landscape-symmetries}
    The directed landscape $\cL$ satisfies the following properties.
    \begin{enumerate}[label={\rm(\roman*)}]
        \item\textup{(Independent increments).}\label{landscale-independent-increments}
        For any finite collection of pairwise disjoint time intervals $(s_1,t_1),\dots,(s_k,t_k)$, the random functions
        \[
            \cL(s_1,\smallbullet; t_1, \smallbullet), \dots, \cL(s_k,\smallbullet; t_k, \smallbullet)
        \]
        are independent.
        \item\textup{(Metric composition).}\label{landscape-metric-composition}
        Almost surely, for any $(s,x;t,y)\in\Rup$ and any $r\in (s,t)$, we have
        \begin{align*}
            \cL(s,x;t,y) = \sup_{z\in \R}\left(\cL(s,x;r,z) + \cL(r,z;t,y)\right).
        \end{align*}
        In particular, we have the \emph{reverse triangle inequality} $\cL(s,x;t,y) \ge \cL(s,x;r,z)+\cL(r,z;t,y)$ for all $r\in (s,t)$ and all $z\in\R$.
        \item\textup{(Translation invariance).}\label{landscape-translation}
        For all $(r,w)\in\R^2$, we have the following equality in distribution of $C(\Rup)$-valued random variables:
        \begin{align*}
            \cL(s,x;t,y) \law \cL(s+r, x+w; t+r, y+w).
        \end{align*}

        \item\textup{(Scaling invariance).}\label{landscape-scaling}
        For all $\lambda>0$, we have the
        following equality in distribution of $C(\Rup)$-valued random variables:
        \begin{align*}
            \cL(s,x;t,y) \law \lambda^{-1/3}\cL(\lambda s, \lambda^{2/3}x; \lambda t, \lambda^{2/3}y).
        \end{align*}

        \item\textup{(Shear invariance).}\label{landscape-shear}        
        For all $\nu \in \R$, we have the following equality in distribution of $C(\Rup)$-valued random variables:
        \begin{align*}
            \cL(s,x;t,y)
            \law \cL(s, x+\nu s; t, y+\nu t)
            + \frac{(x-y-\nu(t-s))^2 - (x-y)^2}{t-s}.
        \end{align*}
        As a special case, for all fixed $s<t$ and $x\in\R$, we have the following equalities in distribution of $C(\R)$-valued random variables:
        \begin{align*}
            \cL(s,\smallbullet;t,x)
            \law 
            \cL(s,x;t,\smallbullet) 
            \law \cL(s, 0; t, \smallbullet)
            + \frac{2x}{t-s}\smallbullet
            \;
            - \frac{x^2}{t-s}.
        \end{align*}
    \end{enumerate}
\end{prop}

We next record two regularity estimates for the directed landscape, which are specializations of results of \cite{dov}.
These will later be combined with \cref{t:KPZ-to-landscape} and the Portmanteau lemma to derive positive-temperature analogues.

\begin{prop}[Uniform pointwise bound, {\cite[Proposition 10.6]{dov}}] \label{p:dov-pointwise}
    For all $L\ge 1$, there exists a random constant $\Cptwise>0$ 
    such that the following is true.
    First, we have
    \begin{align*}
        \P(\Cptwise > m) \le C L^{12} e^{-c m^{3/2}}
        \qquad\text{for all $m\ge 0$},
    \end{align*}
    where $C,c>0$ are universal constants.
    Second, for all $(s,x;t,y)\in [-2L,2L]^4\cap\Rup$, we have
    \begin{align*}
        \left|\cL(s,x;t,y) + \frac{(y-x)^2}{t-s}
        \right|
        \le \Cptwise
        (t-s)^{1/3}\log^{4/3}\left(\frac{8L}{t-s}\right).
    \end{align*}
\end{prop}

\begin{prop}[Spatial H\"older continuity, {\cite[Proposition 10.5]{dov}}]\label{p:dov-holder-general}
    For all $L\ge 1$ and $\epsilon\in(0,1)$,
    there exists a random constant $\CHolder>0$ such that the following holds.
    First, we have
    \begin{align*}
        \P(\CHolder > m) \le CL^{10} \epsilon^{-6} e^{-cm^{3/2}}
        \qquad\text{for all $m\ge 0$},
    \end{align*}
    where $C,c>0$ are universal constants.
    Second, for all $x,x',y,y'\in [-2L,2L]$ and all $-2L<s<t<2L$ with $t-s \ge \epsilon$, we have
    \begin{multline*}
        \left|
            \cL(s,x;t,y) - \cL(s,x';t,y')
            + \frac{(y-x)^2 - (y'-x')^2}{t-s}
        \right|\\
        \le \CHolder \norm{(x,y)-(x',y')}_2^{1/2} \log^{1/2}\left(\frac{8L}{\norm{(x,y)-(x',y')}_2}\right).
    \end{multline*} 
\end{prop}

For later convenience, we record a special case of \cref{p:dov-holder-general}:

\begin{cor}\label{c:dov-holder}
    Fix $L\ge 1$ and $t_0\in (0,1]$, and fix $\gamma\in(0,\frac12)$.
    For all sufficiently small $\e>0$ (depending only on $L,t_0,\gamma$), all $x\in[-L,L]$, and all {$r\in [\e^\gamma, t_0]$}, we have
    \begin{align*}
        \P\left(
            \sup_{\substack{z_1,z_2 \in [-L,L]\\|z_1-z_2| \le 8\e\log^4(1/\e)}}    
            \left|\cL(0,x;r,z_1) - \cL(0,x;r,z_2)\right|
            \ge
            2^{1/3}\e^{1/2}\log^4(1/\e)
        \right) \le  Ce^{-c\log^{9/4}(1/\e)},
    \end{align*}
    where $C,c>0$ depend only on $L,\gamma$.
    Similarly, for all $r\in[0, t_0-\e^\gamma]$,
    \begin{align*}
        \P\left(
            \sup_{\substack{z_1,z_2 \in [-L,L]\\|z_1-z_2| \le 8\e\log^4(1/\e)}}    
            \left|\cL(r,z_1; t_0,x) - \cL(r,z_2;t_0,x)\right|
            \ge
            2^{1/3}\e^{1/2}\log^4(1/\e)
        \right) \le  Ce^{-c\log^{9/4}(1/\e)}.
    \end{align*}
\end{cor}
\begin{proof}
    We only prove the first estimate; the second follows by a symmetric argument.
    Fix $r\in [\e^\gamma,t_0]$.
    Note that 
    \begin{align*}
        \sup_{\substack{z_1,z_2 \in [-L,L]\\|z_1-z_2| \le 8\e\log^4(1/\e)}}
        \left|
            \frac{(z_1-x)^2 - (z_2-x)^2}{r}
        \right|
        &\ls L\e^{1-\gamma}\log^4(1/\e).
    \end{align*}
    Since $1-\gamma > 1/2$, we deduce the following probability bound for all sufficiently small $\e>0$ (depending only on $L,\gamma$):
    \begin{multline*}
        \P\left(
            \sup_{\substack{z_1,z_2 \in [-L,L]\\|z_1-z_2| \le 8\e\log^4(1/\e)}}   
            \left|\cL(0,x;r,z_1) - \cL(0,x;r,z_2)\right|
            \ge
            2^{1/3}\e^{1/2}\log^4(1/\e)\right)\\
        \le
        \P\left(
            \sup_{\substack{z_1,z_2 \in [-L,L]\\|z_1-z_2| \le 8\e\log^4(1/\e)}}    
            \left|\cL(0,x;r,z_1) - \cL(0,x;r,z_2) 
            + \frac{(z_1-x)^2-(z_2-x)^2}{r}\right|
            \ge
            \e^{1/2}\log^4(1/\e)
        \right).
    \end{multline*}
    By \cref{p:dov-holder-general}, there exist $C,c>0$ (depending on $L$) such that the RHS is upper bounded by $C \e^{-6\gamma}\,e^{-c \log^{9/4}(1/\e)}$ for all sufficiently small $\e>0$.
    This completes the proof.
\end{proof}

The next lemma, a crude fourth moment bound for the CDRP free energy, is a quick corollary of the latter's stretched exponential tails established in  \cite{CG20a,CG20b}.
\begin{lemm}[Fourth moment of CDRP free energy]\label{l:fourth-moment}
    For any $x,y\in\R$, we have
    \begin{align*}
        \sup_{n\ge 1}
        \E\left[
            \left|
                \frac{1}{n^{1/3}}\left(\cH^{\xi}_1(0,xn^{2/3};n,yn^{2/3})
                + \frac{n}{24}\right)
            \right|^4
        \right]
        <\infty,
    \end{align*}
    where $\cH^\xi$ was defined in \eqref{e:KPZ-cole-hopf}.
\end{lemm}
\begin{proof}
    By translation and shear invariance (\cref{t:AJRS-Z-scaling}\ref{property-translation}--\ref{property-shear}) and the inequality $|a-b|^4 \le 2^{3}(|a|^4 + |b|^4)$,
    \begin{align*}
        \E\left[
            \left|
                \frac{1}{n^{1/3}}\left(\cH^{\xi}_1(0,xn^{2/3};n,yn^{2/3})
                + \frac{n}{24}\right)
            \right|^4
        \right]
        &=\E\left[
            \left|
                \frac{1}{n^{1/3}}\left(\cH^{\xi}_1(0,0;n,0)
                + \frac{n}{24}
                -\frac{(y-x)^2 n^{1/3}}{2}
                \right)
            \right|^4
        \right]\\
        &\ls
        \E\left[
            \left|
                \frac{1}{n^{1/3}}\left(\cH^{\xi}_1(0,0;n,0)
                + \frac{n}{24}\right)
            \right|^4
        \right]
        + (y-x)^8.
    \end{align*}
    By \cite[Theorem 1.1]{CG20b} and \cite[Theorem 1.11]{CG20a}, there exist constants $C,c>0$ such that for all $n\ge 1$ and all $m\ge 0$,
    \begin{align*}
        &\P\left(
            \frac{1}{n^{1/3}}\left(\cH^{\xi}_1(0,0;n,0)
            + \frac{n}{24}\right)
            > m
        \right)
        \le Ce^{-c m^{3/2}}, \qquad\text{and}\\
        &\P\left(
            \frac{1}{n^{1/3}}\left(\cH^{\xi}_1(0,0;n,0)
                + \frac{n}{24}\right)
            < -m
        \right)
        \le Ce^{-cm^{5/2}}.
    \end{align*}
    This immediately implies the lemma.
\end{proof}

\subsection{Strong Brownian resemblance of energy profiles}\label{ss:locbrown}
As indicated in \cref{ss:idea}, estimating probabilities of polymers passing through small intervals can be reduced to Brownian computations, thanks to the local Brownianity of the Airy/KPZ line ensembles.
This has been quantified in increasingly strong forms in \cite{CH,Ham22,CHH,duncan2}; we will rely on the main result of \cite{CHH}, which we now formulate (stronger statements were proved in \cite{duncan2}, and these would have worked equally well for our purposes).
We use the shorthand
\begin{align}
    \cL_{1}(z) \coloneqq  \cL(0,0;1,z),
    \qquad z\in\R.
\end{align}
By \cite[Definition 10.1]{dov}, the process $\cL_{1}$ is equal in distribution to the top curve of the Airy line ensemble \cite{CH}.
The process $\cL_{1}$ is also known as the \emph{parabolic Airy$_2$ process}, as mentioned in \cref{ss:preface}.

\begin{prop}[Local Brownianity, {\cite[Theorem 1.1]{CHH}}]\label{p:CHH-local-brownian}
    Let $\cL_{1}$ and $\cL_{1}'$ be independent parabolic Airy$_2$ processes, and let $B$ and $B'$ be independent Brownian motions of variance $2$ started from $B(-1)=B'(-1)=0$.
    There exist $C,C'>0$
    such that for any Borel set $S\in\cB(C([-1,1]))^{\otimes 2}$, writing $p \coloneqq  \P((B,B') \in S)$, we have
    \begin{align*}
        \P\Bigl(
            \bigl(\cL_{1}(\smallbullet) - \cL_{1}(-1),\;\; \cL_{1}'(\smallbullet) - \cL_{1}'(-1)\bigr)
            \in S
        \Bigr) 
        \le 
        C p e^{C'\log^{5/6}(2/p)},
    \end{align*}
    where if $p=0$ then the RHS is $0$ by convention. 
\end{prop}
\begin{proof}
    By \cite[Theorem 1.1]{CHH}, for any $S'\in \cB(C([-1,1]))$, writing $p'\coloneqq \P(B\in S')$, we have $\P(\cL_{1}(\cdot) - \cL_{1}(-1) \in S') \le Cp'e^{C'\log^{5/6}(1/p')}$.
    This implies a similar estimate for the product measure, possibly with larger $C,C'$.
    We omit the details. 
\end{proof}

The next lemma is a crude anticoncentration bound for the maximum of a Brownian motion with drift.

\begin{lemm}[Maximum of Brownian motion with drift]\label{l:BM-maximum-drift}
    Let $W$ be a Brownian motion of variance $4$ started from $W(0)=0$.
    Then for all $T>0$, all $\kappa>0$, and all $\mu\in\R$,
    \begin{align*}
        \P\left(
            \max_{z\in[0,T]}\left(W(z) + \mu z\right) \le \kappa
        \right)
        \le C\left(\frac{1}{\sqrt{T}} + |\mu|\right)\kappa,
    \end{align*}    
    where $C>0$ is a universal constant.
\end{lemm}
\begin{proof}
    By Brownian scaling, it suffices to consider the case where $T=1$ and $W$ has variance $1$.
    By the reflection principle for Brownian motion with drift
    (e.g. \cite[Part II, Section 2, Table 1.1.4]{BS96} or \cite[Theorem 5]{Low23}),
    \begin{align*}
        \P\left(
            \max_{z\in[0, 1]} \left(W(z) + \mu z\right) \le \kappa
        \right)
        &=
        1 - \P(W(1) > \kappa - \mu) - e^{2\mu\kappa} \P(W(1) > \kappa + \mu)\\
        &=
        \Phi(\kappa - \mu)
        - e^{2\mu\kappa}
        \Phi(-\kappa - \mu),
    \end{align*}
    where $\Phi$ is the standard Gaussian distribution function.
    Using that $e^{2\mu\kappa} \ge 1 + 2\mu\kappa$, we get
    \begin{align*}
        \Phi(\kappa - \mu)
        - e^{2\mu\kappa}
        \Phi(-\kappa - \mu)
        &\le
        \bigl(\Phi(\kappa - \mu)
        - 
        \Phi(-\kappa - \mu)\bigr)
        -
        2\mu\kappa
        \Phi(-\kappa - \mu)\\
        &\le \frac{1}{\sqrt{2\pi}}\cdot 2\kappa
        + 2|\mu|\kappa,
    \end{align*}
    by bounding the Gaussian density pointwise.
    This completes the proof.
\end{proof}

We now combine \cref{p:CHH-local-brownian} and \cref{l:BM-maximum-drift} to establish an anticoncentration result for the directed landscape.
As mentioned, this implies the delocalization of the geodesic at a given time (geodesics are defined later in \cref{def:geodesics}).

\begin{prop}\label{p:geodesic-deloc}
    Fix $L\ge 1$ and $\ell>0$ and $t_0\in(0, 1]$ and $\gamma\in(0,\frac13)$.
    There exist $\e_0=\e_0(L,\ell,\gamma)>0$ and $C=C(L,\gamma)>0$ such that the following holds.

    For all $0 < \e < \min\{\e_0, (t_0/4)^{1/\gamma}\}$, all $x,y\in [-L,L]$, and all $r\in[\e^\gamma, t_0-\e^\gamma]$, we have
    \begin{multline*}
        \P\left(
            \sup_{|z|\le \ell}
            \bigl[\cL(0,x;r,z) + \cL(r,z;t_0,y)\bigr]
            \le
            \sup_{|z|\le 2\e\log^4(1/\e)}
            \bigl[
                \cL(0,x;r,z) + \cL(r,z;t_0,y)
            \bigr]
            +
            \e^{1/2} \log^5(1/\e)
        \right)\\
        \le 
        C \e^{1-3\gamma}.
    \end{multline*}
    Note that $[\e^\gamma, t_0-\e^\gamma]$ is non-empty, since by hypothesis $\e^\gamma <  t_0/4$.
\end{prop}
\begin{proof}
    We first use the spatial H\"older regularity of $\cL$ to replace the supremum over $|z|\le 2\e\log^4(1/\e)$ with a one-point quantity.
    By a union bound and \cref{c:dov-holder}, we have
    \begin{multline*}
        \P\left(
            \sup_{|z|\le 2\e\log^4(1/\e)}
            \bigl[\cL(0,x; r, z)
            + \cL(r,z;t_0,y)\bigr]
            \ge
            \bigl[\cL(0,x;r,0) + \cL(r,0;t_0,y)\bigr]
            + \e^{1/2}\log^4(1/\e)
        \right)\\
        \begin{aligned}
            &\le
            \P\left(
                \sup_{|z|\le 2\e\log^4(1/\e)}
                \cL(0,x; r, z)
                \ge
                \cL(0,x;r,0)
                + 
                \frac12
                \e^{1/2}\log^4(1/\e)
            \right)\\
            &\quad+
            \P\left(
                \sup_{|z|\le 2\e\log^4(1/\e)}
                \cL(r,z; t_0, y)
                \ge
                \cL(r,0;t_0,y)
                + \frac12
                \e^{1/2}\log^4(1/\e)
            \right)\\
            &\le 2Ce^{-c\log^{9/4}(1/\e)},
        \end{aligned}
    \end{multline*}
    where $C,c>0$ depend on $L,\gamma$.
    Therefore, denoting by $\sE$ the event in the statement of the proposition,
    by a union bound and the crude bound $\e^{1/2}\log^4(1/\e) + \e^{1/2}\log^5(1/\e) \le \e^{1/2}\log^6(1/\e)$, we have
    \begin{multline}\label{e:420}
        \P(\sE)
        \le \P\left(
            \sup_{|z|\le \ell}\bigl[
                \cL(0,x;r,z)
                +
                \cL(r,z;t_0,y)
            \bigr]
            \le
            \bigl[\cL(0,x;r,0) + 
            \cL(r,0;t_0,y)\bigr]
            +
            \e^{1/2}\log^6(1/\e)
        \right)\\
        + 2Ce^{-c\log^{9/4}(1/\e)}.
    \end{multline} 
    We denote by $\sF$ the event inside the probability on the RHS of \eqref{e:420}, and turn towards bounding it.

    By \cref{p:landscape-symmetries}\ref{landscape-shear}, we have the following two distributional equalities
    \begin{align*}
        \cL(0,x;r,\smallbullet) 
        &\law \cL(0,0;r,\smallbullet) + \frac{2x}{r}\smallbullet\; -\frac{x^2}{r}
        \\
        \cL(r,\smallbullet;t_0,x) 
        &\law 
        \cL(r,0;t_0, \smallbullet)
        + \frac{2x}{t_0-r}\smallbullet\; - \frac{x^2}{t_0-r},
    \end{align*}
    and in particular, it also holds that 
    \begin{equation}\label{e:404}
        \begin{aligned}
            \cL(0,x;r,\smallbullet) - 
        \cL(0,x;r,0)
        &\law
        \left[
            \cL(0,0;r,\smallbullet) + \frac{2x}{r}\smallbullet
        \right]
        -\cL(0,0;r,0),
             && 
        \\
        \cL(r,\smallbullet;t_0,x)  - 
        \cL(r,0;t_0,x) 
        &\law
        \left[
            \cL(r,0;t_0, \smallbullet)
            + \frac{2x}{t_0-r}\smallbullet
        \right]
        -\cL(r,0;t_0,0).
        \end{aligned}
    \end{equation}
    We denote
    \begin{align*}
        \mu \coloneqq  \frac{2x}{r} + \frac{2y}{t_0-r}.
    \end{align*}
    Since $r\in[\e^\gamma, t_0 - \e^\gamma]$, we have $|\mu| \ls \e^{-\gamma}$.
    Further, by \cref{p:CHH-local-brownian}, the sum of the two processes in \eqref{e:404} locally resembles a Brownian motion with drift $\mu$.
    The latter resembles a driftless Brownian motion on the interval $[-\e^{2\gamma}, \e^{2\gamma}]$, since $\e^{2\gamma}|\mu| \ls \sqrt{\e^{2\gamma}}$.
    Therefore, assuming $\e_0$ is small enough that $\e_0^{2\gamma} < \ell$, we are led to make the following estimate:
    \begin{align}\label{e:402}
        \P(\sF)
        &= \P\left(
            \sup_{|z|\le \ell}\bigl[
                \cL(0,x;r,z)
                +
                \cL(r,z;t_0,y)
            \bigr]
            -
            \bigl[\cL(0,x;r,0) + 
            \cL(r,0;t_0,y)\bigr]
            \le
            \e^{1/2}\log^6(1/\e)
        \right)\nonumber\\
        &\le
         \P\left(
            \sup_{|z|\le \e^{2\gamma}}\bigl[
                \cL(0,x;r,z)
                +
                \cL(r,z;t_0,y)
            \bigr]
            -
            \bigl[\cL(0,x;r,0) + 
            \cL(r,0;t_0,y)\bigr]
            \le
            \e^{1/2}\log^6(1/\e)
        \right)\nonumber\\
        \overset{\eqref{e:404}}&{=}
            \mathbb{P}\Biggl( 
            \sup_{|z|\le \e^{2\gamma}}
            \bigl[
                \cL(0,0; r ,z) 
                + \cL(r,0;t_0,z) 
                + \mu z
            \bigr]
            -\bigl[
                \cL(0,0;r, 0)
                + \cL(r,0; t_0, 0)
            \bigr]  
            \le \e^{1/2}\log^6(1/\e)
        \Biggr)\nonumber
        \\
        &= 
        \P\Biggl(
            \sup_{|z|\le \e^{2\gamma}}
            \left[
            r^{1/3}\cL_{1}\left(
            \frac{z}{r^{2/3}}
            \right)
            + (t_0-r)^{1/3}\cL_{1}'\left(
                \frac{z}{(t_0-r)^{2/3}}
            \right)
            + \mu z
        \right]
        \nonumber\\
        &\qquad\qquad\qquad\qquad\qquad\qquad\qquad
        -
        \left[
            r^{1/3}\cL_{1}(0)
            + (t_0-r)^{1/3}\cL_{1}'(0)
        \right]
        \le \e^{1/2}\log^6(1/\e)
        \Biggr),
    \end{align}
    where $\cL_{1},\cL_{1}'$ are independent parabolic Airy$_2$ processes (recall that $\cL(0,0;r,\smallbullet)$ and $\cL(r, \smallbullet; t_0,0)$ are independent by \cref{p:landscape-symmetries}\ref{landscale-independent-increments}),
    and where in the last line we used the scaling property of the directed landscape (\cref{p:landscape-symmetries}\ref{landscape-scaling}).

    Note that since $r\in [\e^\gamma, t_0-\e^\gamma]$, the event inside the probability in \eqref{e:402} depends only on $\{\cL_{1}(z)-\cL_{1}(-1)\}_{z\in [-1,1]}$ and $\{\cL_{1}'(z)-\cL_{1}'(-1)\}_{z\in [-1,1]}$ (in fact, it depends only on their restrictions to $z\in [-\e^{4\gamma/3}, \e^{4\gamma/3}]$).
    Therefore,  letting $B$ and $B'$ be independent Brownian motions of variance $2$ started from $B(-1)=B'(-1)=0$, and writing
    \begin{multline*}
        p\coloneqq 
        \P\Biggl(
            \sup_{|z|\le \e^{2\gamma}}
            \left[
            r^{1/3}B\left(
            \frac{z}{r^{2/3}}
            \right)
            + (t_0-r)^{1/3}B'\left(
                \frac{z}{(t_0-r)^{2/3}}
            \right)
            + \mu z
        \right]\\
        -
        \left[
            r^{1/3}B(0)
            + (t_0-r)^{1/3}B'(0)
        \right]
        \le \e^{1/2}\log^6(1/\e)
        \Biggr),
    \end{multline*}
    we can apply \cref{p:CHH-local-brownian} to obtain
    \begin{align}\label{e:406}
        \P(\sF)
        \le C p e^{C'\log^{5/6}(2/p)}.
    \end{align}
    It remains to bound $p$. By Brownian scaling, 
    \begin{align*}
        p &= \P\left(
            \sup_{|z|\le\e^{2\gamma}}
            \left[
                B(z)
                + B'(z) + \mu z
            \right]
            - \left[
               B(0) + B'(0)
            \right]
            \le \e^{1/2}\log^6(1/\e)
        \right)\\
        &= \P\left(
            \sup_{|z|\le\e^{2\gamma}}
            \bigl[W(z) + \mu z\bigr]
            - W(0)
            \le \e^{1/2}\log^6(1/\e)
        \right),
    \end{align*}
    where $W \coloneqq  B+B'$ is a Brownian motion of variance $4$ started from $W(-1)=0$ (by independence of $B,B'$).
    Write $W_\mu(z) \coloneqq  W(z)+\mu z$ and $\kappa \coloneqq  \e^{1/2} \log^6(1/\e)$.
    Then by independence of Brownian increments,
    \begin{align*}
        p
        &= \P\left(
            \sup_{|z|\le \e^{2\gamma}}
            W_\mu(z)
            - W_\mu(0)
            \le 
            \kappa
        \right)\\
        &=
        \P\left(
            \sup_{z\in[-\e^{2\gamma}, 0]}
            W_\mu(z) - W_\mu(0)
            \le 
            \kappa
        \right)
        \P\left(
            \sup_{z\in[0, \e^{2\gamma}]}
            W_\mu(z) - W_\mu(0)
            \le 
            \kappa
        \right)\\
        &=
        \P\left(
            \sup_{z\in[-\e^{2\gamma}, 0]}
            -W_\mu(z)
            \le 
            \kappa
            \,\middle|\, W_\mu(-\e^{2\gamma})=0
        \right)
        \P\left(
            \sup_{z\in[0, \e^{2\gamma}]}
            W_\mu(z)
            \le 
            \kappa
            \,\middle|\, W_\mu(0)=0
        \right)\\
        &=
        \P\left(
            \sup_{z\in[0, \e^{2\gamma}]}
            (W(z) - \mu z)
            \le 
            \kappa
        \right)
        \P\left(
            \sup_{z\in[0, \e^{2\gamma}]}
            (W(z) + \mu z)
            \le 
            \kappa
        \right),
    \end{align*}
    where in the third line we reversed time inside the first probability,
    and where in the last line $W$ denotes a Brownian motion of variance $4$ started from $W(0)=0$.
    Bounding the above two probabilities using \cref{l:BM-maximum-drift}, we get
    \begin{align*}
        p
        &\ls\left(\e^{-\gamma} + |\mu|\right)^2 
        \kappa^2
        \\
        &\ls L^2
        \e^{1 - 2\gamma}\log^{12}(1/\e),
    \end{align*}
    since $\kappa=\e^{1/2}\log^6(1/\e)$ and $|\mu| \le 4L\e^{-\gamma}$.

    Finally, by \eqref{e:420} and \eqref{e:406}, we get that for all sufficiently small $\e>0$ (depending only on $L,\gamma,t_0$),
    \begin{align*}
        \P(\sE)
        &\ls 
        L^{2}
        \e^{1-2\gamma - o(1)}
        + Ce^{-c\log^{9/4}(1/\e)}\\
        &\ls_{L,\gamma} \e^{1-3\gamma}.
    \end{align*}
    This completes the proof of \cref{p:geodesic-deloc}.
\end{proof}

\subsection{Polymer geometry}\label{s:polymer-geometry}

In this subsection we import from the literature several technical results on zero- and positive-temperature polymers that will be used in the proofs of Theorems \ref{t:main} and \ref{t:black-noise}.
Here is an outline:
\begin{itemize}
    \item We define geodesics in the directed landscape in Definitions \ref{def:length-dl} and \ref{def:geodesics}, and record a scaling property of geodesics in \cref{p:DOV-geodesic-holder}.
    \item We record transversal fluctuation estimates for zero-temperature polymers in Propositions \ref{p:DOV-geodesic-holder} and \ref{p:GZ-uniform-TF},
    and for positive-temperature polymers in \cref{p:polymer-TF-DZ}.
    \item We record a result of \cite{DZ24} on the weak convergence of the annealed CDRP measure to the directed landscape geodesic in \cref{t:pathtogeo}.
\end{itemize}

\subsubsection{Zero-temperature polymer (geodesic) geometry }

Let $\cL$ be the directed landscape defined in \cref{d:dl}.
We first recall the definition of geodesics in the directed landscape from \cite[Section 12]{dov}.

\begin{defn}[Length of a path in the directed landscape]\label{def:length-dl}
    For any $s<t$ and any continuous path $\pi:[s,t]\to \R$, we define the \emph{length} of $\pi$ by
    \begin{align*}
        \int_s^t d\cL\circ \pi \coloneqq  
        \inf_{k\in \N}
        \inf_{s=s_0 < s_1 < \cdots < s_k = t}
        \sum_{i=0}^{k-1} \cL(s_i, \pi(s_i); s_{i+1}, \pi(s_{i+1})). 
    \end{align*}
\end{defn}

\begin{defn}[Geodesics in the directed landscape]\label{def:geodesics}
    For $(s,x;t,y)\in\Rup$, a \emph{geodesic}  from $(s,x)$ to $(t,y)$ in the directed landscape
    is a continuous function $\Pi:[s,t]\to \R$ such that $\Pi(s)=x,\Pi(t)=y$, and
    \begin{align}\label{e:def-geo-1}
        \cL(s,x;t,y) = \int_s^t d\cL\circ \Pi.
    \end{align}
    Equivalently, $\Pi$ is a geodesic if for all $k\in\N$ and all $s=s_0<s_1<\cdots<s_k=t$,
    \begin{align}\label{e:def-geo-2}
        \int_s^t d\cL\circ \Pi
        =\sum_{i=0}^{k-1} \cL(s_i,\Pi(s_i); s_{i+1}, \Pi(s_{i+1})).
    \end{align}
    By \cite[Lemma 13.2]{dov}, it holds almost surely that for every $(s,x;t,y)\in\Rup$, there exists a \emph{rightmost} geodesic $\Pi_{(s,x),(t,y)}$ from $(s,x)$ to $(t,y)$, 
    in the sense that $\Pi_{(s,x),(t,y)} \ge \Pi'$ pointwise for any other geodesic $\Pi'$ from $(s,x)$ to $(t,y)$.
\end{defn}

Rightmost geodesics are unique by definition, and we work with them to remove ambiguity when speaking of ``the'' geodesic.
However, uniqueness is irrelevant for our arguments.

\begin{rem}[Directed landscape as an optimization problem over paths]
    Since a.s. there exists a geodesic for every $(s,x;t,y)\in\Rup$, by combining \eqref{e:def-geo-1} with \cref{p:landscape-symmetries}\ref{landscape-metric-composition} and \cref{def:length-dl}, we find that almost surely, for every $(s,x;t,y)\in\Rup$,
    \begin{align}\label{e:DL-as-LPP}
        \cL(s,x;t,y) = \sup_{\pi} \int_s^t d\cL\circ \pi,
    \end{align}
    where the supremum is over all paths $\pi\in C([s,t])$ with $\pi(s)=x$ and $\pi(t)=y$.
    This property will be used throughout the proof of \cref{t:black-noise}.
\end{rem}

We next recall basic scaling and transversal fluctuation properties of geodesics in the directed landscape.

\begin{prop}[Geodesic scaling and modulus of continuity, {\cite[Theorem 1.7]{dov}}]\label{p:DOV-geodesic-holder}
    For every $(s,x;t,y)\in\Rup$, we have the scaling relation
    \begin{align}\label{e:geodesic-KPZ-scaling}
        \left\{\Pi_{(s,x),(t,y)}(s + (t-s)r)\right\}_{r\in[0,1]}
        \law \left\{(t-s)^{2/3}\Pi_{(0,0),(1,0)}(r) + x + (y-x)r\right\}_{r\in[0,1]},
    \end{align}
    where the equality in distribution is as $C([0,1])$-valued random variables.
    For a proof of the above statement, see the discussion below \cite[Theorem 12.1]{dov}.

    Moreover, given $S\ge 1$ and $(s,x;t,y)\in\Rup$ with slope $\left|\frac{y-x}{t-s}\right| \le S$ and $t-s \le 1$,
    there is a random constant $\CTFloc > 0$ with the following two properties.
    First, we have
    \begin{align*}
        \P(\CTFloc > m) \le 
        Ce^{-c m^{3}}
        \qquad\text{for all $m\ge 0$},
    \end{align*}
    where $C,c>0$ are universal constants.
    Second, for all $r,r'\in [s,t]$, we have
    \begin{align}\label{e:geodesic-holder}
        \left|\Pi_{(s,x),(t,y)}(r)-\Pi_{(s,x),(t,y)}(r')\right|
        \le \CTFloc\, |r-r'|^{2/3}\log^{1/3}\left(2\vee \frac{1}{|r-r'|}\right)
        + S|r-r'|.
    \end{align}
    We emphasize that $\CTFloc$ depends on the point $(s,x;t,y)$.
\end{prop}

In addition to the above geodesic modulus of continuity bound, we will also need bounds for the maximum distance that a geodesic wanders from the straight line between its endpoints, with uniform control as the endpoints vary in a compact set. 
For this we use \cite[Lemma 3.11]{GZ22}, stated here in a weaker form tailored to our applications.
\begin{prop}[Geodesic maximum transversal fluctuation, {\cite[Lemma 3.11]{GZ22}}]\label{p:GZ-uniform-TF}
    Given $L\ge 1$, there exists a random constant $\CTF>0$ such that the following is true.
    First, we have 
    \begin{align*}
        \P(\CTF > m) \le Ce^{-cm^{2}}
        \qquad\text{for all $m\ge 0$,}
    \end{align*}
    where $C,c>0$ are constants depending only on $L$.
    Second, for any $(s,x;t,y)\in\Rup\cap[-L,L]^4$, the geodesic $\Pi_{(s,x),(t,y)}$ satisfies
    \begin{align*}
        \sup_{r\in[s,t]}
        \left|
            \Pi_{(s,x),(t,y)}(r) - \frac{x(t-r) + y(r-s)}{t-s}
        \right|
        \le \CTF \,
        (t-s)^{2/3} \log^3\left(2 \vee \frac{1}{t-s}\right).
    \end{align*}
\end{prop}

\subsubsection{Positive-temperature polymer geometry}

We turn now to the positive temperature model, quoting two results of \cite{DZ24} adapted to our setting. The following result bounds the transversal fluctuations of the (quenched) polymer.

\begin{prop}[{Polymer transversal fluctuations, \cite[Corollary 3.5]{DZ24}}]\label{p:polymer-TF-DZ}
    Fix $L\ge 1$ and $T\ge 1$ and 
    $(s,x;t,y)\in\Rup$ with $x,y\in [-L,L]$ and $t-s \le T$.
    There exist $C_1,C_2>0$ depending on $L,T$ such that the following holds.
    
    Fix $\beta\ge (t-s)^{-1/4}$ and let $X\sim \CDRP^{\xi}_\beta(s, x\beta^{2/3}; t,y\beta^{2/3})$.
    Then for all $m\ge 0$,
    \begin{align}\label{e:301}
        \mathbb{P}\left(\mathbb{P}^{\xi}_{\beta,(s,x\beta^{2/3}),(t,y\beta^{2/3})}
        \left(
            \sup_{r\in[s,t]} \left|X(r)\right| > m \beta^{2/3}
        \right)\ge C_{1}e^{-m^{2}/C_{1}}\right)
        \le C_{2}e^{-m^{3}/C_{2}}.
    \end{align} 

    For later convenience, we also record a slight reformulation of the above result.
    For all $\e\in(0,1)$, all $z\in[-\e,\e]$, all $n\ge \e^{-3/2}$, and all $m\ge 0$,
    \begin{align}\label{e:302}
        \P\left(
            \P^\xi_{1,(sn, zn^{2/3}),((s+\e^{3/2})n, zn^{2/3})}
            \left(
            \sup_{r\in[s,s+\e^{3/2}]}\left|
                n^{-2/3}
                X(rn)
            \right|
            \ge  m \e
            \right)
            \ge C_1 e^{-m^2/C_1}
        \right)
        \le C_2 e^{-m^3/C_2}.
    \end{align}

\end{prop}
\begin{proof}
    We will use \cite[Corollary 3.5]{DZ24} which states that
    for all $\epsilon \in (0,1]$ and all 
    $y\in[-L,L]$,
    \begin{align}\label{e:DZ-bound}
        \P\left(
            \P^\xi_{1, (0,0),(\epsilon^{-1},y \epsilon^{-2/3})}
            \left(
                \sup_{r\in[0,1]}
                \left|\epsilon^{2/3} X(\epsilon^{-1}r)\right|
                > m
            \right)
            \ge C_1 e^{-m^2/C_1}
        \right)
        \le C_2 e^{-m^3/C_2}
    \end{align}
    for some $C_1,C_2>0$ depending only on $L$.
    By translation invariance we can assume $s=0$.
    For \eqref{e:302}, let $\wt{z}\coloneqq \e^{-1}z\in[-1,1]$.
    Using \cref{p:CDRP-scaling}\ref{property-CDRP-translation} yields
    \begin{multline*}
        \P^\xi_{1,(0, zn^{2/3}),(\e^{3/2}n, zn^{2/3})}
            \left(
            \sup_{r\in[0,\e^{3/2}]}\left|
                n^{-2/3}
                X(rn)
            \right|
            \ge  m \e
        \right)\\
        \begin{aligned}
            &\law
            \P^\xi_{1,(0,0),(\e^{3/2}n, 0)}
            \left(
                \sup_{r\in[0,\e^{3/2}]}
                \left|
                    n^{-2/3}X(rn) + z
                \right|
                \ge m\e
            \right)\\
            &=
            \P^\xi_{1,(0,0),(\e^{3/2}n, 0)}
            \left(
                \sup_{r\in[0,\e^{3/2}]}
                \left|
                    (\e^{3/2}n)^{-2/3}X(rn) + \wt{z}
                \right|
                \ge m
            \right)\\
            &\le
            \P^\xi_{1,(0,0),(\e^{3/2}n, 0)}
            \left(
                \sup_{r\in[0,\e^{3/2}]}
                \left|
                    (\e^{3/2}n)^{-2/3}X(rn)
                \right|
                \ge m-1
            \right).
        \end{aligned}
    \end{multline*}
    Applying \eqref{e:DZ-bound} with $L=1$ and $\epsilon^{-1} = \e^{3/2}n$ yields \eqref{e:302}.
    We next prove \eqref{e:301}.  We begin with the case $t=1$ and $\beta\ge 1$, and then deduce the general result by rescaling.    We first diffusively rescale to make the temperature $1$, and then translate so the polymer starts at $(0,0)$.
    Namely, using \cref{p:CDRP-scaling}\ref{property-CDRP-translation}--\ref{property-CDRP-scaling}, we have
    \begin{multline*}
        \P^\xi_{\beta, (0,x\beta^{2/3}),(1, y\beta^{2/3})}
            \left(
                \sup_{r\in[0,1]}
                \left|X(r)\right|
                > m \beta^{2/3}
            \right)\\
        \begin{aligned}
        &\law
        \P^\xi_{1, (0,x\beta^{8/3}),(\beta^4, y\beta^{8/3})}
            \left(
                \sup_{r\in[0,1]}
                \left|\beta^{-2}X(\beta^4 r)\right|
                > m \beta^{2/3}
            \right)\\
        &\law
        \P^\xi_{1, (0,0),(\beta^4, (y-x)\beta^{8/3})}
            \left(
                \sup_{r\in[0,1]}
                \left|\beta^{-2}
                    \bigl(X(\beta^4 r)
                    + x\beta^{8/3}\bigr)
                \right|
                > m \beta^{2/3}
            \right)\\
        &\le \P^\xi_{1, (0,0),(\beta^4, (y-x)\beta^{8/3})}
            \left(
                \sup_{r\in[0,1]}
                \left|
                    \beta^{-8/3}
                    X(\beta^4 r)
                \right|
                > m - L
            \right).
        \end{aligned}
    \end{multline*}
    The RHS is small with high probability by \eqref{e:DZ-bound} with $\epsilon = \beta^{-4}$, implying the proposition in the case $t=1$ (say with $T=1$).

    For arbitrary $t>0$ and $\beta \ge t^{-1/4}$, diffusive scaling yields
    \begin{multline}\label{e:beta-t14}
            \P^\xi_{\beta, (0,x\beta^{2/3}),(t, y\beta^{2/3})}
            \left(
                \sup_{r\in[0,t]}
                \left|X(r)\right|
                > m \beta^{2/3}
            \right)\\
        \law
            \P^\xi_{\beta t^{1/4}, (0,x\beta^{2/3}t^{-1/2}),(1, y\beta^{2/3}t^{-1/2})}
            \left(
                \sup_{r\in[0,1]}
                \left|X(r)\right|
                > m \beta^{2/3}t^{-1/2}
            \right),
    \end{multline}
    and we can again apply \eqref{e:DZ-bound} provided $C_1,C_2$ are allowed to depend on $T$.
\end{proof}

Before formulating the next result which is a slight reformulation of \cite[Theorem 1.10]{DZ24}, let us first discuss the validity of the latter.
Indeed, \cite[Theorem 1.10]{DZ24} is formulated as a conditional result relying on the assumption that the KPZ sheet converges to the Airy sheet in a suitable scaling limit (see \cite[Conjecture 1.9]{DZ24} for details).
This has since been proved by \cite{W23},
\footnote{In fact, establishing this convergence was a key step in \cite{W23}'s proof of \cref{t:KPZ-to-landscape}.}
making \cite[Theorem 1.10]{DZ24} unconditional.

We proceed to the aforementioned result of \cite{DZ24}, for which we need the following notation.
For $(s,x;t,y)\in\Rup$ and $\beta>0$, we define the \emph{annealed law} of $X\sim \CDRP_\beta^\xi(s,x;t,y)$ to be the (deterministic) Borel probability measure on $C([s,t])$ given by 
\begin{align*}
    \Pann_{\beta,(s,x),(t,y)} \coloneqq  \E\bigl[\P^\xi_{\beta,(s,x),(t,y)}\bigr].
\end{align*}

\begin{thm}[{Annealed polymer converges to geodesic, \cite[Theorem 1.10]{DZ24}}]\label{t:pathtogeo}
    Fix $x,y\in\R$.
    For $n > 0$, consider the annealed polymer $\overline{X} \sim \Pann_{1,(0,x n^{2/3}), (n,  y n^{2/3})}$, and define the scaled version 
    \begin{align*}
        \overline{X}^{(n)}(r) \coloneqq  n^{-2/3}\, \overline{X}(r n),
        \qquad r\in[0,1].
    \end{align*}
    Restricting $\overline{X}^{(n)}$ to $(0,1)$ yields a random element of $C((0,1))$, where $C((0,1))$ has the topology of uniform convergence on compact sets. 
    As $n\to\infty$, we have the weak convergence
    \begin{align}\label{e:300}
        \overline{X}^{(n)}\Big|_{(0,1)}
        \dto
        2^{1/3} \Pi_{(0,x2^{-1/3}),(1,y2^{-1/3})}
        \Big|_{(0,1)},
    \end{align}
    where $\Pi_{(0,x2^{-1/3}),(1,y2^{-1/3})}$ is the directed landscape geodesic from $(0,x2^{-1/3})$ to $(1,y2^{-1/3})$
    (see \cref{def:geodesics}).
\end{thm}
\begin{proof}
    First, in the case $x=y=0$, \cite[Theorem 1.10]{DZ24} implies that
    \begin{align*}
        \left\{\overline{X}^{(n)}(r)\right\}_{r\in(0,1)}
        \xrightarrow[n\to\infty]{d}
        \left\{\Pi_{(0,0),(2^{1/2},0)}(2^{1/2} r)\right\}_{r\in(0,1)}.
    \end{align*}
    By \eqref{e:geodesic-KPZ-scaling} the RHS is equal in distribution to $2^{1/3}\Pi_{(0,0),(1,0)}\big|_{(0,1)}$, which yields the result for $x=y=0$.

    For general $x,y\in\R$ we use translation and shear invariance of the CDRP to reduce to the case $x=y=0$.
    Let $\overline{X}\sim \Pann_{1,(0,0), (n,0)}$ and let $\overline{X}^{(n)}(r) \coloneqq  n^{-2/3}\, \overline{X}(rn)$ for $r\in[0,1]$.
    Using \cref{p:CDRP-scaling} along with the fact that $\Trans_{-u,-z}$ and $\Shear_{r,-\nu}$ are measure-preserving (\cref{l:Q-measure-preserving}), we obtain that
    \begin{align*}
        \overline{Y}
        \coloneqq  \left\{\overline{X}(r) + xn^{2/3} + r\frac{(y-x)n^{2/3}}{n}\right\}_{r\in[0,n]}
        \sim \Pann_{1,(0,x n^{2/3}),(n,y n^{2/3})}.
    \end{align*}
    Setting $\overline{Y}^{(n)}(r)\coloneqq  n^{-2/3}\,\overline{Y}(rn)$ for $r\in[0,1]$, we get
    \begin{align*}
        \overline{Y}^{(n)}(r) = \overline{X}^{(n)}(r) + x + r(y-x),
        \qquad r\in[0,1].
    \end{align*}
    Since $\overline{X}^{(n)}\big|_{(0,1)} \dto 2^{1/3}\Pi_{(0,0),(1,0)}\big|_{(0,1)}$, 
    it follows that
    \begin{align*}
        \overline{Y}^{(n)}\big|_{(0,1)} \dto \left\{2^{1/3}\Pi_{(0,0),(1,0)}(r) + x + r(y-x)\right\}_{r\in (0,1)},
    \end{align*}
    and the RHS is equal in distribution to the RHS of \eqref{e:300} by \eqref{e:geodesic-KPZ-scaling}.
\end{proof}

\section{Polymer localization and pivotality of thin strips}\label{s:mainproof}

In this short section we formulate the key estimate alluded to earlier in \eqref{efron-stein} in the proof sketch, concerning the pivotality/influence of the white noise inside a thin strip on the polymer free energy (see \cref{t:2main} below).
We then use this estimate to deduce \cref{t:main}.
To orient the reader we first recall the statement of \cref{t:main}.
Recall from \cref{d:narrow-wedge} that for $\beta>0$ and $(s,x;t,y)\in\Rup$, the CDRP free energy profile is defined as
\begin{align}\label{e:501}
    \mathfrak{H}^\xi_\beta(s,x;t,y) \coloneqq  
    \frac{2^{1/3}}{\beta^{4/3}}
    \left(
        \cH_\beta^\xi(s, x2^{1/3}\beta^{2/3}; t, y2^{1/3}\beta^{2/3}) 
        - 2\log\beta 
        + \frac{(t-s)\beta^{4}}{24}
    \right),
\end{align}
where $\cH_\beta^\xi(s,x;t,y) \coloneqq  \log \cZ_\beta^\xi(s,x;t,y)$ (recall \eqref{e:KPZ-cole-hopf} and \cref{def:intro-SHE}), and where $\xi$ is a space-time white noise on $\R^2$. 
By \cref{t:KPZ-to-landscape}, as $\beta\to\infty$ we have $\mf{H}^\xi_{\beta} \dto \cL$ where $\cL$ is the directed landscape.
Our \cref{t:main} states that for $\beta_1,\beta_2\to\infty$ such that $\beta_1 = o(\beta_2)$, we have the joint distributional convergence $(\mf{H}^\xi_{\beta_1}, \mf{H}^\xi_{\beta_2}) \dto (\cL_1,\cL_2)$ where $\cL_1,\cL_2$ are independent directed landscapes.

Next, recall that in \cref{s:resampling-framework-short} we constructed a probability space $(\Xires,\Fres,\Pres)$ supporting a white noise $\xi : \Xires\to\Xi$, such that for every Borel set $B\in\cB(\R^2)$, the space $\Xires$ also supports a white noise $\eta_B$ which can be regarded as obtained from $\xi$ by \emph{resampling $\xi|_B$}.

For $\e,\beta_2>0$ we define the infinite strip
\begin{align}\label{e:strip-unscaled}
    S = S(\e,\beta_2) 
    &\coloneqq  
    \R\times [-\e\beta_2^{2/3}, \e\beta_2^{2/3}].
\end{align}
For $(s,x;t,y)\in\Rup$ and $\beta_1>0$, we define the following partition function and free energy restricted to paths that stay inside $S$:
\begin{align}\label{e:loc}
    \cZ^{\xi}_{\beta_1, S}(s,x;t,y)
    &\coloneqq  \cZ^\xi_{\beta_1}(s,x2^{1/3}\beta_1^{2/3};t,y2^{1/3}\beta_1^{2/3})
    \cdot \P^\xi_{\beta_1,(s,x2^{1/3}\beta_1^{2/3}), (t, y2^{1/3}\beta_1^{2/3})}
    \left(
        \sup_{r\in[s,t]}\left|X(r)\right| \le \e \beta_2^{2/3}
    \right)\nonumber\\[3pt]
    \mf{H}^\xi_{\beta_1, S}(s,x;t,y)
    &\coloneqq  
    \frac{2^{1/3}}{\beta_1^{4/3}}
    \log\left(
        \cZ^{\xi}_{\beta_1, S}(s,x;t,y)
    \right)
    - \frac{2^{4/3}}{\beta_1^{4/3}}\log\beta_1
    + \frac{2^{1/3}}{24}(t-s)\beta_1^{8/3}.
\end{align}

Note that in principle, $\cZ^\xi_{\beta_1,S}(s,x;t,y)$ could equal $0$, making $\mf{H}^\xi_{\beta_1, S}(s,x;t,y)$ ill-defined for any fixed $(s,x;t,y)$ and $\beta_1,\beta_2>0$.
However, as $\beta_1,\beta_2\to\infty$ such that $\beta_1 = o(\beta_2)$, since the polymer at inverse temperature $\beta_1$ typically has transversal fluctuations of size $O(\beta_1^{2/3})$ and hence typically never exits $S(\e,\beta_2)$, 
the restricted free energy $\mf{H}_{\beta_1,S}^\xi$ actually approximates the true free energy $\mf{H}_{\beta_1}^\xi$ uniformly on compact sets.
In fact, for our application, we only need finite-dimensional information as in the following lemma.

\begin{lemm}\label{l:2main}
    Fix $\e>0$.
    For $\beta_2>0$, let $S=S(\e,\beta_2)$ be the strip defined in \eqref{e:strip-unscaled}.
    Then as $\beta_1,\beta_2\to\infty$ such that $\beta_1 = o(\beta_2)$, 
    all finite-dimensional distributions of the random function
    $\mf{H}^{\xi}_{\beta_1, S} -\mf{H}^{\xi}_{\beta_1}$
    converge weakly to $0$ (hence also in probability).
\end{lemm}
\begin{proof}  
Since we only deal with finite-dimensional distributions, 
by a simple union bound it suffices to show that for any $(s,x;t,y)\in\Rup$,
\begin{align*}
    \mf{H}^{\xi}_{\beta_1, S}(s,x;t,y) -\mf{H}^{\xi}_{\beta_1}(s,x;t,y)
    \pto 0.
\end{align*}
By definition of $\mf{H}^\xi_{\beta_1,S}$ and $\mf{H}^\xi_{\beta_1}$ (see \eqref{e:loc} and \eqref{e:501} respectively), and by definition of $S$, we can rewrite the LHS of the above display as
\begin{align*}
    \frac{2^{1/3}}{\beta_1^{4/3}}
    \log 
    \left(
        \P^\xi_{\beta_1,(s,x2^{1/3}\beta_1^{2/3}), (t, y2^{1/3}\beta_1^{2/3})}
        \left(
            \sup_{r\in[s,t]} \left|X(r)\right| \le \e \beta_2^{2/3}
        \right)
    \right).
\end{align*}
Since the prefactor $2^{1/3} / \beta_1^{4/3}$ converges to $0$, we just need to prove that the probability inside the logarithm converges to $1$ in probability, i.e. with high probability the polymer measure concentrates on paths with transversal fluctuations of size $\le \e \beta_2^{2/3}$. 
But this is precisely the content of \cref{p:polymer-TF-DZ}, which in the present setting implies that for all 
$(s,x;t,y)\in\Rup$,
all $\beta_1 \ge (t-s)^{-1/4}$, and all $m\ge 0$,
\begin{align*}
    \P\left(\P^{{\xi}}_{\beta_1,(s,x2^{1/3}\beta_1^{2/3}),(t,y2^{1/3}\beta_1^{2/3})}
    \left(
        \sup_{r\in[s,t]} \left|X(r)\right| \le m \beta_1^{2/3}
    \right) > 1 -  C_{1}e^{-m^{2}/C_{1}}\right)
    > 1 -  C_{2}e^{-m^{3}/C_{2}}
\end{align*}
for some $C_1,C_2>0$ depending on $|x|,|y|, t-s$.
Since $\beta_1\to\infty$, we certainly have $\beta_1 \ge (t-s)^{-1/4}$ eventually.
Setting $m\coloneqq \e(\beta_2/\beta_1)^{2/3}$ in the above display and noting that since $\beta_1=o(\beta_2)$ we have $m\to\infty$, we conclude the lemma.
\end{proof}

We now arrive at one of the main technical results of the paper, which, as already indicated in Section \ref{ss:idea}, is the key input to the proof of \cref{t:main}. 
\begin{thm}[Influence bound]\label{t:2main}
For $\e>0$, let $S=S(\e,\beta_2)$ be the strip defined in \eqref{e:strip-unscaled}.
Let $\xi,\eta_S$ be a pair of coupled white noises as in \eqref{e:def-xi-eta} (so $\eta_S$ is obtained from $\xi$ by resampling $\xi|_{S}$).

There exists $\alpha>0$ such that for all $(s,x;t,y)\in\Rup$ and all sufficiently small $\e>0$,
\begin{align*}
    \limsup_{\beta_2 \to \infty}
    \Eres\left[\left|
        \mf{H}^\xi_{\beta_2}(s,x;t,y)
        - \mf{H}^{\eta_S}_{\beta_2}(s,x;t,y)
    \right|^{2}\right]\le 
    \e^{\alpha}.
\end{align*}
Note that the smallness of $\e$ may depend on $(s,x;t,y)$.
\end{thm}

\begin{rem}
    As indicated in \cref{ss:idea}, in reality the above second moment is of order $\e^{1/2}$.
    An upper bound of this order (up to polylogarithmic corrections) can be proved via a refinement of our arguments, but for simplicity we do not pursue this.
\end{rem}

The proof of \cref{t:2main} is the subject of \cref{s:es}.
We close the present section by using \cref{t:2main} to prove \cref{t:main}.
We first record a simple measure-theoretic lemma:
\begin{lemm}\label{l:abs12}
    Fix $k\in\N$.
    Let $\{X_n\}_{n\in\N}$ and $\{Y_n\}_{n\in\N}$ be sequences of $k$-dimensional random vectors on the same probability space, such that the pair $(X_n,Y_n)$ jointly converges weakly to some $(X,Y)$ as $n\to\infty$.
    Suppose that for every $\e>0$, there exist sequences of random vectors $\{W_n\}_{n\in\N}$ and $\{Z_n\}_{n\in\N}$ such that $Z_n$ and $W_n$ are independent for every $n$, and such that
    \begin{align*}
        \limsup_{n\to\infty} \E\left[
            \norm{X_n - Z_n}_2^2
        \right] 
        \le \e,
        \qquad\text{and}\qquad
        \norm{Y_n - W_n}_2
        \pto 0
        \text{ as $n\to\infty$.}
    \end{align*}
    Then $X$ and $Y$ are independent. 
\end{lemm}
\begin{proof}
    By a diagonal argument, there exist an increasing subsequence $\{n_{i}\}_{i\ge 1}$ 
    and pairs of independent random vectors $Z_i,W_i$ such that 
    $\E\left[\lVert X_{n_{i}}-Z_{i}\rVert_2^2\right] \le 2^{-i}$ and $\P\left(\lVert Y_{n_{i}} - W_i\rVert_2 \ge 2^{-i}\right)\le 2^{-i}$ for all $i\ge 1$.
    It follows (say by Borel--Cantelli) that $(Z_{i}, W_{i})$ converges weakly to $(X,Y)$.
    Since $Z_i$ and $W_i$ are independent for each $i\ge1$, it follows that $X$ and $Y$ are independent, for instance because their joint characteristic function factorizes.
\end{proof}

The proof of \cref{t:main} is now straightforward.
\begin{proof}[Proof of \cref{t:main}] 
    By \cref{t:KPZ-to-landscape}, each of $\mf{H}^{\xi}_{\beta_{1}}$ and $\mf{H}^{\xi}_{\beta_{2}}$ weakly converges to the directed landscape as $\beta_{1},\beta_{2}\to\infty$.
    Therefore, the family of pairs $\{(\mf{H}^{\xi}_{\beta_{1}}, \; \mf{H}^{\xi}_{\beta_{2}})\}_{\beta_1,\beta_2\ge 1}$ is tight due to their marginal tightness.
    So by Prokhorov's theorem, to prove that the pair $(\mf{H}^{\xi}_{\beta_{1}}, \; \mf{H}^{\xi}_{\beta_{2}})$ jointly weakly converges to a pair of independent directed landscapes when $\beta_1,\beta_2\to\infty$ with $\beta_1=o(\beta_2)$, we just need to show that any subsequential joint weak limit must tensorize.
    In fact, since the Borel $\sigma$-algebra on $C(\Rup)$ is generated by finite-dimensional cylinders, we just need to show that all finite-dimensional distributions are asymptotically independent.

    For $\e>0$, let  $S=S(\e,\beta_2)$ be the strip defined in \eqref{e:strip-unscaled}.
    Let $(\xi,\eta_S)$ be a pair of coupled white noises as in \cref{t:2main}.
    Fix $k\in\N$ and $\bfu = (u_1,\dots,u_k)\in(\Rup)^k$.
    Consider, for $\beta\in\{\beta_1,\beta_2\}$, the following finite-dimensional marginals:
    \begin{align*}
        &\mf{H}_{\beta}^\xi(\bfu) \coloneqq  \left(
            \mf{H}_{\beta}^\xi(u_1),\dots,
            \mf{H}_{\beta}^\xi(u_k)
        \right),\\
        &\mf{H}_{\beta}^{\eta_S}(\bfu) \coloneqq  \left(
            \mf{H}_{\beta}^{\eta_S}(u_1),\dots,
            \mf{H}_{\beta}^{\eta_S}(u_k)
        \right).
    \end{align*}
    We similarly define 
    \begin{align*}
        \mf{H}_{\beta_1, S}^\xi(\bfu)
        \coloneqq  \left(
            \mf{H}_{\beta_1, S}^\xi(u_1),\dots, \mf{H}_{\beta_1, S}^\xi(u_k)
        \right),
    \end{align*}
    where $\mf{H}_{\beta_1,S}^\xi$ is defined in \eqref{e:loc}.
    By \cref{l:2main}, as $\beta_1,\beta_2\to\infty$ with $\beta_1=o(\beta_2)$, we have
    \begin{align*}
        \mf{H}_{\beta_1, S}^\xi(\bfu) - \mf{H}_{\beta_1}^\xi(\bfu) \pto 0.
    \end{align*}
    On the other hand, by Theorem \ref{t:2main},  for all sufficiently small $\e>0$, we have
    \begin{align*}
        \limsup_{\beta_2\to\infty}\Eres\left[
            \norm{\mf{H}_{\beta_2}^\xi(\bfu)
            - 
            \mf{H}_{\beta_2}^{\eta_S}(\bfu)
            }_2^2
        \right]
        &\le
        \sum_{i=1}^k
        \limsup_{\beta_2\to\infty}
        \Eres\left[
            \left|\mf{H}_{\beta_2}^\xi(u_i)
            - 
            \mf{H}_{\beta_2}^{\eta_S}(u_i)
            \right|^2
        \right]
        \le k \e^\alpha.
    \end{align*}
    Finally, \cref{p:Z-measurable-real} implies
    that $\xi\mapsto \mf{H}^{\xi}_{\beta_1,S}(\bfu)$ is an $\cF_S$-measurable random vector,
    \footnote{With the notation of \eqref{e:def-ZB}, we have $\cZ^\xi_{\beta_1,S}(s,x;t,y) = \cZ_{B(a,b)}^{C([s,t])}(s,x2^{1/3}\beta_1^{2/3};t,y2^{1/3}\beta_1^{2/3})$ with $B(a,b)=[s,t]\times[-\e\beta_2^{2/3}, \e\beta_2^{2/3}]$.}
    which by \cref{l:Bc-meas-resampling-invariant}\ref{resample-independent} implies that
    \begin{align*}
        \mf{H}^{\xi}_{\beta_1, S}(\bfu)
        \qquad\text{and}\qquad
        \mf{H}^{\eta_S}_{\beta_2}(\bfu)
        \qquad\text{are independent under $\Pres$}.
    \end{align*}
    Lemma \ref{l:abs12} now finishes the proof. 
\end{proof}

\section{Key influence estimate}
\label{s:es}

In this section we prove Theorem \ref{t:2main}. 

\subsection{Changing coordinates}

To aid exposition, we start by reformulating \cref{t:2main} by diffusively rescaling (using \cref{t:AJRS-Z-scaling}) to match the scaling of the discrete last passage percolation model discussed in \cref{ss:idea}.
Namely, we will convert a polymer of length $O(1)$ and inverse temperature $\beta\to\infty$ into a polymer of length $n\to\infty$ and inverse temperature $1$.
Here $\beta$ and $n$ are related by $\beta = n^{1/4}$.

For $n>0$ and $(s,x;t,y)\in\Rup$, write
\begin{equation}\label{e:def-Hn}
    \begin{split}
        \cZ^{\xi,(n)}(s,x;t,y) &\coloneqq  
        \cZ^{\xi}_1(sn, xn^{2/3}; tn, yn^{2/3}),\\
        \cH^{\xi,(n)}(s,x;t,y) &\coloneqq  
        \cH^{\xi}_1(sn, xn^{2/3}; tn, yn^{2/3}),
    \end{split}
\end{equation}
where $\cZ^{\xi}_1$ and $\cH^{\xi}_1$ were defined in \cref{def:intro-SHE} and \eqref{e:KPZ-cole-hopf} (the subscript $1$ is the inverse temperature $\beta=1$).

For $n>0$ and $\e>0$, define the strip
\begin{align}\label{e:new-strip}
    B = B(\e,n) \coloneqq  \R \times [-\e n^{2/3}, \e n^{2/3}].
\end{align}
Let $(\xi,\eta_B)$ be a pair of coupled white noises as in \eqref{e:def-xi-eta} (so $\eta_B$ is obtained from $\xi$ by resampling $\xi|_B$).
We will next show that \cref{t:2main} is equivalent to the following theorem (depicted in Figure \ref{fig:es-overview}):
\begin{thm}[Reformulation of \cref{t:2main}]\label{t:2main-reformulated}
    Let $B=B(\e,n)$ be the strip defined in \eqref{e:new-strip}.
    There exists $\alpha>0$ such that for all $x,y\in\R$ and all sufficiently small $\e>0$, 
    \begin{align*}
        \limsup_{n\to\infty}
        \frac{1}{n^{2/3}}
        \Eres\left[
            \left|
                \cH^{\xi,(n)}(0,x;1,y)
                -
                \cH^{\eta_{B},(n)}(0,x;1,y)
            \right|^2
        \right]
        \le \e^\alpha.
    \end{align*}
\end{thm}

The equivalence of Theorems \ref{t:2main} and \ref{t:2main-reformulated} amounts to straightforward tedious calculations, and the reader may prefer to skip ahead to \cref{s:proof-of-reformulated-ES} where the proof of \cref{t:2main-reformulated} begins.

\addtocontents{toc}{\SkipTocEntry}
\subsection*{Equivalence of \cref{t:2main} and \cref{t:2main-reformulated}}

Recall that \cref{t:2main}
says that for $S=S(\e,\beta) \coloneqq  \R \times [-\e\beta^{2/3}, \e\beta^{2/3}]$, we have 
the following for any $(s,x;t,y)\in\Rup$ and all sufficiently small $\e>0$:
\begin{align*}
    \limsup_{\beta \to \infty}
    \Eres\left[\left|
        \mf{H}^\xi_{\beta}(s,x;t,y)
        - \mf{H}^{\eta_S}_{\beta}(s,x;t,y)
    \right|^{2}\right]\le 
    \e^{\alpha},
\end{align*}
where $\mf{H}^\xi_\beta(s,x;t,y)$ was defined in \cref{d:narrow-wedge}.
We will prove the claimed equivalence using translation invariance and diffusive scaling, as recorded in  \cref{t:AJRS-Z-scaling}.
Briefly, the parameter $n$ in \cref{t:2main-reformulated} will arise as $n = \beta^4$.
Another scaling argument will then reduce to the case $t-s=1$.
Note that we must keep track of the effect of translation and scaling on the \emph{joint distribution} of $(\mf{H}^\xi_{\beta},\; \mf{H}^{\eta_S}_{\beta})$.
We will do this using \cref{l:automorphism-action-resampled}.

From now on we fix $(s,x;t,y)\in\Rup$ and $\beta>0$.
We start by translating to make $s=0$.
Let $\transaff_{s,0}:\R^2\to\R^2$ be the affine translation map from \eqref{e:affine-maps}, and let $\Trans_{s,0}:(\Xi,\cF,\P)\to(\Xi,\cF,\P)$ be the corresponding measure-preserving automorphism from \eqref{e:automorphism-induced}.
By \cref{l:automorphism-action-resampled}, we have the joint distributional equality of white noises:
\begin{align*}
    (\Trans_{s,0}\circ\xi,\; \Trans_{s,0}\circ\eta_S)
    \law
    (\xi,\; \eta_{\transaff_{s,0}^{-1}(S)}).
\end{align*}
To help parse the RHS,
by \eqref{e:affine-maps} the map $\transaff_{s,0}^{-1}$ is just the affine map $(r,z)\mapsto (r-s,z)$.
Since $S$ is translation invariant in the time coordinate, we have $\transaff_{s,0}^{-1}(S)=S$.
Combining the above distributional equality with \cref{t:AJRS-Z-scaling}\ref{property-translation} then gives
\begin{align*}
    \left(
        \mf{H}_{\beta}^{\xi}(s,x;t,y),
        \;
        \mf{H}_{\beta}^{\eta_{S}}(s,x;t,y)
    \right)
    \law
    \left(
        \mf{H}_{\beta}^{\xi}(0,x;t-s,y),
        \;
        \mf{H}_{\beta}^{\eta_{S}}(0,x;t-s,y)
    \right),
\end{align*}
and hence we can assume that $s=0$ in \cref{t:2main}.

We next apply diffusive scaling to reduce to the case $\beta=1$.
Let $\dilaff_{\beta^{-4},\beta^{-2}}$ be the diffusive scaling map from \eqref{e:affine-maps}, and let $\Dil_{\beta^{-4},\beta^{-2}}$ be the corresponding measure-preserving automorphism from \eqref{e:automorphism-induced}.
A straightforward calculation using \eqref{e:affine-maps} shows that $\dilaff_{\beta^{-4},\beta^{-2}}^{-1}(S) = \R \times [-\e\beta^{8/3}, \e\beta^{8/3}]$.
Denote $\wt{\eta}\coloneqq \eta_{\R \times [-\e\beta^{8/3}, \e\beta^{8/3}]}$ (meaning $\wt{\eta}$ is obtained from $\xi$ by resampling $\xi|_{\R \times [-\e\beta^{8/3}, \e\beta^{8/3}]}$).
By the same argument as above, using \cref{l:automorphism-action-resampled} and  \cref{t:AJRS-Z-scaling}\ref{property-scaling} with $\lambda=\beta^2$, 
we have the joint distributional equality
\begin{multline*}
    \left(
        \cZ_{\beta}^{\xi}(0,x2^{1/3}\beta^{2/3};t,y2^{1/3}\beta^{2/3}),\;
        \cZ_{\beta}^{\eta_S}(0,x2^{1/3}\beta^{2/3};t,y2^{1/3}\beta^{2/3})
    \right)\\
    \law
    \left(
        \beta^2\cZ_{1}^{\xi}(0,x2^{1/3}\beta^{8/3};t\beta^{4},y2^{1/3}\beta^{8/3}),
        \;
        \beta^2
        \cZ_{1}^{\wt{\eta}}(0,x2^{1/3}\beta^{8/3};t\beta^{4},y2^{1/3}\beta^{8/3})
    \right).
\end{multline*}
Write $n\coloneqq \beta^4$.
On the RHS of the above display, the resampled region $\R\times [-\e\beta^{8/3}, \e\beta^{8/3}]$ is exactly the strip $B=B(\e,n)$ defined in \eqref{e:new-strip}.
Using this and the notation \eqref{e:def-Hn}, we can rewrite the above distributional equality as
\begin{multline}\label{e:604}
    \left(
        \cZ_{\beta}^{\xi}(0,x2^{1/3}\beta^{2/3};t,y2^{1/3}\beta^{2/3}),\;
        \cZ_{\beta}^{\eta_S}(0,x2^{1/3}\beta^{2/3};t,y2^{1/3}\beta^{2/3})
    \right)\\
    \law
    \left(
        n^{1/2}
        \cZ^{\xi,(n)}(0,x2^{1/3};t,y2^{1/3}),
        \;
        n^{1/2}
        \cZ^{\eta_B,(n)}(0,x2^{1/3};t,y2^{1/3})
    \right).
\end{multline}
Therefore, by the definition of $\mf{H}^\xi_\beta(0,x;t,y)$ (see \eqref{e:501}), we get that
\begin{align*}
    \mf{H}^{\xi}_{\beta}(0,x;t,y)-
    \mf{H}^{\eta_S}_{\beta}(0,x;t,y)
    &\law
    \frac{2^{1/3}}{n^{1/3}}
    \left(
        \cH^{\xi,(n)}(0,x2^{1/3};t,y2^{1/3})
        -
        \cH^{\eta_{B},(n)}(0,x2^{1/3};t,y2^{1/3})
    \right),
\end{align*}
and hence
\begin{multline*}
    \limsup_{\beta\to\infty}
    \Eres\left[
        \left|
            \mf{H}^{\xi}_{\beta}(0,x;t,y)-
            \mf{H}^{\eta_S}_{\beta}(0,x;t,y)
        \right|^2
    \right]
    \\
    = 
    \limsup_{n\to\infty}
        \frac{2^{2/3}}{n^{2/3}}
        \Eres\left[
            \left|
                \cH^{\xi,(n)}(0,x2^{1/3};t,y2^{1/3})
                -
                \cH^{\eta_{B},(n)}(0,x2^{1/3};t,y2^{1/3})
            \right|^2
        \right].
\end{multline*}
This almost implies the equivalence of Theorems \ref{t:2main} and \ref{t:2main-reformulated}---all that remains is to reduce to the case $t=1$.
Assume \cref{t:2main-reformulated} holds for $t=1$.
For arbitrary $t>0$, by the change of variables $m\coloneqq tn$ in \eqref{e:def-Hn}, we have
\begin{multline*}
    \frac{1}{n^{1/3}}\left|
        \cH^{\xi,(n)}(0,x;t,y)
            -
        \cH^{\eta_{B},(n)}(0,x;t,y)
    \right|\\
    =
    \frac{t^{1/3}}{m^{1/3}}
    \left|\cH^{\xi,(m)}(0, xt^{-2/3};1,yt^{-2/3})
            -
    \cH^{\eta_{B},(m)}(0,xt^{-2/3};1,yt^{-2/3})\right|.
\end{multline*}
Similarly, we can write the strip $B=B(\e,n)$ in terms of $m$ as
\begin{align*}
    B = \R \times [-(\e t^{-2/3}) m^{2/3}, (\e t^{-2/3}) m^{2/3}].
\end{align*}
So for all sufficiently small $\e>0$,
\begin{multline*}
    \limsup_{n\to\infty}
    \frac{1}{n^{2/3}}
    \Eres\left[
        \left|
            \cH^{\xi,(n)}(0,x;t,y)
            -
            \cH^{\eta_{B},(n)}(0,x;t,y)
        \right|^2
    \right]\\
    \begin{aligned}
        &=
        t^{2/3}
        \limsup_{m\to\infty}
        \frac{1}{m^{2/3}}
        \Eres\left[
            \left|
                \cH^{\xi,(m)}(0, xt^{-2/3};1,yt^{-2/3})
                        -
                \cH^{\eta_{B},(m)}(0,xt^{-2/3};1,yt^{-2/3})
            \right|^2
        \right]\\
        &\le t^{2/3} (\e t^{-2/3})^\alpha,
    \end{aligned}
\end{multline*}
which is upper bounded by $\e^{\alpha/2}$ for all sufficiently small $\e>0$.
This completes the proof.
\qed
\\

Before continuing, we restate the key input \cref{t:KPZ-to-landscape} from \cite{W23} using our new notation.
\begin{thm}[Reformulation of {\cref{t:KPZ-to-landscape}}]\label{t:KPZ-to-landscape-reformulated}
    We have the following weak convergence
    of $C(\Rup)$-valued random variables:
    \begin{align}
        \frac{2^{1/3}}{n^{1/3}}
        \left(
            \cH^{\xi,(n)}(s,x;t,y)
            +
            \frac{(t-s)n}{24}
        \right) 
        \xrightarrow[n\to\infty]{d} \cL(s,2^{-1/3}x; t, 2^{-1/3}y),
    \end{align}
    where $\cH^{\xi,(n)}$ was defined in \eqref{e:def-Hn}, and where $\cL$ is the directed landscape (see \cref{d:dl}).
\end{thm}

\subsection{Proof of \cref{t:2main-reformulated}}\label{s:proof-of-reformulated-ES}

From now on we fix $x,y\in\R$.
Throughout this section, all constants $C,C',c,c'$, etc. and all ``sufficiently small/large'' statements may depend on $x,y$, which we suppress from the notation.

\begin{figure}[tb]
    \centering
    \includegraphics[width=0.8\textwidth]{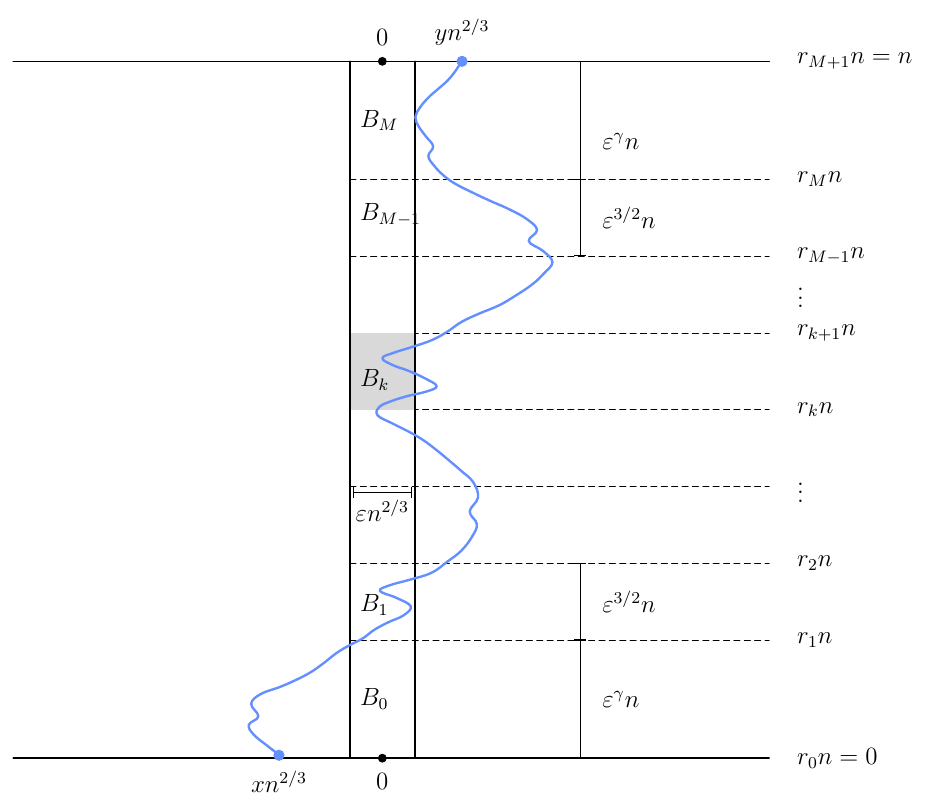}
    \caption{
        The vertical strip is $B = [0,n]\times [-\e n^{2/3}, \e n^{2/3}]$.
        As in \eqref{e:def-boxes-Bk}, $B$ is subdivided into boxes $B_0,B_1,\dots,B_M$, where $B_1,\dots,B_{M-1}$ have height $\le \e^{3/2}n$, and $B_0,B_M$ have height $\asymp \e^\gamma n$ for some $\gamma\in (0,\frac16)$.
        \cref{t:2main-reformulated} asserts a bound for the $L^2$-influence of the white noise in $B$ on  the CDRP free energy $\cH^{\xi,(n)}$.
        By the Efron--Stein inequality \eqref{e:efron-stein-bound}, the influence of $B$ is bounded by the sum of the influences of $B_0,B_1,\dots,B_M$.
        \\[0.8\baselineskip]
        By KPZ energetic fluctuations, resampling the white noise in $B_k=[r_kn, r_{k+1}n]\times [-\e n^{2/3},\e n^{2/3}]$ cannot change the free energy by more than $O((r_{k+1}-r_k)^{1/3}n^{1/3})$.
        In fact, resampling $B_k$ typically changes the free energy negligibly, except on the rare event that before or after resampling, i.e. in the white noise environments $\xi$ or $\eta_k$, the quenched polymer has a reasonable probability of passing through $B_k$.
        }
    \label{fig:es-overview}
\end{figure}

We will prove \cref{t:2main-reformulated} using 
the Efron--Stein inequality (\cref{p:efron-stein-ineq}).
Note that although we are resampling the infinite strip $B=\R\times [-\e n^{2/3}, \e n^{2/3}]$, by directedness $\cH^{\xi,(n)}(0,x;1,y)$
is $\cF_{[0,n]\times\R}$-measurable (\cref{t:AJRAS-Z-existence}\ref{property-adapted}), so it is equivalent to resample $B\cap([0,n]\times\R) = [0,n]\times [-\e n^{2/3},\e n^{2/3}]$.
More precisely, by Lemmas \ref{l:FAFBisFAB} and \ref{l:conditional-variance-formula} we have 
\begin{align*}
    \Eres\biggl[
        \left|
            \cH^{\xi,(n)}
            (0,x;1,y)
            - \cH^{\eta_B,(n)}
            (0,x;1,y)
        \right|^2
    \biggr]
    &=\Eres\biggl[
        \left|
            \cH^{\xi,(n)}
            (0,x;1,y)
            - \cH^{\eta_{B\cap ([0,n]\times\R)},(n)}
            (0,x;1,y)
        \right|^2
    \biggr].
\end{align*}
We therefore redefine $B\coloneqq [0,n]\times [-\e n^{2/3},\e n^{2/3}]$.
We fix $\gamma\in(0, \,\frac16)$ (the precise value will not affect the argument).
Consider the following $\e^{3/2}$-mesh of times (see Figure \ref{fig:es-overview}):
\begin{align}\label{e:resampling-heights}
    \{r_1 < r_2 < \cdots < r_M\} \coloneqq  (\e^{3/2}\Z) \cap [\e^\gamma, 1-\e^\gamma].
\end{align}
We also define 
$r_0 \coloneqq  0$ and $r_{M+1} \coloneqq  1$.
We consider the sequence of boxes 
\begin{align}\label{e:def-boxes-Bk}
    B_k \coloneqq  [r_k n, r_{k+1} n] \times [-\e n^{2/3}, \e n^{2/3}],
    \qquad 
    k\in\lb 0, M\rb,
\end{align}
as depicted in Figure \ref{fig:es-overview}.
This defines the sequence of boxes that we will resample to carry out the argument outlined in Section \ref{ss:idea}. 
While that outline indicated an upper bound of order $\e^{1/2}$ (see \eqref{efron-stein}), here, to simplify the proof at the cost of obtaining a worse bound of $\e^{\alpha}$, we take the two extremal boxes to be of height $\e^{\gamma} n \gg \e^{3/2} n$.
For $k\in\lb0,M\rb$, let $\eta_k\coloneqq \eta_{B_k}$ be the white noise obtained from $\xi$ by resampling $\xi|_{B_k}$ as in \eqref{e:def-xi-eta} (so $(\xi,\eta_k)$ is a pair of coupled white noises under $\Pres$).
By the Efron--Stein inequality (\cref{p:efron-stein-ineq}), we have
\begin{align}\label{e:efron-stein-bound}
    \Eres\left[
            \left|
                \cH^{\xi,(n)}
                (0,x;1,y)
                -
                \cH^{\eta_{B},(n)}
                (0,x;1,y)
            \right|^2
    \right]
    &\le
    \sum_{k=0}^M
    \Eres\left[
            \left|
                \cH^{\xi,(n)}
                (0,x;1,y)
                -
                \cH^{\eta_k,(n)}
                (0,x;1,y)
            \right|^2
        \right].
\end{align}

The rest of this section is devoted to estimating each summand on the RHS of \eqref{e:efron-stein-bound}.
Our analysis will be divided across four subsections:
\begin{itemize}
    \item In \cref{s:favorable-events}, we define a number of high-probability events on the white noise which will be used throughout the analysis. This primarily involves statements saying that the polymer and its free energy satisfy KPZ universality predictions in appropriate quantitative forms.  
    \item In \cref{s:edge-resample}, we estimate the ``edge'' terms in \eqref{e:efron-stein-bound} corresponding to $k=0$ and $k=M$ which will be a straightforward consequence of free energy regularity estimates.
    \item In \cref{s:bulk-resample}, we estimate the ``bulk'' terms in \eqref{e:efron-stein-bound} corresponding to $k\in\lb 1,M-1\rb$.
    In other words, we estimate the ``pivotality'' or ``influence'' of each box $B_k$, i.e. the effect resampling it has on the free energy. 
    This amounts to bounding the contribution to the full partition function made by paths that visit $B_k$.
    \item In \cref{s:combining-edge-bulk}, we combine the results of Sections \ref{s:edge-resample} and \ref{s:bulk-resample} to complete the proof of \cref{t:2main-reformulated}. 
\end{itemize}

\subsection{Favorable events}\label{s:favorable-events}
For most conclusions in this section, we will use the convergence of the KPZ equation and the polymer to the directed landscape and the geodesic therein, respectively, in the zero temperature limit,
together with the Portmanteau lemma.

\subsubsection{Free energy restricted to paths with typical transversal fluctuations}

For the polymer $X\sim \P^\xi_{1,(0,xn^{2/3}), (n, yn^{2/3})}$, 
we denote by $\P^{\xi,(n)}$ the law of $X^{(n)}=\{X^{(n)}(r)\}_{r\in[0,1]}$, where $X^{(n)}(r) \coloneqq  n^{-2/3} X(rn)$.
This scaling ensures that $X^{(n)}$ converges to the directed landscape geodesic (at least in an annealed sense, see \cref{t:pathtogeo}).
The superscript $(n)$ notation is thus consistent with \eqref{e:def-Hn}, in light of \cref{t:KPZ-to-landscape-reformulated}.
For future reference, we note that by \eqref{e:CDRP-fdd}, the multi-point densities of $X^{(n)}$ can be expressed in terms of $\cZ^{\xi,(n)}$ as follows: for all $k\in\N$ and all $0<s_1<\cdots<s_k<1$,
\begin{align}\label{e:scaled-CDRP-fdd}
    \P^{\xi,(n)}
    \left(
        X^{(n)}(s_1)\in dz_1, \dots, X^{(n)}(s_k)\in dz_k
    \right)
    &= n^{2k/3}\, \frac{\prod_{i=0}^k \cZ^{\xi,(n)}(s_i,z_i; s_{i+1},z_{i+1})}{\cZ^{\xi,(n)}(0,x;1,y)}
    dz_1\dots dz_k,
\end{align}
where $(s_0,z_0)\coloneqq (0,x)$ and $(s_{k+1},z_{k+1})\coloneqq (1,y)$,
and where the factor $n^{2k/3}$ is the Jacobian coming from scaling each spatial coordinate $z_i$ by $n^{2/3}$. 

For $k\in\lb 1, M-1\rb$, define the path space event controlling transversal fluctuations at scale $\e^{3/2}$:
\begin{align}\label{e:def-TF-event}
    \XTF_k \coloneqq  \left\{
        X \in C([0,1]) : 
        \sup_{r\in[r_k, r_{k+1}]}\left|X(r) - X(r_{k})\right|
        < \e \log(1/\e)
    \right\}.
\end{align}
Similarly, for $k\in\{0,M\}$, define
\begin{equation}\label{e:def-TF-edge}
    \begin{split}
        \XTF_0 &\coloneqq  \left\{
        X \in C([0,1]) : 
        \left|X(r_1) - x\right|
        < \e^{2\gamma/3}\log(1/\e)
        \right\},\\
        \XTF_M &\coloneqq  \left\{
            X \in C([0,1]) : 
            \left|y - X(r_M)\right|
            < \e^{2\gamma/3} \log(1/\e)
        \right\}.
    \end{split}
\end{equation}
Note that $\XTF_k$ is an open subset of $C((0,1))$ with respect to the topology of uniform convergence on compact sets, for every $k\in\lb 0,M\rb$.
We define the following event on which the rescaled polymer $X^{(n)} \sim \P^{\xi,(n)}$ has typical transversal fluctuations:
\begin{align}\label{e:def-HReg}
    \HReg_{\e,n} \coloneqq  
    \bigcap_{k=0}^M
    \left\{ 
        \xi\in\Xi : 
        \P^{\xi,(n)}(\XTF_k) 
        \ge 1 - e^{-\log^{3/2}(1/\e)}
    \right\}.
\end{align}
For $k\in\lb 0, M\rb$,  define the following restricted partition function and free energy:
\begin{equation}\label{e:def-Zreg-Hreg}
    \begin{split}
    \cZ^{\xi,(n)}_{\reg} 
    = \cZ^{\xi,(n)}_{\reg}(0,x;1,y)
    &\coloneqq  
    \cZ^{\xi,(n)}(0,x;1,y)
    \cdot
    \P^{\xi,(n)}
    (\XTF_k),\\
    \cH^{\xi,(n)}_{\reg}
    = \cH^{\xi,(n)}_{\reg}(0,x;1,y)
    &\coloneqq  \log\left(\cZ^{\xi,(n)}_{\reg}\right),
    \end{split}
\end{equation}
with the convention that $\log 0 \coloneqq  -\infty$.
By strict positivity of $\cZ^{\xi,(n)}(0,x;1,y)$ (\cref{t:AJRAS-Z-existence}\ref{property-strict-positivity}), for all $\xi\in\HReg_{\e,n}$ we have $\cZ^{\xi,(n)}_{\reg}>0$ and hence $\cH^{\xi,(n)}_{\reg}>-\infty$.

The next lemma implies that  $\cH^{\xi,(n)}$ is well-approximated by $\cH^{\xi,(n)}_{\reg}$.

\begin{lemm}[Transversal fluctuation upper bound]\label{l:regcomp}
    There exists $c>0$ such that for all sufficiently small $\e>0$,
    \begin{align*}
        \limsup_{n\to\infty}
        \P(\neg\HReg_{\e,n})
        \le
        e^{-c\log^2(1/\e)}. 
    \end{align*} 
\end{lemm}
\begin{proof} 
By a union bound and Markov's inequality,
\begin{align*}
    \P(\neg \HReg_{\e,n})
    &\le 
    \sum_{k=0}^M
    \P\left(
        \P^{\xi,(n)}(\neg \XTF_k) 
        > e^{-\log^{3/2}(1/\e)}
    \right)\\
    &\le
    e^{\log^{3/2}(1/\e)}
    \sum_{k=0}^M
    \E\left[
        \P^{\xi,(n)}(\neg \XTF_k)
    \right].
\end{align*}
We estimate each summand individually.
The cases $k\in\lb 1,M-1\rb$ and $k\in\{0,M\}$ must be treated separately.
This casework is needed because \cref{t:pathtogeo} is formulated using the topology of $C((0,1))$, not $C([0,1])$, which alone does not imply the rarity of large oscillations near $r=0$ and $r=1$.
However, for $k\in\{0,M\}$ we only need to control the location of the polymer at times $r_1>0$ and $r_M<1$.

Fix $k\in\lb 1,M-1\rb$.
By \cref{t:pathtogeo}, under the annealed measure $\E\bigl[\P^{\xi,(n)}\bigr]$, the path $X^{(n)}\big|_{(0,1)}$ converges weakly to $2^{1/3}\Pi\big|_{(0,1)}$, where $\Pi=\Pi_{(0,x2^{-1/3}),(1,y2^{-1/3})}$ is the directed landscape geodesic from $(0,x2^{-1/3})$ to $(1,y2^{-1/3})$ (defined in \cref{def:geodesics}).
Since $\neg\XTF_k$ is a closed subset of $C((0,1))$, 
it follows by the Portmanteau lemma that
\begin{align*}
    \limsup_{n\to\infty}
    \E\left[
        \P^{\xi,(n)}(\neg \XTF_k)
    \right]
        &\le
        \P\left(
                \sup_{r\in[r_k, r_{k+1}]}
                \left|\Pi(r) - \Pi(r_k)\right|
                \ge 2^{-1/3}\e\log(1/\e)
            \right)\\
        &\le Ce^{-c \log^{2}(1/\e)},
\end{align*}
where the last inequality is by \cref{p:DOV-geodesic-holder} (viz. \eqref{e:geodesic-holder}).

The case $k\in\{0,M\}$ follows from a similar argument as above.
Using that $r_1\ls \e^\gamma$,
we get the following estimate for $k=0$:
\begin{align*}
    \limsup_{n\to\infty}\E\left[
        \P^{\xi,(n)}(\neg \XTF_0) 
    \right]
    &\le
        \P\left(
            \left|\Pi(r_1) - 2^{-1/3}x\right|
            \ge 2^{-1/3}\e^{2\gamma/3}\log(1/\e)
        \right)
    \le C e^{-c\log^2(1/\e)},
\end{align*}
where $C,c>0$ may depend on $\gamma$.
A symmetric argument yields the same bound for $k=M$.
Finally, since $M \ls \e^{-3/2}$, we get 
\begin{align*}
    \limsup_{n\to\infty}
    \P(\neg \HReg_{\e,n})
    \le C' \e^{-3/2}\cdot e^{\log^{3/2}(1/\e)}
    \cdot e^{-c\log^{2}(1/\e)}
    \le e^{-c' \log^2(1/\e)}
\end{align*}
for all sufficiently small $\e$.
\end{proof}

\begin{figure}[htb]
    \centering
    \includegraphics[width=\textwidth]{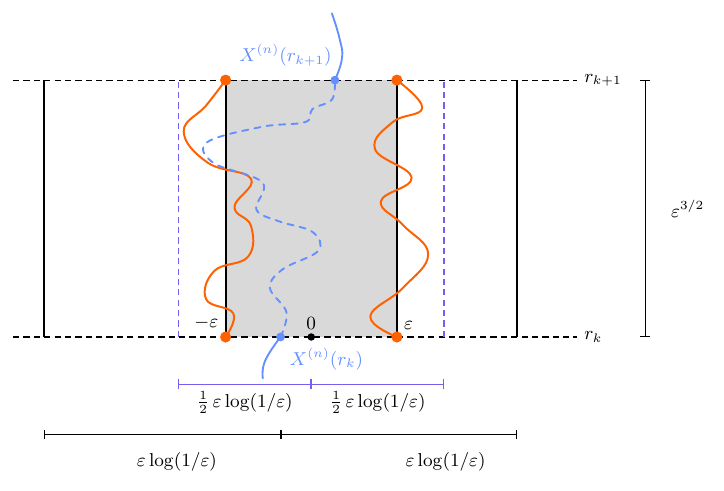}
    \caption{The proof of \cref{l:cond-reg}.
    Shaded gray is the box $[r_k,r_{k+1}] \times [-\e,\e]$.
    The KPZ-scaled polymer $X^{(n)}$ is conditioned to land in $[-\e,\e]$ at times $r_k$ and $r_{k+1}$, and we seek to bound the conditional probability of it wandering further than $\e\log(1/\e)$ from $X^{(n)}(r_k)$ during $[r_k,r_{k+1}]$ (i.e. crossing the vertical solid black lines). 
    By the polymer's Markov property and a union bound, it suffices to bound the probability that an independent  polymer (dashed blue) with endpoints in $[-\e,\e]$ exits $[-\frac12\e\log(1/\e), \frac12\e\log(1/\e)]$ (vertical dashed purple lines) during $[r_k,r_{k+1}]$.
    In fact, by polymer ordering, it suffices to consider the case where the endpoints are either both $-\e$ or $\e$ (drawn in orange): the blue polymer is sandwiched between the orange polymers (in a stochastic dominance sense).
    Finally, the transversal fluctuations of the orange polymers are controlled by \cref{p:polymer-TF-DZ}.
    }
    \label{fig:es-conditional-TF}
\end{figure}

We will also need to control the polymer's transversal fluctuations conditional on its values at the endpoints of a given time interval.
The following lemma accomplishes this.
\begin{lemm}[Conditional transversal fluctuations]\label{l:cond-reg} 
    Define the event 
    \begin{align*}
        \CondReg_{\e,n}
        &\coloneqq  
        \bigcap_{k = 1}^{M-1}
        \left\{
            \xi \in \Xi :
            \P^{\xi,(n)}\left(
                \XTF_k \,\middle|\,
                |X^{(n)}(r_k)|,
                |X^{(n)}(r_{k+1})| \le \e
            \right)
            \ge 1-e^{-\log^{3/2}(1/\e)}
        \right\}.
    \end{align*}
    For all sufficiently small $\e>0$,
    \begin{align*}
        \limsup_{n\to\infty}\P\left(\neg\CondReg_{\e,n}\right) \le e^{-\log^2(1/\e)}.
    \end{align*}
\end{lemm}
\begin{proof}
    Fix $k\in\lb 1,M-1\rb$.
    By the triangle inequality and a union bound,
    \begin{multline*}
        \P^{\xi,(n)}\left(
            \neg\XTF_k\,\middle|\,
            |X^{(n)}(r_k)|,|X^{(n)}(r_{k+1})|\le \e
        \right)\\
        \begin{aligned}
            &\le 
            \P^{\xi,(n)}\left(
                \sup_{r\in[r_k, r_{k+1}]}
                |X^{(n)}(r)| \ge 
                \frac12\e\log(1/\e)
                \,\middle|\,
                |X^{(n)}(r_k)|,|X^{(n)}(r_{k+1})|\le \e
            \right)\\
            &\le
            \P^{\xi,(n)}\left(
                \sup_{r\in[r_k, r_{k+1}]}
                X^{(n)}(r) \ge 
                \frac12\e\log(1/\e)
                \,\middle|\,
                |X^{(n)}(r_k)|,|X^{(n)}(r_{k+1})|\le \e 
            \right)\\
            &\qquad + \P^{\xi,(n)}\left(
                \inf_{r\in[r_k, r_{k+1}]}
                X^{(n)}(r)  \le 
                -\frac12\e \log(1/\e)
                \,\middle|\,
                |X^{(n)}(r_k)|,|X^{(n)}(r_{k+1})|\le \e 
            \right).
        \end{aligned}
    \end{multline*}
    The above is illustrated in Figure \ref{fig:es-conditional-TF} (see the vertical dashed purple lines).
    By polymer ordering (\cref{p:polymer-ordering}), the two terms on the RHS are respectively maximized when $X^{(n)}(r_k)=X^{(n)}(r_{k+1})=\e$ and when $X^{(n)}(r_k)=X^{(n)}(r_{k+1})=-\e$ (see the orange polymers in Figure \ref{fig:es-conditional-TF}).
    More precisely, by the polymer's Markov property (\cref{p:markov-property-CDRP}) and \cref{p:polymer-ordering}, 
    the RHS of the above display is upper bounded by
    \begin{multline*}
        \P^{\xi}_{1,(r_kn, \e n^{2/3}), (r_{k+1}n, \e n^{2/3})}
        \left(
            \sup_{r\in[r_k, r_{k+1}]}
            X^{(n)}(r) \ge 
            \frac12\e\log(1/\e)
        \right)\\
        + 
        \P^{\xi}_{1,(r_kn, -\e n^{2/3}), (r_{k+1}n, -\e n^{2/3})}
        \left(
            \inf_{r\in[r_k, r_{k+1}]}
                X^{(n)}(r)  \le 
            -\frac12\e \log(1/\e)
        \right).
    \end{multline*}
    While we could have used comparison to geodesics in this case as well, we choose to upper bound the above quantities using \cref{p:polymer-TF-DZ} (viz. \eqref{e:302}), which implies the existence of $C_1,C_2>0$ such that for all $|z|\le \e$, all $n\ge \e^{-3/2}$, and all $m>0$,
    \begin{align*}
        \P\left(
            \P^\xi_{1,(r_kn,zn^{2/3}),(r_{k+1}n,zn^{2/3})}
            \left(
                \sup_{r\in[r_k, r_{k+1}]}
                |X^{(n)}(r)|
                \ge m\e
            \right)
            \ge C_1 e^{-m^2/C_1}
        \right)
        \le C_2 e^{-m^3/C_2}.
    \end{align*}
    Thus by the above estimates and a union bound, we have
    \begin{multline*}
        \P\left(
            \P^{\xi,(n)}\left(
            \neg\XTF_k\,\middle|\,
            |X^{(n)}(r_k)|,|X^{(n)}(r_{k+1})|\le \e
            \right)
            \ge 
            e^{-\log^{3/2}(1/\e)}
        \right)
        \\
        \begin{aligned}
        &\le
        2\max_{z\in \{-\e,\e\}}
        \P\left(
            \P^\xi_{1,(r_kn,zn^{2/3}),(r_{k+1}n,zn^{2/3})}
            \left(
                \sup_{r\in[r_k,r_{k+1}]}|X^{(n)}(r)|
                > \frac12 \e  \log(1/\e)
            \right)
            \ge 
            \frac12 e^{-\log^{3/2}(1/\e)}
        \right)\\
        &\le
        2\max_{z\in \{-\e,\e\}}
        \P\left(
            \P^\xi_{1,(r_kn,zn^{2/3}),(r_{k+1}n,zn^{2/3})}
            \left(
                \sup_{r\in[r_k,r_{k+1}]}|X^{(n)}(r)|
                > \frac12 \e  \log(1/\e)
            \right)
            \ge 
            C_1e^{-\frac{1}{C_1}\log^{2}(1/\e)}
        \right)\\
        &\le 2C_2 e^{-\frac{1}{C_2}\log^3(1/\e)}
        \end{aligned}
    \end{multline*} 
    for all sufficiently small $\e>0$ and all $n\ge \e^{-3/2}$.
    Finally, since $M\ls \e^{-3/2}$ , a union bound over $k\in\lb 1,M-1\rb$ yields
    \begin{align*}
        \limsup_{n\to\infty} \P(\neg \CondReg_{\e,n}) \ls \e^{-3/2}e^{-c\log^3(1/\e)},
    \end{align*}
    which implies the lemma.
\end{proof}

\subsubsection{Free energy fluctuations and H\"older continuity}

We will need to control the fluctuations of the free energy accumulated in a small box. Towards this we define the following event:
\begin{multline}\label{e:def-onept-event}
    \HOnePt_{\e,n}
    \coloneqq  
    \bigcap_{k=1}^{M-1}
    \left\{
        \xi : 
        \sup_{|z|,|z'| \le 3\e \log(1/\e)}
        \left|
            \cH^{\xi,(n)}(r_k,z;r_{k+1},z')
            + \frac{\e^{3/2}n}{24}
        \right|
        <
        \e^{1/2} n^{1/3}
        \log^{4}(1/\e)
    \right\}\\ 
    \begin{aligned}
    &\cap
    \left\{
        \xi : 
        \sup_{z:|z-x| \le 2 
        \e^{2\gamma/3}
        \log(1/\e)}
        \left|
            \cH^{\xi,(n)}(0,x; r_1,z)
            + \frac{r_1 n}{24}
        \right|
        <
        \e^{\gamma/3}n^{1/3} 
        \log^{4}(1/\e)
    \right\}\\
    &\cap
    \left\{
        \xi : 
        \sup_{z:|z-y| \le 2
        \e^{2\gamma/3}
        \log(1/\e)}
        \left|
            \cH^{\xi,(n)}(r_M,z;1,y)
            + \frac{(1-r_M)n}{24}
        \right|
        < 
        \e^{\gamma/3}n^{1/3} 
        \log^{4}(1/\e)
    \right\}.
    \end{aligned}
\end{multline}
Note that the above bound is exactly that dictated by the KPZ fluctuation exponent $1/3$, up to the logarithmic factor needed for high-probability statements. 

\begin{lemm}[Uniform pointwise bound]\label{l:onept-high-prob}
    For all sufficiently small $\e>0$,
    \begin{align*}
        \limsup_{n\to\infty}\P(\neg\HOnePt_{\e,n})
        \le e^{-\log^2(1/\e)}.
    \end{align*}
\end{lemm}
\begin{proof}
    We will apply \cref{p:dov-pointwise} and the Portmanteau lemma.
    First, notice that 
    \begin{align*}
        \sup_{|z|,|z'|\le 3\e\log(1/\e)} \frac{(z-z')^2}{\e^{3/2}}
        = 36\e^{1/2}\log^2(1/\e).
    \end{align*}
    Similarly, since $r_1\gs \e^\gamma$ and $1-r_M\gs \e^\gamma$, we have
    \begin{align*}
        &\sup_{z:|z-x| \le 2 \e^{2\gamma/3}\log(1/\e)}
        \frac{(z-x)^2}{r_1}
        \ls r_1^{1/3}\log^2(1/\e) &&\text{and}\\
        &\sup_{z:|z-y| \le 2 \e^{2\gamma/3}\log(1/\e)}
        \frac{(z-y)^2}{1-r_M}
        \ls (1-r_M)^{1/3}\log^2(1/\e).
    \end{align*}
    Therefore, denoting $L\coloneqq 2\max\{|x|,|y|\}+1$ and $r_* \coloneqq  \min\{r_{k+1}-r_k : k\in \lb 0, M\rb\}$, we find that for all $n\ge 1$ and all sufficiently small $\e>0$, we have the inclusions of events
    \begin{multline*}
        \neg\HOnePt_{\e,n}\\
        \begin{aligned}
        &\subset
        \left\{
            \sup_{\substack{(r,z;r',z')\in[-2L,2L]^4 \cap \Rup\\r'-r \ge r_*}}
            \left|
                \frac{\frac{2^{1/3}}{n^{1/3}}\left(\cH^{\xi,(n)}(r,z;r',z') + \frac{(r'-r)n}{24}\right) + \frac{2^{-2/3}(z'-z)^2}{r'-r}}{(r'-r)^{1/3} \log^{4/3}\left(\frac{8}{r'-r}\right)}
            \right|
            \begin{aligned}[t]
                \ge &\log^{4 - 4/3}(1/\e) \\
            &- C\log^{2 - 4/3}(1/\e)
            \end{aligned}
        \right\}\\
        &\subset
        \left\{
            \sup_{\substack{(r,z;r',z')\in[-2L,2L]^4 \cap \Rup\\r'-r \ge r_*}}
            \left|
                \frac{\frac{2^{1/3}}{n^{1/3}}\left(\cH^{\xi,(n)}(r,z;r',z') + \frac{(r'-r)n}{24}\right) + \frac{2^{-2/3}(z'-z)^2}{r'-r}}{(r'-r)^{1/3} \log^{4/3}\left(\frac{8}{r'-r}\right)}
            \right|
            \ge \log^{2}(1/\e)
        \right\}.
        \end{aligned}
    \end{multline*}
    Since the event on the RHS is a closed subset of $C(\Rup)$, it follows by the weak convergence result \cref{t:KPZ-to-landscape-reformulated}, the Portmanteau lemma, and \cref{p:dov-pointwise} that
    \begin{multline*}
        \limsup_{n\to\infty} \P(\neg\HOnePt_{\e,n})
        \\
        \le \P\left(
            \sup_{\substack{(r,z;r',z')\in[-2L,2L]^4 \cap \Rup\\r'-r \ge r_*}}
            \left|
                \frac{\cL(r,2^{-1/3}z;r',2^{-1/3}z') + \frac{2^{-2/3}(z'-z)^2}{r'-r}}{(r'-r)^{1/3} \log^{4/3}\left(\frac{8}{r'-r}\right)}
            \right|
            \ge \log^{2}(1/\e)
        \right)\\
        \le Ce^{-c\log^3(1/\e)}.
    \end{multline*}
    Finally, we decrease the power of $\log(1/\e)$ to absorb $C,c$.
\end{proof}

We next formulate a H\"older continuity estimate for the free energy profile, which will be used to compare extrema of the free energy on small nested intervals.
We define the following event:
\begin{multline}\label{e:def-holder-event}
    \HHol_{\e,n}\\
    \begin{aligned}
    &\coloneqq  
    \bigcap_{k=1}^{M}
    \left\{
        \xi :
        \sup_{\substack{z_1,z_2 \in [-1,1]\\|z_1-z_2| \le 8\e\log(1/\e)}}
        \left|
            \cH^{\xi,(n)}(0,x;r_k,z_1)
            -
            \cH^{\xi,(n)}(0,x;r_k,z_2)
        \right|
        < \e^{1/2}n^{1/3}\log^{4}(1/\e)
    \right\}\\
    &\cap
    \bigcap_{k=1}^{M}
    \left\{
        \xi:
        \sup_{\substack{z_1,z_2 \in [-1,1]\\|z_1-z_2| \le 8\e\log(1/\e)}}
        \left|
            \cH^{\xi,(n)}(r_k,z_1;1,y)
            -
            \cH^{\xi,(n)}(r_k, z_2; 1,y)
        \right|
        < \e^{1/2}n^{1/3}\log^{4}(1/\e)
    \right\}.
    \end{aligned}
\end{multline}

\begin{lemm}\label{l:holder-high-prob}
    For all sufficiently small $\e>0$,
    \begin{align*}
        \limsup_{n\to\infty} \P(\neg \HHol_{\e,n})
        \le e^{-\log^2(1/\e)}.
    \end{align*}
\end{lemm}
\begin{proof}
    We apply \cref{c:dov-holder} and the Portmanteau lemma.
    For notational simplicity we only estimate the increments between points of the form $(0,x;r_k,z)$; the increments between points of the form $(r_k,z;1,y)$ can be handled by a symmetric argument and a union bound.
    
    Denote
    \begin{align*}
        \mathsf{A}_{n,k} \coloneqq  
        \left\{
        \xi :
        \sup_{\substack{z_1,z_2 \in [-1,1]\\|z_1-z_2| \le 8\e\log(1/\e)}}
        \left|
            \cH^{\xi,(n)}(0,x;r_k,z_1)
            -
            \cH^{\xi,(n)}(0,x;r_k,z_2)
        \right|
        \ge \e^{1/2}n^{1/3}\log^{4}(1/\e)
        \right\}.
    \end{align*}
    Since $\mathsf{A}_{n,k}$ is the event that $\cH^{\xi,(n)}$ belongs to a particular closed subset of $C(\Rup)$, it follows by \cref{t:KPZ-to-landscape-reformulated} and the Portmanteau lemma that (writing $\hat{x}\coloneqq 2^{-1/3}x$ and $\hat{z}_i\coloneqq 2^{-1/3}z_i$)
    \begin{multline*}
        \limsup_{n\to\infty}
        \P\left(\bigcup_{k=1}^{M}\mathsf{A}_{n,k}\right)\\
        \begin{aligned}
            &\le 
        \P\left(
            \bigcup_{k=1}^{M}
            \left\{
                \sup_{\substack{z_1,z_2 \in [-1,1]\\|z_1-z_2| \le 8\e\log(1/\e)}}
                \left|
                    \cL(0,\hat{x};r_k,\hat{z}_1)
                    -
                    \cL(0,\hat{x};r_k,\hat{z}_2)
                \right|
                \ge 2^{1/3}\e^{1/2}\log^{4}(1/\e)
            \right\}
        \right)\\
        &\le
        CM e^{-c\log^{17/8}(1/\e)},
        \end{aligned}
    \end{multline*}
    where the last inequality is by \cref{c:dov-holder} and a union bound over $k$.
    Since $M\ls \e^{-3/2}$, we can upper bound the RHS by $e^{-\log^2(1/\e)}$ for all sufficiently small $\e>0$.
\end{proof}

Define the event 
\begin{align}\label{e:def-Good-event}
    \Good_{\e,n} \coloneqq  
    \HReg_{\e,n} 
    \cap \CondReg_{\e,n}
    \cap \HOnePt_{\e,n} 
    \cap \HHol_{\e,n}.
\end{align}

The following lemma lets us focus on analyzing the free energy restricted to paths with typical transversal fluctuations by showing that the contribution to the second moment from the complement event is negligible.
This is by an application of Cauchy--Schwarz and the fourth moment bound \cref{l:fourth-moment}.
\begin{lemm}[Restricting to typical paths and $\Good$]\label{l:restrict-to-Good}
    There exists $c>0$ such that for all sufficiently small $\e>0$ and all $k\in \lb 0,M\rb$,  we have
    \begin{multline*}
        \limsup_{n\to\infty} \frac{1}{n^{2/3}}
        \Eres\left[
            \left|
                \cH^{\xi,(n)}(0,x;1,y) - \cH^{\eta_k,(n)}(0,x;1,y)
            \right|^2
        \right]\\
        \le
        2
        \limsup_{n\to\infty} \frac{1}{n^{2/3}}
        \Eres\left[
            \left|
                \cH^{\xi,(n)}_{\reg} - \cH^{\eta_k,(n)}_{\reg}
            \right|^2
            \1_{\xi,\eta_k \in \Good_{\e,n}}
        \right]
        + e^{-c\log^2(1/\e)},
    \end{multline*}
    where $\cH^{\xi,(n)}_{\reg}$ is defined in \eqref{e:def-Zreg-Hreg}.
\end{lemm}
\begin{proof}
    We suppress the dependence on $(0,x;1,y)$ from the notation.
    Suppose $\xi,\eta_k \in \Good_{\e,n}$.
    Since $\Good_{\e,n}\subset \HReg_{\e,n}$ (defined in \eqref{e:def-HReg}), we have 
    \begin{align*}
        \left|
            \cH^{\xi,(n)}
            - \cH^{\eta_k,(n)}
        \right|
        &\le 
        \left|
            \cH^{\xi,(n)}_{\reg}
            -\cH^{\eta_k,(n)}_{\reg}
        \right|
        - \log\left(1-e^{-\log^{3/2}(1/\e)}\right).
\end{align*}
Then using that $(a-b)^2 \le 2(a^2 + b^2)$ for all $a,b\in\R$, we get
\begin{multline*}
    \limsup_{n\to\infty}
    \frac{1}{n^{2/3}}\Eres\left[
        \left|
            \cH^{\xi,(n)}
            - \cH^{\eta_k,(n)}
        \right|^2
        \1_{\xi,\eta_k\in\Good_{\e,n}}
    \right]\\
    \le
    2
    \limsup_{n\to\infty}
    \frac{1}{n^{2/3}}
    \Eres\left[
        \left|
            \cH^{\xi,(n)}_{\reg}
            -\cH^{\eta_k,(n)}_{\reg}
        \right|^2
        \1_{\xi,\eta_k\in\Good_{\e,n}}
    \right].
    \end{multline*}

    On the other hand, using $(a-b)^2 \le 2(a^2 + b^2)$, the fact that $\xi\law \eta_k$, and Cauchy--Schwarz, we have
    \begin{multline*}
        \frac{1}{n^{2/3}}\Eres\left[
            \left|
                \cH^{\xi,(n)}
                - \cH^{\eta_k,(n)}
            \right|^2
            \1_{\{\xi \not\in \Good_{\e,n}\}\cup\{\eta_k \not\in \Good_{\e,n}\}}
        \right]\\
        \begin{aligned}
        &=
        \frac{1}{n^{2/3}}\Eres\left[
            \left|
                \left(
                    \cH^{\xi,(n)} + \frac{n}{24}
                \right)
                - \left(
                    \cH^{\eta_k,(n)} + \frac{n}{24}
                \right)
            \right|^2
            \1_{\{\xi \not\in \Good_{\e,n}\}\cup\{\eta_k \not\in \Good_{\e,n}\}}
        \right]\\
        &\le 
        \frac{1}{n^{2/3}}\Eres\left[
            4\left|
                \cH^{\xi,(n)}
                + \frac{n}{24}
            \right|^2
            \1_{\{\xi \not\in \Good_{\e,n}\}\cup\{\eta_k \not\in \Good_{\e,n}\}}
        \right]\\
        &\le
        \sqrt{
        16\,
        \Eres\left[
            \left|
                \frac{1}{n^{1/3}}
                \left(
                    \cH^{\xi,(n)}
                    + \frac{n}{24}
                \right)
            \right|^4
        \right]
        \cdot 2\P\left(\neg \Good_{\e,n}\right).
        }
    \end{aligned}
    \end{multline*}
    By \cref{l:fourth-moment}, the above fourth moment is bounded uniformly in $n\ge 1$.
    And by a union bound and Lemmas \ref{l:regcomp}, \ref{l:cond-reg}, \ref{l:onept-high-prob}, and \ref{l:holder-high-prob},
    we have $\limsup_{n\to\infty} \P(\neg \Good_{\e,n}) \le e^{-c\log^2(1/\e)}$ for all sufficiently small $\e>0$.
    This completes the proof.
\end{proof}

We now turn towards bounding the terms in \eqref{e:efron-stein-bound}.

\subsection{Edge influence bound}\label{s:edge-resample}

We first analyze the effect of resampling the initial and final boxes $B_0,B_M$ which are given by $B_0 = [0, r_1n]\times [-\e n^{2/3}, \e n^{2/3}]$ and $B_M = [r_Mn, n]\times [-\e n^{2/3}, \e n^{2/3}]$ (recall \eqref{e:def-boxes-Bk}).
Note that $r_1 \asymp 1-r_M \asymp \e^\gamma$.

We refer to the boxes $B_0,B_M$ as ``edge boxes,'' and the boxes $B_1,B_2,\dots,B_{M-1}$ as ``bulk boxes.''
The next lemma upper bounds the change in free energy caused by resampling an edge box.
We first briefly explain the estimate.
In each of the white noise environments $\xi,\eta_0$, by KPZ fluctuations the free energy accumulated inside the infinite strip $[0,r_1n]\times\R$ fluctuates at scale $(r_1n)^{1/3} \asymp \e^{\gamma/3} n^{1/3}$,
and hence resampling $B_0\subset [0,r_1n]\times\R$ cannot change the free energy by more than $O(\e^{\gamma/3}n^{1/3})$.

\begin{lemm}[Edge resampling bound]\label{l:edge-resampling}
    The following holds for all sufficiently small $\e>0$ and all $n\ge 1$.
    Fix $k\in\{0,M\}$.
    For coupled white noises $(\xi,\eta_k)$ as defined below \eqref{e:def-boxes-Bk},
    it holds $\Pres$-a.s. that
    \begin{align*}
        \left|
            \cH^{\xi,(n)}_{\reg}
            -\cH^{\eta_k,(n)}_{\reg}
        \right|\1_{\xi,\eta_k\in \Good_{\e,n}}
        \le 
        2 \e^{\gamma/3} n^{1/3} \log^{4}(1/\e),
    \end{align*}
    where $\Good_{\e,n}$ was defined in \eqref{e:def-Good-event}.
    In fact, of the four events in \eqref{e:def-Good-event}, only $\HOnePt_{\e,n}$ is relevant for the above bound.
\end{lemm}
\begin{proof}
    We only treat the case $k=0$, as the case $k=M$ follows from a symmetric argument.
    We define the closed interval
    \begin{align*}
        I \coloneqq  \left[
            x - \e^{2\gamma/3}\log(1/\e),\;\;
            x + \e^{2\gamma/3}\log(1/\e)
        \right].
    \end{align*}
    Note that by \eqref{e:scaled-CDRP-fdd} and the definition of $\XTF_0$ (see \eqref{e:def-TF-edge}), we have
    \begin{align*}
        \cZ^{\xi,(n)}_{\,\XTF_0}
        &= n^{2/3}\int_I \cZ^{\xi,(n)}(0,x;r_1,z)
        \cZ^{\xi,(n)}(r_1,z; 1,y)\,dz.
    \end{align*}

    Fix $\xi\in\Good_{\e,n}$.
    Since $\xi\in\HOnePt_{\e,n}$ (see \eqref{e:def-onept-event}), we have 
    \begin{align*}
        \sup_{z\in I}
        \left|
            \cH^{\xi,(n)}(0,x; r_1,z) + \frac{r_1 n}{24}
        \right|
        \le \e^{\gamma/3} n^{1/3} \log^{4}(1/\e).
    \end{align*}
    Equivalently, for all $z\in I$,
    \begin{align}\label{e:650}
        \exp\left(
            -\frac{r_1n}{24} - \e^{\gamma/3} n^{1/3} \log^{4}(1/\e)
        \right)
        \le
        \cZ^{\xi,(n)}(0,x;r_1,z)
        \le \exp\left(
            -\frac{r_1 n}{24} + \e^{\gamma/3} n^{1/3} \log^{4}(1/\e)
        \right).
    \end{align}
    We therefore have the upper and lower bounds
    \begin{align*}
        \cZ^{\xi,(n)}_{\,\XTF_0}
        &\le
        \exp\left(
            -\frac{r_1n}{24} + \e^{\gamma/3}n^{1/3}\log^{4}(1/\e)
        \right)
        n^{2/3}
        \int_{I} 
        \cZ^{\xi,(n)}(r_1,z; 1,y)
        \,dz,\\
        \cZ^{\xi,(n)}_{\,\XTF_0}
        &\ge
        \exp\left(
            -\frac{r_1n}{24} - \e^{\gamma/3}n^{1/3}\log^{4}(1/\e)
        \right)
        n^{2/3}\int_{I} 
        \cZ^{\xi,(n)}(r_1,z; 1,y)
        \,dz.
    \end{align*}
    As a function of $\xi$, the integral on the RHS of the above two inequalities is $\cF_{B_0^c}$-measurable, by \cref{t:AJRAS-Z-existence}\ref{property-adapted} and the fact that integration over $I$ defines a continuous functional on $C(\R)$ (cf. the proof of \cref{l:polymer-measure-is-measurable}).
    The value of the integral is therefore unaffected by resampling $\xi|_{B_0}$; more precisely, by \cref{l:Bc-meas-resampling-invariant}\ref{resample-invariant}, we have
    \begin{align*}
        \int_{I} 
        \cZ^{\xi,(n)}(r_1,z; 1,y)
        \,dz
        =
        \int_{I} 
        \cZ^{\eta_0,(n)}(r_1,z; 1,y)
        \,dz
        \qquad\text{$\Pres$-a.s.}
    \end{align*}
    The lemma follows from the previous two displays.
\end{proof}

Lemma \ref{l:edge-resampling} immediately implies:
\begin{lemm}[Edge influence bound]\label{l:edge-L2-resampling}
    For all sufficiently small $\e>0$ and all $k\in\{0,M\}$,
    \begin{align*}
        \limsup_{n\to\infty} \frac{1}{n^{2/3}}\Eres\left[
            \left|
                \cH^{\xi,(n)}_{\reg} - \cH^{\eta_k, (n)}_{\reg}
            \right|^2
            \1_{\xi,\eta_k \in \Good_{\e,n}}
        \right]
        &\le 4\e^{2\gamma/3}\log^8(1/\e).
    \end{align*}
\end{lemm}

\subsection{Bulk influence bound}\label{s:bulk-resample}

We proceed now to estimating the effect of resampling bulk boxes $B_k$ for $k\in\lb 1,M-1\rb$.
Recall that $B_k = [r_k n,r_{k+1} n]\times[-\e n^{2/3}, \e n^{2/3}]$, with $r_{k+1}-r_k=\e^{3/2}$ and $r_k\ge \e^\gamma $ and $r_{k+1}\le 1-\e^\gamma$, where $\gamma\in(0,\frac16)$ is a fixed parameter (see \eqref{e:resampling-heights}, \eqref{e:def-boxes-Bk}).

Towards applying the Efron--Stein inequality \eqref{e:efron-stein-bound}, we wish to isolate the contribution to the free energy made by paths that enter $B_k$, for any $k\in \lb 1, M-1\rb$.
The results of Section \ref{s:restricting-polymer} show that the contribution to the partition function made by paths that do not intersect $B_k$ is a measurable function of the white noise outside of $B_{k}$, and is therefore unchanged when we resample the noise over $B_{k}$. 
Accordingly, we introduce the following notation:
\begin{align*}
    \cZ^{\xi, (n)}_{\reg}(B_k)
    &\coloneqq 
    \cZ^{\xi,(n)}(0,x;1,y)
    \cdot
    \P^{\xi,(n)}\left(
        \XTF_k
        \cap
        \left\{
            \inf_{r\in[r_k,\,r_{k+1}]}|X^{(n)}(r)|
            \le \e
        \right\}
    \right),\\[3pt]
    \cZ^{\xi, (n)}_{\reg}(B_k^c)
    &\coloneqq 
    \cZ^{\xi,(n)}(0,x;1,y)
    \cdot
    \P^{\xi,(n)}\left(
        \XTF_k
        \cap
        \left\{
            \inf_{r\in[r_k,\,r_{k+1}]}|X^{(n)}(r)|
            > \e
        \right\}
    \right).
\end{align*}
So $\cZ^{\xi,(n)}_{\reg} = \cZ^{\xi, (n)}_{\reg}(B_k) + \cZ^{\xi, (n)}_{\reg}(B_k^c)$.

We first prove an a priori bound for how resampling $\xi|_{B_k}$ affects $\cZ^{\xi,(n)}_{\reg}(B_k)$.
The statement and proof are similar to those of \cref{l:edge-resampling}, as is the underlying reasoning (see the discussion above \cref{l:edge-resampling}, as well as Figure \ref{fig:es-overview}).

\begin{lemm}\label{l:a-priori-change-B}
    There exists $C>0$ such that the following holds for all sufficiently small $\e>0$ and all $n\ge 1$.

    Let $\Good_{\e,n}$ be as in \eqref{e:def-Good-event}.
    Fix $k\in\lb 1, M-1\rb$.
    For all $\xi\in \Good_{\e,n}$, we have $\cZ^{\xi,(n)}_{\reg}(B_k)>0$.
    Moreover, for $(\xi,\eta_k)$ coupled white noises as defined below \eqref{e:def-boxes-Bk}, 
    it holds $\Pres$-a.s. that
    \begin{align*}
        \left|
            \log\left(
            \frac{\cZ^{\xi,(n)}_{\reg}(B_k)}{\cZ^{\eta_k,(n)}_{\reg}(B_k)}
        \right)
        \right|
        \1_{\xi,\eta_k\in\Good_{\e,n}}
        \le 
        C\e^{1/2} n^{1/3} \log^{4}(1/\e)
        + C\log\log(1/\e).
    \end{align*}
\end{lemm}
\begin{proof}
    We use an argument similar to that of \cref{l:edge-resampling}. Fix $\xi\in\Good_{\e,n}$.
    Arguing as in  \eqref{e:650}, since $\xi\in\Good_{\e,n}\subset\HOnePt_{\e,n}$ (see \eqref{e:def-onept-event}), we have that for all $|z|,|z'| \le 3\e \log(1/\e)$,
    \begin{align}\label{e:640}
        \exp\left(
            -\frac{\e^{3/2}n}{24}
            - \e^{1/2}n^{1/3}\log^{4}(1/\e)
        \right)
        \le
        \cZ^{\xi,(n)}(r_k,z;r_{k+1},z')
        \le \exp\left(
            -\frac{\e^{3/2}n}{24}
            + \e^{1/2}n^{1/3}\log^{4}(1/\e)
        \right).
    \end{align}
    And since $\xi\in\Good_{\e,n}\subset\HHol_{\e,n}$ (see \eqref{e:def-holder-event}), we have the following two bounds:
    \begin{equation}\label{e:610}
        \begin{split}
        \sup_{|z|\le 3\e \log(1/\e)}\left|\log\left(
            \frac{\cZ^{\xi,(n)}(0,x;r_k,z)}{\cZ^{\xi,(n)}(0,x;r_k,0)}
        \right)\right|
        &\le \e^{1/2}n^{1/3}\log^{4}(1/\e),
        \\
        \sup_{|z|\le 3\e\log(1/\e)}\left|\log\left(
            \frac{\cZ^{\xi,(n)}(r_{k+1},z; 1,y)}{\cZ^{\xi,(n)}(r_{k+1},0; 1,y)}
        \right)\right|
        &\le \e^{1/2}n^{1/3}\log^{4}(1/\e).
        \end{split}
    \end{equation}
    With the above estimates recorded we proceed to the argument.

    We first lower bound $\cZ^{\xi,(n)}_{\reg}(B_k)$.
    We have the probability lower bound
    \begin{align*}
        \P^{\xi,(n)}\left(
            \XTF_k
            \cap
            \left\{
                \inf_{r\in[r_k,\,r_{k+1}]}|X^{(n)}(r)|
                \le \e
            \right\}
        \right)
        &\ge
        \P^{\xi,(n)}\left(
            \XTF_k \,\middle|\, 
            |X^{(n)}(r_k)|, |X^{(n)}(r_{k+1})| \le \e
        \right)\\
        &\qquad \cdot
        \P^{\xi,(n)}\left(
            |X^{(n)}(r_k)|, |X^{(n)}(r_{k+1})| \le \e
        \right)\\
        &\ge
        \frac12
        \P^{\xi,(n)}\left(
            |X^{(n)}(r_k)|, |X^{(n)}(r_{k+1})| \le \e
        \right),
    \end{align*}
    where the last inequality is because $\xi\in \Good_{\e,n}\subset \CondReg_{\e,n}$ (see \cref{l:cond-reg}) and $1-e^{-\log^{3/2}(1/\e)} \ge \frac12$ for sufficiently small $\e$.
    Note that since $k\in\lb 1,M-1\rb$ we have $r_k>0$ and $r_{k+1}<1$, 
    which allows us to use the multi-point density formula  \eqref{e:scaled-CDRP-fdd} to rewrite the probability on the RHS as
    \begin{multline*}
        \P^{\xi,(n)}\left(
            |X^{(n)}(r_k)|, |X^{(n)}(r_{k+1})| \le \e
        \right)
        \\
        \begin{aligned}
        &=
        \iint_{|z|,|z'| \le \e}
        n^{4/3}
        \frac{\cZ^{\xi,(n)}(0,x;r_k,z)\cZ^{\xi,(n)}(r_k,z;r_{k+1},z')
        \cZ^{\xi,(n)}(r_{k+1},z';1,y)}{\cZ^{\xi,(n)}(0,x;1,y)}
        \,dz\,dz'\\
        \overset{\eqref{e:640},\eqref{e:610}}&{\ge}
        4n^{4/3}\e^2
        \exp\left(
            -\frac{\e^{3/2}n}{24} - 3\e^{1/2}n^{1/3}\log^{4}(1/\e)
        \right)
        \frac{\cZ^{\xi,(n)}(0,x;r_k,0)
        \cZ^{\xi,(n)}(r_{k+1},0;1,y)}{\cZ^{\xi,(n)}(0,x;1,y)}.
        \end{aligned}
    \end{multline*} 
    Multiplying by $\cZ^{\xi,(n)}(0,x;1,y)$ yields
    \begin{equation}\label{e:641}
        \begin{split}
            \cZ^{\xi,(n)}_{\reg}(B_k)
        &\ge
        2 n^{4/3}\e^2
        \exp\left(
            -\frac{\e^{3/2}n}{24} - 3\e^{1/2}n^{1/3}\log^{4}(1/\e)
        \right)\\
        &\qquad \cdot
        \cZ^{\xi,(n)}(0,x;r_k,0)
        \cZ^{\xi,(n)}(r_{k+1},0;1,y).
        \end{split}
    \end{equation}
    By \cref{t:AJRAS-Z-existence}\ref{property-strict-positivity}
    this lower bound is strictly positive, which proves the first claim.

    We upper bound $\cZ^{\xi,(n)}_{\reg}(B_k)$ using a similar argument.
    By definition of $\XTF_k$ (see \eqref{e:def-TF-event}) we have
    \begin{multline}\label{e:upper-bound-ZB}
        \P^{\xi,(n)}\left(
            \XTF_k
            \cap
            \left\{
                \inf_{r\in[r_k,\,r_{k+1}]}|X^{(n)}(r)|
                \le \e
            \right\}
        \right)\\
        \begin{aligned}
            \overset{\eqref{e:def-TF-event}}&{\le} \P^{\xi,(n)}\left(
                |X^{(n)}(r_k)|, |X^{(n)}(r_{k+1})|
                \le \e + 2\e \log(1/\e)
            \right)\\
            \overset{\eqref{e:scaled-CDRP-fdd}}&{\le}
            \iint_{|z|,|z'| \le 3\e \log(1/\e)}
            n^{4/3} 
            \frac{\cZ^{\xi,(n)}(0,x; r_k ,z)\cZ^{\xi,(n)}(r_k,z;r_{k+1},z')
            \cZ^{\xi,(n)}(r_{k+1},z';1,y)}{\cZ^{\xi,(n)}(0,x;1,y)}
            \,dz\,dz'\\
            \overset{\eqref{e:640},\eqref{e:610}}&{\le}
            36 n^{4/3} \e^2 \log^2(1/\e)
            \exp\left(
                -\frac{\e^{3/2}n}{24}
                + 3\e^{1/2}n^{1/3}\log^{4}(1/\e)
            \right)\\
            &\qquad\qquad \cdot
            \frac{\cZ^{\xi,(n)}(0,x; r_k ,0)
            \cZ^{\xi,(n)}(r_{k+1},0;1,y)}{\cZ^{\xi,(n)}(0,x;1,y)},
        \end{aligned}
    \end{multline}
    hence
    \begin{equation}\label{e:642}
        \begin{split}
            \cZ^{\xi,(n)}_{\reg}(B_k)
        &\le 
        36 n^{4/3} \e^2 \log^2(1/\e)
        \exp\left(
            -\frac{\e^{3/2}n}{24}
            +
            3\e^{1/2} n^{1/3} \log^{4}(1/\e)
        \right)\\
        &\qquad \cdot \cZ^{\xi,(n)}(0,x; r_k ,0)
            \cZ^{\xi,(n)}(r_{k+1},0;1,y).
        \end{split}
    \end{equation}

    Since $\xi\mapsto \cZ^{\xi,(n)}(0,x; r_k ,0)\cZ^{\xi,(n)}(r_{k+1},0;1,y)$ is a measurable function of the white noise inside $B_k^c$ (\cref{t:AJRAS-Z-existence}\ref{property-adapted}), it follows by \cref{l:Bc-meas-resampling-invariant}\ref{resample-invariant} that $\cZ^{\xi,(n)}(0,x; r_k ,0)\cZ^{\xi,(n)}(r_{k+1},0;1,y) = \cZ^{\eta_k,(n)}(0,x; r_k ,0)\cZ^{\eta_k,(n)}(r_{k+1},0;1,y)$,
    and hence the lower and upper bounds \eqref{e:641} and \eqref{e:642} hold a.s. with $\eta_k$ in place of $\xi$ on the LHS.
    We deduce that for all sufficiently small $\e>0$, all $n\ge 1$, and $\Pres$-a.e. $\xi,\eta_k\in \Good_{\e,n}$,
    \begin{align*}
        \left|
            \log\left(
                \frac{\cZ^{\xi,(n)}_{\reg}(B_k)}{\cZ^{\eta_k,(n)}_{\reg}(B_k)}
            \right)
        \right|
        &\le 6\e^{1/2}n^{1/3}\log^{4}(1/\e)
        + \log\left(\frac{36\log^2(1/\e)}
        {2}
        \right)\\
        &\le C \e^{1/2} n^{1/3}\log^{4}(1/\e)
        + C\log\log(1/\e)
    \end{align*}
    for some universal constant $C>0$.
\end{proof}

We can easily convert \cref{l:a-priori-change-B} into a bound on $\left|\cH^{\xi,(n)}_{\reg} - \cH^{\eta_k,(n)}_{\reg}
\right|$.

\begin{lemm}[Worst-case bulk resampling bound]\label{l:a-priori-change-full}
    There exists $C>0$ such that for all sufficiently small $\e>0$, all $k\in\lb 1, M-1\rb$, and  all $n\ge 1$,
    it holds $\Pres$-a.s. that
    \begin{align*}
        \left|
                \cH^{\xi,(n)}_{\reg} - \cH^{\eta_k,(n)}_{\reg}
        \right|
        \1_{\xi,\eta_k \in \Good_{\e,n}}
        \le 
        C \e^{1/2} n^{1/3} \log^{4}(1/\e)
        + C\log\log(1/\e),
    \end{align*}
    where $\Good_{\e,n}$ was defined in \eqref{e:def-Good-event}.
\end{lemm}
\begin{proof}
    By \cref{l:meas3}, the map $\xi\mapsto \cZ^{\xi,(n)}_{\reg}(B_k^c)$ is $\cF_{B^c_k}$-measurable.
    Thus by \cref{l:Bc-meas-resampling-invariant}, we have $\cZ^{\xi,(n)}_{\reg}(B_k^c) = \cZ^{\eta_k,(n)}_{\reg}(B_k^c)$ almost surely, and in particular 
    \begin{align}\label{e:611}
        \left|
            \cH^{\xi,(n)}_{\reg} - \cH^{\eta_k,(n)}_{\reg}
        \right|
        \1_{\xi,\eta_k\in\Good_{\e,n}}
        &=
        \left|
            \log\left(
                \frac{\cZ^{\xi,(n)}_{\reg}(B_k)
                +
                \cZ^{\xi,(n)}_{\reg}(B_k^c)
                }
            {\cZ^{\eta_k,(n)}_{\reg}(B_k) +
                \cZ^{\xi,(n)}_{\reg}(B_k^c)}
            \right)
        \right|
        \1_{\xi,\eta_k\in\Good_{\e,n}}.
    \end{align}
    By \cref{l:a-priori-change-B}, for $\xi,\eta_k\in\Good_{\e,n}$ we have $\cZ^{\xi,(n)}_{\reg}(B_k)>0$ and $\cZ^{\eta_k,(n)}_{\reg}(B_k)>0$, and by \cref{t:AJRAS-Z-existence}\ref{property-strict-positivity} we have $\cZ^{\xi,(n)}_{\reg}(B_k^c)\ge 0$.
    Therefore, using that $\left|\log\left(\frac{a+b}{a'+b}\right)\right| \le \left|\log \frac{a}{a'}\right|$ for all $a,a'>0$ and $b\ge 0$,
    we obtain
    \begin{align*}
        \left|
            \cH^{\xi,(n)}_{\reg} - \cH^{\eta_k,(n)}_{\reg}
        \right|\1_{\xi,\eta_k\in\Good_{\e,n}}
        \le
        \left|
            \log\left(
                \frac{\cZ^{\xi,(n)}_{\reg}(B_k)}
                {\cZ^{\eta_k,(n)}_{\reg}(B_k)}
            \right)
        \right|
        \1_{\xi,\eta_k\in\Good_{\e,n}}.
    \end{align*}
    The lemma now follows from \cref{l:a-priori-change-B}.
\end{proof}

Note that the bound in \cref{l:a-priori-change-full} is typically quite lossy, as by KPZ fluctuations, the polymer's location at time $r_kn$ is roughly uniformly distributed on an interval of width $\approx r_k^{2/3}n^{2/3}$, which is much larger than the width of $B_k$ (the latter being $2\e n^{2/3}$),
and on the event that the polymer typically does not enter $B_k$, resampling $B_k$ essentially has no effect on the free energy.
We will make this precise relying on the local Brownianity of the free energy $\cH$.
To start, we introduce the following event, which says that at time $r_k$, the maximum of $\cH^{\xi,(n)}(0,x;r_k,z) + \cH^{\xi,(n)}(r_k,z;1,y)$ on a macroscopic interval is much larger than the maximum of the same function on a small subinterval.
\begin{multline}\label{e:def-argmax-event}
    \Max_{\e,n}^{k} \coloneqq  
    \Biggl\{
        \xi : \sup_{|z|\le 1}\left[
            \cH^{\xi,(n)}(0,x;r_k,z)
            + \cH^{\xi,(n)}(r_k,z; 1,y)
        \right]\\
        >
        \sup_{|z|\le 2\e \log(1/\e)}\left[
            \cH^{\xi,(n)}(0,x;r_k,z)
            + \cH^{\xi,(n)}(r_k,z; 1,y)
        \right]
        + \e^{1/2}n^{1/3} \log^{5}(1/\e)
    \Biggr\}.
\end{multline}
\begin{lemm}[Polymer delocalization]\label{l:Max-probability-bound}
    For all sufficiently small $\e$ and all $k\in\lb 1,M-1\rb$, we have
    \begin{align*}
        \limsup_{n\to\infty}
        \P\left(\neg \Max_{\e,n}^{k}\right)
        \ls \e^{1 - 3\gamma}.
    \end{align*}
\end{lemm}
\begin{proof}
    Fix $k\in\lb 1,M-1\rb$ (so $r_k\in(\e^{3/2}\Z)\cap[\e^\gamma, 1-\e^\gamma]$).
    Define the following subset of $C(\Rup)$:
    \begin{multline*}
        \sS_{\e,n,k}
        \coloneqq 
        \bigg\{
        F \in C(\Rup): 
            \sup_{|z|\le 1}
            \bigl[F(0,x;r_k,z) + F(r_k,z;1,y)\bigr]\\
            \le
            \sup_{|z|\le 2\e\log(1/\e)}
            \bigl[
                F(0,x;r_k,z) + F(r_k,z;1,y)
            \bigr]
            +
            \e^{1/2} \log^5(1/\e)
        \bigg\}.
    \end{multline*}
    Then  $\P(\neg \Max_{\e,n}^{k}) = \P(\frac{1}{n^{1/3}}\cH^{\xi,(n)} \in \sS_{\e,n,k})$.
    Now, observe that $\sS_{\e,n,k}$ is invariant under affine shifts, in the sense that
    for any $F\in C(\Rup)$ and $\mu\in\R$, we have
    \begin{align*}
        F\in \sS_{\e,n,k}
        \qquad\quad
        \textup{iff}
        \qquad\quad
        F(s,z;t,z') + \mu(t-s)
        \in \sS_{\e,n,k},
    \end{align*}
    where $(s,z;t,z')$ denotes a generic element of $\Rup$.
    In particular,
    \begin{align*}
        \P\left(\neg \Max_{\e,n}^{k}\right)
        =\P\left(
            \frac{1}{n^{1/3}}\left(\cH^{\xi,(n)}(s,z;t,z') + \frac{(t-s)n}{24}\right)
            \in\sS_{\e,n,k}
        \right).
    \end{align*}
    Since $\sS_{\e,n,k}$ is a closed subset of $C(\Rup)$, it follows by  \cref{t:KPZ-to-landscape-reformulated} and the Portmanteau lemma that 
    \begin{align*}
        \limsup_{n\to\infty}
        \P\left(
            \neg \Max_{\e,n}^{k}
        \right)
        \le \P\left(
            2^{-1/3}\cL(s,2^{-1/3}z; t, 2^{-1/3}z') \in \sS_{\e,n,k}\right),
    \end{align*}
    and by \cref{p:geodesic-deloc} the RHS is $O(\e^{1-3\gamma})$.
\end{proof}

The next lemma implies that with high probability in $\xi$, the polymer typically does not pass through $B_k$. 
\begin{lemm}[The polymer typically avoids $B_k$]\label{l:brownian-comparison-result}
    There exists $c>0$ such that for all sufficiently small $\e>0$,  
    all $k\in\lb 1, M-1\rb$,
    all $n\ge 1$,
    and all $\xi\in \Good_{\e,n}\cap \Max_{\e,n}^{k}$, we have
    \begin{align*}
        \cZ^{\xi,(n)}_{\reg}(B_k^c)
        \ge \left(
            \exp\left(c\e^{1/2} n^{1/3} \log^5(1/\e)\right) - 1
        \right)
        \cZ^{\xi,(n)}_{\reg}(B_k),
    \end{align*}
    where $\Good_{\e,n}$ was defined in \eqref{e:def-Good-event}.
\end{lemm}
\begin{proof}
    By rearranging, it is equivalent to show that for $\xi\in\Good_{\e,n}\cap\Max_{\e,n}^{k}$, we have
    \begin{align*}
        \cZ^{\xi,(n)}_{\reg} \ge \exp\left(c\e^{1/2} n^{1/3} \log^5(1/\e)\right) \cZ^{\xi,(n)}_{\reg}(B_k).
    \end{align*}
    We will lower bound $\cZ^{\xi,(n)}_{\reg}$ and upper bound $\cZ^{\xi,(n)}_{\reg}(B_k)$.

    Let $I\subset [-1,1]$ be any subinterval of length $|I| = 8\e \log(1/\e)$ such that
    \begin{align*}
        \max_{z\in I}
        \left(
            \cH^{\xi,(n)}(0,x;r_k,z)
            + 
            \cH^{\xi,(n)}(r_k,z; 1,y)
        \right)
        =
        \max_{z\in[-1,1]}        \left(
            \cH^{\xi,(n)}(0,x;r_k,z)
            + 
            \cH^{\xi,(n)}(r_k,z; 1,y)
        \right).
    \end{align*}
    Fix $\xi\in\Good_{\e,n}\cap\Max_{\e,n}^{k}$.
    By the multi-point density formula \eqref{e:scaled-CDRP-fdd} and the positivity of $\cZ$ (\cref{t:AJRAS-Z-existence}\ref{property-strict-positivity}), we have 
    \begin{align*}
        \cZ^{\xi,(n)}_{\reg}
        &=
        n^{2/3}\,
        \P^{\xi,(n)}(\XTF_k)
        \int_\R 
        \cZ^{\xi,(n)}(0,x; r_k, z)\cZ^{\xi,(n)}(r_k, z; 1,y)\,dz\\
        &\ge 
        n^{2/3}\,
        \P^{\xi,(n)}(\XTF_k)\int_I 
        \cZ^{\xi,(n)}(0,x; r_k, z)\cZ^{\xi,(n)}(r_k, z; 1,y)\,dz\\
        &\ge n^{2/3}\,
        \P^{\xi,(n)}(\XTF_k)\,
        |I|
        \exp\left(-2\e^{1/2}n^{1/3}\log^{4}(1/\e)\right)
        \max_{z\in [-1,1]}
            \cZ^{\xi,(n)}(0,x;r_k,z)
            \cZ^{\xi,(n)}(r_k,z; 1,y),
    \end{align*}
    where in the last line we used that $\xi\in \Good_{\e,n} \subset \HHol_{\e,n}$ (see \eqref{e:def-holder-event}) to compare the integrand to its maximum on $I$,  along with the fact that the maximum on $[-1,1]$ is attained on $I$.
    Also, since $\xi\in \Good_{\e,n}\subset\HReg_{\e,n}$, by the definition \eqref{e:def-HReg} of the latter, we can lower bound $\P^{\xi,(n)}(\XTF_k)\ge 1-e^{-\log^{3/2}(1/\e)} \ge \frac12$.
    Further, since $\xi\in \Max_{\e,n}^{k}$ (see \eqref{e:def-argmax-event}), the RHS is lower bounded by
    \begin{align*}
        &\ge
         n^{2/3} \frac{|I|}{2}
        \,\exp\left(-2\e^{1/2}n^{1/3}\log^{4}(1/\e)\right)\\
        &\qquad\cdot
        \exp\left(
            \e^{1/2} n^{1/3} \log^5(1/\e)
        \right)
        \max_{|z| \le 2\e \log(1/\e)}
            \cZ^{\xi,(n)}(0,x;r_k,z)
            \cZ^{\xi,(n)}(r_k,z; 1,y)\\
        &\ge
        n^{2/3} \frac{|I|}{2}
        \exp\left(
            c\e^{1/2} n^{1/3} \log^5(1/\e)
        \right)
        \max_{|z| \le 2\e \log(1/\e)}
            \cZ^{\xi,(n)}(0,x;r_k,z)
            \cZ^{\xi,(n)}(r_k,z; 1,y),
    \end{align*}
    where the second inequality holds for all sufficiently small $\e>0$ uniformly in $n\ge 1$.

    To upper bound $\cZ^{\xi,(n)}_{\reg}(B_k)$, using similar reasoning to \eqref{e:upper-bound-ZB}, we have
    \begin{multline*}
        \P^{\xi,(n)}\left(
            \XTF_k \cap \left\{
                \inf_{r\in[r_k, r_{k+1}]}
                |X^{(n)}(r)| \le \e
            \right\}
        \right)
        \le 
        \P^{\xi,(n)}\left(
            |X^{(n)}(r_k)|
            \le \e + \e \log(1/\e)
        \right)\\
        \overset{\eqref{e:scaled-CDRP-fdd}}{\le}
        \int_{|z| \le 2\e\log(1/\e)}
        n^{2/3}
        \frac{\cZ^{\xi,(n)}(0,x; r_k, z)\cZ^{\xi,(n)}(r_k,z;1,y)}{\cZ^{\xi,(n)}(0,x;1,y)}\,dz,
    \end{multline*}
    and hence
    \begin{align*}
        \cZ^{\xi,(n)}_{\reg}(B_k)
        &\le 
        4 n^{2/3}  \e \log(1/\e)
        \max_{|z|\le 2 \e \log(1/\e)}
        \cZ^{\xi,(n)}(0,x; r_k, z)\cZ^{\xi,(n)}(r_k,z;1,y).
    \end{align*}
    
    Combining the above estimates and using that $|I|=8\e\log(1/\e)$, we get
    \begin{align*}
        \cZ^{\xi,(n)}_{\reg}
        &\ge 
        n^{2/3}
        \frac{|I|}{2}
        \exp\left(
            c \e^{1/2} n^{1/3} \log^5(1/\e)
        \right)
        \frac{1}{4n^{2/3} \e\log(1/\e)}
        \cZ^{\xi,(n)}_{\reg}(B_k)\\
        &=
        \exp\left(
            c
            \e^{1/2} n^{1/3} \log^5(1/\e)
        \right)
        \cZ^{\xi,(n)}_{\reg}(B_k)
    \end{align*}
    as desired.
\end{proof}

We now deduce an $L^2$ influence bound.
\begin{lemm}[Bulk influence bound]\label{l:bulk-resampling}
    There exists $C>0$ such that for all sufficiently small $\e>0$ 
    and all $k\in \lb 1,M-1\rb$, we have 
    \begin{align*}
        \limsup_{n\to\infty}
        \frac{1}{n^{2/3}}
        \Eres\left[
            \left|
                \cH^{\xi,(n)}_{\reg}
                -
                \cH^{\eta_k,(n)}_{\reg}
            \right|^2
            \1_{\xi,\eta_k\in\Good_{\e,n}}
        \right]
        \le 
        C\e^{2 - 3\gamma} \log^{8}(1/\e).
    \end{align*}
\end{lemm}
We remark that since $\gamma\in(0,\frac16)$ was arbitrary, the above exponent $2-3\gamma$ can be made arbitrarily close to the optimal value of $2$ (see the discussion below \eqref{efron-stein}).
\begin{proof}[Proof of \cref{l:bulk-resampling}]
    We decompose the $L^2$ distance as follows:
    \begin{align}
       \Eres\left[
            \left|
                \cH^{\xi,(n)}_{\reg}
                -
                \cH^{\eta_k,(n)}_{\reg}
            \right|^2
            \1_{\xi,\eta_k\in\Good_{\e,n}}
        \right]
        &=
        \Eres\left[
            \left|
                \cH^{\xi,(n)}_{\reg}
                -
                \cH^{\eta_k,(n)}_{\reg}
            \right|^2
            \1_{\xi,\eta_k\in\Good_{\e,n} \cap \Max_{\e,n}^{k}}
        \right]\label{e:goodmax}\\
        &+
        \Eres\left[
            \left|
                \cH^{\xi,(n)}_{\reg}
                -
                \cH^{\eta_k,(n)}_{\reg}
            \right|^2
            \1_{\xi,\eta_k\in\Good_{\e,n}}
            \1_{\neg\{\xi,\eta_k \in \Max_{\e,n}^{k}\}}
        \right]. \label{e:good-nomax}
    \end{align}

    We first estimate \eqref{e:goodmax}.
    Fix $\xi,\eta_k\in\Good_{\e,n}\cap\Max_{\e,n}^{k}$.
    Since $\xi,\eta_k\in\Good_{\e,n}$, the argument used in \eqref{e:611} yields
    \begin{align*}
        \left|\cH^{\xi,(n)}_{\reg} - \cH^{\eta_k,(n)}_{\reg}\right|
        &= \left|
            \log\left(
            \frac{\cZ^{\xi,(n)}_{\reg}(B_k) + \cZ^{\xi,(n)}_{\reg}(B_k^c)}
            {\cZ^{\eta_k,(n)}_{\reg}(B_k)+\cZ^{\xi,(n)}_{\reg}(B_k^c)}
            \right)
        \right|.
    \end{align*}
    Suppose $\cH^{\xi,(n)}_{\reg} - \cH^{\eta_k,(n)}_{\reg} \ge 0$.
    Then by the above display, the non-negativity of $\cZ^{\eta_k,(n)}_{\reg}$, and \cref{l:brownian-comparison-result},
    \begin{align*}
        \left|\cH^{\xi,(n)}_{\reg} - \cH^{\eta_k,(n)}_{\reg} \right|
        &\le
        \log\left(
            \frac{\cZ^{\xi,(n)}_{\reg}(B_k) + \cZ^{\xi,(n)}_{\reg}(B_k^c)}
            {\cZ^{\xi,(n)}_{\reg}(B_k^c)}
        \right)\\
        &\le \log\left(
            1  + \frac{1}{\exp(c\e^{1/2}n^{1/3}\log^5(1/\e)) - 1}
        \right).
    \end{align*}
    In the case $\cH^{\xi,(n)}_{\reg} - \cH^{\eta_k,(n)}_{\reg} < 0$, a symmetric argument yields the same upper bound.
    We obtain
    \begin{multline*}
        \limsup_{n\to\infty}
        \frac{1}{n^{2/3}}
        \Eres\left[
            \left|
                \cH^{\xi,(n)}_{\reg}
                -
                \cH^{\eta_k,(n)}_{\reg}
            \right|^2
            \1_{\xi,\eta_k\in\Good_{\e,n} \cap \Max_{\e,n}^{k}}
        \right]
        \\
        \begin{aligned}
            &\le
            \limsup_{n\to\infty}
            \frac{1}{n^{2/3}}
            \log^2\left(
                1+\frac{1}{\exp(c\e^{1/2}n^{1/3}\log^5(1/\e))-1}
            \right)
            = 0.
        \end{aligned}
    \end{multline*}

    We next estimate \eqref{e:good-nomax}.
    Using \cref{l:a-priori-change-full} to bound $\left|\cH^{\xi,(n)}_{\reg}-\cH^{\eta_k,(n)}_{\reg}\right|$ on $\Good_{\e,n}$, followed by the probability bound \cref{l:Max-probability-bound}, we get 
    \begin{multline*}
        \limsup_{n\to\infty}\frac{1}{n^{2/3}}
        \Eres\left[
            \left|
                \cH^{\xi,(n)}_{\reg}
                -
                \cH^{\eta_k,(n)}_{\reg}
            \right|^2
            \1_{\xi,\eta_k\in\Good_{\e,n}}
            \1_{\neg\{\xi,\eta_k \in \Max_{\e,n}^{k}\}}
        \right]\\
        \begin{aligned}
            &\le
            C\e\log^{8}(1/\e) \cdot 
            \limsup_{n\to\infty} 2\P\left(\neg \Max_{\e,n}^{k}\right)\\
            &\le C' \e^{2- 3\gamma}\log^{8}(1/\e),
        \end{aligned}
    \end{multline*}
    where $C,C'>0$ do not depend on $k$.
    This completes the proof of \cref{l:bulk-resampling}.
\end{proof}

\subsection{Concluding the proof}\label{s:combining-edge-bulk}

\cref{t:2main-reformulated} is an immediate corollary of the Efron--Stein inequality, Lemma \ref{l:restrict-to-Good}, and the bulk/edge influence estimates proved in Lemmas \ref{l:edge-L2-resampling} and \ref{l:bulk-resampling}.

\begin{proof}[Proof of \cref{t:2main-reformulated}]
    By the Efron--Stein bound \eqref{e:efron-stein-bound} and \cref{l:restrict-to-Good},  we have
    \begin{multline*}
        \limsup_{n\to\infty}\frac{1}{n^{2/3}}
        \Eres\left[
            \left|\cH^{\xi,(n)} - \cH^{\eta_B,(n)}\right|^2
        \right]\\
        \begin{aligned}
        \overset{\eqref{e:efron-stein-bound}}&{\le}
        \sum_{k=0}^M
        \limsup_{n\to\infty}\frac{1}{n^{2/3}}
        \Eres\left[
            \left|\cH^{\xi,(n)} - \cH^{\eta_k,(n)}\right|^2
        \right]\\
        &\le
        (M+1)e^{-c\log^2(1/\e)}
        +
        2
        \sum_{k=0}^M
        \limsup_{n\to\infty}\frac{1}{n^{2/3}}
        \Eres\left[
            \left|\cH^{\xi,(n)}_{\reg} - \cH^{\eta_k,(n)}_{\reg}\right|^2
            \1_{\xi,\eta_k\in\Good_{\e,n}}
        \right].
        \end{aligned}
    \end{multline*}
    To bound the sum on the RHS, we use Lemma \ref{l:edge-L2-resampling} to control the edge terms ($k\in\{0,M\}$), and Lemma \ref{l:bulk-resampling} to control the bulk terms ($k\in\lb 1,M-1\rb$):
    \begin{multline*}
        \sum_{k=0}^M
        \limsup_{n\to\infty}\frac{1}{n^{2/3}}
        \Eres\left[
            \left|\cH^{\xi,(n)}_{\reg} - \cH^{\eta_k,(n)}_{\reg}\right|^2
            \1_{\xi,\eta_k\in\Good_{\e,n}}
        \right]\\
        \begin{aligned}
            &= 
        \sum_{k\in\{0,M\}}
        \limsup_{n\to\infty}\frac{1}{n^{2/3}}
        \Eres\left[
            \left|\cH^{\xi,(n)}_{\reg} - \cH^{\eta_k,(n)}_{\reg}\right|^2
            \1_{\xi,\eta_k\in\Good_{\e,n}}
        \right]\\
        &\qquad+
        \sum_{k=1}^{M-1}
        \limsup_{n\to\infty}\frac{1}{n^{2/3}}
        \Eres\left[
            \left|\cH^{\xi,(n)}_{\reg} - \cH^{\eta_k,(n)}_{\reg}\right|^2
            \1_{\xi,\eta_k\in\Good_{\e,n}}
        \right]\\
        &\ls 
        \e^{2\gamma/3} \log^8(1/\e)
        + M \e^{2-3\gamma}\log^{8}(1/\e).
        \end{aligned}
    \end{multline*} 
    Finally, since $M\ls \e^{-3/2}$, we get 
    \begin{align*}
        \limsup_{n\to\infty}\frac{1}{n^{2/3}}
        \Eres\left[
            \left|\cH^{\xi,(n)} - \cH^{\eta_B,(n)}\right|^2
        \right]
        &\ls
        \e^{-3/2} e^{-c\log^2(1/\e)}
        +
        \e^{2\gamma/3} \log^8(1/\e)
        + \e^{\frac12 - 3\gamma} \log^{8}(1/\e)\\
        &\ls \e^{\alpha}
    \end{align*}
    for any $0< \alpha < \min\{\frac{2\gamma}{3},\,\frac12 - 3\gamma\}$ (recall $0<\gamma < 1/6$).
    This finishes the proof of \cref{t:2main-reformulated}.
\end{proof}

The argument implemented in the preceding section can now be adapted to prove \cref{t:black-noise}.

\section{The directed landscape is a two-dimensional black noise}\label{s:black-noise}

We start by recalling the abstract definition of a noise from \cite[Section 6]{EF16} (see also \cite[Definition 3d1 and Section 11b]{Tsir}) which we have deferred until now.

\begin{defn}[Noise]\label{d:noise}
    For $d\in\N$, a \emph{$d$-dimensional noise} consists of a probability space $(\Omega,\cF,\P)$, a collection of sub-$\sigma$-algebras $\{\cF_R\}_{R\in\Rect}$ where $\Rect$ is the set of open $d$-dimensional rectangles in $\R^d$, and a measurable $\R^d$-action $\{\theta_h\}_{h\in\R^d}$ on $(\Omega,\cF,\P)$,
    such that the following properties hold.
    \begin{enumerate}[label={\rm(N\arabic*)},ref={\rm N\arabic*}]
        \item\label{noise-generation}
        \textup{(Exhaustive filtration).} $\cF = \bigvee_{R\in\Rect} \cF_R$.
        \item\label{noise-translation} 
        \textup{(Translation covariance).}
        For all $h\in\R^d$, we have $\theta_h(\cF_R) = \cF_{R+h}$.
        \item\label{noise-indep} 
        \textup{(Independence).}
        For all $R_1,R_2\in\Rect$ with $R_1\cap R_2=\varnothing$, the $\sigma$-algebras $\cF_{R_1}$ and $\cF_{R_2}$ are independent.
        \item\label{noise-join}
        \textup{(Join).}
        For all $R_1,R_2,R_3\in\Rect$ with $R_1\cap R_2=\varnothing$ and $\overline{R_1\cup R_2}=\overline{R_3}$, we have $\cF_{R_1}\vee \cF_{R_2} = \cF_{R_3}$ up to null sets.
    \end{enumerate}
\end{defn}

Note that the above notation for $\sigma$-algebras was already used in the context of white noise in Sections \ref{s:meas-short}--\ref{s:es}.
Since the above $\sigma$-algebras play the same role in this general setting as those in the white noise context (indeed, white noise satisfies the above axioms), we chose to overload the notation to make parallels between the two settings more apparent.
We will not need any white noise inputs in the present section.

\begin{defn}[Black noise]\label{def:black-noise}
    Let $(\Omega,\cF,\{\cF_R\}_{R\in\Rect}, \P,\{\theta_h\}_{h\in\R^d})$ be a $d$-dimensional noise.
    A random variable $X\in L^2(\Omega,\cF,\P)$ is called \emph{linear} if 
    for all $R_1,R_2,R_3\in\Rect$ with $R_1\cap R_2=\varnothing$ and $\overline{R_1\cup R_2}=\overline{R_3}$, we have
    \begin{align*}
        \E[X|\cF_{R_3}] = \E[X|\cF_{R_1}]+\E[X|\cF_{R_2}].
    \end{align*}
    The noise is \emph{black} if the only linear random variable is zero.
    Note that a noise being black depends on the $\sigma$-algebras $\{\cF_R\}_{R\in\Rect}$ only up to $\P$-null sets.
\end{defn}

\begin{rem}[Unpacking {\cref{def:black-noise}}]\label{r:informal-def-black-noise}
We presented the above definition for its brevity, but it is a bit opaque, and unfortunately a complete presentation of a more conceptual definition (such as \cite[Definition 7a1]{Tsir} or \cite[Definition 6d1]{Tsi04}) would require significant preparation.
As a compromise, we now supplement \cref{def:black-noise} with a brief informal discussion on the notion of black noise.
For a much more comprehensive overview, see \cite{Tsir,Tsi04}.

It is helpful to first consider a non-example of a black noise, namely Gaussian white noise on $\R^d$.
As alluded to several times already, Gaussian white noise is stationary under Ornstein--Uhlenbeck (OU) dynamics, making the latter a natural means of perturbing disorder.
The generator of the OU semigroup is an unbounded self-adjoint operator on the space $L^2(\Omega,\cF,\P)$ of functions of Gaussian white noise, with spectrum $\{0,-1,-2,\dots\}$.
The linear random variables are exactly the eigenfunctions with eigenvalue $-1$.
A fundamental result (Wiener chaos decomposition) says that $L^2(\Omega,\cF,\P)$ admits an orthogonal decomposition into eigenspaces of the OU generator.
Moreover, the $0$-eigenspace (constants) and $-1$-eigenspace (linear random variables) together generate all the other eigenspaces, in a certain sense.
This was formulated in different language in \cref{p:chaos-expansion} (see also \cite[Section 1.4]{nualart}).

It turns out that for a general noise as in \cref{d:noise}, there is a natural analogue of OU dynamics.
Just like for Gaussian white noise, in this general setting, the OU generator's eigenvalues are non-positive integers, the $-1$-eigenspace is the set of linear random variables, and the $0$-eigenspace (constants) and $-1$-eigenspace generate all the other eigenspaces (see \cite[Theorem 6a4]{Tsir} for details).
However, unlike the Gaussian case, the OU eigenfunctions generally do not span the whole $L^2(\Omega,\cF,\P)$---if they do, the noise is called \emph{classical}.
Thus for a non-classical noise, there is a nontrivial subspace of $L^2(\Omega,\cF,\P)$ orthogonal to all the eigenspaces---informally, it is the eigenspace corresponding to eigenvalue $-\infty$.
\footnote{In reality, the generator is ill-defined in the non-classical case; only the OU semigroup generalizes to the setting of \cref{d:noise}. However, the language of generators is convenient for exposition.}
A \emph{black noise} is a noise with no classical part: every non-constant observable lies in the $-\infty$-eigenspace of the OU generator.
This is indeed equivalent to \cref{def:black-noise}: if the only linear random variable is zero, or equivalently if the $-1$-eigenspace is trivial, then the $-n$-eigenspace must be trivial for every $1\le n<\infty$.

This presents a natural opportunity to indicate heuristically how black noise connects to the notion of noise sensitivity.
Let $L$ be the OU generator of a black noise, and let $t\mapsto P_t = e^{tL}$ be its semigroup.
Every mean-zero observable $f$ lies in the $-\infty$-eigenspace of $L$, and thus satisfies $P_t f = e^{t(-\infty)}f = 0$ for every $t>0$.
That is, any arbitrarily small perturbation of the noise causes the value of $f$ to become uncorrelated with its initial state.
For a much more comprehensive discussion along these lines, see \cite[Section 5]{Tsi04} and \cite{GS}.
\end{rem}

We now specialize to the directed landscape, recalling first the setup from \cref{s:iop-black-noise}.
Let $\Omega \coloneqq  C(\Rup)$ and $\cF\coloneqq \sigma(\cL(u):u\in\Rup) = \cB(C(\Rup))$, and let $\P$ be the law of the directed landscape $\cL$ (see \cref{d:dl}).
We realize the directed landscape as the canonical process on $(\Omega,\cF,\P)$, i.e. $\cL:\Omega\to\Omega$ is the identity map.
Let $\Rect = \{(s,t)\times(a,b)\}_{s<t, a<b}$ be the collection of open rectangles in $\R^2$.
It will be convenient later to also consider infinite strips and rectangles containing some of their boundary edges; let $\overline{\Rect}$ be the collection of all rectangles of the form $I\times (a,b)$ where $-\infty \le a<b\le \infty$ and $I\in\{(s,t),[s,t],(s,t],[s,t)\}$ for $-\infty<s<t<\infty$.
Like in \cref{s:iop-black-noise}, for any $R\in\overline{\Rect}$ and any $(s,x),(t,y)\in R$ with $s<t$, we define (by analogy with \eqref{e:DL-as-LPP}) the \emph{restricted length}
\begin{align}\label{e:restricted-length}
    \cL_R(s,x;t,y) \coloneqq  
    \sup\left\{
        \int_s^t d\cL\circ \pi : 
        \pi\in C([s,t]),
        \quad
        \pi(s)=x, \pi(t)=y,
        \quad (r,\pi(r))\in R,\;\forall r\in[s,t]
    \right\},
\end{align}
where $\int_s^t d\cL\circ \pi$ is defined in \cref{def:length-dl}.
As mentioned in \cref{s:iop-black-noise}, $\cL_R$ can be viewed as the ``internal metric'' induced on $R$ by $\cL$.
We will later prove (see \cref{l:restricted-L-properties}\ref{property-LR-continuous} and \eqref{e:852S})
that $\cL_R$ defines a measurable map $(\Omega,\cF) \to C(\Rup(R))$, where $\Rup(R)\coloneqq \{(s,x;t,y)\in\Rup : (s,x),(t,y)\in R\}$, and where $C(\Rup(R))$ has the topology of uniform convergence on compact sets.
This permits us to define the $\sigma$-algebra $\cF_R\subset\cF$ generated by $\cL_R$:
\begin{align}\label{e:restricted-sigma-alg}
    \cF_R \coloneqq  
    \sigma(\cL_R)
    =
    \sigma\bigl(
        \cL_R(u) : u \in \Rup(R)
    \bigr).
\end{align}
Finally, we have a natural $\R^2$-action $\{\theta_h\}_{h\in\R^2}$ on $(\Omega,\cF,\P)$ given by translation: for $h=(h_1,h_2)\in\R^2$ and $\omega\in\Omega=C(\Rup)$, define $\theta_h\omega$ by 
\begin{align}\label{e:translation-action-DL}
    (\theta_h\omega)(s,x;t,y) \coloneqq  \omega(s - h_1, x - h_2; t - h_1, y - h_2)
    \qquad\text{for } (s,x;t,y)\in\Rup.
\end{align}

As indicated in \cref{s:iop-black-noise}, our \cref{t:black-noise} follows immediately from the following \cref{t:2d-noise}.

\begin{thm}[The directed landscape is a 2D noise]\label{t:2d-noise}
    Let $(\Omega,\cF,\P)$ be the canonical probability space of the directed landscape $\cL$, let $\{\cF_R\}_{R\in\Rect}$ be as in \eqref{e:restricted-sigma-alg}, and let $\{\theta_h\}_{h\in\R^2}$ be as in \eqref{e:translation-action-DL}.
    The quintuple $(\Omega,\cF, \{\cF_R\}_{R\in\Rect}, \P, \{\theta_h\}_{h\in\R^2})$ is a two-dimensional noise.
\end{thm}

\begin{proof}[Proof of \cref{t:black-noise}]
    By \cite[Theorem 1.10]{HP24}, the temporal filtration $\{\cF_{[s,t]\times\R}\}_{s<t}$ 
    is a one-dimensional black noise (with $\theta$ in \eqref{e:translation-action-DL} restricted to the subgroup $\{(h_1,0)\}\subset \R^2$).
    Although this filtration is defined in terms of the closed strip $[s,t]\times \R$, since $\cL_{[s,t]\times\R}(s',\smallbullet;t',\smallbullet)=\cL(s',\smallbullet;t',\smallbullet)$ for any $s\le s' < t' \le t$, it follows by continuity of $\cL$ that  $\cF_{[s,t]\times\R}=\cF_{(s,t)\times\R}$ up to null sets,
    and hence $\{\cF_{(s,t)\times \R}\}_{s<t}$ is also a one-dimensional black noise.
    Combining this with \cref{t:2d-noise} and a general result of \cite[Theorem 1.13 and Section 1.6]{TsirBA} showing that a two-dimensional noise with a one-dimensional black noise marginal must be a two-dimensional black noise, we immediately conclude the result.
    For the reader's convenience, we include below an alternative argument to the same end.

    We must show that for the two-dimensional noise of \cref{t:2d-noise}, the only linear random variable is zero (see \cref{def:black-noise}).
    Let $X$ be a linear random variable. 
    Then for any $s<r<t$ and $n\ge 1$,
    \begin{align}\label{e:linear}
        \E[X|\cF_{(s,t)\times (-n,n)}] = 
        \E[X|\cF_{(s,r)\times(-n,n)}] + \E[X|\cF_{(r,t)\times(-n,n)}].
    \end{align}
    We claim that
    \begin{align*}
        \bigvee_{n=1}^\infty \cF_{(s,t)\times (-n,n)}
        = \cF_{(s,t)\times \R}
        \qquad \text{up to null sets.}
    \end{align*}
    The inclusion $\subset$ follows from the fact that $\cL_{(s,t)\times(-n,n)}$ is a measurable function of the restriction of $\cL$ to $\Rup((s,t)\times\R)$ (proven in \eqref{e:852S}), 
    and the inclusion $\supset$ follows by noting that for any fixed $u\in \Rup((s,t)\times\R)$, we have $\cL_{(s,t)\times (-n,n)}(u) \uparrow \cL(u)$ as $n\to\infty$ because the $\cL(u)$-geodesic stays inside $(s,t)\times (-n,n)$ for large enough $n$.  
    Taking $n\to\infty$ in \eqref{e:linear} and using martingale convergence yields
    \begin{align*}
        \E[X|\cF_{(s,t)\times \R}] = 
        \E[X|\cF_{(s,r)\times\R}] + \E[X|\cF_{(r,t)\times\R}].
    \end{align*}
    So $X$ is also linear with respect to the one-dimensional noise $\{\cF_{(s,t)\times\R}\}_{s<t}$.
    The latter is black, hence $X=0$.
\end{proof}

The rest of this section is devoted to proving \cref{t:2d-noise}.
We need to verify that the properties \ref{noise-generation}--\ref{noise-join} hold for $(\Omega,\cF,\{\cF_R\}_{R\in\Rect},\P, \{\theta_h\}_{h\in\R^2})$.
As indicated in \cref{s:iop-black-noise}, of the four properties, the hardest to check is \ref{noise-join}.
Indeed, we will first quickly prove \ref{noise-generation}--\ref{noise-indep}, with the rest of the section devoted to establishing \ref{noise-join}.

\subsection{Proofs of Properties \ref{noise-generation}--\ref{noise-indep}}\label{s:first-three-properties}

\hphantom{a} 

\addtocontents{toc}{\SkipTocEntry}
\subsection*{Proof of Property \ref{noise-generation} ($\cF$ is generated by $\{\cF_R\}_{R\in\Rect}$)}

Since by definition $\cF = \sigma(\cL(u):u\in\Rup)$,
it suffices to show that for any fixed $u\in\Rup$, the random variable $\cL(u)$ is measurable with respect to $\bigvee_{n=1}^\infty \cF_{R_n}$ for $R_n\coloneqq (-n,n)\times(-n,n)$.
Observe that $\cL_{R_n}(u) \le \cL_{R_{n+1}}(u)$ for all $n\ge 1$.
Moreover, almost surely, the graph of the $\cL(u)$-geodesic is contained in $R_n$ for all sufficiently large $n$, and hence $\cL_{R_n}(u) = \cL(u)$.
In particular, $\cL_{R_n}(u) \uparrow \cL(u)$ almost surely.
Since $\cL_{R_n}(u)$ is $\cF_{R_n}$-measurable, the increasing limit $\cL(u)$ must be measurable with respect to $\bigvee_{n=1}^\infty \cF_{R_n}$.
\qed

\addtocontents{toc}{\SkipTocEntry}
\subsection*{Proof of Property \ref{noise-translation} (translation invariance)}

Property \ref{noise-translation} is immediate from the translation invariance of the directed landscape (\cref{p:landscape-symmetries}\ref{landscape-translation}).
\qed

\addtocontents{toc}{\SkipTocEntry}
\subsection*{Proof of Property \ref{noise-indep} (independence across disjoint rectangles)}

We need the following technical lemma as input to our arguments.
This will also be used in the proof of Property \ref{noise-join}.
As before, for $R\in\Rect$ we write
\begin{align*}
    \Rup(R)\coloneqq \{(s,x;t,y)\in\Rup : (s,x),(t,y)\in R\}.
\end{align*}

\begin{figure}[tbhp]
    \centering
    \includegraphics[width=0.8\textwidth]{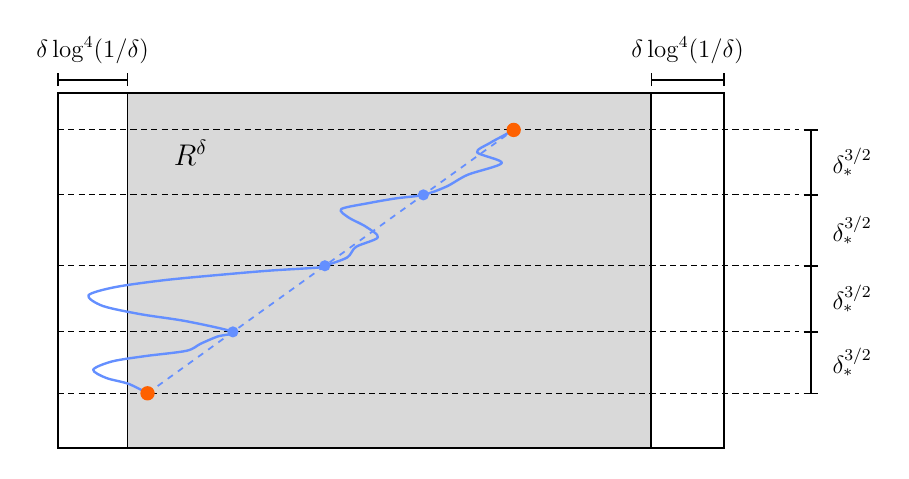}
    \caption{
        Depicted are the rectangles $R^\d\subset R$, with $R^\d$ shaded and $R\setminus R^\d$ unshaded.
        We fix any $u\in\Rup(R^\d)$ and draw its endpoints in orange.
        By geodesic transversal fluctuation estimates (\cref{p:GZ-uniform-TF}), there exists a random $\d_*>0$ such that all geodesics of height $\le \d_*^{3/2}$ with endpoints inside $R^\d$ do not exit $R$. 
        Thus by concatenating $O(\d_*^{-3/2})$-many such geodesics (blue), we obtain a piecewise-geodesic path that is contained in $R$ and connects the endpoints of $u$.
        }
    \label{fig:noise-thin-cylinders}
\end{figure}

\begin{lemm}[Properties of the restricted length $\cL_R$]\label{l:restricted-L-properties}
    Almost surely, the following properties hold for every $R = (s_0,t_0) \times (x_0,y_0) \in \Rect$.
    \begin{enumerate}[label={\rm(\roman*)}]
        \item\label{property-LR-convolution} 
        \textup{(Metric composition).}
        For every $(s,x;t,y)\in\Rup(R)$ and every $s<r<t$, we have 
        \begin{align*}
            \cL_R(s,x;t,y) = \sup_{z\in(x_0,y_0)}\left(
                \cL_R(s,x;r,z) + \cL_R(r,z; t,y)
            \right).
        \end{align*}
        \item\label{property-LR-finite}\textup{(Finiteness).}
        For every $(s,x;t,y)\in\Rup(R)$, we have $|\cL_R(s,x;t,y)|<\infty$.
        \item\label{property-LR-continuous}
        \textup{(Continuity).}
        $\cL_R$ is a continuous function $\Rup(R) \to \R$.
    \end{enumerate} 
\end{lemm}
\begin{proof}

    \ref{property-LR-convolution}:
    Immediate from the definition \eqref{e:restricted-length}.

    \ref{property-LR-finite}:
    By \eqref{e:DL-as-LPP} and \eqref{e:restricted-length}, we have  $\cL_R \le \cL < \infty$ pointwise.
    So we just need to prove that $\cL_R > -\infty$.
    As depicted in Figure \ref{fig:noise-thin-cylinders}, the proof proceeds by showing that for any $u\in\Rup(R)$, there exists a path with endpoints given by $u$ formed by concatenating finitely many \emph{unrestricted} $\cL$-geodesics that are all contained in $R$. 
    We will be brief here, but a detailed argument to this effect appears later in \cref{l:moments-wh-L} (arguments of this style have also appeared before, see for instance \cite[Proposition 4.5]{BGZ21}).
    To start, for any $\d>0$, write $R^\d \coloneqq  (s_0,t_0) \times (x_0 + \d\log^4(1/\d),\; y_0 - \d\log^4(1/\d))$.
    Since $R$ is an open rectangle,  we have $u\in \Rup(R^\d)$ for all sufficiently small $\d>0$.

    For any $(s,x;t,y)\in\Rup(R^\d)$ with $t-s \le \d^{3/2}$, 
    since the spatial points $x,y$ lie at distance $> \d\log^4(1/\d)$ from the boundary of $(x_0,y_0)$, 
    geodesic transversal fluctuation estimates imply that typically
    the geodesic $\Pi_{(s,x),(t,y)}$ does not exit $R$ during its journey, i.e.
    \begin{align*}
        (r, \Pi_{(s,x),(t,y)}(r)) \in R\qquad\text{for all }r\in[s,t].
    \end{align*}
    In fact, by \cref{p:GZ-uniform-TF}, the above holds for every such $(s,x;t,y)$ simultaneously with probability tending to $1$ as $\d\to 0$.
    Said differently, \cref{p:GZ-uniform-TF} implies that there exists a random $\d_*\in(0, \d/2)$ such that for all $(s,x;t,y)\in\Rup(R^\d)$ with $t-s \le \d_*^{3/2}$, the geodesic $\Pi_{(s,x),(t,y)}$ does not exit $R$ (see Figure \ref{fig:noise-thin-cylinders}).
    In particular, choosing any sequence of $O(\d_*^{-3/2})$-many points in $(x_0 + \d\log^4(1/\d),\; y_0 - \d\log^4(1/\d))$, and connecting them with geodesics with each geodesic running for time $O(\d_*^{3/2})$,
    we obtain a path $\pi$ that is contained in $R$, with endpoints $u$, and whose length  $\int d\cL\circ \pi$ can be written as $\cL(u_1)+\cL(u_2)+\cdots+\cL(u_n)$ where $n=O(\d_*^{-3/2})$,
    and where the big-$O$s conceal deterministic constants that depend only on $R$.
    Since $\cL>-\infty$ pointwise, the claim follows.

    In fact, it will be useful to upgrade the above to a uniform lower bound on $\cL_R(u')$ for all points $u'$ in a small neighborhood of $u$.
    This follows by noticing that for any $u'$ sufficiently close to $u$, 
    say $\norm{u'-u}<\d_*^{100}$,
    repeating the above construction produces a piecewise-geodesic path with endpoints $u'$ that is contained in $R$ and 
    comprised of the \emph{same} number $n$ of geodesic segments.
    Further, by \cref{p:dov-pointwise}, there is a random constant $\fC>0$ such that
    for every $u'$, every geodesic segment involved in the corresponding path has length at least $-\fC\,\d_*^{-3/2}$.
    Since the number of such segments is $n = O(\d_*^{-3/2})$, summing this length lower bound shows that for any sufficiently small open neighborhood $U$ of $u$, we have
    \begin{align}\label{e:unif-lower-cty}
        \inf_{u'\in U} \cL_R(u') \ge -C\,\fC\,\frac{1}{\d_*^3},
    \end{align}
    where $C=C(R)>0$ is deterministic.

    \ref{property-LR-continuous}: 
    We will use a similar idea to \ref{property-LR-finite}:
    any restricted path that is competitive in \eqref{e:restricted-length} must locally look like a geodesic,
    and hence near the starting point and ending point, one can couple to an (unrestricted) directed landscape and use its continuity.
    It will be useful to refer to the \emph{restricted geodesic}, i.e. the path maximizing \eqref{e:restricted-length}.
    Although restricted geodesics can be shown to exist, we will not do so.
    Thus, in the argument below, ``the restricted geodesic'' is only a shorthand for ``any path that could possibly attain the maximum in \eqref{e:restricted-length}.''

    Let $\d_*>0$ be the random constant defined above, fix $u=(s,x;t,y)\in\Rup(R)$,
    and let $U$ be an open neighborhood of $u$ small enough that the uniform lower bound \eqref{e:unif-lower-cty} holds.
    We claim that for all $u'=(s,x;t',y')\in U$ (note the starting points coincide with those of $u$), we have
    \begin{align}\label{e:approximate-convolution}
            \cL_R(u') =
            \sup_{z\in (x_0+\d_*^{2}, y_0-\d_*^2)}
            \Bigl(
                \cL_R(s,x; t-\delta_*^{10}, z)
                +
                \cL_R(t-\delta_*^{10},z; t',y')
            \Bigr).
    \end{align}
    In words, this says that at time $t-\delta_*^{10}$, the restricted geodesic
    must intersect the interval $(x_0+\d_*^{2}, y_0-\d_*^{2})$.
    To see this, note that by Proposition \ref{p:dov-pointwise}, the maximum $\cL$-length of \emph{any} path with endpoints $u'$ that passes through a point $(t-\delta_*^{10}, z)$ with $z\not\in  (x_0+\d_*^{2}, y_0-\d_*^2)$ is upper bounded by $-\fC\,\delta_*^{-8}$,  which is much less than the lower bound of $\cL_R(u') \gs -\fC\,\d_*^{-3}$ asserted in \eqref{e:unif-lower-cty}. 
    Next, by the aforementioned transversal fluctuation bound \cref{p:GZ-uniform-TF}, the $\cL$-geodesic between any $(t',y')$ and any such $(t-\delta_*^{10},z)$ must stay within distance $\delta_*^{6}$ of the straight line joining its endpoints, and hence must be entirely contained in $R$.
    Therefore, the restricted geodesic must coincide with some $\cL$-geodesic  on the time interval $[t-\d_*^{10}, t']$.
    This upgrades \eqref{e:approximate-convolution} to the following identity:
    \begin{align}\label{e:approximate-convolution1}
            \cL_R(u') =
            \sup_{z\in (x_0+\d_*^{2}, y_0-\d_*^2)}
            \Bigl(
                \cL_R(s,x; t-\delta_*^{10}, z)
                +
                \cL(t-\delta_*^{10},z; t',y')
            \Bigr),
        \end{align}
    where in the second term $\cL_R$ has been replaced by $\cL$.
    Continuity of $\cL_R(s,x;t,y)$ in the last two coordinates now follows from the continuity of $\cL$.
    Namely, given $\e>0$, there exists a neighborhood $U'\subset U$ of $u$ such that
    for all $u'=(s,x;t',y')\in U'$,
    \begin{align*}
        \sup_{z\in (x_0+ \d_*^2, y_0 - \d_*^2)}
        \left|
        \cL(t-\delta_*^{10}, z; t,y)
        -\cL(t-\delta_*^{10},z; t',y')
        \right|
        <\e.
    \end{align*}
    A symmetric argument for the first two coordinates $(s,x)$ along with the triangle inequality implies continuity in all four coordinates.    
\end{proof}

With the above lemma in hand, we now prove Property \ref{noise-indep} using a spatial mixing property of $\cL$ as sketched in \cref{s:iop-black-noise}.

\begin{figure}[tbhp]
    \centering \hspace{-1cm}
    \includegraphics[width=0.9\textwidth]{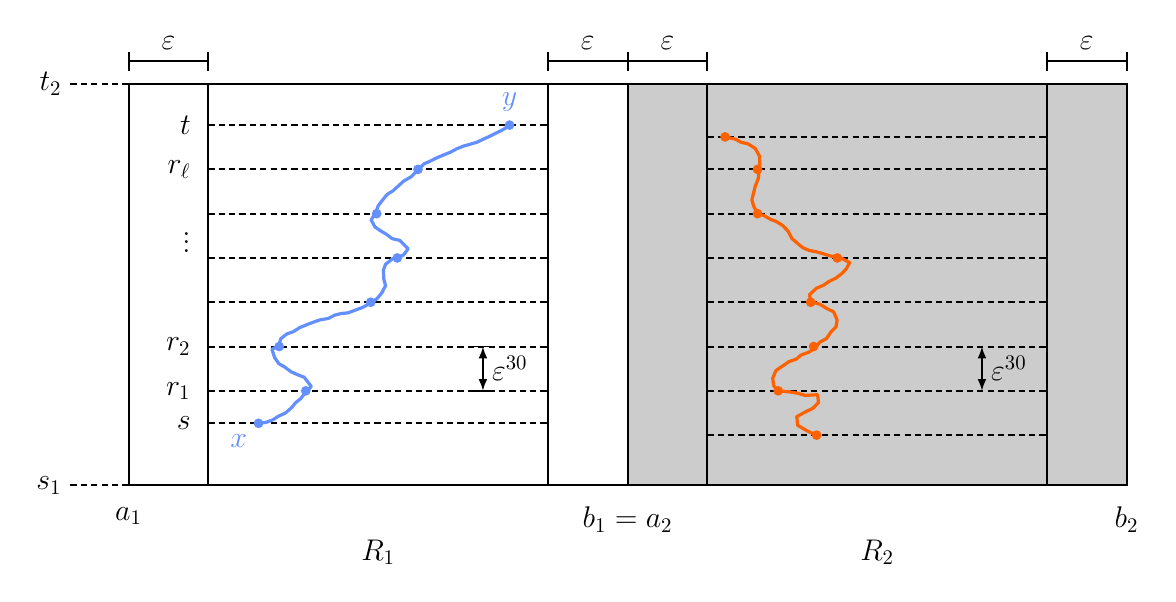}
    \caption{
        The proof of Property \ref{noise-indep} (disjoint independence).
        We consider disjoint open rectangles $R_1$ (unshaded) and $R_2$ (shaded), and aim to prove independence of $\cL_{R_1}$ and $\cL_{R_2}$.
        Depicted is the case where $R_1,R_2$ share a face and have the same bottom/top coordinates.
        By combining the metric composition property of $\cL_{R_i}$ with $\cL$-geodesic transversal fluctuation estimates, 
        we can approximate $\cL_{R_i}$-lengths by concatenating $\cL$-geodesics (blue and orange)  of height $\le \e^{30}$ with endpoints in $(a_i+\e,b_i-\e)$, analogous to Figure \ref{fig:noise-thin-cylinders}.
        Thus it suffices to establish asymptotic independence of such $\cL$-lengths as $\e\to 0$.
        Since $R_1,R_2$ are disjoint, the horizontal distance between strips of the form $[r_j, r_{j+1}]\times(a_1+\e,b_1-\e)$ and $[r_j, r_{j+1}]\times(a_2+\e,b_2-\e)$ is at least $2\e$.
        The desired independence then follows by a spatial mixing property of $\cL$.
        }
    \label{fig:noise-independence}
\end{figure}

\begin{proof}[Proof of Property \ref{noise-indep}]
    
    For $i\in\{1,2\}$ let $R_i = (s_i,t_i)\times(a_i,b_i)\in \Rect$, with $R_1\cap R_2=\varnothing$.
    If $(s_1,t_1)\cap(s_2,t_2)=\varnothing$, then independence of $\cF_{R_1},\cF_{R_2}$ follows from the independence of the directed landscape's temporal increments (\cref{p:landscape-symmetries}\ref{landscale-independent-increments}).
    Assume instead $(s_1,t_1)\cap(s_2,t_2)\ne\varnothing$, so $(a_1,b_1)\cap(a_2,b_2)=\varnothing$.
    We assume without loss of generality that $b_1\le a_2$.
    For notational simplicity we will also assume that $(s_1,t_1)=(s_2,t_2)$.
    The general case follows from a straightforward modification of the argument.

    Our goal is to show that $\cL_{R_1}$ and $\cL_{R_2}$ are independent.
    It will become apparent through the course of the argument that $\cL_{R_i}$ is a measurable map $\Omega \to C(\Rup(R_i))$ (see \eqref{e:852S} below).
    Therefore, it suffices to prove independence of all finite-dimensional distributions of $\cL_{R_1}$ and $\cL_{R_2}$.

    Fix any $s_1<s<t<t_1$ and $\e>0$, and define the set of mesh times (depicted in Figure \ref{fig:noise-independence})
    \begin{align*}
        \{s=r_0<r_1<\cdots< r_{\ell+1}=t\}
        &\coloneqq 
        ((\e^{30}\Z)\cap (s,t))\cup \{s,t\}.
    \end{align*}
    As is probably apparent, the choice of $30$ as the exponent in the mesh size is arbitrary; any large enough number would suffice.
    Notice that $\ell \le (t_1-s_1) \e^{-30}$, and $r_{j+1}-r_j\le \e^{30}$ for all $j\in\lb 0,\ell\rb$.
    For any $(s,x;t,y)\in \Rup(R_i)$,  by the definition of $\cL_R$ in \eqref{e:restricted-length}, it follows by considering paths restricted to $[s,t]\times(a_i+\e,b_i-\e)$, that we have the pointwise convergence
    \begin{align}\label{e:852}
        \cL_{R_i}
        (s,x;t,y) 
        = 
        \lim_{\e\to 0}
        \sup_{x_1,\dots,x_\ell\in (a_i+\e,b_i-\e)}\sum_{j=0}^{\ell}
        \cL_{R_i}(r_j,x_j;r_{j+1},x_{j+1}),
    \end{align}
    where $x_0\coloneqq x$ and $x_{\ell+1}\coloneqq y$. 
    Further, again by the transversal fluctuation bound \cref{p:GZ-uniform-TF}, for all small enough $\e$, the $\cL$-geodesic between any such $(r_j,x_{j})$ and $(r_{j+1},x_{j+1})$ stays inside $R_i$. 
    Let us call this event $\sA_{i,\e}$, 
    the probability of which can be seen to be at least $1 - Ce^{-c\e^{-10}}$.
    The preceding discussion implies that on the event $\sA_{i,\e}$, we have
    \begin{equation}\label{restricequal}
        \cL_{R_i}(r_j,x_j;r_{j+1},x_{j+1})=\cL(r_j,x_j;r_{j+1},x_{j+1})
        \qquad
        \forall j\in\lb 0,\ell\rb,\;
        \forall x_j,x_{j+1} \in (a_i+\e, b_i-\e).
    \end{equation}
    It follows by Borel--Cantelli that almost surely, for all sufficiently small $\e>0$, we can replace $\cL_{R_i}$ by $\cL$ on  the RHS of \eqref{e:852}  to obtain
    (cf. the caption of Figure \ref{fig:noise-independence})
    \begin{align}\label{e:852S}
        \cL_{R_i}(s,x;t,y) 
        = 
        \lim_{\e\to 0}
        \sup_{x_1,\dots,x_\ell\in (a_i+\e,b_i-\e)}\sum_{j=0}^{\ell}
        \cL(r_j,x_j;r_{j+1},x_{j+1}).
    \end{align}
    For future use, we introduce notation for the RHS above:
    for any function $F:\Rup \to \R$, define $\cM_{i,\e}(F) : \Rup(R_i) \to [-\infty,\infty]$ by
    \begin{align}\label{e:853}
        \cM_{i,\e}(F)(s,x;t,y) \coloneqq  
        \sup_{x_1,\dots,x_\ell\in (a_i+\e,b_i-\e)}\sum_{j=0}^{\ell}
        F(r_j,x_j;r_{j+1},x_{j+1}).
    \end{align}
    Note that above, $\ell$ and $r_0,r_1,\dots,r_{\ell+1}$ depend on $(s,x;t,y)$.

    We next use a spatial mixing property of the directed landscape, as stated in \cite[Lemma 3.2]{HP24}, to decouple $\cL_{R_1}$ and $\cL_{R_2}$ on small scales:
    \begin{lemm}[{\cite[Lemma 3.2]{HP24}}]\label{l:DL-mixing}
        Fix $\e>0$.
        There exists a coupling of three directed landscapes $(\cL^\e,\cL_1^\e,\cL_2^\e)$ such that $\cL_1^\e$ and $\cL_2^\e$  are independent, and such that there exists a random variable $X^\e>0$ satisfying
        \begin{align*}
            \cL^\e(s,x;t,y) &= \cL_1^\e(s,x;t,y)\quad \text{for all } s_1 \le s < t \le t_1 \text{ with } t-s \le \e^{30} \text{ and all } x,y < b_1-X^\e,\quad\text{and}\\
            \cL^\e(s,x;t,y) &= \cL_2^\e(s,x;t,y)\quad \text{for all } s_1 \le s < t \le t_1 \text{ with } t-s \le \e^{30} \text{ and all } x,y > b_1 + X^\e.
        \end{align*}
        Moreover, $\P(X^\e > m \e^{19}) \le Ce^{-cm^3}$ for all $m\ge 0$, where $C,c>0$ depend only on $s_1,t_1,b_1$.
    \end{lemm}
    \begin{proof}[Proof of \cref{l:DL-mixing}]
        In fact, \cite[Lemma 3.2]{HP24} is stated as follows:
        \footnote{Strictly speaking, \cite[Lemma 3.2]{HP24} is stated for $r=0$ and $\e=1$, but KPZ scaling (\cref{p:landscape-symmetries}\ref{landscape-scaling}) implies the version stated here.}
        for any fixed $r\in\R$, 
        there exists a coupling $(\cL^\e,\cL^\e_1,\cL^\e_2)$ with $\cL^\e_1,\cL^\e_2$ independent, and a random variable $Y^\e>0$ 
        such that
        \begin{align*}
            \cL^\e(s,x;t,y) &= \cL_1^\e(s,x;t,y)\qquad \text{for all } r \le s < t \le r+\e^{30} \text{ and all } x,y < b_1 - Y^\e,\quad\text{and}\\
            \cL^\e(s,x;t,y) &= \cL_2^\e(s,x;t,y)\qquad \text{for all } r \le s < t \le r+\e^{30} \text{ and all } x,y > b_1 + Y^\e,
        \end{align*}
        and $\P(Y^\e > m \e^{20}) \le Ce^{-cm^{3}}$.
        This is proved in \cite{HP24} by considering a pre-limiting model driven by i.i.d. noise (Exponential LPP), and using geodesic transversal fluctuation bounds to show the following.
        For any two geodesics $\Pi_1$ and $\Pi_2$ across the strip $[r,r+\e^{30}]\times \R$ such that the endpoints of $\Pi_1$ and $\Pi_2$ lie at distance $\gg \e^{20}$ to the left and right of the vertical line $[r,r+\e^{30}]\times\{b_1\}$ respectively, they explore disjoint---hence independent---parts of the noise with high probability.
        \cref{l:DL-mixing} follows from the same proof, along with a union bound to control the transversal fluctuations of geodesics across each of $O(\e^{-30})$-many strips of height $\e^{30}$ covering $[s_1,t_1]\times\R$, which is easily accounted for by reducing the exponent in the tail bound from $20$ to $19$.
    \end{proof}

    By \eqref{e:852S}--\eqref{e:853}, we have
    \begin{align*}
        \left|
            \cL_{R_i}(\cL^\e)(s,x;t,y)
            - 
            \cM_{i,\e}(\cL^\e)(s,x;t,y)
        \right|
        \xrightarrow[\e\to0]{p} 0,
    \end{align*}
    where we wrote $\cL_{R_i}(\cL^\e)$ to emphasize that $\cL_{R_i}$ is a measurable function of a directed landscape $\cL^\e$.
    On the other hand, by \cref{l:DL-mixing},
    we have $\P(X^\e < \e) \to 1$ as $\e\to 0$.
    In particular, since $(a_2+\e)-(b_1-\e)\ge 2\e$ and $r_{j+1}-r_j \le \e^{30}$, it follows from \cref{l:DL-mixing} that 
    \begin{align*}
            \left|
                \cM_{i,\e}(\cL^\e)(s,x;t,y) - \cM_{i,\e}(\cL^\e_i)(s,x;t,y)
            \right|
        \xrightarrow[\e\to0]{p} 0.
    \end{align*}
    We conclude that
    \begin{align}\label{e:8530}
        \left|
            \cL_{R_i}(\cL^\e)(s,x;t,y)
            - 
            \cM_{i,\e}(\cL^\e_i)(s,x;t,y)
        \right|
        \xrightarrow[\e\to0]{p} 0.
    \end{align}
    Since $\cL_1^\e$ and $\cL_2^\e$ are independent, this essentially finishes the proof.
    Below we record the remaining straightforward details for completeness.
    We repeat the above construction for a finite collection of points $u_1,\dots,u_p\in \Rup(R_1)$ and $v_{1},\dots,v_q\in\Rup(R_2)$.
    By \eqref{e:8530} and a union bound, we have 
    \begin{multline*}
        \biggl\lVert\left(
            \bigl\{\cL_{R_1}(\cL^{\e})(u_j)\bigr\}_{j\in\lb 1,p\rb},\;
            \bigl\{
                \cL_{R_2}(\cL^{\e})(v_j)
            \bigr\}_{j\in\lb 1,q\rb}
        \right)\\
        - \left(
            \bigl\{\cM_{1,\e}(\cL_1^{\e})(u_j)\bigr\}_{j\in\lb 1,p\rb},\;
            \bigl\{
                \cM_{2,\e}(\cL_2^{\e})(v_j)
            \bigr\}_{j\in\lb 1,q\rb}
        \right)
        \biggr\rVert_2
        \xrightarrow[\e\to0]{\;p\;}0.
    \end{multline*}
    Since $\cL^{\e}\law \cL$ for all $\e$, passing to a subsequence $\e_k\to 0$ yields
    \begin{align*}
        \left(
            \bigl\{\cM_{1,\e_k}(\cL_1^{\e_k})(u_j)\bigr\}_{j\in\lb 1,p\rb},
            \bigl\{
                \cM_{2,\e_k}(\cL_2^{\e_k})(v_j)
            \bigr\}_{j\in\lb 1,q\rb}
        \right)
        \xrightarrow[\e_k \to 0]{\;d\;} 
        \left(
            \bigl\{\cL_{R_1}(u_j)\bigr\}_{j\in\lb 1,p\rb},
            \bigl\{
                \cL_{R_2}(v_j)
            \bigr\}_{j\in\lb 1,q\rb}
        \right).
    \end{align*}
    But $\cL_1^{\e_k},\cL_2^{\e_k}$ are independent for every $k$, so the limiting pair is also independent.
    This completes the proof of Property \ref{noise-indep}.
\end{proof}

The rest of this section is devoted to proving Property \ref{noise-join}.
\subsection{Proof of Property \ref{noise-join} (join property)}\label{ss:join-property}

We need to show that for any $R_1,R_2,R_3\in\Rect$ with $R_1\cap R_2=\varnothing$ and $\overline{R_1\cup R_2}=\overline{R_3}$, we have
\begin{align*}
    \cF_{R_1}\vee \cF_{R_2} = \cF_{R_3}
    \qquad\text{up to null sets}.
\end{align*}
We first show that $\cF_{R_1} \subset  \cF_{R_3}$ and that $\cF_{R_2} \subset  \cF_{R_3}$ follows analogously. This implies $\cF_{R_1}\vee \cF_{R_2}\subset \cF_{R_3}$. The proof of the other inclusion will then occupy the rest of the section. 

That $\cF_{R_1} \subset  \cF_{R_3}$ follows immediately from noticing that  \eqref{e:852S} implies that 
    \begin{align}\label{e:852SS}
        \cL_{R_1}(\cL)(s,x;t,y) 
        = 
        \lim_{\e\to 0}
        \sup_{x_1,\dots,x_\ell\in (a_1+\e,b_1-\e)}\sum_{j=0}^{\ell}
        \cL_{R_3}(\cL)(r_j,x_j;r_{j+1},x_{j+1}),
    \end{align}
    where we have used \eqref{restricequal} to replace $\cL$ by $\cL_{R_3}$ and that all points $(r_j,x_j;r_{j+1},x_{j+1})$ as in \eqref{e:852S} are elements of $\Rup(R_1)$ and hence of $\Rup(R_3)$ since $R_1\subset R_3$.

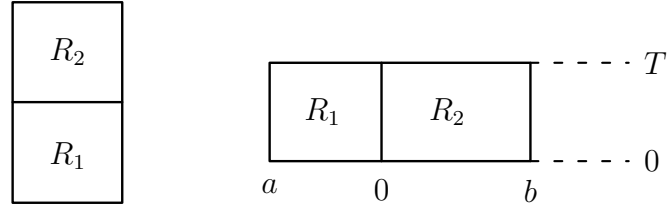
\begin{figure}[hbt]
\centering
\begin{tikzpicture}[
    line cap=round,
    line join=round,
    every node/.style={font=\large}
]


\draw[line width=0.9pt] (0,0) rectangle (1.45,2.65);
\draw[line width=0.9pt] (0,1.33) -- (1.45,1.33);

\node at (0.73,0.62) {$R_1$};
\node at (0.73,2.00) {$R_2$};


\begin{scope}[xshift=3.4cm]

    \draw[line width=0.9pt] (0,0.55) rectangle (3.45,1.85);

    \draw[line width=0.9pt] (1.48,0.55) -- (1.48,1.85);

    \node at (0.70,1.20) {$R_1$};
    \node at (2.35,1.20) {$R_2$};

    \node[below=3pt] at (0,0.55) {$a$};
    \node[below=3pt] at (1.48,0.55) {$0$};
    \node[below=3pt] at (3.45,0.55) {$b$};

    \draw[line width=0.8pt, dash pattern=on 3pt off 6pt]
        (3.45,1.85) -- (4.75,1.85);
    \draw[line width=0.8pt, dash pattern=on 3pt off 6pt]
        (3.45,0.55) -- (4.75,0.55);

    \node[right=2pt] at (4.75,1.85) {$T$};
    \node[right=2pt] at (4.75,0.55) {$0$};

\end{scope}

\end{tikzpicture}
\caption{
        For disjoint open rectangles $R_1,R_2\in\Rect$ such that $\overline{R_1\cup R_2}$ is a closed rectangle, their closures share a side in one of the above two ways (up to relabeling).
        }
    \label{fig:rectangle-cases}
\end{figure}

To show the other inclusion, it suffices to show that for any point $u\in \Rup(R_3)$, the random variable $\cL_{R_3}(u)$ has a $\cF_{R_1}\vee \cF_{R_2}$-measurable version.
Note that the closed rectangles $\overline{R_1}$ and $\overline{R_2}$ must share a side, since $R_1\cap R_2=\varnothing$ and $\overline{R_1\cup R_2}$ is a closed rectangle.
As indicated in \cref{s:iop-black-noise}, we will split our analysis into two cases, depending on whether the shared side is parallel to the temporal axis or the spatial axis.
Namely, write $R_i = (s_i,t_i)\times(x_i,y_i)$ for $i\in\{1,2,3\}$.
By relabeling, we may assume that $s_1\le s_2, t_1\le t_2, x_1\le x_2, y_1\le y_2$.
Then one of the following two cases occurs (see Figure \ref{fig:rectangle-cases}):
\begin{enumerate}[label={\rm(\arabic*)}]
    \item\textup{(Temporally adjacent boxes).}\label{case-temporal-adjacent} $t_1=s_2$ and $(x_1,y_1)=(x_2,y_2)$, and hence $R_3 = (s_1,t_2)\times(x_1,y_2)$;
    \item\textup{(Spatially adjacent boxes).}\label{case-spatial-adjacent} $(s_1,t_1) = (s_2,t_2)$ and $y_1=x_2$, and hence $R_3 = (s_1,t_1) \times (x_1,y_2)$.
\end{enumerate}

Case \ref{case-temporal-adjacent} is straightforward:

\addtocontents{toc}{\SkipTocEntry}
\subsection*{Proof of Property \ref{noise-join}, Case \ref{case-temporal-adjacent}}

    By translation invariance, it suffices to consider the case $(x_1,y_1) = (-a,a)$ for some $a>0$ and $s_1=0$.
    So $R_1 = (0,t_1)\times(-a,a)$ and $R_2=(t_1,t_2)\times(-a,a)$ and $R_3 = (0,t_2)\times(-a,a)$.

    Fix $u=(s,x;t,y)\in\Rup(R_3)$.
    We need to show that $\cL_{R_3}(u)$ has a $\cF_{R_1}\vee\cF_{R_2}$-measurable version.
    If $u\in \Rup(R_1)$ or $u\in\Rup(R_2)$ then this is trivial, since by directedness the set of admissible paths in \eqref{e:restricted-length} is the same as for $R_3$.
    So we only need to consider the situation where the endpoints of $u$ straddle or lie on the boundary edge $R_3\setminus (R_1\cup R_2) = \{t_1\}\times(-a,a)$, that is,  $s\le t_1 \le t$.
    For simplicity we only handle the case where they straddle the boundary ($s<t_1<t$), as the other case will follow from essentially the same argument.

    Fix $\e_0>0$ small enough that $t_1-\e_0 > s$.
    For all $z\in(-a,a)$ and all $\e\in(0,\e_0]$, we have by directedness
    \begin{align*}
        \cL_{R_1}(s,x;t_1-\e, z)
        = \cL_{R_3}(s,x;t_1-\e, z).
    \end{align*}
    By \cref{l:restricted-L-properties}\ref{property-LR-continuous},  the function $\cL_{R_3}(s,x;\smallbullet,\smallbullet)$ is uniformly continuous on $[t_1-\e_0, t_1]\times [-a+\d,a-\d]$ for any $\d>0$,
    hence
    \begin{align*}
        \lim_{\e\to0}
        \sup_{z\in[-a+\d,a-\d]}\left|
            \cL_{R_1}(s,x;t_1-\e,z)
            - \cL_{R_3}(s,x;t_1,z)
        \right| = 0.
    \end{align*}
    A symmetric argument shows that for any $(t,y)\in R_2$,
    \begin{align*}
        \lim_{\e\to0} 
        \sup_{z\in[-a+\d,a-\d]}\left|
            \cL_{R_2}(t_1+\e,z;t,y)
            - \cL_{R_3}(t_1,z;t,y)
        \right| = 0.
    \end{align*}
    From this we get the identity
    \begin{multline*}
        \sup_{z \in [-a+\d,a-\d]}
        \bigl(
            \cL_{R_3}(s,x;t_1,z) + \cL_{R_3}(t_1,z;t,y)
        \bigr)
        \\
        =
        \lim_{\e\to0}
        \sup_{z \in [-a+\d,a-\d]}
        \bigl(
            \cL_{R_1}(s,x;t_1-\e,z) + \cL_{R_2}(t_1+\e,z;t,y)
        \bigr),
    \end{multline*}
    and the RHS clearly has a $\cF_{R_1}\vee\cF_{R_2}$-measurable version for any fixed $\d>0$ (taking the supremum is a measurable operation since $\cL_{R_i}$ is continuous).
    Taking $\d\to 0$ and using \cref{l:restricted-L-properties}\ref{property-LR-convolution} and the fact that $R_3$ is an open set, we get
    \begin{align*}
        \cL_{R_3}(s,x;t,y)
        &=\lim_{\d\to0}
        \sup_{z \in [-a+\d,a-\d]}
        \bigl(
            \cL_{R_3}(s,x;t_1,z) + \cL_{R_3}(t_1,z;t,y)
        \bigr)\\
        &=
        \lim_{\d\to 0}
        \left(
            \lim_{\e\to0}
            \sup_{z \in [-a+\d,a-\d]}
            \bigl(
                \cL_{R_1}(s,x;t_1-\e,z) + \cL_{R_2}(t_1+\e,z;t,y)
            \bigr)
        \right),
    \end{align*}
    and hence $\cL_{R_3}(s,x;t,y)$ has a $\cF_{R_1}\vee\cF_{R_2}$-measurable version.
    This completes the proof.
    \qed \\

\addtocontents{toc}{\SkipTocEntry}
\subsection*{Proof of Property \ref{noise-join}, Case \ref{case-spatial-adjacent}}

The proof of Property \ref{noise-join}, Case \ref{case-spatial-adjacent} is rather long and will occupy the rest of the section. 
While the strategy was outlined in \cref{s:iop-black-noise}, we begin with a roadmap to aid in navigating the different steps.

By translation invariance, it suffices to consider the case $s_1=0$ and $y_1=0$.
In other words, we consider boxes of the form
\begin{align*}
    R_1 &= (0,T)\times(a,0),\\
    R_2 &= (0,T) \times (0,b),\\
    R_3 &= (0,T)\times(a,b)
\end{align*}
for $T>0$ and $a<0<b$, as depicted on the RHS of Figure \ref{fig:rectangle-cases}.
For the rest of this section we fix any $u_*=(s_*,x_*;t_*,y_*)\in\Rup(R_3)$.
It suffices to show that  $\cL_{R_3}(u_*)$ has a $\cF_{R_1}\vee\cF_{R_2}$-measurable version.\\

\begin{itemize}
    \item In \cref{s:noise-conditional-variance-reduction} we reduce to proving a certain conditional variance bound (see \cref{p:cond-var-delta}).
    \item
    Proving the above conditional variance estimate occupies the rest of the section starting from \cref{s:proof-of-CE-estimate}.
    This consists of several steps:
    \begin{itemize}
        \item In \cref{s:discrete-LPP-proxy} we approximate $\cL$ with a discrete last passage percolation (LPP) model. 
        This is not an essential step, but it will make the problem finitary, which will be technically convenient and helpful for presenting the arguments.
        \item In \cref{s:discrete-lpp-reduction} we
        introduce a discrete proxy for the data in $\cF_{R_1}\vee\cF_{R_2}$, complementing the discrete LPP model of \cref{s:discrete-LPP-proxy}.
        After that, we present a detailed outline of the proof of \cref{p:cond-var-delta}, which centers around an Efron--Stein resampling argument for the discrete LPP model resembling that in \cref{ss:idea}.
        \item In Sections \ref{s:mesh-step-1}--\ref{s:mesh-step-2} we complete the proof of \cref{p:cond-var-delta}, following the strategy outlined in \cref{s:discrete-lpp-reduction}.
    \end{itemize}
\end{itemize}

As mentioned, our first aim is to reduce the proof of Property \ref{noise-join}, Case \ref{case-spatial-adjacent}
to  a certain conditional variance estimate.
At a high level, recalling the discussion at the end of \cref{s:iop-black-noise}, 
this will consist of reducing to the toy problem of proving Property \ref{noise-join} for infinite strips (discussed below \eqref{semi1}) by subdividing $R_3$ into short rectangles as in Figure \ref{fig:noise-delta-mesh} below, so as to access the Airy line ensemble estimates required for an Efron--Stein argument like that used to prove \cref{t:main}.
The reduction presented below will only take a couple of pages, but completely explaining how it advances the above goal is difficult at this point and we will not attempt to do so.
We will instead revisit this after the ideas involved in the proof have been laid out in more detail (in \cref{s:discrete-lpp-reduction}).

\subsection{Reduction to a conditional variance estimate}\label{s:noise-conditional-variance-reduction}

\begin{figure}[tbh]
    \centering
    \hspace{2cm}\includegraphics[width=0.8\textwidth]{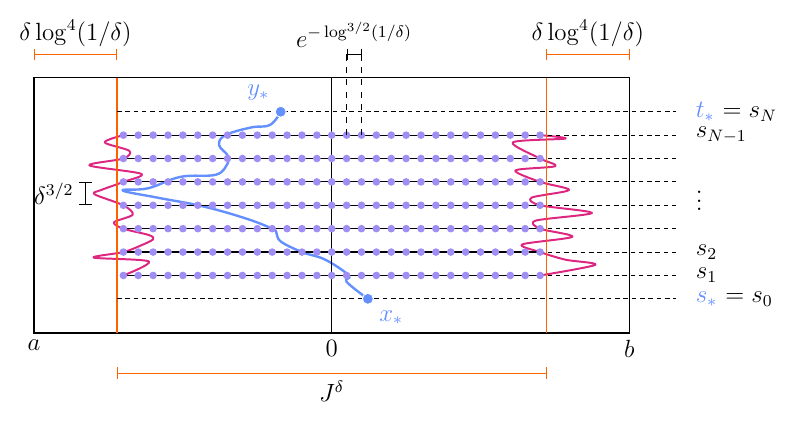}
    \caption{
        The large black rectangles are $R_1=(0,T)\times (a,0)$ and $R_2=(0,T)\times(0,b)$.
        We fix a point $u_*=(s_*,x_*;t_*,y_*)\in\Rup(R_3)$, drawn in blue, where $R_3 = (0,T)\times (a,b)$.
        The orange interval is $J^\d = (a+\d\log^4(1/\d), b-\d\log^4(1/\d))$.
        The horizontal solid black lines are separated from each other by height $\d^{3/2}$; they are of the form $\{s_i\}\times J^\d$ for $i\in\lb 1,N-1\rb$.
        The purple mesh points on the black lines are separated from their horizontally-adjacent neighbors by distance $e^{-\log^{3/2}(1/\d)}$.
        Drawn in blue is a piecewise-geodesic path between $e^{-\log^{3/2}(1/\d)}$-mesh points at times $s_0,s_1,\dots,s_{N}$, which is contained in $J^\d$ at times $s_0,\dots,s_N$.
        The pink paths are also piecewise-geodesic.
        Since $s_{i+1}-s_i=\d^{3/2}$, each geodesic segment has transversal fluctuations $O(\d\log^3(1/\d))$  by \cref{p:GZ-uniform-TF}.
        Since these geodesics start and end inside $J^\d$, which lies at distance $\d\log^4(1/\d) \gg \d\log^3(1/\d)$ from the boundary of the box $R_3$, none of them has a large enough transversal fluctuation to exit $R_3$.
    }
    \label{fig:noise-delta-mesh}
\end{figure}

We fix a small $\d\in(0,\frac12)$ such that $t_* - s_*=N\d^{3/2}$ for some $N\in\N$, and divide the time interval $[s_*,t_*]$ into thin strips of height $\d^{3/2}$.
For $i\in\lb0,N\rb$, denote $s_i\coloneqq  s_* + i\d^{3/2}$ (illustrated in Figure \ref{fig:noise-delta-mesh}).
We will eventually take $\d\to0$, and we assume from now on that $\d<1/2$ and that $b-a > 2\d\log^4(1/\d)$.
This lets us define the open interval
\begin{align}\label{e:def-Jdelta}
    J^\d &\coloneqq  \left(a+\d\log^4(1/\d),\;\; b-\d\log^4(1/\d)\right),
\end{align}
which is drawn at the bottom of Figure \ref{fig:noise-delta-mesh}.

The key technical estimate in this section is the following conditional variance bound.

\begin{prop}\label{p:cond-var-delta}
    There exists $C,c>0$ depending only on $R_1,R_2$, such that for all sufficiently small $\d$,
    \begin{align*}
        \max_{i\in\lb 0, N-1\rb}
        \sup_{z,w\in J^\d}
        \E\left[
            \Var\left(
                \cL(s_i,z ;s_{i+1}, w)\,\middle|\,\cF_{R_1}\vee\cF_{R_2}
            \right)
        \right]
        \le Ce^{-c\log^2(1/\d)}.
    \end{align*}
\end{prop}

Note that \cref{p:cond-var-delta} is different than the bound \eqref{equiv} discussed in Section \ref{s:iop-black-noise}: here there is no $\e$-buffer around the interface between $R_1$ and $R_2$.
An estimate along the latter lines will be a key ingredient in the proof of \cref{p:cond-var-delta}, yielding the result upon taking $\e\to 0$. The estimate above is in fact an approximate version of  Property \ref{noise-join}. The approximation stems from the fact that the observable in \cref{p:cond-var-delta} is $\cL$, not $\cL_{R_3}$.
Property \ref{noise-join} implies that the conditional variance of $\cL_{R_3}$ given $\cF_{R_1}\vee\cF_{R_2}$ is zero; the bound in \cref{p:cond-var-delta} is controlled by the failure probability for coupling $\cL$ and $\cL_{R_3}$.

Postponing the proof of \cref{p:cond-var-delta} to the upcoming \cref{s:proof-of-CE-estimate}, we now complete our analysis of Case \ref{case-spatial-adjacent}, thereby also completing the proof of \cref{t:2d-noise}.

The basic strategy is as follows.
By \cref{p:cond-var-delta}, the information in $\cF_{R_1}\vee\cF_{R_2}$ determines the $\cL$-lengths across the short rectangles $[s_i,s_{i+1}]\times J^\d$, up to a super-polynomially small error.
Also, the shortness of these rectangles lets us couple the $\cL$-lengths to their $\cL_{R_3}$ counterparts up to super-polynomial errors.
Thus since there are only polynomially-many rectangles, the metric composition property for $\cL$ is the same as that for $\cL_{R_3}$ with high probability.
This implies that $\cF_{R_1}\vee\cF_{R_2}$ also determines $\cL_{R_3}$ up to a super-polynomial error, which upon taking $\d\to0$ completes the proof in Case \ref{case-spatial-adjacent}.

We define the spatial mesh (depicted in Figure \ref{fig:noise-delta-mesh} as black points)
\begin{align}
    M^\d &\coloneqq 
    J^\d\cap 
    (e^{-\log^{3/2}(1/\d)}\Z).
\end{align}
For $z\in J^\d$, we denote by $\overline{z}$ the closest point to $z$ in $M^\d$ (we fix an arbitrary deterministic rule for breaking ties).
We define a random function $f^\d:\lb 0,N-1\rb \times J^\d \times J^\d \to\R$ by
\begin{align*}
    f^\d(i,z,w) \coloneqq  \E\left[
        \cL(s_i, \overline{z}; s_{i+1}, \overline{w})\,\middle|\,\cF_{R_1}\vee \cF_{R_2}
    \right].
\end{align*}
{So $f^\d(i,\cdot,\cdot)$ is piecewise constant,} with $f^\d(i,z,w)=f^\d(i,\overline{z},\overline{w})$ for all $z,w\in J^\d$.
By a union bound over pairs of mesh points $\overline{z},\overline{w}\in M^\d$, Chebyshev's inequality, and the conditional variance bound in \cref{p:cond-var-delta}, we have
\begin{multline}\label{e:866}
    \P\left(
        \max_{(i,\overline{z},\overline{w})\in\lb 0,N-1\rb\times M^\d\times M^\d}
        \left|
            \cL(s_i,\overline{z};s_{i+1},\overline{w})
            -
            f^\d(i,\overline{z},\overline{w})
        \right|
        > e^{-\log^{3/2}(1/\d)}
    \right)\\
    \begin{aligned}
        &\le  e^{C\log^{3/2}(1/\d)}  \cdot
    \max_{(i,\overline{z},\overline{w})\in\lb 0,N-1\rb\times M^\d\times M^\d}
    \P\left(
        \left|
            \cL(s_i,\overline{z};s_{i+1},\overline{w})
            -
            f^\d(i,\overline{z},\overline{w})
        \right|
        > e^{-\log^{3/2}(1/\d)}
    \right)\\
    &\le C' e^{-c'\log^2(1/\d)},
    \end{aligned}
\end{multline}
where $C,C',c'>0$ depend only on $R_1,R_2$.

By the H\"older regularity of $\cL$ recorded in \cref{p:dov-holder-general}, since $s_{i+1}-s_i = \d^{3/2}$ for all $i$, there exist $C,c,c'>0$ depending only on $R_1,R_2$ such that 
    \begin{multline}\label{e:holder-for-restricted}
        \P\left(
            \max_{i\in \lb 0,N-1\rb}
            \sup_{\substack{(z,w),(z',w')\in J^\d\times J^\d\\|z-z'|\vee |w-w'| \le e^{-\log^{3/2}(1/\d)}}}
            \left|
                \cL(s_i,z;s_{i+1},w)
                -\cL(s_i,z';s_{i+1},w')
            \right|
            > e^{-\frac{1}{10}\log^{3/2}(1/\d)}
        \right)\\
        \le
        C\d^{-9}\exp\left(
            -c e^{c'\log^{3/2}(1/\d)}
        \right).
    \end{multline}

Combining \eqref{e:866} and \eqref{e:holder-for-restricted}, we find that with probability at least $1-Ce^{-c\log^2(1/\d)}$,
the following estimates hold for all $(i,z,w)\in\lb 0,N-1\rb\times J^\d\times J^\d$:
\begin{align*}
    \left|
        \cL(s_i,z;s_{i+1},w) - f^\d(i,z,w)
    \right|
    &\le
    \left|
        \cL(s_i,z;s_{i+1},w)
        -\cL(s_i,\overline{z};s_{i+1},\overline{w})
    \right|
    +
    \left|
        \cL(s_i,\overline{z};s_{i+1},\overline{w}) - f^\d(i,\overline{z},\overline{w})
    \right|\\
    &\le
    e^{-\frac{1}{10} \log^{3/2}(1/\d)}
    + e^{-\log^{3/2}(1/\d)}\\
    &\le 2e^{-\frac{1}{10} \log^{3/2}(1/\d)}.
\end{align*}
Therefore, defining
\begin{align}\label{e:880}
    \cL^\d(u_*) \coloneqq  
    \sup_{z_1,\dots,z_{N-1} \in J^\d}
    \sum_{i=0}^{N-1}
    \cL(s_i,z_i;s_{i+1},z_{i+1}),
\end{align}
where $z_0\coloneqq x_*$ and $z_{N}\coloneqq y_*$, we get that
\begin{align}\label{e:881}
    \P\left(
        \left|
        \cL^\d(u_*) 
        - 
        \sup_{z_1,\dots,z_{N-1}\in J^\d}
        \sum_{i=0}^{N-1} 
        f^\d(i,z_i, z_{i+1})
        \right|
        \le e^{-\frac{1}{20} \log^{3/2}(1/\d)}
    \right)
    \ge 1 - Ce^{-c\log^2(1/\d)},
\end{align}
where we used the inequality $|\sup_{\bfz} g(\bfz) - \sup_{\bfz} h(\bfz)| \le \sup_{\bfz}|g(\bfz)-h(\bfz)|$, and 
where $C,c>0$ depend only on $R_1,R_2$.
Here we used that $N = (t_*-s_*)\d^{-3/2}$ grows only polynomially in $\d^{-1}$.
Note that by definition, $\sup_{z_1,\dots,z_{N-1}\in J^\d}\sum_{i=0}^{N-1} f^\d(i,z_i, z_{i+1})$ is  $\cF_{R_1}\vee\cF_{R_2}$-measurable.
To help parse this, note that $\cL^\d(u_*)$ is just a maximization over paths of the form depicted in Figure \ref{fig:noise-delta-mesh} in blue.

Finally, we will show that 
\begin{align}\label{e:85222}
    \cL^\d(u_*) \xrightarrow{\d\to 0} \cL_{R_3}(u_*)
    \qquad\text{almost surely.}
\end{align}
Postponing the justification of this to the next paragraph, we now complete the proof.
By the argument used to prove \eqref{e:852}, we have
\begin{align}\label{e:9100}
    \cL_{R_3}(u_*)
    = \lim_{\d\to 0}
    \sup_{z_1,\dots,z_{N-1}\in J^\d} \sum_{i=0}^{N-1} \cL_{R_3}(s_i,z_i;s_{i+1},z_{i+1})
    \qquad\text{almost surely.}
\end{align}
Combining this with \eqref{e:881}--\eqref{e:85222}, we obtain
\begin{align*}
    \sup_{z_1,\dots,z_{N-1}\in J^{\d}}
        \sum_{i=0}^{N-1} 
        f^{\d}(i,z_i, z_{i+1})
    \xrightarrow{\d\to0} 
    \cL_{R_3}(u_*)
    \qquad\text{in probability.}
\end{align*}
The LHS is by definition $\cF_{R_1}\vee\cF_{R_2}$-measurable for any fixed $\d>0$, and therefore the limit $\cL_{R_3}(u_*)$ has a $\cF_{R_1}\vee\cF_{R_2}$-measurable version.
This implies Property \ref{noise-join} in Case \ref{case-spatial-adjacent}.

It remains to prove \eqref{e:85222}.
The argument is essentially the same as that in \eqref{e:852S}.
In light of \eqref{e:9100}, it suffices to prove that almost surely, for all sufficiently small $\d>0$,
\begin{align}\label{e:910}
    \cL^\d(u_*) = \sup_{z_1,\dots,z_{N-1}\in J^\d} \sum_{i=0}^{N-1} \cL_{R_3}(s_i,z_i;s_{i+1},z_{i+1}).
\end{align}
For this, note that by \cref{p:GZ-uniform-TF}, there is a random constant $\CTF>0$ whose law depends only on $R_1,R_2$, 
such that for every $\d$,
\begin{align*}
    \sup_{i\in \lb 0,N-1\rb}
    \sup_{z,w\in J^\d}
    \sup_{r\in[s_i,s_{i+1}]}
    \left|
        \Pi_{(s_i,z),(s_{i+1},w)}(r) - \frac{z(s_{i+1}-r) + w(r-s_i)}{\d^{3/2}}
    \right|
    \le \CTF\, \d \log^3(1/\d)
\end{align*}
where we used that $s_{i+1}-s_i = \d^{3/2}$ for all $i\in\lb 0, N-1\rb$.
Recall that the distance between $\partial J^\d$ and $\partial (a,b)$ is $\d\log^4(1/\d)$, where $\partial$ means boundary.
Therefore, for $\d$ small enough such that $\CTF < \log(1/\d)$, every geodesic $\Pi_{(s_i,z),(s_{i+1},w)}$ for $i\in \lb0,N-1\rb$ and $z,w\in J^\d$ has its graph contained in $R_3 = (0,T)\times (a,b)$ (see the segments comprising the blue and pink paths in Figure \ref{fig:noise-delta-mesh}).
This implies that
\begin{align*}
    \cL(s_i,z; s_{i+1},w) = \cL_{R_3}(s_i,z; s_{i+1},w)
    \qquad\forall i\in\lb 0,N-1\rb, \quad\forall z,w\in J^\d,
\end{align*}
which implies \eqref{e:910} and thus proves \eqref{e:85222}. \qed
\\

We have proved Property \ref{noise-join} in Case \ref{case-spatial-adjacent},  modulo \cref{p:cond-var-delta}.
We now turn towards proving the latter.

\subsection{Proof of Proposition \ref{p:cond-var-delta} (conditional variance at scale $\d$)}\label{s:proof-of-CE-estimate}

\begin{figure}[ht]
    \centering
    \includegraphics[width=\linewidth]{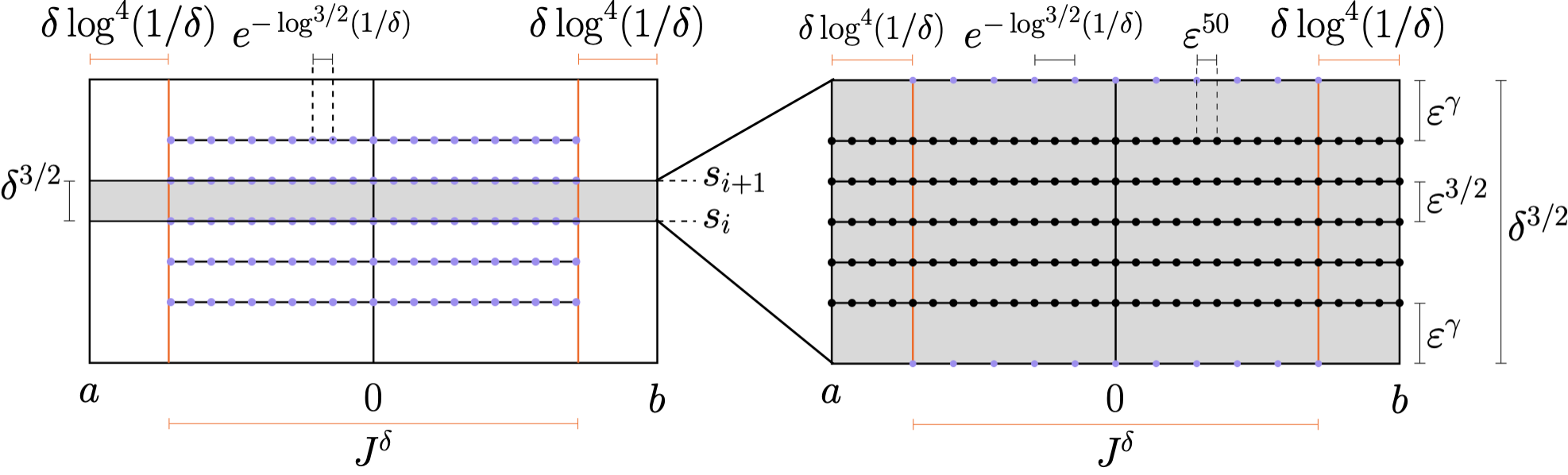}
    \caption{
        The two levels of discretization.\\
        \textbf{Left:} The scale-$\d$ discretization scheme from Figure \ref{fig:noise-delta-mesh}.
        Shaded is a short strip of the form $[s_i,s_{i+1}]\times J^\d$. 
        \cref{p:cond-var-delta} asserts a certain conditional variance bound for $\cL$-lengths across the shaded strip.\\
        \textbf{Right:}
        We zoom in on the shaded strip in the left figure (not drawn to scale).
        To prove \cref{p:cond-var-delta}, we introduce a further space-time discretization indexed by a parameter $\e\ll \d$ (detailed in \cref{s:discrete-LPP-proxy}).
        The black mesh points are spaced from their horizontal neighbors by distance $\e^{50}$.
        The horizontal lines are separated by distance $\e^{3/2}$, except for the first/last line which are separated from the top/bottom of the strip by distance $\e^\gamma\gg \e^{3/2}$.
        Note that $\e^{50}$-mesh points also populate the top and bottom edges of the strip, but these are suppressed for legibility.
    }
    \label{fig:noise-delta-mesh-zoom-in}
\end{figure}

We first make a reduction to simplify notation.
Consider the short rectangles
\begin{align*}
    S_1 &\coloneqq  [0,\d^{3/2}]\times (a,0),\\
    S_2 &\coloneqq  [0,\d^{3/2}]\times(0,b),\\
    S_3 &\coloneqq  [0,\d^{3/2}]\times(a,b).
\end{align*}
By Property \ref{noise-join} in Case \ref{case-temporal-adjacent} (proved at the start of \cref{ss:join-property}), 
for any $j\in \lb 0,N-1\rb$, setting $h=(s_j,0)$, where $s_j\coloneqq  s_* + j\d^{3/2}$, we have $\cF_{R_1} = \cF_{(0,s_j)\times(a,0)}\vee \cF_{S_1+h} \vee \cF_{(s_{j+1},T)\times(a,0)}$ up to null sets (because by continuity, $\cF_{S_1+h}=\cF_{\mathrm{Int}(S_1+h)}$ up to null sets).
By Property \ref{noise-indep} (proved in \cref{s:first-three-properties}), the $\sigma$-algebras $\cF_{(0,s_j)\times(a,0)}, \cF_{S_1+h}, \cF_{(s_{j+1},T)\times(a,0)}$ on the RHS are independent.
Moreover, by directedness, $\cL(s_j,z;s_{j+1},w)$ is independent of $\cF_{(0,s_j)\times(a,0)}\vee \cF_{(s_{j+1},T)\times(a,0)}$ for any $z,w$ (\cref{p:landscape-symmetries}\ref{landscale-independent-increments}).
The analogous statements hold for $R_2,S_2$ by symmetry.
Therefore, by translation invariance, \cref{p:cond-var-delta} is equivalent to:
\begin{align}\label{e:886}
    \sup_{z,w\in J^\d}
    \E\left[
        \Var(\cL(0,z;\d^{3/2}, w)\mid \cF_{S_1}\vee\cF_{S_2})
    \right]
    \le C e^{-c\log^2(1/\d)}.
\end{align}
We turn towards this now.
Recall that the proof strategy outlined in Section \ref{ss:idea} was in the context of a discrete LPP problem.
Working with a discrete model will indeed be convenient for exposition, and will facilitate certain manipulations involving conditional variances.
Therefore, we will prove \eqref{e:886} by approximating $\cL(0,z;\d^{3/2},w)$ with a discrete last passage percolation (LPP) model.
We remark that our particular discretization scheme is ad hoc, and others would have worked equally well for our purposes (see \cref{r:other-discretizations} below).
Finally, before delving into the details, we present in Figure \ref{fig:noise-delta-mesh-zoom-in} a rough depiction of the two layers of discretization involved in our arguments.
Recall that this was already hinted at in Sections \ref{ss:idea} and \ref{s:iop-black-noise}.
The height-$\d^{3/2}$ mesh will be used to essentially couple to the half-infinite strip problem discussed in Section \ref{s:iop-black-noise}, and the upcoming height-$\e^{3/2}$ mesh to carry out the proof strategy involving the Efron--Stein inequality.

\subsection{Discrete LPP proxy}\label{s:discrete-LPP-proxy}

We fix any constant $\gamma\in(0,\frac{1}{6})$ (the exact value is irrelevant to our arguments), which will play the same role as the constant $\gamma$ appearing in \cref{s:es} in the proof of \cref{t:2main}.
{For the rest of \cref{s:black-noise}, all multiplicative constants $C,C',C'',c,c',c''$, etc. may depend on $S_1,S_2,\gamma$, but not on any other parameters unless stated otherwise.
So for example, we write $A=O(B)$ or $A\ls B$ if $|A|\le CB$ for some $C=C(S_1,S_2,\gamma)>0$,
and $A=O_p(B)$ or $A \ls_p B$ if $|A|\le CB$ for some $C=C(S_1,S_2,\gamma,p)>0$, for some parameter $p$.
}

Throughout the coming discussion, we take some $\e>0$ small enough that $\d^{3/2} =  2\e^\gamma + (n-1)\e^{3/2}$ for some $n\in\N$ (note that $n$ is unrelated to Section \ref{s:noise-conditional-variance-reduction}'s $N$).
{We will eventually choose $\e$ to depend on $\delta$ in a super-polynomial way:  $\e \asymp e^{-c\log^{2}(1/\d)}$ for some $c>0$ that will be specified later.
However, for most of the section we will develop estimates with $\e$ as a general parameter assuming only that it is smaller than some large power of $\d$, say $\e\le \d^{100/\gamma}$.
} 

\begin{figure}[tbhp]
    \centering
    \includegraphics[width=\linewidth]{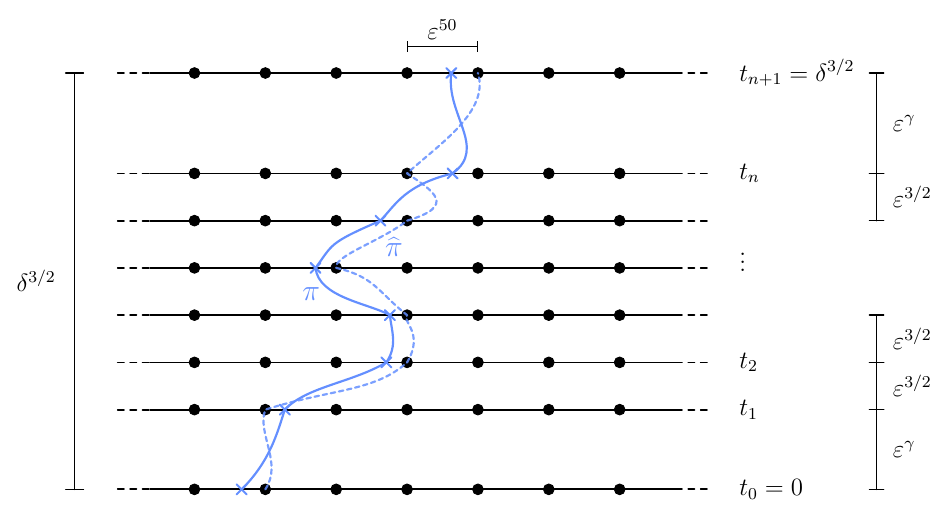}
    \caption{Part of the strip $[0,\d^{3/2}]\times \R$.
    The horizontal lines correspond to the times $t_0,t_1,\dots,t_{n+1}$ defined in \eqref{e:eps-mesh-times}.
    Each horizontal line carries a copy of $\e^{50}\Z$ depicted as black points (figure not to scale).\\
    In solid blue is a directed path $\pi$ from the bottom of the strip to the top.
    Its locations at times $t_0,t_1,\dots,t_{n+1}$ are marked with Xs. 
    These Xs are rounded to $\e^{50}\Z$ as in \cref{def:rounding}.
    The ``rounding'' of $\pi$, denoted $\wh{\pi}$, is the piecewise-geodesic path (dashed blue) interpolating the rounded Xs.
    }
    \label{fig:noise-rounding-eps-times}
\end{figure}

We define the times
\begin{equation}\label{e:eps-mesh-times}
    \begin{split}
        t_0 &\coloneqq  0,\\
        t_k &\coloneqq  \e^\gamma + (k-1)\e^{3/2} \quad\text{for } k\in \lb 1,n\rb,\\
        t_{n+1} &\coloneqq \d^{3/2}. 
    \end{split}
\end{equation}
So $t_1-t_0 = \e^\gamma = t_{n+1}-t_n$, and $t_{k+1}-t_k=\e^{3/2}$ for $k\in\lb 1,n-1\rb$.
These times are depicted on the RHS of Figure \ref{fig:noise-rounding-eps-times}.

For future reference, we note that if $\e\le \d^{100/\gamma}$, then for all $k\in\lb 0,n\rb$, we have
\begin{align}\label{e:985-exponent}
    \frac{t_{k+1} - t_k}{\d^{3/2}} 
    \le 
    (t_{k+1}-t_k)^{0.985},
\end{align}
which can be seen by considering the cases $t_{k+1}-t_k = \e^{\gamma}, \e^{3/2}$ and using that $\gamma < 1/6$.

\begin{defn}[Mesh paths]\label{def:mesh-path}
    We say that a continuous path $\pi:[0,\d^{3/2}]\to \R$ is a \emph{mesh path} if
    \begin{itemize}
        \item $\pi(t_k) \in \e^{50}\Z$ for all $k\in\lb 0,n+1\rb$,
        \item $|\pi(t_{k+1})-\pi(t_k)| < (t_{k+1}-t_k)^{2/3}\log(1/\e)$ for all $k\in\lb 0,n\rb$, and
        \item  
        $\pi|_{[t_k,t_{k+1}]} = \Pi_{(t_k,\pi(t_k)), (t_{k+1}, \pi(t_{k+1}))}$ for all $k\in \lb 0,n\rb$.
    \end{itemize}
    In other words, a mesh path is obtained by selecting a mesh point $v_k\in \e^{50}\Z$ corresponding to time $t_k$ for every $k\in\lb 0,n+1\rb$, connecting $(t_k,v_k)$ and $(t_{k+1},v_{k+1})$ with a directed landscape geodesic, and then concatenating these geodesics across $k$ (as the reader may have already anticipated, the exponent 50 is an arbitrary large number and any other such choice would have sufficed too).
    Note that we impose a hard constraint on a mesh path's local transversal fluctuations.
    \footnote{
        This is just for technical convenience, and should be viewed as a mild constraint:
        we will be interested in the longest mesh path, which is approximately given by rounding the $\cL$-geodesic to $\e^{50}\Z$, and the $\cL$-geodesic satisfies the stated local transversal fluctuation bounds with high probability (Proposition \ref{p:DOV-geodesic-holder}).
    }
    See Figure \ref{fig:noise-rounding-eps-times} for a depiction of the spatial mesh $\e^{50}\Z$ and times $t_0,\dots,t_{n+1}$ involved in this definition.
    Since any mesh path $\pi$ is piecewise geodesic, its length is given by
    \begin{align*}
        \int_0^{\d^{3/2}} d\cL\circ \pi 
        = \sum_{k=0}^n \cL(t_k,\pi(t_k); t_{k+1},\pi(t_{k+1})).
    \end{align*}
\end{defn}

\begin{defn}[Rounding points and paths]\label{def:rounding}
    For $z\in \R$, we denote by $\wh{z}$ the unique point in $\e^{50}\Z$ satisfying $|\wh{z}| \le |z|$ and $|z-\wh{z}|<\e^{50}$.
    This convention ensures that if $z\in J^{\d}$, then $\wh{z}\in J^{\d}$.

    For any $0\le k_1 < k_2 \le n+1$ and any continuous path $\pi : [t_{k_1}, t_{k_2}]\to \R$, we define a path $\wh{\pi} : [t_{k_1}, t_{k_2}] \to \R$ by setting
    \begin{align*}
        \wh{\pi}\big|_{[t_{k}, t_{k+1}]} \coloneqq  \Pi_{(t_k, \wh{\pi(t_k)}), (t_{k+1}, \wh{\pi(t_{k+1})})}
        \qquad\text{for } k\in \lb k_1,k_2-1\rb.
    \end{align*}
    In words, we round $\pi(t_{k_1}), \pi(t_{k_1+1}),\dots,\pi(t_{k_2})$ to $\e^{50}\Z$, then interpolate between these points with $\cL$-geodesics (see also Figure \ref{fig:noise-rounding-eps-times}).
    Note that if $\pi$ has geodesic-like local transversal fluctuations, then $\wh{\pi}$ is a mesh path (\cref{def:mesh-path}).
    Also, if $\pi$ is itself a mesh path, then $\pi=\wh{\pi}$.

    In our applications,  $\pi$ will itself be piecewise geodesic, so we think of $\wh{\pi}$ as obtained by \emph{rounding} $\pi$ to $\e^{50}\Z$.
\end{defn}

\begin{defn}[Discrete LPP model]\label{def:discrete-LPP}
    For $z,w\in\R$, we define 
    \begin{align*}
        \wh{\cL}^\e(0,z;\d^{3/2},w) \coloneqq  
        \max_{\pi:\wh{z}\to\wh{w}} \int_0^{\d^{3/2}} d\cL\circ \pi,
    \end{align*}
    where the maximum is over all mesh paths $\pi$ with $\pi(0)=\wh{z}$ and $\pi(\d^{3/2})=\wh{w}$. 
    Note that the maximum is indeed attained, as there are only finitely many mesh paths connecting $\wh{z}$ to $\wh{w}$.
    Writing $u=(0,z;\d^{3/2},w)$, we define the \emph{mesh geodesic} $\Gamma_u$ to be the mesh path attaining the above maximum (with ties broken in any deterministic measurable way).
\end{defn}

\begin{rem}[Other discretization schemes]\label{r:other-discretizations}
    The spatial mesh of width $\e^{50}$ used in Definitions \ref{def:mesh-path}--\ref{def:discrete-LPP} was chosen somewhat arbitrarily; any sufficiently large power of $\e$ will suffice for our purposes, and we did not attempt to optimize this.
    In fact, it is not hard to extend our arguments to an LPP proxy given by replacing the spatial mesh $\e^{50}\Z$ with a countable dense subset of $\R$.
    We opted for the more finitary description above to aid exposition. 
\end{rem}

We now prove three regularity estimates that will allow us to compare $\cL$ and $\wh{\cL}^\e$, and to bound moments of $\wh{\cL}^\e$.
While the conclusions are intuitively clear, the proofs are somewhat technical and the reader may find it helpful to skip ahead to \cref{s:discrete-lpp-reduction} on first reading and return here as needed.
The first lemma, a consequence of the H\"older regularity of the directed landscape, says that rounding a geodesic $\Pi_{(t_{k_1},z), (t_{k_2},w)}$ to $\e^{50}\Z$ at times $t_{k_1},t_{k_1+1},\dots,t_{k_2}$ changes its length by $O(\e^{10})$ 
(in reality this is $O(\e^{25 - 3/2})$ up to polylog factors).

\begin{lemm}[Rounding geodesics to $\e^{50}\Z$]\label{l:tail-and-comparison-mesh}
    Let $t_0,t_1,\dots,t_{n+1}$ be as in \eqref{e:eps-mesh-times}.
    For all sufficiently small $\e>0$ and all $m\ge 0$, we have
    \begin{align*}
        \P\left(
            \max_{0\le k_1 < k_2 \le n+1}\,
            \sup_{z,w\in J^{\d}}
            \left|\cL(t_{k_1},z; t_{k_2},w) 
            - 
            \int_{t_{k_1}}^{t_{k_2}} d\cL \circ \wh{\Pi}_{(t_{k_1},z),(t_{k_2},w)}
            \right|
            > m \e^{10}
        \right)
        \le Ce^{-cm^{6/7}},
    \end{align*}
    where $\wh{\Pi}_{(t_{k_1},z),(t_{k_2},w)}$ is the path obtained by rounding the geodesic $\Pi_{(t_{k_1},z),(t_{k_2},w)}$ to $\e^{50}\Z$ as in \cref{def:rounding}.
\end{lemm}
\begin{proof}
    We first briefly outline the argument.
    We will use the spatial H\"older continuity of $\cL$  to assert that rounding $\Pi=\Pi_{(t_{k_1},z),(t_{k_2},w)}$ to $\e^{50}\Z$ cannot change its length by more than $O(\e^{10})$.
    However, $\cL$ is only locally H\"older continuous, with a random H\"older norm on each compact set.
    Since $\Pi$ exits any fixed compact set with positive probability, we must track how $\cL$'s local H\"older norm grows as the compact set increases, and balance this with the rarity of $\Pi$ having a large transversal fluctuation.

    For $L\ge 1$, let $\CHolder^{L} > 0$ be the corresponding random constant from \cref{p:dov-holder-general} with $\epsilon \coloneqq  \e^{3/2}$.
    This means that $\CHolder^{L}$ controls the spatial modulus of continuity of $\cL$ restricted to points $(s,x;t,y)\in\Rup\cap [-2L,2L]^4$ with time separation $t-s\ge \e^{3/2}$.
    By \cref{p:dov-holder-general}, we have $\P(\CHolder^{L} > m) \le C L^{10} \e^{-9} e^{-cm^{3/2}}$ for all $m\ge0$, where $C,c>0$ are universal constants.
    Moreover, we can arrange that $\CHolder^{L} \ge \CHolder^{L'}$ for $L\ge L'$. 

    Fix $0 \le k_1 < k_2 \le n+1$. 
    Since $t_{k_2}-t_{k_1} \ge \e^{3/2}$, it follows by \cref{p:dov-holder-general} that
    \begin{multline*}
        \sup_{\substack{z,z',w,w'\in [-2L,2L]\\|z-z'|, |w-w'| \le \e^{50}}}
        \left|\cL(t_{k_1}, z; t_{k_2}, w) - \cL(t_{k_1}, z'; t_{k_2}, w') + \frac{(z-w)^2 - (z'-w')^2}{t_{k_2}-t_{k_1}}\right|\\
        \begin{aligned}
            &\le C\, 
            \CHolder^{L}\, \e^{25} \log^{1/2}\left(\frac{8L}{\e^{50}}\right)\\
            &\le C'\,\CHolder^{L}\,L^{1/2}\e^{24},
        \end{aligned}
    \end{multline*}
    where $C,C'>0$ are universal constants.
    Next, considering the quadratic term on the LHS, since $t_{k_2}-t_{k_1}\ge \e^{3/2}$, we have
    \begin{align*}
        \sup_{\substack{z,z',w,w'\in [-2L,2L]\\|z-z'|, |w-w'| \le \e^{50}}}
        \left|
            \frac{(z-w)^2 - (z'-w')^2}{t_{k_2}-t_{k_1}}
        \right|
        \le \frac{CL\e^{50}}{\e^{3/2}}
        \le C L \e^{40},
    \end{align*}
    where $C>0$ is a universal constant.
    Applying this in the previous display, we get
    \begin{align}\label{e:940}
        \sup_{\substack{z,z',w,w'\in [-2L,2L]\\|z-z'|, |w-w'| \le \e^{50}}}
        \left|\cL(t_{k_1}, z; t_{k_2}, w) - \cL(t_{k_1}, z'; t_{k_2}, w')\right|
        \le C''(\CHolder^{L}+1)L\,\e^{24}.
    \end{align}

    Let $\CTF$ be the random constant controlling geodesic transversal fluctuations defined in \cref{p:GZ-uniform-TF} with $L=b-a+1$.
    By \cref{p:GZ-uniform-TF}, there is a universal constant $C>0$ such that for every $0\le k_1<k_2 \le n+1$, every $z,w\in (a,b)$, and every geodesic $\Pi_{(t_{k_1},z),(t_{k_2},w)}$, we have
    \begin{align}\label{e:941}
        \sup_{r\in [t_{k_1}, t_{k_2}]}\left|\Pi_{(t_{k_1},z),(t_{k_2},w)}(r)\right| \le C\,\CTF + b-a 
        =: \fL.
    \end{align}
    By \cref{p:GZ-uniform-TF}, we have $\P(\fL > m) \le Ce^{-cm^2}$ for some $C,c>0$ depending only on $S_1,S_2$.

    Fix $0 \le k_1 < k_2 \le n+1$ and $z,w\in (a,b)$.
    Write $u\coloneqq (t_{k_1}, z;t_{k_2},w)$ and let $\Pi_u$ be the geodesic corresponding to $\cL(u)$.
    By the definition of geodesics (see \eqref{e:def-geo-2}), the triangle inequality, and \eqref{e:940}--\eqref{e:941}, we have
    \begin{align*}
        \left|
            \cL(u) - 
            \int_{t_{k_1}}^{t_{k_2}}d\cL\circ \wh{\Pi}_u
        \right|
            &\le
            \sum_{k=k_1}^{k_2-1}
            \left|
                \cL(t_k, \Pi_u(t_k); t_{k+1}, \Pi_u(t_{k+1}))
                - \cL(t_k, \wh{\Pi_u(t_k)}; t_{k+1}, \wh{\Pi_u(t_{k+1})})
            \right|\\
            &\le
            C''(k_2-k_1)(\CHolder^{\fL} + 1)\fL\,\e^{24}\\
            &\le C''(\CHolder^{\fL} + 1)\fL\,\e^{20},
    \end{align*}
    since $k_2-k_1 \le n+1 \le \e^{-3/2}$.
    Note that the above estimate holds for all $z,w\in (a,b)$ and $0\le k_1<k_2\le n+1$ simultaneously.
    Fixing $L\ge 1$ (to be specified momentarily), we get by a union bound
    \begin{multline*}
        \P\left(
            \max_{0\le k_1<k_2\le n+1}
            \sup_u 
            \left|
            \cL(u) - 
            \int_{t_{k_1}}^{t_{k_2}}d\cL\circ \wh{\Pi}_u
            \right|
            > m\e^{10}
        \right)\\
        \le
        \P\left(
            \max_{0\le k_1<k_2\le n+1}
            \sup_u
            \left|
                \cL(u) - 
                \int_{t_{k_1}}^{t_{k_2}}d\cL\circ \wh{\Pi}_u
            \right|
            > m\e^{10},
            \quad
            \fL \le L
        \right)
        +
        \P(\fL > L)\\
        \begin{aligned}
            &\le
            \P\Bigl(C''(\CHolder^{L} + 1)L\,\e^{20} > m\e^{10}\Bigr)
            +
            \P(\fL > L)\\
            &\le C L^{10} \e^{-9} \exp\left(
                - c\left(\frac{m}{L\e^{10}}\right)^{3/2}
            \right)
            + C' e^{-c' L^2}.
        \end{aligned}
    \end{multline*}
    Assuming $m\ge 1$, straightforward algebra shows that the exponents match when $L = (m\e^{-10})^{3/7}$, and plugging this in yields an upper bound of
    \begin{align*}
        C'' \,
        \frac{m^{30/7}}{\e^{363/7}}
        \exp\left(
            -c'' \frac{m^{6/7}}{\e^{60/7}}
        \right).
    \end{align*}
    Finally, we absorb the $m^{30/7}/\e^{363/7}$ prefactor into the constant $c''$, drop the $\e^{60/7}$ from the exponent's denominator to get the simpler-looking upper bound stated in \cref{l:tail-and-comparison-mesh}, and increase $C''$ to extend to $m<1$.
\end{proof}

The next lemma compares $\cL$ and $\wh{\cL}^\e$.
Let us briefly explain its statement.
Recall that the mesh paths in the definition of $\wh{\cL}^\e$ (see \cref{def:discrete-LPP}) are subject to a hard local transversal fluctuation constraint (see \cref{def:mesh-path}).
On the (high-probability) event that the $\cL$-geodesic has typical local transversal fluctuation behavior, \cref{l:tail-and-comparison-mesh} ensures that rounding the $\cL$-geodesic to $\e^{50}$ yields a mesh path with comparable length.
This comparison breaks down when the $\cL$-geodesic has a large local transversal fluctuation; this is the source of the probability bound stated below. 

\begin{lemm}[Comparing $\cL$ to $\wh{\cL}^\e$]\label{l:moments-tails-wh-L}
    Let $\wh{\cL}^\e$ be as in \cref{def:discrete-LPP}.
    There exist $C,c>0$ such that for all $\e\le \d^{100/\gamma}$, we have
    \begin{align*}
        \inf_{z,w\in J^\d}
        \P\left(
            \left|\cL(0,z;\d^{3/2},w)
            -
            \wh{\cL}^\e(0,z;\d^{3/2},w)\right|
            \le \e^9
        \right)
        &\ge 1-Ce^{-c\log^2(1/\e)}.
    \end{align*}
\end{lemm}
\begin{proof}
    Fix $u=(0,z;\d^{3/2},w)$ with $z,w\in J^\d$, and write $\wh{u}\coloneqq (0,\wh{z};\d^{3/2},\wh{w})$ for its rounding to $\e^{50}\Z$.
    Let $\CTFloc=\CTFloc(\wh{u})>0$ be the random constant controlling the modulus of continuity of the geodesic $\Pi_{\wh{u}}$, as defined in \cref{p:DOV-geodesic-holder} with $S = \frac{b-a+1}{\d^{3/2}}$.

    Let $t_0,\dots,t_{n+1}$ be as in \eqref{e:eps-mesh-times}.
    We have the local transversal fluctuation bound
    \begin{align*}
        \max_{k\in\lb 0,n\rb}
        \left|
            \Pi_{\wh{u}}(t_{k+1}) - \Pi_{\wh{u}}(t_k)
        \right|
        \overset{\eqref{e:geodesic-holder}}&{\le} \CTFloc(\wh{u})\, (t_{k+1}-t_k)^{2/3}\log^{1/3}\left(\frac{1}{t_{k+1}-t_k}\right)
        + \frac{b-a+1}{\d^{3/2}}(t_{k+1}-t_k)\\
        \overset{\eqref{e:985-exponent}}&{\le}
        \CTFloc(\wh{u})\, (t_{k+1}-t_k)^{2/3}\log^{1/3}\left(\frac{1}{\e^{3/2}}\right)
        + (b-a+1)(t_{k+1}-t_k)^{0.985},
    \end{align*}
    where in the second line we used that $t_{k+1}-t_k \ge \e^{3/2}$ for every $k\in\lb 0,n\rb$.
    Since $0.985 > 2/3$ and $t_{k+1}-t_k \le \e^\gamma$, we can choose a deterministic constant $c_0=c_0(S_1,S_2,\gamma)>0$ such that on the event
    \begin{align*}
        \sA_{u} \coloneqq  \{\CTFloc(\wh{u}) < c_0 \log^{2/3}(1/\e)\},
    \end{align*}
    we have
    \begin{align}\label{e:929}
        \max_{k\in\lb 0,n\rb}
        \left|
            \Pi_{\wh{u}}(t_{k+1}) - \Pi_{\wh{u}}(t_k)
        \right|
        \le \frac12 (t_{k+1}-t_k)^{2/3} \log(1/\e).
    \end{align}
    And by \cref{p:DOV-geodesic-holder}, we have
    \begin{align}\label{e:930}
        \sup_u \P(\neg\sA_{u}) \le  Ce^{-c\log^2(1/\e)}.
    \end{align}

    Continuing from \eqref{e:929}, on the event $\sA_u$, rounding $\Pi_{\wh{u}}$ to $\e^{50}\Z$ as in \cref{def:rounding} yields a path $\wh{\Pi}_{\wh{u}}$ satisfying
    \begin{align*}
        \max_{k\in\lb 0,n\rb}\left|\wh{\Pi}_{\wh{u}}(t_{k+1})-\wh{\Pi}_{\wh{u}}(t_k)\right| 
        < (t_{k+1}-t_k)^{2/3} \log(1/\e),
    \end{align*}
    and hence $\wh{\Pi}_{\wh{u}}$ is a mesh path from $\wh{z}$ to $\wh{w}$ (see \cref{def:mesh-path}).
    Since $\wh{\cL}^\e(u)$ is the maximum length among \emph{all} such mesh paths (see \cref{def:discrete-LPP}), we get
    \begin{align}\label{e:931}
        \wh{\cL}^\e(u) 
        &\ge \int_0^{\d^{3/2}} d\cL\circ \wh{\Pi}_{\wh{u}}
        \qquad\text{on the event $\sA_u$.}
    \end{align}
    By \eqref{e:DL-as-LPP}, we have the deterministic upper bound $\wh{\cL}^\e(u)\le \cL(\wh{u})$.
    Therefore, by \eqref{e:931},
    \begin{align}\label{e:932}
        \sup_u
        \P\left(
            \left|\cL(u) - \wh{\cL}^\e(u)\right| > \e^9
        \right)
            &\le
            \sup_u\P\left(
                \cL(\wh{u}) - \wh{\cL}^\e(u) > \frac{\e^9}{2}
            \right)
            +
            \sup_u\P\left(
                \Bigl|\cL(u) - \cL(\wh{u})\Bigr| > \frac{\e^9}{2}
            \right).
    \end{align}
    We bound the two terms on the RHS.
    For the first term, applying \eqref{e:931} followed by \cref{l:tail-and-comparison-mesh} and \eqref{e:930}, we get 
    \begin{align*}
        \sup_u\P\left(
            \cL(\wh{u}) - \wh{\cL}^\e(u) > \frac{\e^9}{2}
        \right)
            &\le
            \sup_u \P\left(
                    \cL(\wh{u}) -
                    \int_{0}^{\d^{3/2}}d\cL \circ \wh{\Pi}_{\wh{u}}
                     > \frac{\e^9}{2}
                \right)
            + \sup_u\P(\neg\sA_u)\\
                &\le Ce^{-c\e^{-6/7}} + C'e^{-c'\log^2(1/\e)}\\
                &\le C''e^{-c''\log^2(1/\e)}.
    \end{align*}
    For the second term on the RHS of \eqref{e:932}, we use the spatial H\"older regularity of $\cL$ (\cref{p:dov-holder-general}) to assert that $\left|\cL(u) - \cL(\wh{u})\right| \ls \e^{25}\log^{1/2}(1/\e)$ with high probability uniformly in $u$.
    More precisely, by \cref{p:dov-holder-general}, we have
    \begin{align*}
        \sup_u\P\left(
                \Bigl|\cL(u) - \cL(\wh{u})\Bigr| > \frac{\e^9}{2}
            \right)
        &\le C \d^{-9} e^{-c\e^{-21}}\\
        &\le C e^{-c'\e^{-20}},
    \end{align*}
    where we used that $(25 - 9)\cdot \frac{3}{2} > 21$ and that 
    $\d^{-9} \le \e^{-9\gamma/100} < \e^{-3/200}$.
    This finishes the proof.
\end{proof}

The next lemma bounds moments of $\wh{\cL}^\e$.
Since $\wh{\cL}^\e$ is upper bounded by $\cL$, the main task is to control its lower tail.
{To do this we consider two cases.
If the $\cL$-geodesic has typical local transversal fluctuations, we can round it to $\e^{50}\Z$ as in \cref{def:rounding} to produce a mesh path without significantly changing its passage time.
On the rare event where this fails, we instead construct a mesh path by concatenating many short geodesics, and lower bound their passage time.}

\begin{lemm}[Moment bounds for discretized LPP model]\label{l:moments-wh-L}
    Let $\wh{\cL}^\e$ be as in Definition \ref{def:discrete-LPP}.
    For all $p\ge 1$ and all $\e \le \delta^{100/\gamma}$, we have      
    \begin{align*}
        \sup_{z,w\in J^{\d}}
        \E\left[
            \left|\wh{\cL}^\e(0,z;\d^{3/2},w)
                + \frac{(z-w)^2}{\d^{3/2}}
            \right|^p
            \right]
        \ls_{p}
        \d^{p/2} \log^{4p/3}(1/\d).
    \end{align*}
\end{lemm}
\begin{proof}
    By the inequality $|x+y|^p \le 2^{p-1}(|x|^p + |y|^p)$, it suffices to separately bound moments of the positive and negative parts.

    For the positive part, by \eqref{e:DL-as-LPP}
    and \cref{p:landscape-symmetries}\ref{landscape-scaling}--\ref{landscape-shear},
    we have 
    \begin{align*}
        \wh{\cL}^\e(0,z;\d^{3/2},w)
        \overset{\eqref{e:DL-as-LPP}}{\le}
        \cL(0,\wh{z};\d^{3/2},\wh{w})
        &\law
        \d^{1/2}\cL(0,0;1,0) - \frac{(\wh{w}-\wh{z})^2}{\d^{3/2}}.
    \end{align*}
    Adding the parabolic term, we get
    \begin{align*}
        \wh{\cL}^\e(0,z;\d^{3/2},w)
        + \frac{(z-w)^2}{\d^{3/2}}
        &\preceq_{\mrm{sd}}
        \d^{1/2}\cL(0,0;1,0)
        + O\left(\frac{\e^{50}}{\d^{3/2}}\right)\\
        &\preceq_{\mrm{sd}}
        \d^{1/2}\cL(0,0;1,0)
        + O(\d^{100}),
    \end{align*}
    since $\e^{50} \le \d^{5000/\gamma} < \d^{30000}$.
    In particular,
    \begin{align}\label{e:positive-part}
        \sup_{z,w\in J^\d}
        \E\left[
            \max\left\{\wh{\cL}^\e(0,z;\d^{3/2},w) + \frac{(z-w)^2}{\d^{3/2}}, \;0\right\}^p
        \right]
        &\ls_{p} 
        \d^{p/2}\,
        \E\left[\bigl|\cL(0,0;1,0)\bigr|^p\right]
        + \d^{100p}\nonumber\\
        &\ls_{p} \d^{p/2},
    \end{align}
    because all moments of $\cL(0,0;1,0)$ are finite (it follows the Tracy--Widom GUE distribution, see the proof of \cite[Lemma 10.4]{dov}).

    We turn next to the negative part.
    We will explicitly construct a mesh path whose length satisfies an appropriate lower bound.
    As mentioned, a construction to this effect was outlined earlier in the proof of  \cref{l:restricted-L-properties}\ref{property-LR-finite}.

    Write $u\coloneqq (0,z;\d^{3/2},w)$ and define the event
    \begin{align*}
        \sA_u \coloneqq  \left\{
            \left|\cL(u) - \wh{\cL}^\e(u)\right| \le \e^{9}
        \right\}.
    \end{align*}
    Since $\e\le \d^{100/\gamma}$, we have $\P(\sA_u)\ge 1 - Ce^{-c\log^2(1/\e)}$ uniformly in $u$.
    Let $\Cptwise$ be the random constant controlling the pointwise values of $\cL$ 
    defined in \cref{p:dov-pointwise} with $L=b-a+1$.
    On the event $\sA_u$,
    we have by \cref{p:dov-pointwise} that
    \begin{align*}
        \wh{\cL}^\e(u) + \frac{(z-w)^2}{\d^{3/2}}
        &\ge \cL(u) + \frac{(z-w)^2}{\d^{3/2}} - \e^9\\
        &\ge -\Cptwise\,\d^{1/2}\log^{4/3}\left(\frac{8(b-a+1)}{\d^{3/2}}\right)
        - \e^9\\
        &\gs -(\Cptwise+1)\d^{1/2}\log^{4/3}(1/\d)
    \end{align*}
    and hence 
    \begin{align}\label{e:968}
        \E\left[
            \max\left\{-\wh{\cL}^\e(u) - \frac{(z-w)^2}{\d^{3/2}},\; 0\right\}^p
            \1_{\sA_u}
        \right]
        &\ls_{p}
        \d^{p/2}\log^{4p/3}(1/\d)
        \E\left[
            (\Cptwise+1)^p
        \right]\nonumber\\
        &\ls_{p} \d^{p/2}\log^{4p/3}(1/\d),
    \end{align}
    by the tail bound in \cref{p:dov-pointwise}.

    It remains to lower bound $\wh{\cL}^\e(u)$ on $\neg\sA_u$.
    Set $\mu \coloneqq  \frac{w-z}{\d^{3/2}}$, and set
    \begin{align*}
        z_k &\coloneqq  \wh{z+t_k \mu},
        \qquad k\in\lb 0,n+1\rb,
    \end{align*}
    where the notation is from \cref{def:rounding}.
    We have the following oscillation estimate:
    \begin{equation}\label{e:965}
        \begin{split}
            \left|z_{k+1}-z_k\right| 
            &\le 
            |\mu|(t_{k+1}-t_k) + 2\e^{50}\\
            &\le 
            \frac{b-a+2}{\d^{3/2}}(t_{k+1}-t_k)\\
            \overset{\eqref{e:985-exponent}}&{\le}
            (b-a+2)(t_{k+1}-t_k)^{0.985}.
        \end{split}
    \end{equation}
    By the first line of \eqref{e:965}, we have
    \begin{align}\label{e:difference-parabolas}
        \left|
            \sum_{k=0}^n
            \frac{(z_{k+1}-z_k)^2}{t_{k+1}-t_k}
            - \frac{(z-w)^2}{\d^{3/2}}
        \right|
        &\le 2(n+1)\e^{50}
        \le \e^{40},
    \end{align}
    since $(n+1)\le \e^{-3/2}$.
    Further, by the third line of \eqref{e:965}, for small enough $\d$ (hence small enough $\e$),
    concatenating the geodesics $\Pi_{(t_k,z_k),(t_{k+1},z_{k+1})}$ across $k\in\lb 0,n\rb$ yields a mesh path $\pi:[0,\d^{3/2}]\to\R$ with $\pi(0)=\wh{z}$ and $\pi(\d^{3/2})=\wh{w}$.
    Using \cref{p:dov-pointwise} to lower bound the length of each geodesic segment, we get the following lower bound for the length of $\pi$ (and hence for $\wh{\cL}^\e(u)$):
    \begin{align*}
        \wh{\cL}^\e(u) + \frac{(z-w)^2}{\d^{3/2}}
        &\ge
        \int_0^{\d^{3/2}} d\cL\circ \pi
        + \frac{(z-w)^2}{\d^{3/2}}\\
        &=
        \sum_{k=0}^n \cL(t_k,z_k; t_{k+1},z_{k+1})
        + \frac{(z-w)^2}{\d^{3/2}}
        \\
        &\ge \sum_{k=0}^n \left[
            - \frac{(z_{k+1}-z_k)^2}{t_{k+1}-t_k}
            - \Cptwise\,(t_{k+1}-t_k)^{1/3}\log^{4/3}\left(\frac{8(b-a+1)}{t_{k+1}-t_k}\right)
        \right]
        + \frac{(z-w)^2}{\d^{3/2}}
        \\
        \overset{\eqref{e:difference-parabolas}}&{\ge}
        \sum_{k=0}^n \left[
            - \Cptwise\,\e^{\gamma/3}\log^{4/3}\left(\frac{8(b-a+1)}{\e^{3/2}}\right)
        \right]
        - \e^{40}
        \\
        &=
        -(n+1)\,\Cptwise\,\e^{\gamma/3}\log^{4/3}\left(\frac{8(b-a+1)}{\e^{3/2}}\right)
        - \e^{40}
        \\
        &\gs
        -(\Cptwise+1)\,\frac{1}{\e^{3/2}},
    \end{align*}
    where we used that $n+1 \le \e^{-3/2}$ and that $\e \le \d$.
    From this and \cref{p:dov-pointwise} we get
    \begin{align*}
        \E\left[
            \max\left\{-\wh{\cL}^\e(u)-\frac{(z-w)^2}{\d^{3/2}},\; 0\right\}^p
        \right]
        &\ls_{p} \frac{1}{\e^{3p/2}}\E[(\Cptwise + 1)^p]\\
        &\ls_{p} \frac{1}{\e^{3p/2}}.
    \end{align*}
    Therefore, by Cauchy--Schwarz and the aforementioned probability bound for $\sA_u$, we have
    \begin{align*}
        \E\left[
            \max\left\{-\wh{\cL}^\e(u)-\frac{(z-w)^2}{\d^{3/2}},\; 0\right\}^p
            \1_{\neg \sA_u}
        \right]
        &\le \E\left[
            \max\left\{-\wh{\cL}^\e(u)-\frac{(z-w)^2}{\d^{3/2}},\; 0\right\}^{2p}
        \right]^{1/2}
        \P(\neg \sA_u)^{1/2}\\
        &\ls_{p} \frac{1}{\e^{3p/2}}\cdot e^{-c\log^2(1/\e)}\\
        &\le e^{-c'\log^2(1/\e)}\\
        &\le e^{-c'\log^2(1/\d)}.
    \end{align*}
    for some $c'=c'(p)>0$.
    Combining this with \eqref{e:968} gives
    \begin{align*}
        \E\left[
            \max\left\{-\wh{\cL}^\e(u)-\frac{(z-w)^2}{\d^{3/2}},\;0\right\}^p
        \right]
        &\ls_{p}
        \d^{p/2}\log^{4p/3}(1/\d) + e^{-c'\log^2(1/\d)}\\
        &\ls_{p} \d^{p/2}\log^{4p/3}(1/\d).
    \end{align*}
    \cref{l:moments-wh-L} now follows by combining the above with the bound for the positive part \eqref{e:positive-part}.
\end{proof}

\subsection{Using the discrete LPP proxy to prove \cref{p:cond-var-delta}}\label{s:discrete-lpp-reduction}

We assume from now on that $\e$ is small enough that $[-\e\log^4(1/\e),\e\log^4(1/\e)] \subset J^{\d}$, where $J^\d$ is defined in \eqref{e:def-Jdelta}.  
Recall from above \eqref{e:eps-mesh-times} that $\e$ satisfies $\d^{3/2} = 2\e^\gamma + (n-1)\e^{3/2}$ for some $n\in\N$.
At the cost of being repetitive, we recall from the discussion at the beginning of Section \ref{s:discrete-LPP-proxy} that we will eventually take $\e$ to be $e^{-c\log^{2}(1/\d)}$ for some $c>0$ (to be specified later), and that for the moment we will develop estimates with $\e$ as a general parameter satisfying $\e \le \delta^{100/\gamma}$.

Write $\d' \coloneqq  \d/2$.
We will assume from now on that $\d$ is small enough that $J^{\d}\subset J^{\d'}$, and such that
$\mathrm{dist}(J^{\d}, \partial J^{\d'}) > \frac{1}{10} \d\log^4(1/\d)$ (equivalently, $\d\log^4(1/\d) > \frac{10}{9}\d'\log^4(1/\d')$),
where $\mathrm{dist}(\cdot,\cdot)$ is the usual distance between subsets of $\R$, and where $J^\d$ is defined in \eqref{e:def-Jdelta}.

We partition $\e^{50}\Z$ into four sets:
\begin{equation}\label{e:mesh-sets}
    \begin{split}
        M^{\near}_{S_1}
        &\coloneqq  (\e^{50}\Z)\cap (a+\d'\log^4(1/\d'), -\e\log^4(1/\e)),\\
        M^{\near}_{S_2}
        &\coloneqq  (\e^{50}\Z)\cap (\e\log^4(1/\e), b - \d'\log^4(1/\d')),\\
        M^{\far}
        &\coloneqq  
        (\e^{50}\Z)\cap (
            \R\setminus J^{\d'}
        ),\\
        M^\bdry
        &\coloneqq  (\e^{50}\Z)\cap [-\e\log^4(1/\e), \e\log^4(1/\e)].
    \end{split}
\end{equation}
We then define the following random vectors (depicted below in Figure \ref{fig:noise-eps-mesh} via highlighted paths):
\begin{equation}\label{e:def-Y-vectors}
    \begin{split}
    \bL^{\near} &\coloneqq  
    \bigcup_{k=1}^{n-1}\left\{\cL(t_k,v;t_{k+1},v') : 
    v,v'\in M^{\near}_{S_1}\cup M^{\near}_{S_2},
    \quad |v-v'| < \e\log(1/\e)
    \right\},\\
    \bL^{\far} &\coloneqq  
    \bigcup_{k=1}^{n-1}
    \left\{\cL(t_k,v;t_{k+1},v') : v,v'\in
    \e^{50}\Z,
    \quad \{v,v'\}\cap M^\far \ne\varnothing,
    \quad |v-v'| < \e\log(1/\e)
    \right\},\\
    \bL^{\bdry} &\coloneqq  
    \bigcup_{k=1}^{n-1}
    \left\{\cL(t_k,v;t_{k+1},v') : v,v'\in
    \e^{50}\Z,
    \quad 
    \{v,v'\}\cap M^{\bdry} \ne \varnothing,
    \quad
    |v-v'| < \e\log(1/\e)
    \right\}\\
    &\qquad\cup
    \left\{
        \cL(0,v;t_{1},v') : v\in (\e^{50}\Z)\cap J^\d,\quad v'\in \e^{50}\Z,\quad |v-v'|< \e^{2\gamma/3}\log(1/\e)
    \right\}\\
    &\qquad\cup
    \left\{
        \cL(t_{n},v;t_{n+1},v') : v'\in (\e^{50}\Z)\cap J^\d,\quad v\in \e^{50}\Z,\quad |v-v'|< \e^{2\gamma/3}\log(1/\e)
    \right\}.
    \end{split}
\end{equation}
Here we are abusing the notation $\{\}$ and $\cup$: the above collections are indeed vectors, i.e. ordered lists, and $\cup$ denotes concatenation of ordered lists.
We also define the restricted length counterpart of $\bL^\near$:
\begin{equation}\label{e:bL-restricted}
    \begin{split}
        \bL^\near_{S_1} &\coloneqq 
        \bigcup_{k=1}^{n-1}\left\{\cL_{S_1}(t_k,v;t_{k+1},v') : 
        v,v'\in M^{\near}_{S_1},
        \quad |v-v'| < \e\log(1/\e)
        \right\},\\
        \bL^\near_{S_2} &\coloneqq  
        \bigcup_{k=1}^{n-1}\left\{\cL_{S_2}(t_k,v;t_{k+1},v') : 
        v,v'\in M^{\near}_{S_2},
        \quad |v-v'| < \e\log(1/\e)
        \right\},\\
        \bL^\near_{S_1\cup S_2}
        &\coloneqq  \bL^\near_{S_1}\cup \bL^\near_{S_2}.
    \end{split}
\end{equation}
Note that for $v\in M^\near_{S_1}$ and $v'\in M^\near_{S_2}$ we have $|v-v'|> 2\e\log^4(1/\e)$, and therefore $\bL^\near$ and $\bL^\near_{S_1\cup S_2}$ are vectors of the same (finite) length.

\begin{figure}[thbp]
    \centering
    \includegraphics[width=\textwidth]{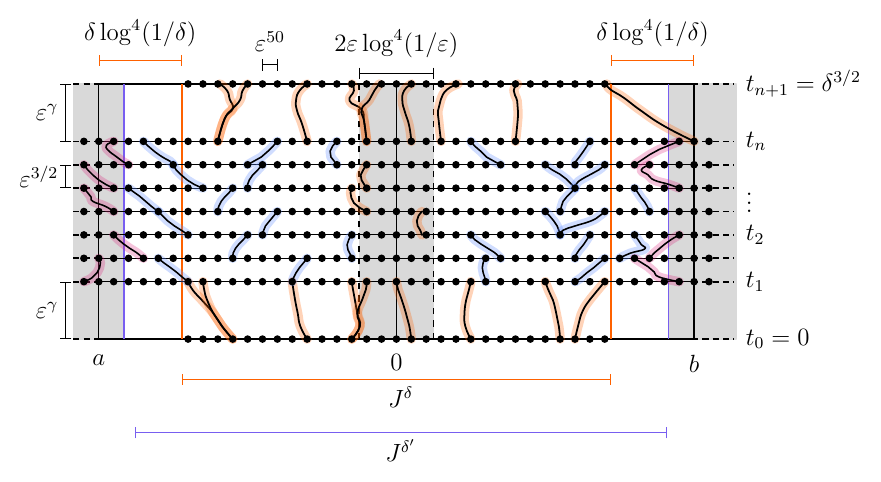}
    \caption{
        As in Figure \ref{fig:noise-delta-mesh-zoom-in} (right), we have zoomed in on the short horizontal strip $S_3 \coloneqq  [0,\d^{3/2}]\times (a,b)$ (large black rectangle).
        In the text, $S_1$ denotes the portion of the strip to the left of $0$ and $S_2$ denotes the portion to the right.
        The orange interval $J^\d$ and the purple interval $J^{\d'}$ are defined as in \eqref{e:def-Jdelta} (recall $\d'\coloneqq \d/2)$.
        The narrow shaded interval at the center is $[-\e\log^4(1/\e), \e\log^4(1/\e)]$, which acts as a cushion around the interface $[0, \d^{3/2}]\times \{0\}$ between $S_1$ and $S_2$.
        We divide the time interval $[0,\d^{3/2}]$ into times $t_0,t_1,\dots,t_{n+1}$ as in \eqref{e:eps-mesh-times}, depicted on the right.
        The black mesh points are separated from their horizontally-adjacent neighbors by distance $\e^{50}$.
        Note that the mesh points at times $t_1,\dots,t_n$ extend infinitely to the left and right, into the shaded regions beyond the strips $S_1,S_2$.\\
        \textbf{Highlighted paths:}
        The short black paths are directed landscape geodesics between various pairs of temporally-adjacent mesh points.
        They are highlighted according to which of the vectors $\bL^\near,\bL^\far,\bL^\bdry$ their lengths belong to:
        \blue for $\bL^\near$, \pink for $\bL^\far$, and \orange for $\bL^\bdry$.\\
        \textbf{Local transversal fluctuation constraint:}
        For each of the geodesics, the horizontal displacement of its  endpoints is required to be less than $(t_{k+1}-t_k)^{2/3}\log(1/\e)$.
        \\
        \textbf{Restricted lengths:}
        Like $\bL^\near$, the vectors $\bL^\near_{S_1},\bL^\near_{S_2}, \bL^\near_{S_1\cup S_2}$ are also indexed by the endpoints of the \blue paths, but the entries of the latter are defined with respect to the restricted lengths $\cL_{S_1}$ and $\cL_{S_2}$.
        }
    \label{fig:noise-eps-mesh}
\end{figure}

With the above notation in place, we turn to outlining our proof of \cref{p:cond-var-delta}, which as argued in \eqref{e:886} is equivalent to
\begin{align}\label{recall12}
    \sup_u \E[\Var(\cL(u)\mid \cF_{S_1}\vee \cF_{S_2})] \le Ce^{-c\log^2(1/\d)},
\end{align}
where the supremum is over $u=(0,z;\d^{3/2},w)$ with $z,w\in J^\d$.

The proof broadly consists of two steps.
We now present a high-level schematic outline before proceeding to the details in the upcoming subsections.

\addtocontents{toc}{\SkipTocEntry}
\subsection*{Step 1: Coupling to unrestricted lengths at scale $\e$} The first step allows us to work with the discrete model. Further, the conditioning on the restricted sigma-algebra $\cF_{S_1}\vee \cF_{S_2}$ will be replaced by an unrestricted sigma algebra. That is, we will prove that for all small $\d>0$ and $\e\le \d^{100/\gamma}$, we have
\begin{align}\label{e:mesh-step-1}
    \sup_u\E[\Var(\cL(u)\mid \cF_{S_1}\vee \cF_{S_2})]
    \le
    \sup_u\E\left[\Var(\wh{\cL}^\e(u)\mid \bL^\near)\right]
    + O(\e^9),
\end{align}
where the supremum is over $u=(0,z;\d^{3/2},w)$ with $z,w\in J^\d$.
In particular, choosing $\e = e^{-c\log^2(1/\d)}$ will yield an error term of the form on the RHS of \eqref{recall12}.

Importantly, in \eqref{e:mesh-step-1}, the RHS conditional variance does not involve any restricted quantities and only involves (a discrete proxy for) $\cL$.
The unrestrictedness of the expression is essential for our later argument, which leverages Airy line ensemble estimates that are only available in the unrestricted setting (recall the proof sketch in \cref{ss:idea}).
The replacement of the conditioning $\sigma$-algebras is justified by the fact that by geodesic transversal fluctuation estimates, $\bL^\near_{S_1\cup S_2}=\bL^\near$ with high probability.

\textbf{Step 1} is carried out in \cref{s:mesh-step-1}.

\addtocontents{toc}{\SkipTocEntry}
\subsection*{Step 2: Bounding influence using geodesic transversal fluctuations and the Efron--Stein inequality}

We first essentially reduce to the problem considered in \cref{s:iop-black-noise} of proving the join property for two half-infinite strips and to proving a conditional variance estimate analogous to \eqref{equiv}.

\begin{figure}[tbp]
    \centering
    \begin{subfigure}[c]{0.95\textwidth}
        \includegraphics[width=\linewidth]{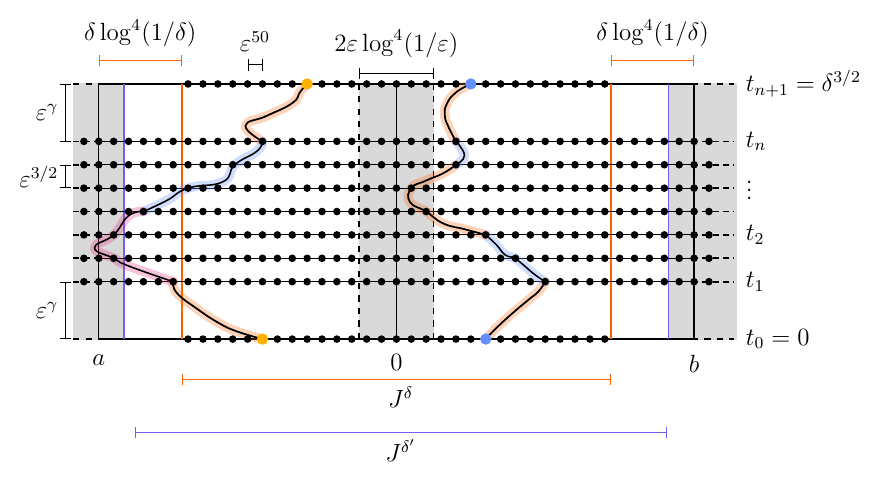}
    \end{subfigure}\\
    \vfill
    \begin{subfigure}{0.6\textwidth}
        \includegraphics[width=\linewidth]{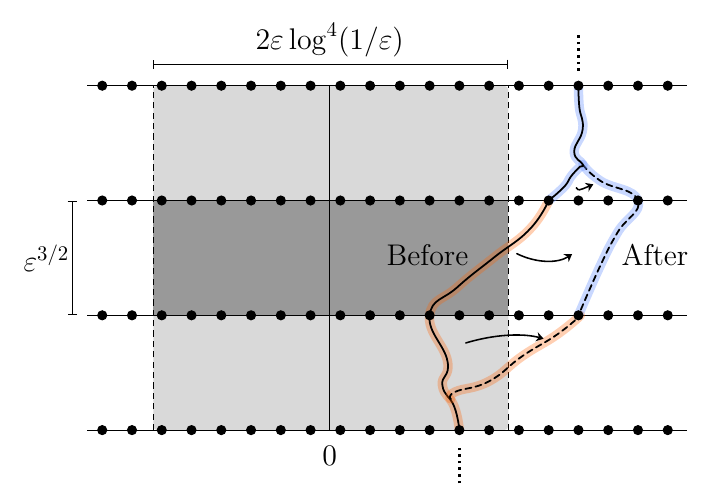}
    \end{subfigure}
    \caption{
        Resampling in the discrete LPP model.\\
        \textbf{Top:}
        Two mesh geodesics with different endpoints are drawn (black paths with orange/blue endpoints) with the segments highlighted as in Figure \ref{fig:noise-eps-mesh}.
        We seek to estimate the influence on the mesh geodesic lengths of all the \orange and \pink segments in Figure \ref{fig:noise-eps-mesh}, i.e. the segments with at least one endpoint in the shaded region (in the actual argument, we also estimate the influence of all segments crossing the unshaded strips $[t_0,t_1]\times\R$ and $[t_n,t_{n+1}]\times\R$).
        By transversal fluctuation bounds, mesh geodesics typically never enter the left or right shaded regions, and hence do not contain any \pink segments; consequently \pink segments have negligible influence.
        The influence of the \orange segments is estimated via the Efron--Stein inequality.\\
        \textbf{Bottom:} Close-up of the shaded thin strip in the top panel.  The solid path is a portion of the mesh geodesic before resampling the noise (lengths of \orange segments intersecting the darker box as well as the $\pink$ segments inside the same horizontal strip).
        After resampling the box, the mesh geodesic locally reroutes, changing course to follow the dashed path.
    }
    \label{fig:noise-resample}
\end{figure}

Viewing $\wh{\cL}^\e(u)$ as a function of $(\bL^\near,\bL^\far,\bL^\bdry)$, the conditional variance on the RHS of \eqref{e:mesh-step-1} is (up to constants) just the squared $L^2$-norm of the change in the value of $\wh{\cL}^\e(u)$ caused by resampling $(\bL^\far,\bL^\bdry)$ (i.e. resampling the lengths of the \orange and \pink paths in Figure \ref{fig:noise-eps-mesh}).
However, by geodesic transversal fluctuation bounds, with high probability the mesh geodesic does not exit $[0,\d^{3/2}]\times J^{\d'}$, i.e. the mesh geodesic's length typically does not involve the entries of $\bL^\far$ (\pink paths).
This allows us to implement the strategy outlined in Section \ref{s:iop-black-noise} in the case of two half-infinite strips, and upper bound the influence of $\bL^{\bdry}$ using the Efron--Stein inequality (as depicted in Figure \ref{fig:noise-resample}). 
In fact, we will prove that there exist constants $C,c,\alpha>0$ 
such that for all small enough $\delta>0$ and $\e\le \delta^{100/\gamma}$, we have
\begin{align}\label{e:noise-ES-unsimplified123}
    \sup_{u}
    \E\left[
        \Var(\wh{\cL}^\e(u)\mid \bL^\near)
    \right]
    &\le C\left(\e^{\alpha}+\frac{1}{\e^{3/2}}
    e^{-c\log^2(1/\d)}\right),
\end{align}
where the supremum is over all $u=(0,z;\d^{3/2},w)$ with $z,w\in J^\d$.

The $e^{-c\log^2(1/\d)}$ term is essentially the influence of the \pink paths stemming from the rare event that the geodesic exits $[0,\d^{3/2}]\times J^{\d'}$,
and the $\e^{\alpha}$ term is the influence of the \orange paths, corresponding exactly to the bound obtained in the positive-temperature setting (see \cref{t:2main}).

Finally, choosing $\e = e^{-c'\log^2(1/\d)}$ 
for small enough $c'>0$,
as mentioned after \eqref{e:mesh-step-1}, 
we will obtain a bound of $\ls e^{-c''\log^2(1/\d)}$.

\textbf{Step 2} is carried out in \cref{s:mesh-step-2}.

\addtocontents{toc}{\SkipTocEntry}
\subsection*{Putting everything together}

Before presenting the details, we quickly show how \eqref{e:mesh-step-1} and \eqref{e:noise-ES-unsimplified123} imply \eqref{recall12} (and hence \cref{p:cond-var-delta}).

\begin{proof}[Proof of \cref{p:cond-var-delta}] 
    Choosing $\e = e^{-c\log^2(1/\d)}$ (up to a $\Theta(1)$ prefactor needed for the divisibility condition $(\d^{3/2} - 2\e^\gamma)/\e^{3/2} \in \N$, see \cref{s:discrete-LPP-proxy})
    for small enough $c>0$, we get
    \begin{align*}
        \textup{LHS of \eqref{recall12}}
        \overset{\eqref{e:mesh-step-1}}&{\le} \sup_u \E\left[\Var(\wh{\cL}^\e(u)\mid \bL^\near)\right] 
        + C e^{-9c\log^2(1/\d)}\\
        \overset{\eqref{e:noise-ES-unsimplified123}}&{\le}
        C' e^{-c'\log^2(1/\d)},
    \end{align*}
    which is exactly \eqref{recall12}.
    As argued above \eqref{e:886}, this implies \cref{p:cond-var-delta}.
\end{proof}

Before implementing the above program, we make the following remark. 

\begin{rem}
A natural approach to 
bounding the LHS of \eqref{e:noise-ES-unsimplified123}
assuming $\bL^{\near},\bL^{\far},\bL^{\bdry}$
were all independent would be to apply another Efron--Stein inequality. 
Namely, recalling that $\wh{\cL}^\e(u)$ is a function of $(\bL^{\near},\bL^{\far},\bL^{\bdry})$, we would have had
\begin{align}\label{approach23}
    \E\left[
        \Var(\wh{\cL}^\e(u)\mid \bL^\near)\right] \le   \E\left[
        \Var(\wh{\cL}^\e(u)\mid \bL^\near, \bL^\bdry)\right]+ \E\left[
        \Var(\wh{\cL}^\e(u)\mid \bL^\near, \bL^\far)
    \right].
\end{align}
This separates the influences of the \pink and \orange paths.
Now, transversal fluctuation estimates for the mesh geodesic may be used to show that the first term is $\ls e^{-c\log^2(1/\d)}$,
while an Efron--Stein argument can be used to deduce that the second term is $\ls \e^{\alpha}$ much like in the infinite half-strip setting considered in \cref{s:iop-black-noise}.
Instead, the proof of  \eqref{e:noise-ES-unsimplified123} combines these two arguments, and only delivers the weaker bound where the first term is replaced by $\frac{1}{\e^{3/2}}e^{-c\log^2(1/\d)}$.
The reason for this is that for \eqref{approach23} to hold, we really need a form of independence between $\bL^{\far}$ and $\bL^{\bdry}$ \emph{conditional} on $\bL^{\near}$. While available mixing estimates (such as \cref{l:DL-mixing} which featured in our earlier arguments) imply a strong form of decay of correlation between $\bL^{\far}$ and $\bL^{\bdry}$, an analogous conditional mixing statement isn't quite available. 
Although we expect such a statement to be true and plan to investigate this in the future,
here we will only prove the weaker bound in \eqref{e:noise-ES-unsimplified123} which already suffices for our applications.
\end{rem}

\subsection{Completing Step 1 (unrestricted lengths at scale $\e$)}\label{s:mesh-step-1}

To prove \eqref{e:mesh-step-1}, we will first make three estimates using abstract properties of conditional variance, and then bound various error terms arising in the process.

Fix $z,w\in J^\d$ and write $u\coloneqq (0,z;\d^{3/2},w)$.
We want to bound $\E[\Var(\cL(u)\mid \cF_{S_1}\vee\cF_{S_2})]$ as in \eqref{e:mesh-step-1}.
We will first replace $\cF_{S_1}\vee\cF_{S_2}$ with $\bL^{\near}$.
After that, it will be easy to replace $\cL(u)$ with $\wh{\cL}^\e(u)$ thanks to their strong similarity guaranteed by \cref{l:moments-tails-wh-L}.

Since $\bL^\near_{S_1\cup S_2}$ is $\cF_{S_1}\vee\cF_{S_2}$-measurable (the former is defined in \eqref{e:bL-restricted}), it follows by a general property of conditional variance (\cref{l:conditioning-on-less}) that
\begin{align}\label{e:8910}
    \E[\Var(\cL(u)\mid \cF_{S_1}\vee\cF_{S_2})]
    &\le
    \E[\Var(\cL(u)\mid \bL^\near_{S_1\cup S_2})].
\end{align}
We next wish to bound \eqref{e:8910} by the conditional variance given the \emph{unrestricted} data $\bL^\near$.
At the end of this subsection we will show, using geodesic transversal fluctuation estimates, that $\bL^\near_{S_1\cup S_2} = \bL^\near$ with high probability.
(It may help to note that $\bL^\near_{S_1\cup S_2} = \bL^\near$ is equivalent to all the \blue-highlighted geodesics in Figure \ref{fig:noise-eps-mesh} being contained in $S_1\cup S_2$.)
In light of this, it is convenient to use the following general conditional variance lemma to replace the conditioning on the RHS of \eqref{e:8910} with a conditioning on $\bL^\near$, at the cost of an error depending on $\P(\bL^\near_{S_1\cup S_2} \ne \bL^\near)$.

\begin{lemm}[Conditioning on similar data]\label{l:cond-var-comparison}
    Let $X$ be a random variable such that $\P(|X|>m) \le Ce^{-cm^{3/2}}$ for some $C,c>0$ and all $m\ge0$.
    Let $Y,Z$ be random $d$-dimensional vectors on the same probability space as $X$, for some $d\ge 1$.
    We have
    \begin{align*}
        \Bigl|
            \E[\Var(X|Y)] - \E[\Var(X|Z)]
        \Bigr|
        \le C'\,\P(Y\ne Z)^{0.9}
    \end{align*} 
    for some $C'>0$ depending only on $C,c$.
    (The exponent $0.9$ can be replaced with any $\kappa<1$.)
\end{lemm}
\vspace{-0.5ex}
\noindent
The proof of \cref{l:cond-var-comparison} is deferred to Appendix \ref{s:conditional-variance}.
\vspace{2ex}

Before continuing with the argument, we make a brief notational change to facilitate later moment estimates.
Since we are only interested in one-point central moments, we may translate $\cL$ by any function without changing the problem.
In particular, it will be convenient to add a parabolic term to $\cL$ to leverage stationarity.
Accordingly, we define
\begin{align*}
    \cK(0,z;\d^{3/2},w) &\coloneqq  \cL(0,z;\d^{3/2},w) + \frac{(z-w)^2}{\d^{3/2}},\\
    \wh{\cK}^\e(0,z;\d^{3/2},w) &\coloneqq  \wh{\cL}^\e(0,z;\d^{3/2},w) + \frac{(z-w)^2}{\d^{3/2}}.
\end{align*}

Continuing, by the above lemma, we have
\begin{equation}\label{e:891}
    \begin{split}
        \E[\Var(\cL(u)\mid \bL^\near_{S_1\cup S_2})]
        &=  \E[\Var(\cK(u)\mid \bL^\near_{S_1\cup S_2})]\\
        &\le
        \E[\Var(\cK(u)\mid \bL^\near)]
        + O\Bigl(\P(\bL^\near_{S_1\cup S_2} \ne \bL^\near)^{0.9}\Bigr),
    \end{split}
\end{equation}
where the big-$O$ constant is universal thanks to the uniform one-point tail estimate
\begin{align}\label{e:TW-tail}
    \sup_u \P(|\cK(u)| > m) 
    = \P(|\cL(0,0;1,0)| > m \d^{-1/2})
    \le Ce^{-cm^{3/2}}
    \qquad\text{for all $m\ge 0$},
\end{align}
a consequence of \cref{p:landscape-symmetries}\ref{landscape-scaling}--\ref{landscape-shear} and the fact that $\cL(0,0;1,0)$ follows the Tracy--Widom GUE distribution.

\begin{rem}
    \cref{l:cond-var-comparison} is one of our main technical motivations for approximating $\cL$ with a discrete LPP model.
    It provides a convenient way to formalize the intuition that $\cL$ restricted to $S_3$ is ``almost $\cF_{S_1}\vee\cF_{S_2}$-measurable.''
\end{rem}

Continuing from \eqref{e:891}, and looking towards \eqref{e:mesh-step-1}, we approximate $\cK(u)$ by its discrete counterpart $\wh{\cK}^\e(u)$.
By the property of conditional variance recorded in \cref{l:difference-conditional-variances}, we have
\begin{multline*}
    \E\left[
        \Var(\cK(u)\mid \bL^\near)
    \right]
    \le 
    \E\left[
        \Var(\wh{\cK}^\e(u)\mid \bL^\near)
    \right]\\
    +
    \left(
        \E\left[\cK(u)^2\right]^{1/2}
        +\E\left[\wh{\cK}^\e(u)^2\right]^{1/2}
    \right)
    \E\left[\left(\cK(u) - \wh{\cK}^\e(u)\right)^2\right]^{1/2}.
\end{multline*}
We now bound the second term on the RHS.
For the factor $\E\left[\left(\cK(u) - \wh{\cK}^\e(u)\right)^2\right]^{1/2}$, we have
\begin{align*}
    \E\left[\left(\cK(u) - \wh{\cK}^\e(u)\right)^2\right]^{1/2}
    &\le
    \e^9
    +
    \E\left[\left(\cK(u) - \wh{\cK}^\e(u)\right)^4\right]^{1/4}
    \P\left(\left|\cK(u) - \wh{\cK}^\e(u)\right|>\e^9\right)^{1/4}\\
    &\le \e^9
    + \left(
        \E\left[\cK(u)^4\right]^{1/4}
        +
        \E\left[
            \wh{\cK}^\e(u)^4
        \right]^{1/4}
    \right)
    Ce^{-c\log^2(1/\e)},
\end{align*}
where the first line is by Cauchy--Schwarz, and the second is by the $L^4$-norm triangle inequality and the probability bound \cref{l:moments-tails-wh-L} (note that $\cK-\wh{\cK}^\e = \cL-\wh{\cL}^\e$).
By \eqref{e:TW-tail} and \cref{l:moments-wh-L},
the second and fourth moments of $\cK(u)$ and $\wh{\cK}^\e(u)$ are all bounded uniformly in $u,\d,\e$.
Combining this with the above two displays, we obtain 
\begin{align*}
    \E\left[
        \Var(\cK(u)\mid \bL^\near)
    \right]
    &\le
    \E\left[
        \Var(\wh{\cK}^\e(u)\mid \bL^\near)
    \right]
    + 
    O(\e^{9}).
\end{align*}

Applying the above in \eqref{e:891} and replacing $\cK,\wh{\cK}^\e$ with their translates $\cL,\wh{\cL}^\e$, we obtain
\begin{align*}
    \E[\Var(\cL(u)\mid \cF_{S_1}\vee \cF_{S_2})]
    &\le
    \E\left[
        \Var(\wh{\cL}^\e(u)\mid \bL^\near)
    \right]
    + 
    O\Bigl(
        \P(\bL^\near_{S_1\cup S_2} \ne \bL^\near)^{0.9} 
        + \e^{9}
    \Bigr).
\end{align*}
Finally, we will show in the upcoming \cref{l:mesh-restricted-equals-full-simpler} that $\P(\bL^\near_{S_1\cup S_2} \ne \bL^\near) = O(e^{-c\log^{2}(1/\e)})$.
Given this, by the above display, we get
\begin{align}\label{e:8920}
    \E[\Var(\cL(u)\mid \cF_{S_1}\vee \cF_{S_2})]
    &\le
    \E\left[
        \Var(\wh{\cL}^\e(u)\mid \bL^\near)
    \right]
    + O(\e^9),
\end{align}
which is exactly \textbf{Step 1}'s \eqref{e:mesh-step-1}.
Note in particular that by choosing $\e$ of the form $e^{-c\log^2(1/\d)}$, we get an error term of the form $O(e^{-c\log^2(1/\d)})$.

We close this subsection with the following lemma that was mentioned just above.

\begin{lemm}\label{l:mesh-restricted-equals-full-simpler}
    There exist $C,c>0$ such that for all small $\delta>0$ and all $\e \le \delta^{100/\gamma}$,
    we have $\P\left(\bL^\near_{S_1\cup S_2} \ne \bL^\near\right) \le Ce^{-c\log^{2}(1/\e)}$.
\end{lemm}
\begin{proof}
    We first estimate geodesic transversal fluctuations via a similar argument as in the reduction from \cref{p:cond-var-delta} to Property \ref{noise-join}, see below \eqref{e:910}.
    Let $\CTF>0$ be the random constant from \cref{p:GZ-uniform-TF} corresponding to $L=b-a+1$.
    Let $c_0>0$ be such that $c_0 \log(1/\e)  \cdot \log^3(1/\e^{3/2}) < \log^4(1/\e)$.
    Then by \cref{p:GZ-uniform-TF}, on the event $\sA \coloneqq  \{\CTF \le c_0\log(1/\e)\}$, 
    we have that 
    \begin{align*}
        \max_{k\in \lb 1,n-1\rb}
        \max_{\substack{(v,v') \in (M^\near_{S_1})^2 \cup (M^\near_{S_2})^2\\ |v-v'|<\e\log(1/\e)}}
        \;\sup_{r\in[t_k, t_{k+1}]}
        \left|
            \Pi_{(t_k, v), (t_{k+1},v')}(r)
            - \frac{v(t_{k+1}-r) + v'(r-t_k)}{t_{k+1}-t_k}
        \right|
        < \e\log^4(1/\e)
    \end{align*}
    because $t_{k+1}-t_k = \e^{3/2}$ (see \eqref{e:eps-mesh-times}).
    Since the distance from $M^\near_{S_1}\cup M^\near_{S_2}$ to the boundary $\partial((a,b)\setminus\{0\})$ is at least $\e\log^4(1/\e)$ (see \eqref{e:mesh-sets} and Figure \ref{fig:noise-eps-mesh}),
    the above display implies that every geodesic on the LHS spends its entire journey inside $[0,\d^{3/2}]\times ((a,b)\setminus\{0\})$.
    In particular, all the restricted lengths coincide with their unrestricted counterparts: for all $k\in\lb 1,n-1\rb$ and all $v,v'\in M^\near_{S_1}$, we have by \eqref{e:restricted-length} that
    \begin{align*}
        \cL(t_k,v;t_{k+1},v')
        = \cL_{S_1}(t_k,v;t_{k+1},v'),
    \end{align*}
    and the same holds for $v,v'\in M^\near_{S_2}$ and $\cL_{S_2}$.
    In other words, we have $\sA \subset \{\bL^\near_{S_1\cup S_2} = \bL^\near\}$, by definition of $\bL^\near_{S_1\cup S_2},\bL^\near$ (see \eqref{e:def-Y-vectors}--\eqref{e:bL-restricted}).
    It remains to note that $\P(\neg \sA) \le Ce^{-c\log^2(1/\e)}$ by \cref{p:GZ-uniform-TF}.
\end{proof}

\subsection{Step 2 (mesh geodesic transversal fluctuations and Efron--Stein)}\label{s:mesh-step-2}

In this subsection we prove the following bound (recorded already in \eqref{e:noise-ES-unsimplified123}).

\begin{prop}[Influence of $\bL^\bdry$ and $\bL^\far$]\label{p:half-strip-ES}
    Let $\gamma\in(0,\frac{1}{6})$ be the parameter fixed at the beginning of \cref{s:discrete-LPP-proxy}.
    There exist $C,c,\alpha>0$ depending only on $S_1,S_2,\gamma$
    such that for all small enough $\delta>0$ and $\e\le \delta^{100/\gamma}$, we have
    \begin{align}\label{e:noise-ES-unsimplified}
        \sup_{u}
        \E\left[
            \Var(\wh{\cL}^\e(u)\mid \bL^\near)
        \right]
        &\le C\left(\e^{\alpha}+\frac{1}{\e^{3/2}}
        e^{-c\log^2(1/\d)}\right),
    \end{align}
    where the supremum is over all $u=(0,z;\d^{3/2},w)$ with $z,w\in J^\d$.
\end{prop}

As indicated above \eqref{e:noise-ES-unsimplified123}, a key ingredient in the proof of \cref{p:half-strip-ES} will be mesh geodesic transversal fluctuation bounds.
We begin by preparing the latter.

\addtocontents{toc}{\SkipTocEntry}
\subsection*{Mesh geodesic transversal fluctuations}

Let $\Gamma_u$ be the mesh geodesic for $\wh{\cL}^\e(u)$ as in \cref{def:discrete-LPP}.
We view $\Gamma_u$ as a $C([0,\d^{3/2}])$-valued measurable function of the data $(\bL^\near, \bL^\far, \bL^\bdry)$.
Define the event
\begin{equation}\label{e:inside-event}
    \Inside_u \coloneqq  \left\{
        (\bL^\near, \bL^\far, \bL^\bdry) : 
        \Gamma_u(r) \in J^{\d'}\quad\text{for all } r\in [0,\d^{3/2}]
    \right\}.
\end{equation}
\begin{lemm}[Mesh geodesic typically stays inside]\label{l:no-dependence-on-far}
    For all $\e\le \d^{100/\gamma}$, we have
    \[
        \sup_{u}\P(\neg \Inside_u) \le C e^{-c\log^2(1/\d)},
    \]
    where the supremum is over $u=(0,z;\d^{3/2},w)$ with $z,w\in J^\d$.
\end{lemm}

To prove \cref{l:no-dependence-on-far}, we need to bound the transversal fluctuations of mesh geodesics (see \cref{def:discrete-LPP}).
Since mesh geodesics are not really geodesics, we cannot appeal to known bounds for geodesic transversal fluctuations (such as Propositions \ref{p:DOV-geodesic-holder} and \ref{p:GZ-uniform-TF}, see also the pink paths in Figure \ref{fig:noise-delta-mesh}).
Instead, we will rely on the fact that with high probability, \emph{any} path with a large transversal fluctuation has a significantly shorter length than the geodesic with the same endpoints (a weaker form of this was used earlier in the proof of \cref{l:restricted-L-properties}\ref{property-LR-continuous}).
This together with a statement that the mesh geodesic has length close to that of the true geodesic will finish the argument.

\begin{lemm}[Paths with large transversal fluctuations are short]\label{l:large-TF-is-uncompetitive}
    For $z,w\in J^\d$, define the set of paths
    \begin{align*}
        \mathrm{BigTF}^\d(z;w) \coloneqq  
        \left\{
            \pi \in C([0,\d^{3/2}]) : \pi(0)=z,\quad \pi(\d^{3/2})=w,\quad 
            \exists r\in[0,\d^{3/2}]\text{ such that }
            \pi(r) \not\in J^{\d'}
        \right\}.
    \end{align*}
    Define also the event
    \begin{align}\label{e:def-big-tf-uncompetitive}
        \NoBigTF \coloneqq  
        \left\{
            \inf_{z,w\in J^\d}
            \left(
                \cL(0,z;\d^{3/2},w) - 
                \sup_{\pi\in \mathrm{BigTF}^\d(z;w)}
                \int_0^{\d^{3/2}} d\cL\circ \pi
            \right)
            \ge 
            \d^{1/2}\log^7(1/\d)
        \right\}.
    \end{align}
    There exist $C,c>0$ such that for all sufficiently small $\d$,
    \begin{align*}
        \P\left( \NoBigTF \right)
        \ge 1- Ce^{-c\log^2(1/\d)}.
    \end{align*}
\end{lemm}
\begin{proof}
    Fix $\pi\in\mathrm{BigTF}^\d(z;w)$, and let $r_0\coloneqq \inf\{r\in[0,\d^{3/2}]: \pi(r)\not\in J^{\d'}\}$ be the first time $\pi$ exits $J^{\d'}$.
    By the reverse triangle inequality for $\cL$ (see \cref{p:landscape-symmetries}\ref{landscape-metric-composition}), we have
    \begin{align*}
        \int_0^{\d^{3/2}}d\cL\circ \pi \le 
        \cL(0,z;r_0,\pi(r_0)) + \cL(r_0,\pi(r_0); \d^{3/2},w),
    \end{align*}
    with equality attained if $\pi$ is the concatenation of two geodesics.
    By continuity we have $\pi(r_0) \in \partial J^{\d'}=\{a+\d'\log^4(1/\d'), b-\d'\log^4(1/\d')\}$, so we obtain
    \begin{align*}
        \sup_{\pi\in \mathrm{BigTF}^\d(z;w)}
            \int_0^{\d^{3/2}} d\cL\circ \pi
        &=
        \sup_{\tau \in (0, \d^{3/2})}
        \max_{\zeta\in \partial J^{\d'}}
        \left(
            \cL(0,z;\tau,\zeta) + \cL(\tau,\zeta; \d^{3/2},w)
        \right).
    \end{align*}

    Next, by \cref{p:dov-pointwise},
    the following holds with probability at least $1-Ce^{-c\log^2(1/\d)}$:
    \begin{align*}
        \sup_{\substack{0\le s < t \le \d^{3/2}\\ x,y\in [a,b]}}
        \left|
            \cL(s,x;t,y) + \frac{(x-y)^2}{t-s}
        \right| 
        \le \d^{1/2}\log^{8/3}(1/\d).
    \end{align*}
    On the above event, we have that
    \begin{multline*}
        \inf_{z,w\in J^\d}
        \left[
        \cL(0,z;\d^{3/2},w)
        -
        \sup_{\tau\in(0,\d^{3/2})}
        \max_{\zeta\in\partial J^{\d'}}
        \left(
            \cL(0,z;\tau,\zeta) + \cL(\tau,\zeta;\d^{3/2},w)
        \right)
        \right]\\
        \begin{aligned}
        &\ge
        \inf_{z,w\in J^\d}
        \left[
            -\frac{(z-w)^2}{\d^{3/2}}
            + \inf_{\tau\in(0,\d^{3/2})}\min_{\zeta\in\partial J^{\d'}}
            \left(
                \frac{(z-\zeta)^2}{\tau} + \frac{(w-\zeta)^2}{\d^{3/2}-\tau}
            \right)
        \right]
        - 3\d^{1/2}\log^{8/3}(1/\d)\\
        &\ge
        \inf_{z,w\in J^\d}
        \min_{\zeta\in\partial J^{\d'}}
        \left(
            \frac{(z-\zeta)^2 + (w-\zeta)^2 - (z-w)^2}{\d^{3/2}}
        \right)
        -3\d^{1/2}\log^{8/3}(1/\d)\\
        &=\inf_{z,w\in J^\d}
        \min_{\zeta\in\partial J^{\d'}}
        \left(
            \frac{2(\zeta-z)(\zeta-w)}{\d^{3/2}}
        \right)
        -3\d^{1/2}\log^{8/3}(1/\d)\\
        &= \frac{2(\d\log^4(1/\d) - \d'\log^4(1/\d'))^2}{\d^{3/2}} 
        -3\d^{1/2}\log^{8/3}(1/\d)\\
        &\ge
        \frac{1}{50}\d^{1/2}\log^8(1/\d)
        -3\d^{1/2}\log^{8/3}(1/\d),
        \end{aligned}
    \end{multline*}
    where in the last line we used that $\d\log^4(1/\d) - \d'\log^4(1/\d') > \frac{1}{10}\d\log^4(1/\d)$ (see the start of \cref{s:discrete-lpp-reduction}).
    It remains to note that the RHS is lower bounded by $\d^{1/2}\log^7(1/\d)$
    for all sufficiently small $\d>0$.
\end{proof}

We are ready to prove \cref{l:no-dependence-on-far}.

\begin{proof}[Proof of \cref{l:no-dependence-on-far}]
    Using the notation of \cref{l:large-TF-is-uncompetitive}, by \eqref{e:inside-event} we have
    \begin{align*}
        \neg \Inside_u 
        =
        \left\{\Gamma_u \in \mathrm{BigTF}^\d(\wh{z};\wh{w})\right\},
    \end{align*}
    where $\Gamma_u$ is the mesh geodesic (\cref{def:discrete-LPP}).
    So on the event $\NoBigTF \cap \neg \Inside_u$ (see \cref{l:large-TF-is-uncompetitive}), we have
    \begin{align*}
        \cL(0,\wh{z};\d^{3/2},\wh{w}) - 
        \wh{\cL}^\e(0,z;\d^{3/2},w)
        \ge 
        \d^{1/2} \log^7(1/\d).
    \end{align*}
    On the other hand, since $\e \le \d^{100/\gamma}$, \cref{l:moments-tails-wh-L} implies that with probability at least $1-Ce^{-c\log^2(1/\e)}$, we have
    \begin{align*}
        \cL(0,\wh{z};\d^{3/2},\wh{w}) - \wh{\cL}^\e(0,z;\d^{3/2},w)
        \le \e^9,
    \end{align*}
    which is incompatible with the estimate in the previous display.
    Finally, using \cref{l:large-TF-is-uncompetitive} to bound the probability of $\neg \NoBigTF$, we obtain
    \begin{align*}
        \sup_u 
        \P(\neg\Inside_u)
        \le
        \P(\neg\NoBigTF) + Ce^{-c\log^2(1/\e)}
        \le C' e^{-c'\log^2(1/\d)}.
    \end{align*}
    This completes the proof of \cref{l:no-dependence-on-far}.
\end{proof}

With the above preparations made, we will now prove Proposition \ref{p:half-strip-ES}.

\addtocontents{toc}{\SkipTocEntry}
\subsection*{Efron--Stein argument for Proposition \ref{p:half-strip-ES}}

As mentioned earlier, we will be resampling the lengths of all the \orange and \pink geodesics in Figure \ref{fig:noise-eps-mesh}.
By the Efron--Stein inequality, it suffices to understand the effect of resampling the lengths on a single time interval $[t_k,t_{k+1}]$, as depicted in Figure \ref{fig:noise-resample}.
Like in the positive-temperature setting (\cref{s:es}), the analysis will differ in the cases $k\in\{0,n\}$ and $k\in\lb 1,n-1\rb$.
Accordingly, we now split $\bL^\far$ and $\bL^\bdry$ into several vectors, one for each time slice.

For $k\in\lb 1,n-1\rb$, write
\begin{align*}
    \LFarComp{k} \coloneqq  
    \left\{\cL(t_k,v;t_{k+1},v') : v,v'\in
    \e^{50}\Z,
    \quad 
    \{v,v'\}\cap M^{\far} \ne \varnothing,
    \quad
    |v-v'| < \e\log(1/\e)
    \right\},
\end{align*}
so that $\bL^\far = \bigcup_{k=1}^{n-1} \LFarComp{k}$ (see \eqref{e:def-Y-vectors}).
Similarly, for $k\in\lb 1,n-1\rb$, write
\begin{align}\label{e:Yedge-components-middle}
    \LBdryComp{k} \coloneqq  
    \left\{\cL(t_k,v;t_{k+1},v') : v,v'\in
    \e^{50}\Z,
    \quad 
    \{v,v'\}\cap M^{\bdry} \ne \varnothing,
    \quad
    |v-v'| < \e\log(1/\e)
    \right\},
\end{align}
and write
\begin{equation}\label{e:Yedge-components}
    \begin{split}
        \LBdryComp{0} &\coloneqq  
    \left\{
        \cL(0,v;t_{1},v') : v\in (\e^{50}\Z)\cap J^\d, \quad v'\in \e^{50}\Z,\quad |v-v'|< \e^{2\gamma/3}\log(1/\e)
    \right\},\\
    \LBdryComp{n} &\coloneqq \left\{
        \cL(t_{n},v;\d^{3/2},v') : v'\in (\e^{50}\Z)\cap J^\d, \quad v\in \e^{50}\Z,\quad |v-v'|< \e^{2\gamma/3}\log(1/\e)
    \right\},
    \end{split}
\end{equation}
so that $\bL^{\bdry} = \bigcup_{k=0}^n \LBdryComp{k}$.

By \cref{def:discrete-LPP}, 
$\wh{\cL}^\e(u)$ is a measurable function of $(\bL^\near,\bL^\far,\bL^\bdry)$,
so we can write
\begin{align}\label{e:def-fu}
    \wh{\cL}^\e(u) = f_u(\bL^\near, \bL^\far, \bL^\bdry)
\end{align}
for some measurable function $f_u : \R^\infty \to \R$.
In other words, $\wh{\cL}^\e(u)$ is a function of the lengths of all the highlighted geodesics in Figure \ref{fig:noise-eps-mesh} (and many others not drawn there).

Note that $\LBdryComp{0}, (\LFarComp{1},\LBdryComp{1}), \dots, (\LFarComp{n-1},\LBdryComp{n-1}), \LBdryComp{n}$ are conditionally independent given $\bL^{\near}$,
due to the independence of $\cL(t_k,\smallbullet;t_{k+1},\smallbullet)$ across $k\in\lb 0,n\rb$ (\cref{p:landscape-symmetries}\ref{landscale-independent-increments}).
This conditional independence permits us to apply the Efron--Stein inequality (\cref{p:efron-stein-classical}) to bound 
\begin{equation}\label{e:efron-stein-noise}
    \begin{split}
        \E\left[
            \Var(\wh{\cL}^\e(u)\mid \bL^\near)
        \right]
        &\le
        \sum_{k=1}^{n-1} \E\left[
            \left|
                f_u(\bL^{\near}, \bL^{\far},\bL^{\bdry})
                - f_u(\bL^{\near}, \bL^{\far,(k)},\bL^{\bdry,(k)})
            \right|^2
        \right]\\
        &\qquad+\E\left[\left|
                f_u(\bL^{\near}, \bL^{\far},\bL^{\bdry})
                - f_u(\bL^{\near}, \bL^{\far},\bL^{\bdry,(0)})
            \right|^2\right]\\
        &\qquad+\E\left[\left|
                f_u(\bL^{\near}, \bL^{\far},\bL^{\bdry})
                - f_u(\bL^{\near}, \bL^{\far},\bL^{\bdry,(n)})
            \right|^2\right],
    \end{split}
\end{equation}
where $(\bL^{\near}, \bL^{\far,(k)},\bL^{\bdry,(k)})$ is obtained from $(\bL^{\near}, \bL^{\far},\bL^{\bdry})$ by replacing $(\LFarComp{k}, \LBdryComp{k})$ with an independent copy $(\LFarCompRe{k}, \LBdryCompRe{k})$ conditional on everything else (as mentioned, this resampling procedure is illustrated in Figure \ref{fig:noise-resample}, bottom).
Similarly, for the $k=0,n$ terms,  $(\bL^{\near}, \bL^{\far},\bL^{\bdry,(k)})$ is obtained from $(\bL^{\near}, \bL^{\far},\bL^{\bdry})$ by replacing $\LBdryComp{k}$ with a conditionally independent copy $\LBdryCompRe{k}$.
Thus, for the $k=0,n$, terms, we are  resampling the entire discretized directed landscape $\wh{\cL}^\e$ on the strip $[t_k,t_{k+1}]\times\R$.

We will bound each term on the RHS of \eqref{e:efron-stein-noise}, focusing separately on the cases $k\in\{0,n\}$ and $k\in\lb 1,n-1\rb$, just like in \cref{s:es}.
Before splitting into cases, we first prove a worst-case upper bound for the change in $\wh{\cL}^\e$ caused by resampling at height $t_k$, applicable for any $k\in\lb 0,n\rb$.
{
To state the result in a unified way, we write
\begin{equation*}
    \begin{split}
        \bL &\coloneqq  (\bL^\near, \bL^\far, \bL^\bdry),\\
        \bL^{(k)} &\coloneqq  (\bL^\near, \bL^{\far,(k)}, \bL^{\bdry,(k)}), \qquad k\in\lb 1,n-1\rb,\\
        \bL^{(k)} &\coloneqq  (\bL^\near, \bL^{\far}, \bL^{\bdry,(k)}), \qquad k\in\{0,n\}.
    \end{split}
\end{equation*}
We recall from \eqref{e:inside-event} the event $\Inside_u$, which says that the mesh geodesic $\Gamma_u$ stays inside $[0,\d^{3/2}]\times J^{\d'}$ throughout its journey:
\begin{equation}\label{e:inside-recall}
    \begin{split}
        \Inside_u &\coloneqq  \left\{
            \bL : 
            \Gamma_u(r) \in J^{\d'}\quad\textup{for all } r\in[0,\d^{3/2}]
        \right\}\\
        \overset{\eqref{e:def-Y-vectors}}&{=}\left\{
            \bL : \textup{$\Gamma_u$ only uses the entries in $\bL^\near \cup \bL^\bdry$}
        \right\}.
    \end{split}
\end{equation}
Next, consider the event that before and after resampling the noise at level $k$, the mesh geodesic stays inside $[0,\d^{3/2}]\times J^{\d'}$:
\begin{equation}\label{e:both-inside}
    \begin{split}
        \BothIn^k_u &\coloneqq  
        \{
        \bL \in \Inside_u
        \}
        \cap \{\bL^{(k)} \in \Inside_u\}\\
        &= 
        \left\{
            \textup{$\Gamma_u$ only uses the entries of $\bL^\near \cup \bL^\bdry$}
        \right\}\\
        &\qquad
        \cap
        \left\{
            \textup{$\Gamma_u^{(k)}$ only uses the entries of $\bL^\near \cup \bL^{\bdry,(k)}$}
        \right\},
    \end{split}
\end{equation}
where $\Gamma_u$ and $\Gamma_u^{(k)}$ are the mesh geodesics corresponding to $\bL$ and $\bL^{(k)}$, respectively.
By Lemma \ref{l:no-dependence-on-far} and a union bound, we have 
\begin{equation}\label{good1243}
    \sup_u \P\left(\neg \BothIn^k_u\right) 
    \le Ce^{-c\log^2(1/\d)}.
\end{equation}
}

Since it suffices for our arguments, we state the following bound only on the above good event.

\begin{lemm}[Worst-case resampling bound]\label{l:uniform-bound-as}
    There exists a random constant $\Cresample>0$ such that the following holds.
    First, we have
    \begin{align}\label{e:Cresample-prob}
        \P(\Cresample > m\log(1/\e)) \le C e^{-cm^{3/2}}
        \qquad\textup{for all $m\ge 0$},
    \end{align}
    where $C,c>0$ depend only on $S_1,S_2$.
    Second, for all $k\in\lb 0,n\rb$ and all $u=(0,z;\d^{3/2},w)$ with $z,w\in J^\d$,
    we have
    \begin{align}\label{e:Cresample-bound}
        \left|
            f_u(\bL)
            - f_u(\bL^{(k)})
        \right|
        \1_{\BothIn^k_u}
        \le 
        2\Cresample\,(t_{k+1}-t_k)^{1/3}\log^2(1/\e).
    \end{align}
    Finally, for all $p\ge1$ and all $k\in\lb 0,n\rb$,
    \begin{align}\label{e:uniform-pth-moment}
        \sup_u
        \E\left[
            \left|
                f_u(\bL)
                - f_u(\bL^{(k)})
        \right|^p
        \right]
        \ls_{p} 
        (t_{k+1}-t_k)^{p/3} \log^{3p}(1/\e)
        + e^{-c\log^2(1/\d)}.
    \end{align}
\end{lemm}
\begin{proof}[Proof of \cref{l:uniform-bound-as}]

    We begin by bounding the $\ell^\infty$-norm of the vector $\LBdryComp{k}$ for each $k\in\lb 0,n\rb$, using the uniform pointwise bound on $\cL$ provided by \cref{p:dov-pointwise}.

    For $v\in J^\d$, if $v'$ satisfies $|v-v'| < \e\log(1/\e)$, then $v' \in (a,b)$.
    Recalling \eqref{e:Yedge-components-middle}--\eqref{e:Yedge-components}, we deduce the following crude bound:
    for any $k\in\lb 0,n\rb$, we have
    \begin{align*}
        \norm{\LBdryComp{k}}_\infty
        &\le
        \max_{\substack{v,v'\in (\e^{50}\Z)\cap (a,b)\\|v-v'|<(t_{k+1}-t_k)^{2/3}\log(1/\e)}}
        \left|\cL(t_k,v;t_{k+1},v')\right|.
    \end{align*}
    Let $\Cptwise$ be the random constant controlling extrema of $\cL$, as defined in \cref{p:dov-pointwise} with $L=b-a+1$.
    By \cref{p:dov-pointwise}, we have
    \begin{multline*}
        \max_{\substack{v,v'\in (\e^{50}\Z)\cap (a,b)\\|v-v'|<(t_{k+1}-t_k)^{2/3}\log(1/\e)}}
        \left|\cL(t_k,v;t_{k+1},v')\right|\\
        \begin{aligned}
            &\le
            \Cptwise\,(t_{k+1}-t_k)^{1/3}
            \log^{4/3}\left(\frac{8(b-a+1)}{t_{k+1}-t_k}\right)
            +
            \frac{((t_{k+1}-t_k)^{2/3}\log(1/\e))^2}{t_{k+1}-t_k}\\
            &\le C(\Cptwise+1)(t_{k+1}-t_k)^{1/3} \log^2(1/\e),
        \end{aligned}
    \end{multline*}
    for some $C>0$.
    By the tail behavior of $\Cptwise$ (see \cref{p:dov-pointwise}), setting $\fC\coloneqq C(\Cptwise + 1)$, we  have $\P(\fC > m) \le C' e^{-cm^{3/2}}$ for all $m\ge 0$.

    Having bounded the $\ell^\infty$-norms $\norm{\LBdryComp{k}}_\infty$, we proceed to the resampling bound.
    Fix $k\in\lb 0,n\rb$.
    By definition we have $\cL^{\bdry}_k \law \widetilde\cL^{\bdry}_k$.
    Therefore, letting $\fC$ be as in the previous paragraph, there exists another random constant $\fC^{(k)}$ with the same distribution as $\fC$ such that
    \begin{align*}
        \norm{\LBdryComp{k}}_\infty\vee \norm{\LBdryCompRe{k}}_\infty
        \le 
        (\fC \vee \fC^{(k)})\,(t_{k+1}-t_k)^{1/3} \log^2(1/\e).
    \end{align*}
    Let $\Cresample \coloneqq  
    \max\{\fC, \fC^{(0)},\dots,\fC^{(n)}\}$.
    By a union bound, we have that for any $m\ge 1$,
    \begin{align*}
        \P\left(\Cresample > m\log(1/\e)\right) &\le (n+2)\P(\fC > m\log(1/\e))\\
        &\le C \e^{-3/2} e^{-cm^{3/2}\log^{3/2}(1/\e)}\\
        &\le C'e^{-c' m^{3/2}},
    \end{align*}
    since $n\le \e^{-3/2}$ is only polynomial in $\e^{-1}$.
    By increasing $C'$ we can extend the above to all $m\ge 0$, which proves \eqref{e:Cresample-prob}.
    It remains to prove \eqref{e:Cresample-bound}.

    Let $\Gamma_u$ be the mesh geodesic (\cref{def:discrete-LPP}) corresponding to $\bL$, 
    and let $\Gamma^{(k)}_u$ be the mesh geodesic corresponding to  $\bL^{(k)}$. Since we are assuming that the event $\BothIn^k_u$ occurs (see \eqref{e:both-inside}),
    it follows that
    $\Gamma_u$ only uses entries of $\bL^\near \cup \bL^\bdry$,
    and $\Gamma_u^{(k)}$ only uses entries of
    $\bL^\near \cup \bL^{\bdry,(k)}$.
    Thus, for any $k\in\lb 0,n\rb$, we have
    \begin{multline*}
        f_u(\bL) - f_u(\bL^{(k)})\\
        \begin{aligned}
            &\le
            \Bigl(\cL(t_k,\Gamma_u(t_k);t_{k+1},\Gamma_u(t_{k+1}))
            -\wt{\cL}(t_k,\Gamma_u(t_k);t_{k+1},\Gamma_u(t_{k+1}))\Bigr)\1_{\{\Gamma_u(t_k), \Gamma_u(t_{k+1})\} \cap M^\bdry \ne \varnothing}\\
            &\le 2\left(\norm{\LBdryComp{k}}_\infty\vee \norm{\LBdryCompRe{k}}_\infty\right)\\
            &\le 2 \Cresample\,(t_{k+1}-t_k)^{1/3}\log^2(1/\e),
        \end{aligned}
    \end{multline*}
    where $\wt{\cL}(t_k,\Gamma_u(t_k);t_{k+1},\Gamma_u(t_{k+1}))$ is the entry of $\bL^{(k)}$ corresponding to the points $\Gamma_u(t_k),\Gamma_u(t_{k+1})\in \e^{50}\Z$.
    The same argument with $\bL$ and $\bL^{(k)}$ interchanged completes the proof of  \eqref{e:Cresample-bound}.
    
    \eqref{e:uniform-pth-moment} is an easy consequence of \eqref{e:Cresample-bound}.
    To see this, partition the $p\th$ moment into $\BothIn^k_u$ and $\neg\BothIn^k_u$.
    On $\BothIn^k_u$, apply the almost-sure bound \eqref{e:Cresample-bound} and the tail estimate \eqref{e:Cresample-prob}.
    On $\neg\BothIn^k_u$, use Cauchy--Schwarz along with the probability bound \eqref{good1243} and the moment bound \cref{l:moments-wh-L}.
\end{proof}

\subsubsection{The case $k\in\lb 1,n-1\rb$ via Brownian comparison}

Like in the positive-temperature setting, to effectively bound the $k\in \lb 1,n-1\rb$ terms in \eqref{e:efron-stein-noise}, we must account for the fact that the mesh geodesic is \emph{delocalized} and hence typically does not pass through the pairs of mesh points involved in $\LBdryComp{k}$.
Equivalently, the mesh geodesic typically lies outside $[-\e\log^4(1/\e), \e\log^4(1/\e)]$ at times $t_k$ and $t_{k+1}$ (the shaded strip in the center of Figure \ref{fig:noise-resample}).
This is captured by the following lemma, which is analogous to \cref{l:brownian-comparison-result}.

\begin{lemm}[Mesh geodesic delocalization]\label{l:mesh-geo-deloc}
    For $u=(0,z;\d^{3/2},w)$ with $z,w\in J^\d$, 
    let $\Gamma_u$ be the mesh geodesic (see \cref{def:discrete-LPP}).
    For $k\in\lb1,n-1\rb$, we define the event
    \begin{align*}
        \MeshGeoDeloc^k_u \coloneqq  
        \left\{
            (\bL^\near, \bL^\far, \bL^\bdry) : 
            \left|
                \Gamma_u(t_k)
            \right|
            \wedge
            \left|
                \Gamma_u(t_{k+1})
            \right|
        > \e\log^4(1/\e)
        \right\},
    \end{align*}
    where $t_k,t_{k+1}$ were defined in \eqref{e:eps-mesh-times}.
    Note that the above event depends on $\e$.

    There exists $C>0$ such that for all small $\d>0$ and $\e\le \delta^{100/\gamma}$,
    \begin{align*}
        \sup_u
        \max_{k\in\lb 1,n-1\rb}
        \P\left(\neg \MeshGeoDeloc^k_u\right) \le  C \e^{1-3\gamma},
    \end{align*}
    where the supremum is over all $u=(0,z;\d^{3/2},w)$ with $z,w\in J^\d$.
\end{lemm}
\begin{proof}
    For $u=(0,z;\d^{3/2},w)$, define the event
    \begin{align*}
        \sG_{u} \coloneqq  \left\{
            \left|
                \cL(0,\wh{z};\d^{3/2},\wh{w})-\wh{\cL}^\e(0,\wh{z};\d^{3/2},\wh{w})
            \right| \le \e^9
        \right\}.
    \end{align*}
    We have $\sup_u \P(\neg\sG_{u})\le Ce^{-c\log^2(1/\e)}$ by \cref{l:moments-tails-wh-L}. To prove the lemma, by a union bound it suffices to upper bound
    \begin{align*}
        \sup_{u}
        \max_{k\in\lb 1,n\rb}
        \P\left(
            \left|\Gamma_u(t_k)\right| \le \e\log^4(1/\e)
        \right).
    \end{align*}
    By the above upper bound for $\P(\neg\sG_u)$, we have
    \begin{align*}
        \sup_u\P\left(
            \left|\Gamma_u(t_k)\right| \le \e\log^4(1/\e)
        \right)
        &\le \sup_u\P\left(
            \left\{
                \left|\Gamma_u(t_k)\right| \le \e\log^4(1/\e)
            \right\}
        \cap 
        \sG_u
        \right)
        + Ce^{-c\log^2(1/\e)}.
    \end{align*}
    Fix $u=(0,z;\d^{3/2},w)$ and $k\in\lb 1,n\rb$, and assume that the event inside the probability on the RHS of the above display occurs.
    Then we have 
    \begin{align*}
        \cL(0,\wh{z};\d^{3/2},\wh{w}) - \e^9
        &\le
        \wh{\cL}^\e(0,\wh{z};\d^{3/2},\wh{w})\\
        &=
        \max_{\substack{|v| \le \e\log^4(1/\e),\\v\in \e^{50}\Z}}
        \left[
            \wh{\cL}^\e(0,\wh{z};t_k, v) + \wh{\cL}^\e(t_k,v;\d^{3/2},\wh{w})
        \right]\\
        &\le
        \sup_{|v| \le \e\log^4(1/\e)}
        \left[
            \cL(0,\wh{z};t_k, v) + \cL(t_k,v;\d^{3/2},\wh{w})
        \right],
    \end{align*}
    where in the second line, the last passage times are defined analogously to \cref{def:discrete-LPP} using mesh paths restricted to $[0,t_k]$ and $[t_k,\d^{3/2}]$.
    On the other hand, by the metric composition property of $\cL$ (\cref{p:landscape-symmetries}\ref{landscape-metric-composition}), we have
    \begin{align*}
        \cL(0,\wh{z};\d^{3/2},\wh{w}) \ge \sup_{|v|\le 1} \left[
            \cL(0,\wh{z}; t_k, v) + \cL(t_k,v;\d^{3/2},\wh{w})
        \right].
    \end{align*}
    The above two displays imply that
    \begin{align*}
        \sup_{|v|\le 1} \left[
            \cL(0,\wh{z};t_k, v) + \cL(t_k,v;\d^{3/2},\wh{w})
        \right]
        \le
        \sup_{|v|\le \e\log^4(1/\e)} \left[
            \cL(0,\wh{z};t_k, v) + \cL(t_k,v;\d^{3/2},\wh{w})
        \right]
        + \e^9.
    \end{align*}
    Since $\e \le \d^{100/\gamma}$,  by \cref{p:geodesic-deloc} (applied with $t_0 = \d^{3/2}$), the above occurs with probability $O(\e^{1-3\gamma})$.
    This completes the proof.
\end{proof}

We are ready to prove \cref{p:half-strip-ES}.
To start, we combine Lemmas \ref{l:uniform-bound-as} and \ref{l:mesh-geo-deloc} to bound the terms on the RHS of \eqref{e:efron-stein-noise} corresponding to $k\in\lb 1,n-1\rb$.
To aid in this, 
we define, for $k\in\lb 1,n-1\rb$, the event (analogously to \eqref{e:both-inside})
\begin{equation*}
    \begin{split}
        \BothDeloc_u^k &\coloneqq  
        \left\{(\bL^\near, \bL^\far, \bL^\bdry) \in \MeshGeoDeloc_u^k\right\}\\
        &\qquad \cap \left\{(\bL^\near, \bL^{\far,(k)}, \bL^{\bdry,(k)}) \in \MeshGeoDeloc_u^k\right\}\\
        &= \left\{\Gamma_u \textup{ does not use an entry of } \cL^\bdry_k\right\}
        \cap\left\{\Gamma_u^{(k)} \textup{ does not use an entry of } \wt{\cL}^\bdry_k\right\},
    \end{split}
\end{equation*}
where $\Gamma_u$ and $\Gamma_u^{(k)}$ are the mesh geodesics corresponding respectively to $(\bL^\near,\bL^\far,\bL^\bdry)$ and $(\bL^\near, \bL^{\far,(k)}, \bL^{\bdry,(k)})$.
Then recalling from \eqref{e:both-inside} the event $\BothIn^k_u$, we have
\begin{multline}\label{e:921}
    \E\left[
        \left|
            f_u(\bL^\near, \bL^{\far}, \bL^{\bdry})
            - f_u(\bL^\near, \bL^{\far, (k)}, \bL^{\bdry, (k)})
        \right|^2
    \right]\\
    \begin{aligned}
    &= 
    \E\left[
        \left|
            f_u(\bL^\near, \bL^{\far}, \bL^{\bdry})
            - f_u(\bL^\near, \bL^{\far, (k)}, \bL^{\bdry, (k)})
        \right|^2 \1_{\BothIn_u^k} \1_{\BothDeloc_u^k}
        \right]\\
    &+
    \E\left[
        \left|
            f_u(\bL^\near, \bL^{\far}, \bL^{\bdry})
            - f_u(\bL^\near, \bL^{\far, (k)}, \bL^{\bdry, (k)})
        \right|^2
       \1_{\BothIn_u^k} \1_{\neg \BothDeloc_u^k}
    \right]\\
    &+
    \E\left[
        \left|
            f_u(\bL^\near, \bL^{\far}, \bL^{\bdry})
            - f_u(\bL^\near, \bL^{\far, (k)}, \bL^{\bdry, (k)})
        \right|^2
        \1_{\neg \BothIn_u^k}
    \right].
    \end{aligned}
\end{multline}

We bound the above three terms separately. 
On $\BothIn_u^k\cap\BothDeloc_u^k$, we have 
\[
    f_u(\bL^\near, \bL^{\far}, \bL^{\bdry})
            = f_u(\bL^\near, \bL^{\far, (k)}, \bL^{\bdry, (k)}),
\]
since on the time interval $[t_k,t_{k+1}]$, 
both geodesics $\Gamma_u$ and $\Gamma_u^{(k)}$ only use values from $\cL^{\near}_k$.
So the first term on the RHS of \eqref{e:921} vanishes.

Next, on $\BothIn_u^k\cap\neg\BothDeloc_u^k$, Lemma \ref{l:uniform-bound-as} (viz. \eqref{e:Cresample-bound}) implies
\begin{align*}
    \left|f_u(\bL^\near, \bL^{\far}, \bL^{\bdry}) - f_u(\bL^\near, \bL^{\far, (k)}, \bL^{\bdry, (k)})\right|
    \le 2 \Cresample\,\e^{1/2}\log^2(1/\e).
\end{align*}
Then using that $\P(\neg \BothDeloc_u^k) \le C \e^{1-3\gamma}$ by \cref{l:mesh-geo-deloc} and a union bound, we obtain 
\begin{multline*}
    \E\left[
        \left|
            f_u(\bL^\near, \bL^{\far}, \bL^{\bdry})
            - f_u(\bL^\near, \bL^{\far, (k)}, \bL^{\bdry, (k)})
        \right|^2
       \1_{\BothIn_u^k} \1_{\neg \BothDeloc_u^k}
    \right]\\
    \begin{aligned}
        &\le 
        \E\left[
            \Cresample^2\, \1_{\neg \BothDeloc_u^k}
        \right]
        \cdot
        4\e\log^4(1/\e)\\
        &\le 
        \left(
            C\e^{1-3\gamma}
            \log^4(1/\e)
            + \E\left[
                \Cresample^2 \1_{\Cresample > \log^2(1/\e)}
            \right]
        \right)
        \cdot
        4\e\log^4(1/\e)\\
        \overset{\eqref{e:Cresample-prob}}&{\ls} \e^{2-3\gamma} 
        \log^{8}(1/\e),
    \end{aligned}
\end{multline*}

Finally, for the third term in \eqref{e:921}, we use Cauchy--Schwarz along with the the moment bounds from \cref{l:moments-wh-L} and the fact that $\P(\neg \BothIn_u^k)\le Ce^{-c\log^2(1/\d)}$ by \eqref{good1243}, thereby obtaining
\begin{align*}
    \E\left[
        \left|
            f_u(\bL^\near, \bL^{\far}, \bL^{\bdry})
            - f_u(\bL^\near, \bL^{\far, (k)}, \bL^{\bdry, (k)})
        \right|^2
       \1_{\neg \BothIn_u^k}
    \right]\ls e^{-c\log^2(1/\d)}.
\end{align*}             

Combining the above, we obtain a bound of $O\left(\e^{2-3\gamma} \log^{8}(1/\e) + e^{-c\log^2(1/\d)}\right)$ for each term on the RHS of \eqref{e:efron-stein-noise} corresponding to $k\in \lb 1,n-1\rb$.  
For the $k=0,n$ terms in \eqref{e:efron-stein-noise}, we simply use \eqref{e:uniform-pth-moment} to bound them by $O\left(\e^{2\gamma/3}\log^6(1/\e) + e^{-c\log^2(1/\d)}\right)$.
Putting things together, and recalling that $n\le \e^{-3/2}$, we obtain
\begin{align*}
    \E\left[\Var(\wh{\cL}^\e(u) \mid \bL^\near)\right]
    &\ls 
    \frac{1}{\e^{3/2}}\cdot\left[ \e^{2-3\gamma} \log^{8}(1/\e)+
    e^{-c\log^2(1/\d)}\right]
    + \left[\e^{2\gamma/3} \log^{6}(1/\e) + e^{-c\log^2(1/\d)}\right]
    \\
    &\ls \left(\e^{1/2-3\gamma} + \e^{2\gamma/3}\right)\log^8(1/\e)
    + \frac{1}{\e^{3/2}} e^{-c\log^2(1/\d)},
\end{align*}
uniformly in $u$.
Since $1/2 - 3\gamma > 0$, this finishes the proof of \cref{p:half-strip-ES}, and thus also of \cref{t:2d-noise}.
\qed

\vspace{.3in}
As mentioned after the statement of \cref{t:black-noise}, we end this section with a brief discussion on how a similar strategy may be employed to show that the SHF is a 3D noise, and hence a 3D black noise by the result of \cite{GT25} which showed it to be a 1D black noise in the temporal direction.

\subsection{Towards proving the critical 2D stochastic heat flow (SHF) is a 3D black noise}\label{ss:SHF-black-noise}
The SHF is the universal scaling limit of critical $2+1$-dimensional directed polymer models \cite{CSZ23}, and in light of its 1D black noise property and \cref{t:black-noise} it is natural to conjecture that the SHF is a 3D black noise.
By the same argument used to deduce \cref{t:black-noise} from \cref{t:2d-noise}, it suffices to prove that the SHF is a 3D noise.
While we do not pursue a full proof considering that the focus of the present article is $1+1$-dimensional models and that it is already of significant length, in this subsection we attempt to convey that the 3D noise property of SHF appears to be significantly \emph{easier} to prove than \cref{t:2d-noise}. 

As in the 2D case, only the 3D join property (counterpart of Property \ref{noise-join}) is not an immediate fact and this is what we discuss.
Consider two disjoint open 3D boxes that share a face $S$.
To prove the 3D join property, one must show that the SHF can be reconstructed given its data off of $S$.
If the boxes are temporally adjacent, then this follows from the SHF's semigroup property \cite{T24,CM26}.
The other case is where $S$ lies in a plane orthogonal to one of the spatial coordinate axes.
In this case, the desired reconstruction is equivalent to showing that ``resampling" the SHF on a small neighborhood of $S$ conditionally on the remainder has a negligible effect on the SHF. This in fact follows by dimensionality considerations. To avoid setting up the necessary language to deal with the singularity of the SHF, we will instead illustrate this with the pre-limiting example of the $2+1$-dimensional lattice polymer in the critical temperature scaling window.  
We start by setting up the required notation.

Let $\pi=(\pi(t))_{t\in\Zpos}$ be the simple symmetric random walk on $\Z^2$, let $P$ be its law, and let $E$ denote expectation with respect to $P$.
For $n\ge 1$, let $R_n \coloneqq  E^{\otimes 2}[|\{t\in\lb 0,n-1\rb : \pi_1(t) = \pi_2(t)\}|]$ be the expected number of collisions of two independent random walks $\pi_1,\pi_2$ up to time $n-1$.
Define $\beta_n>0$ through the relation
\begin{align*}
    e^{\beta_n^2} - 1 = \frac{1}{R_n},
\end{align*}
which is essentially the same as setting $\beta_n = \sqrt{\pi/\log n}$, up to lower order terms.
Let $\omega=(\omega(t,x))_{(t,x)\in\Zpos \times \Z^2}$ be a field of i.i.d. standard Gaussian variables.
For integers $0\le s < t$ and $x,y\in\Z^2$, the \emph{point-to-point} polymer partition function is 
\begin{align*}
    Z_n(s,x;t,y) \coloneqq  E\left[
        \exp\left(\beta_n \sum_{r=s}^{t-1} \omega(r,\pi(r)) - \frac{\beta_n^2(t-s)}{2}\right)
        \1_{\pi(t)=y}
        \,\middle|\, \pi(s)=x
    \right].
\end{align*}
Note that by reversibility of random walk, $Z_n(s,x;t,y)=Z_n(s,y;t,x)$.
The \emph{point-to-plane} partition function is
\begin{align*}
    Z_n(s,x;t,\ast) &\coloneqq  
    \sum_{y\in\Z^2} Z_n(s,x;t,y)\\
    &=E\left[
        \exp\left(\beta_n \sum_{r=s}^{t-1} \omega(r,\pi(r)) - \frac{\beta_n^2(t-s)}{2}\right)
        \,\middle|\, \pi(s)=x
    \right].
\end{align*}
Note that $\E[Z_n(s,x;t,\ast)]=1$.
We consider the length-$n$ point-to-plane polymer, with its starting point averaged over a diffusively-scaled ball.
Namely, let $B=B_n$ be the ball of radius $\sqrt{n}$ centered at the origin in $\Z^2$.
Define
\begin{align*}
    Z_n(0,B;n,\ast)
    \coloneqq  \frac{1}{|B|}\sum_{x\in B} Z_n(0,x;n,\ast).
\end{align*}
By \cite[Theorem 1.1]{CSZ23}, as $n\to\infty$, the average $Z_n(0,B;n,\ast)$ converges in distribution to the SHF averaged over the unit ball in $\R^2$ and is a unit-order quantity.
To demonstrate the feasibility of proving the 3D noise property for SHF, we now quickly estimate the effect on $Z_n(0,B;n,\ast)$ of resampling the disorder 
inside the plane $S\coloneqq \{(t,y)\in\Zpos\times\Z^2 : y_2=0\}$.
Let $\wt{Z}_n(0,B;n,\ast)$ be  obtained from $Z_n(0,B;n,\ast)$ by resampling $\omega(t,y)$ for every $(t,y)\in S$, i.e. replacing each of them by an i.i.d, copy $\wt\omega(t,y)$.
We claim that
\begin{align}\label{e:resampling-SHF}
    \E\left[
    \left |
        Z_n(0,B;n,\ast) - \wt{Z}_n(0,B;n,\ast)
    \right |
    \right]
    \ls \frac{1}{n^{1/4 - o(1)}}.
\end{align}
Since $Z_n(0,B;n,\ast)$ and hence also $\wt{Z}_n(0,B;n,\ast)$ are unit-order quantities, the above estimate implies that the influence of $S$ is negligible. 

The argument is based on a simple Efron--Stein bound.
Note however that the above bound is stated in terms of the $L^1$-norm (the reason for this will become clear shortly).

To prove \eqref{e:resampling-SHF}, it will be convenient to notice that by directedness and diffusivity of $\pi$, the partition function $Z_{n}(0,B;n,\ast)$ depends primarily on the disorder variables $\omega(t,y)$ such that $t\in\lb 0,n-1\rb$ and $|y| \lesssim \sqrt{n}$. To take advantage of this, we define the following truncated version of the partition function.
We start with the following set of good trajectories that do not deviate by more than $\sqrt n \log n$ from the origin in the first $n-1$ steps. That is, 
\[
    A_{n}\coloneqq \left\{|\pi(r)|\le \sqrt{n}\log n\quad \text{ for all } r\in\lb 0,n-1\rb\right\}.
\]
Then for $0 \le s < t\le n$ and $x,y\in\Z^2$, we define
\begin{align*}
    \overline{Z}_n(s,x;t,y) \coloneqq  E\left[
        \exp\left(\beta_n \sum_{r=s}^{t-1} \omega(r,\pi(r)) - \frac{\beta_n^2(t-s)}{2}\right)
        \1_{\pi \in A_{n}}\,\1_{\pi(t)=y}
        \,\middle|\, \pi(s)=x
    \right],
\end{align*}
as well as its point-to-plane counterpart
\begin{align*}
    \overline{Z}_n(s,x;t,\ast) &\coloneqq  
    \sum_{y\in\Z^2} \overline{Z}_n(s,x;t,y).
\end{align*}
We analogously define $\overline{Z}_n(0,B;n,\ast)$ to be
$\frac{1}{|B|}\sum_{x\in B} \overline{Z}_n(0,x;n,\ast)$. 
Similarly, let $\wt{\overline{Z}}_n(0,B;n,\ast)$ be the counterpart when the disorder on $S$ is resampled. 
Now, notice that 
\begin{align*}
    \E\left[
        \left|Z_n(0,B;n,\ast)-\overline{Z}_n(0,B;n,\ast)\right|
    \right] 
    &=
    P_B\left(
        \exists r \in \lb 0,n-1\rb : |\pi(r)| > \sqrt{n}\log n
    \right)\\
    &\le Ce^{-c\log^2 n}, 
\end{align*}
where $P_B$ denotes the law of $\pi$ with $\pi(0)$ distributed uniformly in $B$.
(An $L^2$ analogue of the above bound is unfortunately more complicated to prove, which is why we stated \eqref{e:resampling-SHF} in the $L^{1}$ form.)
By the triangle inequality, to prove \eqref{e:resampling-SHF}, it suffices to prove 
\begin{align*}
    \E\left[
    \left|
        \overline{Z}_n(0,B;n,\ast) - \wt{\overline{Z}}_n(0,B;n,\ast)
    \right|
    \right]
    \ls \frac{1}{n^{1/4 - o(1)}}.
\end{align*}
We will in fact prove, using the Efron--Stein inequality, that
\begin{align}\label{e:resampling-SHF2}
    \E\left[
        \left(
            \overline{Z}_n(0,B;n,\ast) - \wt{\overline{Z}}_n(0,B;n,\ast)
        \right)^2
    \right]
    \ls \frac{1}{n^{1/2 - o(1)}}.
\end{align}

For $(t,y)\in S$, let $\wt{\overline{Z}}_n^{(t,y)}(0,B;n,\ast)$ be the version of $\overline{Z}_n^{(t,y)}(0,B;n,\ast)$ obtained by resampling only  $\omega(t,y)$.
By the Efron--Stein inequality,
\begin{align}\label{efron32}
    \E\left[
        \left(
            \overline{Z}_n(0,B;n,\ast) - \wt{\overline{Z}}_n(0,B;n,\ast)
        \right)^2
    \right]
    &\le
 \sum_{t=0}^{n-1}
    \sum_{|y_1|\le \sqrt{n}\log n} 
    \E\left[
        \left(
            \overline{Z}_n(0,B;n,\ast) - \wt{\overline{Z}}_n^{(t,y)}(0,B;n,\ast)
        \right)^2
    \right],
\end{align}
where $y=(y_1,0)$. Note that above we used that $\overline{Z}_n(0,B;n,\ast)$ only depends on $\omega(t,y)$ with $0\le t\le n-1$ and $|y|\le \sqrt{n}\log n$.
We now estimate each term on the RHS.
Using the Chapman--Kolmogorov equation 
(which also holds for $\overline{Z}_n$)
\begin{align*}
    \overline{Z}_n(0,B;n,\ast)
    &= \frac{1}{|B|}\sum_{x\in B}
        \sum_{y\in\Z^2}
        \overline{Z}_n(0,x;t,y)\overline{Z}_n(t,y;n,\ast),
\end{align*}
we can write the increment on the RHS of \eqref{efron32} as
\begin{align*}
    \overline{Z}_n(0,B;n,\ast) - \wt{\overline{Z}}_n^{(t,y)}(0,B;n,\ast)
    &=\left(\frac{1}{|B|}\sum_{x\in B}
    \overline{Z}_n(0,x;t,y) \right)
    \left(
        \overline{Z}_n(t,y;n,\ast)
        -
        \wt{\overline{Z}}_n^{(t,y)}(t,y;n,\ast)
    \right).
\end{align*}

Since the two factors on the RHS are independent, we get
\begin{multline*}
    \E\left[
        \left(
            \overline{Z}_n(0,B;n,\ast) - \wt{\overline{Z}}_n^{(t,y)}(0,B;n,\ast)
        \right)^2
    \right]\\
    \begin{aligned}
    &=
    \E\left[
        \left(\frac{1}{|B|}\sum_{x\in B}
        \overline{Z}_n(0,x;t,y) \right)^2
    \right]\cdot
    \E\left[
        \left(
            \overline{Z}_n(t,y;n,\ast)
            -
            \wt{\overline{Z}}_n^{(t,y)}(t,y;n,\ast)
        \right)^2
    \right]\\
    &\le
   \frac{1}{|B|^2}
    \E\left[
        Z_n(0,0;t,\ast)^2
    \right]
    \cdot
    4\E\left[
        Z_n(t,y;n,\ast)^2
    \right] 
    \\
    &\ls 
    \frac{\log^2 n}{n^2},
    \end{aligned}
\end{multline*}
where in the first inequality we used that $0 \le \overline{Z}_n\le  Z_n$ pointwise and that $\sum_{x\in B}{Z}_n(0,x;t,y)\le \sum_{x\in \Z^2}
{Z}_n(0,x;t,y)$, and the latter, by reversibility of random walk and translation invariance, is equal in distribution to $Z_n(0,0;t,\ast)$.
Further, in the last line, we used that the point-to-plane partition function has second moment $O(\log n)$ by \cite[Proposition A.2]{CSZ19}. Note that above, the key reason why the influence of $\omega(t,y)$ is small is that a tiny fraction of diffusive paths pass through $(t,y)$, and hence they contribute insignificantly to $Z_n(0,0;n,\ast)$. Let us mention that computations of a similar flavor to the above ($L^1$ analogues) have also appeared in \cite{CD25} in the study of enhanced noise sensitivity for the SHF.
Plugging the above bound into \eqref{efron32} yields an influence bound of $O(\frac{\log^3 n}{n^{1/2}})$, which proves \eqref{e:resampling-SHF2}. 
\qed

\begin{rem}
    We end this discussion with an interesting contrast between the strong disorder and intermediate disorder polymers,
    using as examples the zero-temperature model of $1+1$-dimensional discrete LPP and the $2+1$-dimensional critical SHF.

    For LPP in dimension $1+1$, the length-$n$ geodesic, whose total energy is of order $n$ with fluctuations of order $n^{1/3}$, spends time $n^{1/3}$ on a given line orthogonal to the spatial axis in a window of width $n^{2/3}$.
    Thus each line's energetic contribution is the same as the overall energetic fluctuation, and hence cannot be ignored.
    To this end, the argument in Sections \ref{s:noise-conditional-variance-reduction}--\ref{s:mesh-step-2} shows that the contribution from a single line is concentrated, implying that the contributions before and after resampling cancel each other out to leading order.

    In contrast, in the intermediate disorder regime, the $2+1$-dimensional polymer partition function (corresponding to the SHF) is a unit-order quantity, and is essentially a function of all disorder variables in a cylinder of height $n$ and radius $n^{1/2}$.
    As \eqref{e:resampling-SHF2} shows, the total contribution to the squared energy by any hyperplane (one of $n^{1/2}$ many) is at most $O(\frac{1}{n^{1/2-o(1)}})$, which is negligible.
    A matching lower bound, up to sub-polynomial factors, can be obtained by a further Efron--Stein argument resampling the $n^{1/2}$-many hyperplanes one by one. 
    We also expect identical statements to hold for $1+1$-dimensional intermediate disorder polymer models.
    That Efron--Stein arguments yield almost sharp results in the intermediate disorder cases can be attributed to the fact that low-degree (at most logarithmic in the system size) terms dominate the corresponding chaos expansions.
\end{rem}

\section{Mutual singularity via Gaussian multiplicative chaos}\label{s:mutual-singularity}

In this section we prove \cref{t:pathsing}, which says that for $\beta_1\ne \beta_2$, the CDRP measures $\P^\xi_{\beta_1},\P^\xi_{\beta_2}$ are mutually singular almost surely (see \cref{ss:CDRP} for the definition of $\P^\xi_\beta$).
As mentioned in \cref{ss:idea}, the proof relies on the recent  interpretation in \cite{QRV25} of the CDRP measure as a Gaussian multiplicative chaos (GMC) measure on path space.
Here is an outline of the section:
\begin{itemize}
    \item In \cref{ss:fa-preliminaries} we record 
    basic functional-analytic inputs.
    \item In \cref{ss:gmc-construction} we construct the CDRP measure as a GMC, following \cite{QRV25}.
    \item In \cref{ss:energy-level-sets} we characterize the typical ``energy'' of a path sampled from the GMC measure, and deduce \cref{t:pathsing}.
    The analysis rests on the fact that exponentially tilting a Gaussian has the effect of shifting its mean, as captured by Girsanov's theorem.
\end{itemize}

\subsection{Functional analysis preliminaries}\label{ss:fa-preliminaries}

This subsection records elementary functional-analytic inputs needed for our later arguments.
The reader may find it helpful to just glance at \cref{l:nice-L2-basis} and then proceed to \cref{ss:gmc-construction}, returning here as needed.

For $A\subset \R^2$, we denote by $C_c(A)$ the set of continuous real-valued functions on $\R^2$ with compact support contained in $A$.
When $A$ is compact, we equip $C_c(A)$ with the usual supremum norm: $\lVert \varphi\rVert_{\infty} \coloneqq  \sup_{x\in \R^2} |\varphi(x)|$. 
We equip $C_c(\R^2)$ with its canonical LF-topology, i.e. the finest locally convex topology for which the inclusion maps $C_c([-k,k]^2)\hookrightarrow C_c(\R^2)$ are continuous for every integer $k\ge 1$.
Concretely, this means that a sequence $\{\varphi_n\}_{n\ge 1}\subset C_c(\R^2)$ converges to $\varphi\in C_c(\R^2)$ if and only if there exists a compact set $K\subset \R^2$ such that $\varphi,\varphi_1,\varphi_2,\dots$ are all supported in $K$, and $\lim_{n\to\infty}\lVert \varphi_n-\varphi\rVert_{\infty}=0$.
For further background on the topology of $C_c(\R^2)$, see \cite[Chapter 13, Example II]{Treves}.

The following lemma asserts the existence of an orthonormal basis of $L^2(\R^2)$ whose linear span is also dense in $C_c(\R^2)$.
This will be used at various points in our analysis (see especially the crucial technical result \cref{l:L2-norm-occupation-measure}).
\begin{lemm}\label{l:nice-L2-basis}
    There exists a sequence $\{\ee_i\}_{i\ge 1} \subset C_c(\R^2)$ of continuous compactly supported functions with the following two properties:
    \begin{enumerate}
    \item The linear span of $\{\ee_i\}_{i\ge 1}$ is dense in $C_c(\R^2)$.
    \item $\{\ee_i\}_{i\ge 1}$ is an orthonormal basis of $L^2(\R^2)$.
\end{enumerate}
\end{lemm}
\begin{proof}
    Let $\{\varphi_i\}_{i\ge 1}\subset C_c(\R^2)$ be a countable dense subset of $C_c(\R^2)$
    (to construct such a subset, use Stone--Weierstrass to select a countable dense subset of $C_c([-n,n]^2)$ for each $n\ge 1$ and form the union).
    Let $S$ be the linear span of $\{\varphi_i\}_{i\ge 1}$.
    Then $S$ is dense in $C_c(\R^2)$, since $\{\varphi_i\}_{i\ge 1}\subset S$.
    Applying the Gram--Schmidt process to $\{\varphi_i\}_{i\ge 1}$ yields an $L^2(\R^2)$-orthonormal sequence $\{\ee_i\}_{i\ge 1}$ whose span is $S$.
    Since each $\varphi_i$ is continuous and compactly supported, the same is true of the $\ee_i$.

    It remains to show that $\{\ee_i\}_{i\ge 1}$ is an orthonormal basis of $L^2(\R^2)$, or equivalently that $S$ is dense in $L^2(\R^2)$.
    Fix $f\in L^2(\R^2)$ and $\varepsilon>0$.
    Since $C_c(\R^2)$ is dense in $L^2(\R^2)$, there exists $\psi\in C_c(\R^2)$ such that $\lVert f-\psi\rVert_{L^2}<\varepsilon$.
    Choose a sequence $\{\psi_n\}_{n\ge 1}\subset S$ converging to $\psi$ in $C_c(\R^2)$.
    By definition, this means there exists a compact set $K$ containing the supports of $\psi,\psi_1,\psi_2,\dots$, and  $\lim_{n\to\infty}\lVert \psi-\psi_n\rVert_\infty =0$.
    Choose $n$ large enough that $\lVert \psi-\psi_n\rVert_\infty < \varepsilon/\sqrt{|K|}$.
    Then we have
    \begin{align*}
        \lVert \psi-\psi_n\rVert_{L^2}
        \le 
        \lVert \psi-\psi_n\rVert_\infty
        \sqrt{|K|}
        < \varepsilon,
    \end{align*}
    hence
    \begin{align*}
        \lVert f-\psi_n\rVert_{L^2}
        \le \lVert f-\psi\rVert_{L^2}
        + \lVert \psi-\psi_n\rVert_{L^2}
        < 2\varepsilon.
    \end{align*}
    This proves that $S$ is dense in $L^2(\R^2)$.
\end{proof}

\subsection{GMC measure construction}\label{ss:gmc-construction}

From now on we fix an orthonormal basis $\{\ee_i\}_{i\ge 1}$ of $L^2(\R^2)$ whose linear span is dense in $C_c(\R^2)$, as provided by \cref{l:nice-L2-basis}.
We also fix $(s,x;t,y)\in\Rup$ and suppress it from the notation.

Let $X=\{X(r)\}_{r\in [s,t]}$ be a Brownian bridge from $X(s)=x$ to $X(t)=y$,
and let $\gamma^X$ be the occupation measure of its graph:
\begin{align}\label{e:occupation-measure}
    \gamma^X(A) \coloneqq  \int_s^t \1_{(r,X(r))\in A}\,dr,
    \qquad A\in \cB(\R^2).
\end{align}
Note that $\gamma^X$ is singular with respect to two-dimensional Lebesgue measure, since by Fubini the graph of any $X\in C([s,t])$ has zero Lebesgue measure.
We denote
\begin{align*}
    \gamma_i^X \coloneqq  \int_{\R^2} \ee_i\,d\gamma^X
    = \int_s^t \ee_i(r,X(r))\,dr.
\end{align*}
Since each $\ee_i$ is continuous and compactly supported, it follows that $X\mapsto \gamma_i^X$ is a continuous function on path space $C([s,t])$.
This continuity will play a key technical role in our proof of \cref{t:pathsing} (see the proof of \cref{p:energy-level-set}, especially around \eqref{e:portmanteau-consequence}).

Let $\xi$ be  space-time white noise on $\R^2$ (\cref{def:white-noise}), and set
\begin{align*}
    \xi_i \coloneqq  \xi(\ee_i)
    \qquad\text{for } i\ge 1.
\end{align*}
Then $\xi_1,\xi_2,\dots$ are i.i.d. standard Gaussian random variables, since  $\{\ee_i\}_{i\ge 1}$ is an orthonormal basis of $L^2(\R^2)$.
Let $\Pbr=\mrm{P}_{\mrm{BB},(s,x),(t,y)}$ be the Brownian bridge measure on path space $C([s,t])$, and let $\Ebr$ be the corresponding expectation.
For $\beta \ge 0$ we consider prelimiting \emph{un-normalized GMC measures} given by exponential tilts of $\Pbr$: for $n\ge 1$,
\begin{align}\label{e:unnormalized-GMC-prelimit}
    \frac{dM_{\beta,n}^\xi}{d\Pbr}(X) \coloneqq  
    \exp\left(
        \beta \sum_{i=1}^n  \xi_i \gamma_i^X
        - \frac{\beta^2}{2}\sum_{i=1}^n (\gamma_i^X)^2
    \right),
    \qquad X \in C([s,t]).
\end{align}
Note that for fixed $X\in C([s,t])$, the sum $\sum_{i=1}^n \xi_i \gamma_i^X$ is centered Gaussian with variance $\sum_{i=1}^n (\gamma_i^X)^2$, and hence $\E[\frac{d M_{\beta,n}^\xi}{d\Pbr}(X)] = 1$.
It follows by Fubini that the \emph{annealed un-normalized measure} is $\E[M_{\beta,n}^\xi] = \Pbr$.  

The partition function is defined as
\begin{align*}
    Z^{\xi}_{\beta,n} &\coloneqq  M^{\xi}_{\beta,n}(C([s,t]))\\
    &= \int_{C([s,t])} 
    \exp\left(
        \beta \sum_{i=1}^n  \xi_i \gamma_i^X
        - \frac{\beta^2}{2}\sum_{i=1}^n (\gamma_i^X)^2
    \right)
    \,d\Pbr(X).
\end{align*}
Since $\E[M^\xi_{\beta,n}]=\Pbr$, we have $\E[Z^{\xi}_{\beta,n}]=1$,
and hence $Z^\xi_{\beta,n}<\infty$ a.s.
We also clearly have $Z^\xi_{\beta,n}>0$ a.s., and hence
we can define the prelimiting \emph{normalized GMC measures}
\begin{align}\label{e:normalized-GMC-prelimit}
    \frac{1}{Z_{\beta,n}^\xi} M^\xi_{\beta,n}.
\end{align}

Next, following \cite{QRV25}, we will show that, for $\P$-a.e. realization of the white noise $\xi$, the probability measure $\frac{1}{Z_{\beta,n}^\xi} M^\xi_{\beta,n}$ converges weakly to a ($\xi$-dependent) probability measure on $C([s,t])$ as $n\to\infty$.
The limiting measure will be exactly the CDRP measure defined in \cref{ss:CDRP}, and the limiting partition function $Z_\beta^\xi$ will coincide with $p(t-s, y-x)^{-1}\cZ_\beta^\xi(s,x;t,y)$ as in \eqref{wienertilt}, where $\cZ_\beta^\xi$ is the solution to \eqref{e:she-intro} and $p$ is the heat kernel.

The weak convergence of $\frac{1}{Z^\xi_{\beta,n}}M^\xi_{\beta,n}$ will be proved by combining a certain martingale property of the un-normalized measures $M^\xi_{\beta,n}$ with a uniform integrability statement (in fact $L^2$-boundedness).
We begin with the latter:

\begin{lemm}[$L^2$-boundedness, {\cite{QRV25}}]\label{l:L2-bound}
    For any $\beta\ge0$ and any bounded continuous function $f : C([s,t])\to \R$, we have
    \begin{align*}
        \sup_{n\ge 1}\E\left[\left|\int_{C([s,t])}f(X)\,dM^\xi_{\beta,n}(X)\right|^2\right]
        <\infty.
    \end{align*}
\end{lemm}
\begin{proof}
    Since $f$ is bounded, it suffices to bound the second moment of the partition function $Z^\xi_{\beta,n} = M^\xi_{\beta,n}(C([s,t]))$.
    This was done in \cite[Section 7.2]{QRV25} (in fact they bound the $p\th$ moment for every $p\ge 1$), but we sketch the argument for the reader's convenience.
    A straightforward calculation (see \cite[Proposition 39]{QRV25}) yields
    \begin{align*}
        \E\left[
            (Z^\xi_{\beta,n})^2
        \right]
        &= \Ebr^{\otimes 2}\left[
            \exp\left(
                \beta^2 \sum_{i=1}^n \gamma_i^X \gamma_i^{X'}
            \right)
        \right],
    \end{align*}
    where $X,X'$ are independent Brownian bridges.
    \cite{QRV25} proves the convergence
    \begin{align}\label{e:mutual-intersection}
        \Ebr^{\otimes 2}\left[
            \exp\left(
                \beta^2 \sum_{i=1}^n \gamma_i^X \gamma_i^{X'}
            \right)
        \right]
        \xrightarrow{n\to\infty}
        \Ebr^{\otimes 2}\left[
            \exp\left(
                \beta^2 \sum_{i=1}^\infty \gamma_i^X \gamma_i^{X'}
            \right)
        \right]
        < \infty,
    \end{align}
    and in fact the infinite sum on the RHS is the mutual intersection local time $L$ of the processes
    $\{(r,X(r))\}_{r\in[s,t]}$ and $\{(r,X'(r))\}_{r\in[s,t]}$ (see \cite[Lemma 50, Proposition 74, and Remark 75]{QRV25}).
    To see heuristically why $L$ has all exponential moments, note that
    \[
        L = \int_s^t\int_s^t \delta_0(r-r')\delta_0(X(r)-X'(r'))\,dr\,dr'
        = \int_s^t \delta_0(X(r)-X'(r))\,dr.
    \]
    Since $X,X'$ are independent, the RHS above is just the local time in $0$ of the Brownian bridge $r\mapsto X(r)-X'(r)$, which is well-known to have a Gaussian tail.
\end{proof}

As mentioned, the other ingredient besides \cref{l:L2-bound} for the weak convergence of the normalized GMC measures is the martingale structure of the un-normalized measures, as stated below.

\begin{lemm}[Martingale convergence]\label{l:martingale-convergence}
    For any $\beta\ge 0$ and any bounded continuous function $f:C([s,t])\to\R$, the sequence
    \[
        \int_{C([s,t])}f(X)\,dM^\xi_{\beta,n}(X),\qquad n\ge 1
    \]
    forms a martingale with respect to the filtration $\mathcal{F}_n\coloneqq \sigma(\xi_1,\dots,\xi_n)$.
    It is bounded in $L^2(\mathbb{P})$ by \cref{l:L2-bound}, and hence converges a.s. and in $L^2(\mathbb{P})$.
\end{lemm}
\begin{proof}
    For any $m>n\ge 1$, by \eqref{e:unnormalized-GMC-prelimit} and Fubini, we have
    \begin{align*}
        \E\left[
            \int_{C([s,t])}f(X)\,dM^\xi_{\beta,m}(X)
            \middle|\cF_n
        \right]
        &=
        \E\Biggl[
            \int_{C([s,t])}
                f(X)\exp\left(
                    \beta \sum_{i=1}^m  \xi_i \gamma_i^X
                - \frac{\beta^2}{2}\sum_{i=1}^m (\gamma_i^X)^2
                \right)
                d\Pbr(X)
            \Bigg|\cF_n
        \Biggr]\\
        &= \int_{C([s,t])}
        f(X)
        \E\left[
            \exp\left(
                \beta \sum_{i=1}^m  \xi_i \gamma_i^X
            - \frac{\beta^2}{2}\sum_{i=1}^m (\gamma_i^X)^2
            \right)
            \middle| \cF_n
        \right]
        d\Pbr(X)\\
        &=\int_{C([s,t])}
        f(X)
            \exp\left(
                \beta \sum_{i=1}^n  \xi_i \gamma_i^X
                - \frac{\beta^2}{2}\sum_{i=1}^n (\gamma_i^X)^2
            \right)
        d\Pbr(X)\\
        &=
        \int_{C([s,t])}f(X)\,dM^\xi_{\beta,n}(X),
    \end{align*}
    where the third equality is because $\xi_1,\xi_2,\dots$ are i.i.d. standard Gaussian random variables.
\end{proof}

Using the above two lemmas, we now deduce the a.s. weak convergence of GMC measures using a general abstract argument reminiscent of \cite[Section 6]{Berestycki17}.

\begin{lemm}[Quenched weak convergence of GMC measures]\label{l:quenched-weak-convergence}
    Almost surely, the un-normalized measures $M_{\beta,n}^\xi$ converge weakly to a finite measure $M_{\beta}^\xi$.
    The limiting annealed un-normalized measure is $\E[M_\beta^\xi] = \Pbr$.

    In fact, $M_{\beta}^\xi$ is nonzero almost surely, so the normalized measures $\frac{1}{Z_{\beta,n}^\xi}M_{\beta,n}^\xi$ also converge weakly to the probability measure $\frac{1}{Z_{\beta}^\xi}M_{\beta}^\xi$, where $Z_{\beta}^\xi \coloneqq  M_{\beta}^\xi(C([s,t]))$.
\end{lemm}
\begin{proof}
    By combining \cref{l:martingale-convergence} with the fact that $\E[M_{\beta,n}^\xi]=\Pbr$ and the tightness of $\Pbr$, 
    one can show that $\{M_{\beta,n}^\xi\}_{n\ge 1}$ forms a tight sequence almost surely (e.g. using Doob's inequality and Borel--Cantelli).
    Since $C([s,t])$ is separable, there exists a sequence $\{f_k\}_{k\ge 1}$ of continuous bounded functions on $C([s,t])$ that separates finite Borel measures on $C([s,t])$.
    \footnote{That is, for finite Borel measures $\mu$ and $\nu$, we have $\mu=\nu$ iff $\int f_k\,d\mu = \int f_k\, d\nu$ for all $k$.}
    By \cref{l:martingale-convergence}, almost surely, the integrals $\int f_k \,dM^\xi_{\beta,n}$ converge as $n\to\infty$ for every $k\ge 1$ simultaneously.
    Therefore, $\{M_{\beta,n}^\xi\}_{n\ge 1}$ has at most one weak limit point a.s., and hence converges weakly a.s.

    By \cite[Lemma 43]{QRV25}, the strict positivity of the limiting partition function $Z_\beta^\xi$ is subject to a zero-one law: $\P(Z_\beta^\xi > 0)\in\{0,1\}$.
    By \cref{l:martingale-convergence} we have $Z_{\beta,n}^\xi\to Z_{\beta}^\xi$ in $L^2(\P)$,  so in particular $\E[Z_{\beta}^\xi]=1$ and hence $\P(Z_\beta^\xi > 0)=1$.
    A similar argument shows that $\E[M_\beta^\xi]=\Pbr$.
\end{proof}

As promised, the normalized GMC measure constructed above coincides with the CDRP measure:
\begin{lemm}[Coincidence with CDRP, {\cite[Theorem 4]{QRV25}}]\label{l:coincidence-with-CDRP}
    For any $\beta\ge0$, we have
    \begin{align*}
        \frac{1}{Z_\beta^\xi} M_\beta^\xi = 
        \P^\xi_\beta
        \qquad\text{a.s.}
    \end{align*}
    where $M_\beta^\xi,Z_\beta^\xi$ are from \cref{l:quenched-weak-convergence}, and $\P^\xi_\beta = \P^\xi_{\beta,(s,x),(t,y)}$ is the CDRP measure from \cref{t:AJRAS-polymer-existence}.
\end{lemm}
\begin{proof}
    For $\beta=0$ this is trivial, as $\frac{1}{Z^\xi_0}M_0^\xi = \Pbr = \P^\xi_0$.
    For $\beta>0$, this was shown in \cite[Theorem 4]{QRV25} by proving a.s. coincidence of finite-dimensional distributions (see also \cite[Lemma 87]{QRV25}).
\end{proof}

Having expressed the CDRP as a GMC measure, we turn now to analyzing its support.

\subsection{Energy level sets}\label{ss:energy-level-sets}

The upcoming \cref{p:energy-level-set} identifies the typical ``energy'' of the CDRP $X\sim \P^\xi_\beta$, relying on the description of the CDRP measure $\P^\xi_\beta$ as the a.s. weak limit of $\frac{1}{Z^\xi_{\beta,n}}M^\xi_{\beta,n}$ provided by \cref{l:quenched-weak-convergence,l:coincidence-with-CDRP}.
Our argument closely resembles an argument in \cite[Proposition 3.4]{DS11} showing that the $\gamma$-Liouville quantum gravity measure concentrates on $\gamma$-thick points of the Gaussian free field \cite{HMP10}.

The analysis hinges on the following measure-theoretic lemma, which is one of the main reasons we require our $L^2(\R^2)$-basis to satisfy the properties stated in \cref{l:nice-L2-basis}.

\begin{lemm}[CDRP Hamiltonian has infinite variance pointwise]\label{l:L2-norm-occupation-measure}
    Let $\{\ee_i\}_{i\ge 1}\subset C_c(\R^2)$ be an orthonormal basis of $L^2(\R^2)$ as in Lemma \ref{l:nice-L2-basis}.
    Let $\mu$ be a Radon measure on $\R^2$, and for $i\ge 1$ let $\mu_i\coloneqq  \int_{\R^2} \ee_i\,d\mu$ (the integral exists because $\ee_i$ is continuous and compactly supported).
    If $\sum_{i=1}^\infty \mu_i^2 < \infty$, then $\mu$ is absolutely continuous with respect to Lebesgue measure.

    Thus, since for any $X\in C([s,t])$ the occupation measure $\gamma^X$ in \eqref{e:occupation-measure} is singular with respect to Lebesgue measure,  we have
    \begin{align*}
        \lim_{n\to\infty} \VarXi\left(\sum_{i=1}^n \xi_i \gamma_i^X\right)
        =
        \sum_{i=1}^\infty (\gamma_i^X)^2
        =\infty.
    \end{align*}
\end{lemm}
\begin{proof}
    Since $\{\ee_i\}_{i\ge 1}$ is an orthonormal basis of $L^2(\R^2)$, the sum $\psi\coloneqq \sum_{i=1}^\infty \mu_i \ee_i$ converges in $L^2(\R^2)$.
    We will show that $\mu = \psi\,dx$.

    Let $S$ be the linear span of $\{\ee_i\}_{i\ge 1}$.
    For any $\varphi = \sum_{i=1}^n c_i \ee_i\in S$, we have
    \begin{align*}
        \int_{\R^2} \varphi\,d\mu
        = \sum_{i=1}^n c_i \mu_i
        &= \langle \varphi, \psi \rangle_{L^2(\R^2)}
        = \int_{\R^2} \varphi \psi\,dx.
    \end{align*}
    Since $\mu$ is Radon, $\varphi\mapsto \int_{\R^2} \varphi\,d\mu$ is a continuous linear functional on $C_c(\R^2)$,
    and by Cauchy--Schwarz the same is true for $\varphi\mapsto \int_{\R^2}\varphi \psi\,dx$.
    Then since $S$ is dense in $C_c(\R^2)$ by \cref{l:nice-L2-basis}, we get
    \begin{align*}
        \int_{\R^2} \varphi\,d\mu
        = \int_{\R^2} \varphi \psi\,dx
        \qquad\text{for all } \varphi\in C_c(\R^2).
    \end{align*}
    The LHS is a positive linear functional on $C_c(\R^2)$, so the RHS is too.
    By the Riesz--Markov--Kakutani representation theorem (e.g. \cite[Theorem 7.2]{Foll99}), we conclude that $\mu = \psi\,dx$.
\end{proof}

We also record the following elementary property of Gaussian random variables, which follows from a direct calculation (or Girsanov's theorem).
The motivation for this is explained right after.
\begin{lemm}\label{l:girsanov}
    Let $H \sim N(0,\sigma^2)$ be a centered Gaussian random variable of variance $\sigma^2$.
    For any Borel set $S\in\cB(\R)$,
    \begin{align*}
        \E\left[\exp\left(H - \frac{\sigma^2}{2}\right)\1_{H \in S}\right]
        = \P(H + \sigma^2 \in S).
    \end{align*}
    In other words, the law of $H$ under the tilted measure $\exp(H-\frac{\sigma^2}{2})\,d\P(H)$ is $N(\sigma^2,\sigma^2)$.
\end{lemm}

For us, the upshot of \cref{l:girsanov} is that if the variance $\sigma^2$ of $H$ diverges to $\infty$, then under the tilted measure, $H$ concentrates around $\sigma^2$, because $|H-\sigma^2| \ls \sigma \ll \sigma^2$.
In particular, since by \cref{l:L2-norm-occupation-measure} the variance of the approximate CDRP Hamiltonian $H^\xi_n(X) \coloneqq  \sum_{i=1}^n \xi_i \gamma_i^X$ diverges as $n\to\infty$, 
this should imply that for $X$ sampled from the CDRP measure, we have the concentration
\[
    \beta H^\xi_n(X) = \beta^2 \VarXi(H^\xi_n(X))(1+o(1)),
    \qquad\text{or equivalently}\qquad
    \frac{\sum_{i=1}^n \xi_i \gamma_i^X}{\sum_{i=1}^n (\gamma_i^X)^2} = \beta(1+o(1)).
\]
This is indeed the content of our aforementioned energy level set result.
As mentioned, the above reasoning goes back to \cite{DS11} in the context of Liouville quantum gravity.

\begin{prop}[Energy level set]\label{p:energy-level-set}
    There exists a deterministic sequence $\{n_k\}_{k\ge 1}\subset \Z_{\ge 1}$ with $n_k\uparrow\infty$ such that for any $\beta\ge 0$, it holds almost surely that
    \begin{align*}
        \P_\beta^\xi\left(
            \lim_{k\to\infty}
            \frac{\sum_{i=1}^{n_k} \xi_i \gamma^X_i}{\sum_{i=1}^{n_k} (\gamma^X_i)^2} = \beta
        \right)
        = 1.
    \end{align*}
    In particular, for any $\beta_1\ne \beta_2$, the CDRP measures $\P_{\beta_1}^\xi$ and $\P_{\beta_2}^\xi$ are mutually singular almost surely, as claimed in \cref{t:pathsing}.
\end{prop}

\begin{rem}
    We restricted to a subsequence $n_k$ in the above result only to simplify the proof; it is straightforward to refine the argument to prove that $\P^\xi_\beta\bigl(\lim_{n\to\infty} \frac{\sum_{i=1}^n \xi_i \gamma_i^X}{\sum_{i=1}^n (\gamma_i^X)^2} = \beta\bigr)=1$ almost surely.
    However the stated subsequential limit result already implies mutual singularity, so for simplicity we do not pursue this.
\end{rem}

\begin{proof}[Proof of \cref{p:energy-level-set}]
    By \cref{l:coincidence-with-CDRP}, the proposition is equivalent to 
    \begin{align*}
        \frac{1}{Z^\xi_\beta}
        M^\xi_\beta\left(\lim_{k\to\infty}
        \frac{\sum_{i=1}^{n_k} \xi_i \gamma^X_i}{\sum_{i=1}^{n_k} (\gamma^X_i)^2} = \beta\right) = 1,
    \end{align*}
    and this is what we will prove.

    For $n\ge 1$ let
    \begin{align*}
        H^\xi_n(X) \coloneqq  \sum_{i=1}^n \xi_i \gamma_i^X,\qquad \xi\in\Xi, \quad X\in C([s,t]),
    \end{align*}
    where recall from \eqref{e:def-Xi} that $\Xi=\R^{L^2(\R^2)}$ is the space where the white noise $\xi$ takes values.
    Note that for fixed $X\in C([s,t])$, the random variable $\xi\mapsto H^\xi_n(X)$ is centered Gaussian with variance 
    \begin{align*}
        \sigma_n^2(X) \coloneqq  \VarXi(H^\xi_n(X))
        =\sum_{i=1}^n (\gamma_i^X)^2,
    \end{align*}
    with the convention that a centered Gaussian of variance zero is identically $0$.
    By \cref{l:L2-norm-occupation-measure}, there exists a sequence $n_k\to\infty$ such that
    \begin{align*}
        \Pbr\left(\left\{X\in C([s,t]):\sigma_{n_k}^2(X) \ge 2^k\right\}\right) \ge 1-2^{-k}.
    \end{align*}
    Therefore, by Borel--Cantelli, the set of paths
    \begin{align*}
        B \coloneqq  \left\{
            X\in C([s,t]) : 
            \sigma_{n_k}^2(X) \ge 2^k
            \textup{ for all sufficiently large $k$}
        \right\}
    \end{align*}
    has $\Pbr(B)=1$.
    By \cref{l:quenched-weak-convergence} we have $\E[M_\beta^\xi(\neg B)] = \Pbr(\neg B) = 0$,
    so by non-negativity $M_\beta^\xi(\neg B)=0$ a.s.

    We fix a realization of the white noise $\xi$ and define the following measurable subset of path space:
    \begin{align*}
        A_{\beta,k}^\xi \coloneqq  \left\{X \in C([s,t]) : 
            \left|
                H^\xi_{n_k}(X) - \beta \sigma_{n_k}^2(X)
            \right|
            \le k
            \sigma_{n_k}(X)
        \right\}.
    \end{align*}
    Fix $k\ge 1$.
    Since $X\mapsto\gamma_i^X$ is a continuous function on $C([s,t])$, the functions $H^\xi_{n_k}$ and $\sigma_{n_k}^2$ are also continuous, so $A_{\beta,k}^\xi$ is a closed subset of $C([s,t])$.
    Therefore, by \cref{l:quenched-weak-convergence} and the Portmanteau lemma, it holds $\P$-a.s. that
    \begin{align}\label{e:portmanteau-consequence}
        M^\xi_\beta(\neg A_{\beta,k}^\xi)
        \le 
        \liminf_{n\to\infty} 
        M^\xi_{\beta,n}(\neg A_{\beta,k}^\xi).
    \end{align}
    It is straightforward to show that $\xi\mapsto M^\xi_\beta(\neg A_{\beta,k}^\xi)$ and $\xi\mapsto M^\xi_{\beta,n}(\neg A_{\beta,k}^\xi)$ are measurable functions of the white noise.
    Averaging over $\xi$, applying Fatou's lemma, and then applying Fubini, we get
    \begin{multline*}
        \E\Bigl[
            M^\xi_\beta(\neg A_{\beta,k}^\xi)
        \Bigr]
        \le \liminf_{n\to\infty}
        \E\Bigl[
            M^\xi_{\beta,n}(\neg A_{\beta,k}^\xi)
        \Bigr]\\
        \overset{\eqref{e:unnormalized-GMC-prelimit}}{=} \liminf_{n\to\infty}
        \int_{C([s,t])}
        \E\left[
            \exp\left(\beta H_{n}^\xi(X) - \frac{\beta^2 \sigma_{n}^2(X)}{2}\right)
            \1_{|H_{n_k}^\xi(X) - \beta \sigma_{n_k}^2(X)| > k\sigma_{n_k}(X)}
        \right]
        d\Pbr(X)\\
        =
        \int_{C([s,t])}
        \E\left[
            \exp\left(\beta H_{n_k}^\xi(X) - \frac{\beta^2 \sigma_{n_k}^2(X)}{2}\right)
            \1_{|H_{n_k}^\xi(X) - \beta \sigma_{n_k}^2(X)| > k\sigma_{n_k}(X)}
        \right]
        d\Pbr(X),
    \end{multline*}
    where the last line is because $H_{n}^\xi(X) - H_{n_k}^\xi(X)$ is independent of $H_{n_k}^\xi(X)$.
    Continuing, using that $H_{n_k}^\xi(X) \sim N(0,\sigma_{n_k}^2(X))$, \cref{l:girsanov} implies that the RHS is equal to
    \begin{align*}
        \int_{C([s,t])}
        \P\Bigl(
            \left|N(0, \sigma_{n_k}^2(X))\right| > k\sigma_{n_k}(X)
        \Bigr)\,
        d\Pbr(X),
    \end{align*}
    and therefore
    \begin{align*}
        \E\left[
            M^\xi_\beta(\neg A_{\beta,k}^\xi)
        \right]
        \le 2e^{-k^2/2}.
    \end{align*}
    So by Markov's inequality,
    \begin{align*}
        \P\left(
            M^\xi_\beta(\neg A_{\beta,k}^\xi) > e^{-k^2/4}
        \right)
        \le 2e^{-k^2/4}.
    \end{align*}
    Now defining the event
    \begin{align*}
        A^\xi_\beta \coloneqq  
        \liminf_{k\to\infty}A^\xi_{\beta,k}
        =
        \left\{X:\left|H^{\xi}_{n_k}(X) - \beta \sigma_{n_k}^2(X)\right|
            \le k \sigma_{n_k}(X)
            \textup{ for all sufficiently large $k$}\right\}
    \end{align*}
    and applying Borel--Cantelli twice, we get that $M^\xi_\beta(\neg A^\xi_\beta)=0$ almost surely.
    Since $M^\xi_\beta(\neg B)=0$, it follows that
    \begin{align*}
        \frac{1}{Z^\xi_\beta} M^\xi_\beta(A^\xi_\beta \cap B)
        =1.
    \end{align*}
    Finally, for $X\in A^\xi_\beta\cap B$, since $\sigma^2_{n_k}(X) \ge 2^k$ eventually, we have
    \begin{align*}
        \limsup_{k\to\infty} \left|
            \frac{H_{n_k}^\xi(X)}{\sigma_{n_k}^2(X)} - \beta
        \right|
        \le \limsup_{k\to\infty} \frac{k}{2^{k/2}}
        =0
    \end{align*}
    as desired.
\end{proof}

\section{Failure of a spectral approach, and the onset of chaos via coarse-graining}
\label{s:heuristic}

In this section we discuss coarse-graining schemes for the CDRP.
This perspective finds two applications: a heuristic justification of Conjecture \ref{c:main}, and a rigorous proof of energetic stability in part of the sub-critical regime.
However, before getting into coarse-graining, for context we first discuss a different natural spectral approach.
\\

As alluded to in \cref{ss:preface}, there is a natural correspondence between noise sensitivity properties of a function and its expressibility in terms of low-degree polynomials (see also \cref{r:informal-def-black-noise}).
This naturally forges a connection between the study of such properties and harmonic analysis.
For instance, for a square-integrable functional $f$ of independent $\mathrm{Bernoulli}(1/2)$ bits $X_1,X_2,\dots, X_n$,
the fraction of mass in low-frequency modes in its Fourier--Walsh expansion measures ``how much'' of the observable can be described using low-degree polynomials.
Going further, if the bits $X_1,\dots,X_n$ are each resampled independently at rate $1$,
the correlation between the values of $f$ at times $0$ and $t$ can be represented as an exponential moment of the \emph{spectral distribution} of $f$, which is the random variable that takes value $k\in\Zpos$ with probability proportional to the sum of squares of all Fourier--Walsh coefficients of $f$ corresponding to modes of frequency $k$.
For background on Fourier--Walsh
spectra and noise sensitivity, see \cite{GS}.

The above notions have exact analogues in the Gaussian setting:
the Fourier--Walsh expansion corresponds to the Wiener chaos decomposition,
the spectral distribution is defined using squared $L^2$-norms of projections onto homogeneous Wiener chaoses,
and the resampling dynamics corresponds to Ornstein--Uhlenbeck (OU) dynamics.
The relationship between the spectral distribution and decorrelation under OU dynamics is identical to the Bernoulli case.
Thus, a natural first approach to studying chaotic properties of a functional of Gaussian white noise is to study its spectral distribution.

\subsection{Spectral distribution of the polymer partition function}\label{s:spectral-failure}

Gu and Komorowski \cite{GK25} recently initiated the study of the spectral distribution of the stochastic heat equation with flat initial data, viewed as a functional of Gaussian space-time white noise $\xi$.
Equivalently, this is the point-to-line CDRP partition function 
$Z_\beta \coloneqq  \int_\R \cZ^\xi_\beta(0,0;1,x)\,dx$,
where $\cZ^\xi_\beta(0,0;1,x)$ is defined in \cref{def:intro-SHE}.
Writing $S_\beta$ for the spectral distribution of $Z_\beta$ (as defined above),
they proved the following central limit theorem:
\begin{align}\label{e:gk-clt}
    \frac{S_\beta - \frac{1}{2}\beta^4}{\beta^2} \xrightarrow[\beta\to\infty]{d} N(0,1).
\end{align}
(In fact, $S_\beta$ is approximately a Poisson random variable with mean $\frac12\beta^4$.)
In particular, as $\beta\to\infty$, the dominant contribution to the second moment of $Z_\beta$ comes from Wiener chaoses of degree $\frac12\beta^4 \pm O(\beta^2)$.
In \cite{GK25} the above result was stated in the fixed-$\beta$ and long-time regime, i.e. for $\int_\R \cZ^\xi_\beta(0,0;t,x)\,dx$ as $t\to\infty$ (these are equivalent by \cref{t:AJRS-Z-scaling}\ref{property-scaling}).
In those coordinates, they observed that under an OU perturbation of the white noise, the above result implies a noise sensitivity (decorrelation) threshold at perturbation scale $t^{-1}$. 
By similar reasoning, one can use \eqref{e:gk-clt} to conclude that as $\beta_1, \beta_1\to\infty$ with $\beta_2 > \beta_1$,
\begin{equation*}
    \Corr(Z_{\beta_1},Z_{\beta_2}) = 
    \begin{cases}
      1-o(1), & \beta_2^4 - \beta_1^4 \ll \beta_1^2,\\
      o(1), & \beta_2^4 - \beta_1^4 \gg \beta_1^2.
    \end{cases}
\end{equation*}
The above conditions are equivalent to $\beta_2 - \beta_1 \ll \beta_1^{-1}$ and $\beta_2 - \beta_1 \gg \beta_1^{-1}$.
In contrast, Conjecture \ref{c:main} predicts the onset of free energy decorrelation to occur when $\beta_2-\beta_1 \gg \beta^{1/3}_1$. 
This discrepancy reflects the difficulty of studying free energy properties in a low temperature/strong disorder setting, where the partition function is heavy-tailed and its moments and correlation structure reveal little about the free energy's typical behavior.

As an aside, accessing free energy correlations via moments of the partition function is possible given sufficient information about mixed fractional moments near the origin, since
\begin{align*}
    \Cov(\log Z_{\beta_1},\,\log Z_{\beta_2})
    =
    \partial_a\partial_b
    \log \E\bigl[Z_{\beta_1}^a Z_{\beta_2}^b\bigr]
    \Big|_{a=b=0},
\end{align*}
assuming the above is well-defined.
While a non-rigorous analysis along these lines was carried out in \cite{SY02} using the replica Bethe ansatz \cite{K87}, rigorous understanding of non-integer moments remains quite limited (see however \cite{DT21}).

\subsection{Coarse-graining and Conjecture {\ref{c:main}}}
A more promising approach to understanding free energy correlations is via a coarse-graining procedure, as we now heuristically describe.
This heuristic forms the basis of Conjecture \ref{c:main}.

The key observation is that, after an appropriate coarse-graining,
a perturbation of the inverse temperature
can be viewed as a perturbation of the \emph{effective disorder}.
This connects temperature chaos for the CDRP
to disorder chaos for last passage percolation (LPP) models, where the critical perturbation scale for disorder chaos is expected---and in some settings known---to be $n^{-1/3}$ (see \cite{GH20,GH23} and the discussion in \cref{ss:preface}).

\begin{figure}[tbh]
\centering
    \includegraphics[width=\linewidth]{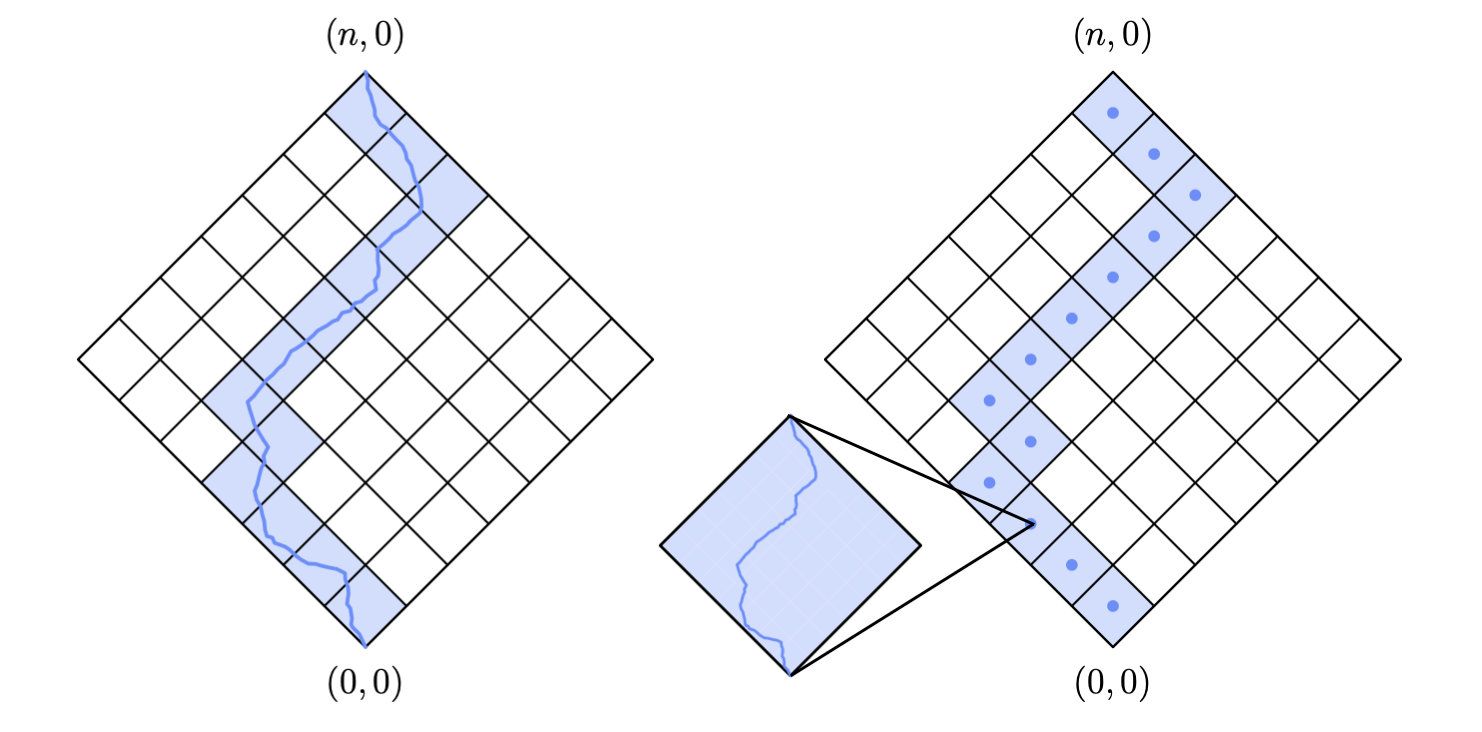}
\caption{
    Coarse-graining of the CDRP.\\
    \textbf{Left:}
    A trajectory of the length-$n$ CDRP is drawn as a continuous blue path.
    The space-time plane is divided into unit blocks, and the blocks the CDRP passes through (shaded blue) form a coarse version of its trajectory.\\  
    \textbf{Right:}
    The CDRP is erased, leaving behind its coarse trajectory.
    The blocks are identified with the vertices of $\Z^2$ (blue dots), and the CDRP's coarse trajectory forms an approximate sample from a lattice directed polymer model, where each block's disorder variable is the length-$1$ CDRP free energy across it (depicted as a short blue path in the inset diamond).
}
    \label{fig:coarse}
\end{figure}

We only discuss the one-point energetic version of Conjecture \ref{c:main} for the CDRP.
Let $\chi_1$ and $\chi_2$ be independent copies of the one-point directed landscape marginal $\cL(0,0;1,0)$. 
Recalling \eqref{energy1}, the predicted one-point energetic decorrelation statement is that as $\beta_1,\beta_2\to\infty$ with $\beta_2 > \beta_1$,
 \begin{align}\label{energy2} 
  \left(
      \mathfrak H^\xi_{\beta_1}(0,0;1,0),\;
      \mathfrak H^\xi_{\beta_2}(0,0;1,0)
  \right)
  \dto
  \begin{cases}
      (\chi_1,\chi_1),
      & \beta_2-\beta_1\ll \beta_1^{1/3}, \\[2mm] 
      (\chi_1,\chi_2),
      & \beta_2-\beta_1\gg \beta_1^{1/3}.
  \end{cases}
\end{align}
To relate this prediction to a more familiar disorder chaos picture, set $n\coloneqq \beta_1^4$ and $1+\varepsilon \coloneqq \frac{\beta_2}{\beta_1}.$
By the CDRP height-temperature scaling relation (\cref{t:AJRS-Z-scaling}\ref{property-scaling}),
the pair on the LHS of \eqref{energy2} has the same law (up to a deterministic centering and normalization) as the pair
\begin{align*}
      \left(\mathfrak H^\xi_{1}(0,0;n,0),\;\mathfrak H^\xi_{1+\varepsilon}(0,0;n,0)\right).
\end{align*}
Thus the problem is equivalent to comparing two long polymers in the same white noise environment, one at inverse temperature $1$ and the other at inverse temperature $1+\varepsilon$.

We now coarse-grain the space-time plane into unit blocks (boxes), arranged in a directed lattice structure as in Figure \ref{fig:coarse}.
For each unit block $B_v$, let $p_v$ and $q_v$ be its lower and upper endpoints, respectively.
Define the local block free energies 
\begin{align*}
    X_v \coloneqq  \log \mathcal Z^\xi_1(p_v;q_v),
    \qquad\qquad
    Y_v \coloneqq  \log \mathcal Z^\xi_{1+\varepsilon}(p_v;q_v).
\end{align*}
The full CDRP partition functions $\cZ^\xi_1(0,0;n,0)$ and $\cZ^{\xi}_{1+\e}(0,0;n,0)$ are not exactly equal in distribution to directed polymer partition functions built from the block variables $\{X_v\}_{v\in\Z^2}$ and $\{Y_v\}_{v\in\Z^2}$: 
the block variables are correlated, and the coarse-graining ignores boundary matching and other local effects.
Nevertheless, as a heuristic proxy, it is natural to compare the two coarse-grained partition functions
\begin{align*}
    \wt{Z}_1(0,0;n,0)
    &\coloneqq 
    \sum_{\pi}
    \exp\left(\sum_{v\in\pi} X_v\right),
    \qquad\qquad
    \wt{Z}_{1+\e}(0,0;n,0)
    \coloneqq 
    \sum_{\pi}
    \exp\left(\sum_{v\in\pi} Y_v\right),
\end{align*}
where the sums are over all directed lattice paths in the coarse-grained lattice.
In this coarse-grained description, both polymers have inverse temperature one and hence are in the \emph{strong disorder} regime.
The change in temperature has instead been converted to a change of effective disorder: the two disorder environments are now $\{X_v\}_{v\in\Z^2}$ and $\{Y_v\}_{v\in\Z^2}.$ 
Thus, the effect of coarse-graining is twofold: 
it transforms an intermediate disorder polymer model into a strong disorder polymer model, and it transforms a temperature perturbation into a disorder perturbation.

Since polymers in the strong disorder regime are expected to share the large-scale features of their zero-temperature counterparts, it is natural to compare this situation with disorder chaos in LPP.
Recall from \cref{ss:preface} that when Brownian LPP undergoes a disorder perturbation given by OU dynamics, the onset of geodesic chaos occurs at time scale $n^{-1/3}$ \cite{GH20,GH23}, and the same $n^{-1/3}$ threshold is expected to govern  geodesic energy decorrelation.
Since running OU dynamics for time $t$ causes the disorder's correlation with its initial state to decay by $1-e^{-t} \sim t$,
the critical time scale $t\asymp n^{-1/3}$ corresponds to a \emph{decay of disorder correlation} of order $n^{-1/3}$.
Therefore, for the coarse-grained model, we expect the onset of energetic decorrelation to also occur when 
\begin{align*}
    \Corr(X_v,Y_v)= 1-\Theta(n^{-1/3}).
\end{align*}

It remains to estimate $\Corr(X_v,Y_v)$.
A Taylor expansion suggests that we may pass from block free energies to block partition functions,
\footnote{
    Since $X_v, Y_v$ are free energies of unit-length CDRPs at bounded inverse temperatures, 
    they are of unit order, and hence $\E[(X_v-Y_v)^2]\asymp \E[(\exp(X_v)-\exp(Y_v))^2]$.
} 
and the correlation of the latter can be shown to be $1 - \Theta(\e^2)$ for small $\e>0$ using the chaos expansion \eqref{e:chaos-Z} or the Feynman--Kac formula \cite{BC95}.
Combining these observations, we find that the critical threshold for the onset of temperature chaos is defined through the relation
\begin{align*}
    1 - \Theta(n^{-1/3}) =  1 - \Theta(\e^2).
\end{align*}
Finally, recalling that $n=\beta_1^4$ and $1+\varepsilon=\frac{\beta_2}{\beta_1}$, we obtain the threshold predicted in Conjecture \ref{c:main}:
as $\beta_1\to\infty$,
\begin{align*}
    \left(\mathfrak{H}^\xi_{\beta_1}(\smallbullet,\smallbullet;\smallbullet,\smallbullet),
    \;\mathfrak{H}^\xi_{\beta_2}(\smallbullet,\smallbullet;\smallbullet,\smallbullet)\right)
    \dto
    \begin{cases} 
        (\cL_1, \cL_1),  & \beta_2-\beta_1\ll \beta_1^{1/3},
        \\[2mm]
        (\cL_1, \cL_2),  &\beta_2-\beta_1\gg \beta_1^{1/3}.
    \end{cases}
\end{align*}

\subsection{Rigorously proving energetic stability for a polynomial gap}\label{ss:polynomial-stability}

We will now briefly indicate how the H\"older regularity estimates of $\cZ^\xi_\smallbullet(\smallbullet,\smallbullet;\smallbullet,\smallbullet)$ proved in \cite{AJRS22} can be combined with a rigorous version of the above heuristic coarse-graining construction to prove energetic stability, i.e. $(\mf{H}^\xi_{\beta_1},\mf{H}^\xi_{\beta_2}) \dto (\cL_1,\cL_1)$, when $\beta_2-\beta_1 \ll \beta_1^{-\alpha}$ for some  $\alpha>0$.
Note that this is significantly weaker than Conjecture \ref{c:main}, which predicts that stability holds up to  $\beta_2-\beta_1 \ll \beta_1^{1/3}$ which in particular diverges as $\beta_1\to\infty$.
However, we choose to present this argument since it seems that the available estimates of regularity from \cite{AJRS22} by themselves only imply stability for a gap that is exponentially decaying in $\beta$.
Thus, coarse-graining offers a way to bootstrap estimates at small scales to stronger estimates at larger scales by leveraging the geometric structure of polymer models, a theme that has featured centrally in the recent works \cite{MR4586220,BSS23}.

As before, it is equivalent to consider polymers of length $n=\beta_1^4$ at inverse temperatures $1$ and $1+\e = \frac{\beta_2}{\beta_1}$, where $\e$ is a suitable power of $n^{-1}$ (to be specified later).
We first analyze the partition function $\cZ_\beta^\xi(0,0;n,0)$ for $\beta\in\{1, 1+\e\}$.
By the Chapman--Kolmogorov equation, we have
\begin{align*}
    \cZ_\beta^\xi(0,0;n,0)
    &= \int_{\R^{n-1}} \prod_{i=0}^{n-1} \cZ_\beta^\xi(i,x_i; i+1, x_{i+1})\,dx_1\dots dx_{n-1},
\end{align*}
where $x_0\coloneqq 0$ and $x_n\coloneqq 0$.
This decomposition plays a similar role as the coarse-grained lattice polymer model above.
While the length-$n$ CDRP is expected to behave like the geodesic, and hence have transversal fluctuations of size $n^{2/3}$, its finite-dimensional distributions are still absolutely continuous to those of Brownian bridge, and we will take advantage of this to obtain some crude bounds.
We fix some $K>0$ (we will later choose $K\gg n^{2/3}$) and define
\begin{align*}
    W_\beta
    &\coloneqq   
    \int_{[-K,K]^{n-1}} \prod_{i=0}^{n-1} \cZ^\xi_\beta(i,x_i;i+1,x_{i+1})\,dx_1\dots dx_{n-1},\\
    R_\beta &\coloneqq  \int_{\R^{n-1}\setminus [-K,K]^{n-1}} \prod_{i=0}^{n-1} \cZ^\xi_\beta(i,x_i;i+1,x_{i+1})\,dx_1\dots dx_{n-1},
\end{align*}
so that $\cZ^\xi_\beta(0,0;n,0) = W_\beta + R_\beta$.
We now show that the remainder $R_\beta$ is small using that 
$\E[\cZ_\beta^\xi(s,x;t,y)] = p(t-s,y-x)$ where $p(r,z)=\frac{1}{\sqrt{2\pi r}}e^{-z^2/2r}$ is the heat kernel, which follows from the chaos expansion \eqref{e:chaos-Z}.
By Fubini's theorem and the independence of $\cZ^\xi_\beta(i,\smallbullet;i+1,\smallbullet)$ across $i$ (\cref{l:independence-disjoint}\ref{independence} and \cref{t:AJRAS-Z-existence}\ref{property-adapted}), we have
\begin{align}\label{e:Rbeta-upper}
    \E[R_\beta] 
        &= \int_{\R^{n-1} \setminus [-K,K]^{n-1}}
            \prod_{i=0}^{n-1} 
            p(1, x_{i+1}-x_i)
            \,dx_1\dots dx_{n-1}\nonumber\\
        &= p(n,0) \cdot \mathrm{P}^{\mathrm{BB}}_{(0,0),(n,0)}\left(\max_{i\in\lb 1,n-1\rb} |X(i)| > K\right)\nonumber\\
        &\ls e^{-c_0 K^2/n},
\end{align}
where $X$ is a Brownian bridge from $X(0)=0$ to $X(n)=0$, and where $c_0>0$ is a universal constant.

The above implies that $\cZ^\xi_\beta(0,0;n,0) =  W_\beta + o(1)$ with high probability. However this is not enough to compare logarithms. Towards this, we will now prove a stronger multiplicative comparison, namely $\cZ^\xi_\beta(0,0;n,0) =  W_\beta(1+o(1))$, by combining the above bound on $R_\beta$ with a lower bound on $\cZ^\xi_\beta(0,0;n,0)$.
For the latter, by \cite[Theorem 1.1]{CG20b}, there exist  constants $C,c,c_1>0$ such that for all $n\ge 1$ and all $\beta\in\{1,1+\e\}$,
\begin{align}\label{e:Z-lower-tail}
    \P\left(
        \cZ^\xi_\beta(0,0;n,0) > e^{-c_1 n}
    \right)
    \ge 1 - Ce^{-cn^2}.
\end{align}
Denote 
\begin{align*}
    \kappa \coloneqq  e^{-c_1 n}.
\end{align*}
We set $K \coloneqq  (5 c_1/c_0)^{1/2}n$, so that by \eqref{e:Rbeta-upper}, we have $\E[R_\beta] \ls \kappa^5$ uniformly in $\beta\in\{1,1+\e\}$.
This along with Markov's inequality easily implies a lower tail bound for $W_\beta$:
\begin{align*}
    \P\left(W_\beta > \frac{\kappa}{2}\right)
    &= \P\left(\cZ^\xi_\beta(0,0;n,0) > \frac{\kappa}{2} + R_\beta\right)\\
    &\ge \P\left(\cZ^\xi_\beta(0,0;n,0) > \kappa\right)
    - \P\left(R_\beta > \frac{\kappa}{2}\right)\\
    &\ge 1 - Ce^{-cn^2}  - C'\kappa^4\\
    &\ge 1-C''e^{-c' n},
\end{align*}
uniformly in $\beta\in\{1,1+\e\}$.
We may now compare $\log \cZ^\xi_\beta(0,0;n,0)$ and $\log W_\beta$ using the inequality $\log(1+x)\le x$ for $x\ge 0$ along with the above probability bound:
\begin{align*}
    \P\left(
        \left|
            \log \cZ^\xi_\beta(0,0;n,0)
            - \log W_\beta
        \right|
        > \kappa
    \right)
    &=
    \P\left(
        \log\left(
            1+
            \frac{R_\beta}{W_\beta}
        \right)
        > \kappa
    \right)\\
    &\le
    \P(R_\beta > \kappa W_\beta)\\
    &\le
    \P\left(
        R_\beta > \frac{\kappa^2}{2}
    \right)
    + \P\left(W_\beta \le \frac{\kappa}{2}\right)\\
    &\ls \kappa^3 + e^{-c'n}\\
    &\ls e^{-c''n},
\end{align*}
where again everything is uniform in $\beta\in\{1,1+\e\}$.
We deduce that
\begin{align}\label{e:100000}
    \left|
        \log \left(\frac{\cZ^\xi_{1+\e}(0,0;n,0)}{\cZ^\xi_{1}(0,0;n,0)}\right)
        -
        \log\left(\frac{W_{1+\e}}{W_1}\right)
    \right|
    \xrightarrow[n\to\infty]{p} 0.
\end{align}

To establish energetic stability, it remains to show that for small enough $\e$, we have $|\log W_{1+\e} - \log W_1| \ll n^{1/3}$ with high probability.
Write
\begin{align*}
    \delta \coloneqq  
    \sup_{\substack{i\in \lb 0,n-1\rb\\x,y\in [-K,K]}}
    \left|
      \frac{\cZ_{1+\e}^\xi(i,x;i+1,y)}{\cZ_{1}^\xi(i,x;i+1,y)} - 1
    \right|,
\end{align*}
where as before $K\coloneqq (5c_1/c_0)^{1/2}n$.
On the event that $\d<1$, we have
\begin{align}\label{e:1000}
    (1-\d)^n \le \frac{W_{1+\e}}{W_1} \le (1+\d)^n.
\end{align}
It remains to choose $\e$ small enough (as a function of $n$) that $\d=o(1/n)$ with high probability.
In fact, we will show that it suffices to take $\e=n^{-\alpha}$ for some $\alpha>0$, using the $1^-$-H\"older continuity of $\beta\mapsto \cZ_\beta^\xi(s,x;t,y)$.
For this, it is convenient to use the normalized partition function $\wh{\cZ}_\beta^\xi(s,x;t,y) \coloneqq  p(t-s,y-x)^{-1}\cZ_\beta^\xi(s,x;t,y)$.
Fix any $\eta\in(0,1)$.
By Cauchy--Schwarz, the H\"older seminorm moment bounds in \cite[Proposition 3.8(i)]{AJRS22},
\footnote{Let us warn the reader that our notation does not match \cite{AJRS22}.
Their $Z_\beta$ is the same as our $\cZ_\beta$, and their $\cZ_\beta(s,x;t,y)$ is the same as our $\wh{\cZ}_\beta(s,x;t,y)$.}
and the negative moment bounds in \cite[Corollary 3.10(ii)]{AJRS22}, 
we have that for all sufficiently large $q\ge 1$,
\begin{multline*}
    \E[\delta^{q}]
    =
    \E\left[
        \sup_{\substack{i\in \lb 0,n-1\rb\\x,y\in [-K,K]}}
        \left|
            \frac{\wh{\cZ}_{1+\e}^\xi(i,x;i+1,y)}{\wh{\cZ}_{1}^\xi(i,x;i+1,y)} - 1
        \right|^q
    \right]\\
    \begin{aligned}
    &\le
    \E\left[
        \sup_{\substack{i\in \lb 0,n-1\rb\\x,y\in [-K,K]}}\frac{1}{\wh{\cZ}_{1}^\xi(i,x;i+1,y)^{2q}}
      \right]^{\frac12}
      \E\left[
        \sup_{\substack{i\in \lb 0,n-1\rb\\x,y\in [-K,K]}}
        \left|
            \wh{\cZ}_{1+\e}^\xi(i,x;i+1,y) - \wh{\cZ}_{1}^\xi(i,x;i+1,y)
        \right|^{2q}
    \right]^{\frac12}\\
    &\le
    \left(n
    \E\left[
        \sup_{x,y\in [-K,K]}\frac{1}{\wh{\cZ}_{1}^\xi(0,x;1,y)^{2q}}
      \right]\right)^{\frac12}
      \left(n
      \E\left[
        \sup_{x,y\in [-K,K]}
        \left|
            \wh{\cZ}_{1+\e}^\xi(0,x;1,y) - \wh{\cZ}_{1}^\xi(0,x;1,y)
        \right|^{2q}
    \right]\right)^{\frac12}\\
    &\ls_{q,\eta}
    n^{\frac12}K^{3q}
    \cdot
    n^{\frac12}K^q\e^{\eta q}\\
    &= (n^{1/q} K^4 \e^\eta)^{q},
    \end{aligned}
\end{multline*}
where in the third line we bounded the supremum over $i\in\lb0,n-1\rb$ by the sum.
Choosing some $q\ge 100$ and applying Markov's inequality, we get that
\begin{align*}
    \delta \le n K^4 \e^\eta
    \qquad\qquad\text{with probability $1-O(n^{-99})$.}
\end{align*}
Recalling that $K=\Theta(n)$, we now assume that $\e \le n^{-7/\eta}$, so that the above bound implies $\d \ls  n^{-2}$.
Then \eqref{e:1000} implies that $\log(W_{1+\e}/W_1) \pto 0$,
which by \eqref{e:100000} is equivalent to
\begin{align*}
    \left|\log \cZ^\xi_{1+\e}(0,0;n,0) - \log \cZ^\xi_1(0,0;n,0)\right| \xrightarrow[n\to\infty]{p}0.
\end{align*}
In particular, the two free energies differ by $o(n^{1/3})$ with high probability, and hence by \cref{t:KPZ-to-landscape-reformulated} they are strongly correlated and jointly converge (after centering and rescaling by $n^{-1/3}$) to the same copy of the directed landscape.
It remains only to convert back to $\beta_1=n^{1/4}$ and $\beta_2=\beta_1(1+\e)$.
Since $\eta\in(0,1)$ was arbitrary, this implies energetic stability in the regime $\beta_2-\beta_1 \ll\beta_1^{-\alpha}$ for any $\alpha> 28 - \frac14$. 
\qed

\appendix
\section{Properties of conditional variance}\label{s:conditional-variance}

This appendix records a few elementary properties of conditional variance.

The first lemma says that conditioning on less information increases conditional variance (in expectation).

\begin{lemm}[Conditioning on less information]\label{l:conditioning-on-less}
    Let $X$ be a square-integrable random variable on a probability space $(\Omega,\cF,\P)$, and let $\cG_1\subset \cG_2\subset \cF$ be any nested sub-$\sigma$-algebras.
    We have
    \begin{align*}
        \E[\Var(X|\cG_2)] \le \E[\Var(X|\cG_1)].
    \end{align*}
\end{lemm}
\begin{proof}
    Denote $\E_{\cG}X\coloneqq \E[X|\cG]$.
    By the law of total variance, the claimed inequality is equivalent to $\Var(\E_{\cG_1}X) \le \Var(\E_{\cG_2}X)$.
    We assume WLOG that $\E[X]=0$.
    Then by the tower property and conditional Jensen's inequality, we have
    \begin{align*}
        \Var(\E_{\cG_1}X)
        =
        \E\left[
            \bigl(\E_{\cG_1}X\bigr)^2
        \right]
        &= \E\left[
            \bigl(\E_{\cG_1}
                \E_{\cG_2}X\bigr)^2
        \right]\\
        &\le 
        \E\left[
            \E_{\cG_1}\bigl(\E_{\cG_2}X\bigr)^2
        \right]\\
        &= \Var(\E_{\cG_2}X)
    \end{align*}
    as desired.
\end{proof}

The next lemma is just the fact that conditional expectation is an orthogonal projection in $L^2$. 

\begin{lemm}\label{l:conditional-variance-variational}
    Let $X$ be a square-integrable random variable and let $Y$ be a random $d$-dimensional vector on the same probability space, for some $d\in\N$.
    Then
    \begin{align*}
        \E[\Var(X|Y)]
        = \E[(X-\E[X|Y])^2]
        = \inf_f \E[(X-f(Y))^2],
    \end{align*}
    where the infimum is over all measurable functions $f:\R^d\to\R$ such that $f(Y)$ is square-integrable.
\end{lemm}

The next lemma bounds differences of conditional variances in terms of $L^2$ norms.

\begin{lemm}\label{l:difference-conditional-variances}
    Let $X_1,X_2$ be square-integrable random variables on a probability space $(\Omega,\cF,\P)$, and let $\cG\subset \cF$ be any sub-$\sigma$-algebra.
    Then
    \begin{align*}
        \Bigl|
            \E[\Var(X_1|\cG)] - \E[\Var(X_2|\cG)]
        \Bigr|
        \le
        \bigl(\E[X_1^2]^{1/2} + \E[X_2^2]^{1/2}\bigr)\cdot\E\bigl[(X_1-X_2)^2\bigr]^{1/2}.
    \end{align*} 
\end{lemm}
\begin{proof}
    For any $V,W \in L^2(\Omega,\cF,\P)$, Cauchy--Schwarz implies $\left|\norm{V}^2-\norm{W}^2\right|\le \norm{V+W}\norm{V-W}$, where $\norm{\smallbullet} \coloneqq  \sqrt{\E[\smallbullet^2]}$ is the $L^2$ norm.
    In particular, writing $P\smallbullet \coloneqq  \E[\smallbullet|\cG]$, we have
    \begin{align*}
        \left|
            \norm{X_1-PX_1}^2 - \norm{X_2-PX_2}^2
        \right|
        &\le
        \bigl(
            \norm{X_1-PX_1} + \norm{X_2-P X_2}
        \bigr)
        \norm{X_1-X_2 - P(X_1-X_2)}\\
        &\le 
        \bigl(\norm{X_1} + \norm{X_2}\bigr)\norm{X_1-X_2},
    \end{align*}
    where we used that $I - P$ is an orthogonal projection, hence a contraction.
\end{proof}

We conclude this appendix with the proof of \cref{l:cond-var-comparison}.

\begin{proof}[Proof of \cref{l:cond-var-comparison}]
    Assume first that $|X|\le K$ is bounded.
    Then $\E[X|Z]=f(Z)$ for some Borel measurable $f:\R^d\to[-K,K]$.
    By \cref{l:conditional-variance-variational},
    \begin{align*}
        \E[\Var(X|Y)]
        &\le \E\left[\left(X-f(Y)\right)^2\right]\\
        &\le
        \E\left[(X-f(Z))^2\right]
        +
        \E\left[(X-f(Y))^2 \1_{Y\ne Z}\right]
        \\
        &\le
        \E[\Var(X|Z)] +
        4K^2\, \P(Y\ne Z).
    \end{align*}
    Repeating this argument with $Y,Z$ interchanged leads to
    \begin{align}\label{e:800}
        \Bigl|
            \E[\Var(X|Y)] - \E[\Var(X|Z)]
        \Bigr|
        \le 4K^2\, \P(Y\ne Z).
    \end{align}

    We now handle the unbounded case.
    Fix $K>0$ (to be specified later), and let $X_K \coloneqq  \max\{-K,\min\{X,K\}\}$ be the truncation of $X$ to $[-K,K]$.
    By \eqref{e:800}, we have
    \begin{align}\label{e:801}
        \Bigl|
            \E[\Var(X_K|Y)] - \E[\Var(X_K|Z)]
        \Bigr|
        \le 4K^2\,\P(Y\ne Z),
    \end{align}
    so we just need to compare $\E[\Var(X_K|U)]$ to $\E[\Var(X|U)]$ for $U\in\{Y,Z\}$.

    By \cref{l:difference-conditional-variances},
    \begin{align*}
        \Bigl|\E[\Var(X|U)] - \E[\Var(X_K|U)]\Bigr|
        &\le \bigl(\E[X^2]^{1/2}+\E[X_K^2]^{1/2}\bigr)
        \cdot\E[(X-X_K)^2]^{1/2}\\
        &\le 2\E[X^2]^{1/2}\cdot\E[(X-X_K)^2]^{1/2},
    \end{align*}
    where we used that $|X_K| \le |X|$.
    For the second factor on the RHS, using that $|X-X_K| \le |X|$, and combining Cauchy--Schwarz with the hypothesized tail behavior of $|X|$, we get
    \begin{align*}
        \E[(X-X_K)^2]^{1/2}
        &=\E\left[
            (X-X_K)^2\1_{|X|>K}
        \right]^{1/2}\\
        &\le \E\left[X^2\1_{|X|>K}\right]^{1/2}\\
        &\le \E[X^4]^{1/4}\cdot
        Ce^{-cK^{3/2}}.
    \end{align*}
    Since $\E[X^2]^{1/2}$ and $\E[X^4]^{1/4}$ are upper bounded by some constant depending only on $C,c$, we conclude via \eqref{e:801} that
    \begin{align*}
        \Bigl|
            \E[\Var(X|Y)] - \E[\Var(X|Z)]
        \Bigr|
        \le 4K^2\,\P(Y\ne Z) + C' e^{-cK^{3/2}}
    \end{align*}
    where $C'>0$ depends only on $C,c$.
    If $\P(Y\ne Z)>0$, then setting $K=\P(Y\ne Z)^{-0.05}$ readily implies the lemma.
    If instead $\P(Y\ne Z)=0$, then we just send $K\to\infty$ in the above display.
\end{proof}

\bibliographystyle{alpha}
\bibliography{bib.bib}

\end{document}